\documentclass[%
a4paper,twoside,11pt,titlepage=false,bibliography=plainheading]{scrbook} 
\usepackage{natbib}
\usepackage[latin1]{inputenc}
\usepackage[italian,english]{babel}

\usepackage{graphicx}
\usepackage{float}
\usepackage{wrapfig}
\usepackage[scaled=.90]{helvet}
\usepackage{textcomp}
\usepackage{framed} 
\usepackage[T1]{fontenc} %
\usepackage{lmodern} 

\usepackage{hhline}
\usepackage{multicol}

\usepackage[paper=a4paper,dvips=false,pdftex=false,vtex=false]{geometry}
\makeatletter{}
\let\orig@Itemize =\itemize
\let\orig@Enumerate =\enumerate
\let\orig@Description =\description
\def\Nospacing{\itemsep=0pt\topsep=0pt\partopsep=0pt%
\parskip=0pt\parsep=0pt}

\def\noitemsep{
\renewenvironment{itemize}{\orig@Itemize\Nospacing}{\endlist}
\renewenvironment{enumerate}{\orig@Enumerate\Nospacing}{\endlist}
\renewenvironment{description}{\orig@Description\Nospacing}%
{\endlist}
}
\def\doitemsep{
\renewenvironment{itemize}{\orig@Itemize}{\endlist}
\renewenvironment{enumerate}{\orig@Enumerate}{\endlist}
\renewenvironment{description}{\orig@Description}{\endlist}
}
\AtBeginDocument{\noitemsep}
\makeatother{}

\usepackage{pstricks}

\definecolor{mmlogored}{rgb}{0.93,0.25,0.34}
\definecolor{mmlogoblue}{rgb}{0.09,0.26,0.30}
\newcommand{\matemm}{%
  \textsf{\color{mmlogoblue}\textbf{mate}\color{mmlogored}matita}}

\definecolor{pdcolor3}{rgb}{.85,.95,.99}

\newlength{\tmplength}
\newsavebox{\tmpbox}
\newsavebox{\tmpboxb}

\definecolor{rossosfondo}{rgb}{.99,.82,.80}
\newenvironment{tutor}[1][Tutor]{%
  \sbox{\tmpboxb}{#1}%
  \begin{lrbox}{\tmpbox}\begin{minipage}{0.87\linewidth}%
      \raggedright%
    }{%
    \end{minipage}%
  \end{lrbox}%
  \begin{trivlist}
  \item
  \rotatebox[origin=mc]{90}{\usebox{\tmpboxb}}%
  ~%
  \psframebox*[framearc=0.2,framesep=6pt,fillcolor=rossosfondo]{\usebox{\tmpbox}}%
  \hspace*{\fill}%
\end{trivlist}
}

\usepackage{soul}
\sethlcolor{pdcolor3}
\let\bambini\hl

\definecolor{blulogo}{rgb}{0,.25,.48}
\definecolor{rossoscuro}{rgb}{0.93,0.25,0.34}
\newcommand{\attivita}[1]{{\fontfamily{pzc}\selectfont{}\color{rossoscuro}#1}\/}

\makeatletter
\newenvironment{consegna}{%
  \MakeFramed {\FrameRestore\@setminipage}}%
{\par\unskip\endMakeFramed}
\newenvironment{dadocente}{%
  \setlength{\fboxrule}{2pt}%
  \MakeFramed {\advance\hsize-30pt \FrameRestore}}%
{\endMakeFramed}
\makeatother

\makeatletter
\newenvironment{studente}[1][Studente]{%
  \sbox{\tmpboxb}{#1}%
  \setlength{\tmplength}{\hsize}%
  \addtolength{\tmplength}{-40pt}%
    \MakeFramed{\advance\hsize-40pt\FrameRestore\@setminipage}%
    }{\par\unskip\endMakeFramed}%
\makeatother

\newcommand{\materiali}{\par\medskip%
  {\color{rossoscuro}\textbf{Materiali consegnati agli alunni}}%
  \nopagebreak[4]%
  \par\smallskip}

\newenvironment{calendario}{\nopagebreak[4]\begin{multicols}{2}%
    \raggedright{}}{\end{multicols}}
\newenvironment{datiincontro}{}{}

\usepackage[bookmarks=false,colorlinks=true,citecolor=mmlogoblue,anchorcolor=mmlogored,urlcolor=mmlogored]{hyperref}

\begin{document}
\frontmatter{}%
\selectlanguage{english}%
\title{Promoting a practice of active student-centred instruction into
  the mathematics classroom:\\ \matemm{}'s ``turnkey laboratory'' kits} %
\author{Marina Cazzola} %
\date{} %
\maketitle

\setcounter{secnumdepth}{0}
\newcommand{\imgdir}{.}

\selectlanguage{english}

In spite of
research documenting the importance of active, student-centered
instruction, we have to regret that it is not widespread
in real teaching practice at school, as teachers too often rely on
self perpetuating ``traditional'' methods.

My experience at the University of Milano-Bicocca, as a lecturer in
the pre-service teacher training university degree course, led me to
promote Problem-Based Learning (PBL) as a methodology particularly
effective for the teaching of mathematics. In fact, I became convinced
that PBL is particularly suitable for developing a deeper
understanding of mathematical facts and methodologies, as well as for
leading the learners to acquire fundamental skills such as critical
thinking and problem solving abilities.

Meanwhile I started collaborating with the Center \matemm{}, an
Interuniversity Research Center for the Communication and Informal
Learning of Mathematics, founded in early 2005, having its origins in
the experience of promoting mathematics by four Italian universities:
Milano, Milano-Bicocca, Pisa and Trento. Among the many activities
carried out by \matemm{} \citep[e.g. see][]{czz2007}, the more
connected to these types of issues are the designing and testing of
problems suitable for conducting PBL sessions in the classroom,
alongside the investigation on the characteristics making a problem
suitable to be treated in a mathematical ``laboratory''.  Furthermore,
one of the fields of research of \matemm{} is the design of
mathematical interactive
exhibitions~\citep[e.g. see][]{czz2005:matemilano}, trying to connect
popularization of science with teaching: from this experience,
\matemm{} has conceived a set of prototypes of \emph{ready to use}
kits to be used to implement PBL laboratory sessions in the classroom.

The aim of this report is to describe the activities proposed in such
kits and document a field trial conducted during the academic year
2009-10.  In that year, the ``Ufficio scolastico regionale per la
Lombardia'' sponsored a refresher course for in-service primary and
secondary school teachers (grades K-8), namely ``Kit di laboratorio:
Tra numeri e forme''. During the course the participants experimented
the prototypes of the kits produced by \matemm{}, experiencing PBL in
their classrooms, under the supervision of university tutors. The
participants' reports are collected in this report.

\section{Methods}
These kinds of activities are located in the framework of PBL.  For a
neat definition of PBL we refer to~\citet{savery06:metabpl} according
to which PBL ``is an instructional (and curricular) learner-centered
approach that empowers learners to conduct research, integrate theory
and practice, and apply knowledge and skills to develop a viable
solution to a defined problem''. Typically a PBL session follows these
steps~\citep[e.g. see][]{czz2008b:madrid}:
\begin{itemize}
\item pupils are given a problem;
\item they discuss the problem and/or work on the problem in small
  groups, collecting information useful to solve the problem;
\item all the pupils gather together to compare findings and/or
  discuss conclusions; new problems could arise from this discussion,
  in this case
\item pupils go back to work on the new problems, and the cycle starts
  again.
\end{itemize}
The next section will be devoted to the description of the problems
proposed in the kits. We can not go into details of the theory,
however we refer to~\citet{hmelosilver04:whatlearn} for a report of
ongoing research on the aspects that make a problem suitable for these
kinds of activities, in order to allow students to draw full benefits.
\matemm{} experienced that, expecially with younger students, a game
can be the starting point for a PBL session, and in this sense, as we
will see, the kits contain many games as well as problems.

In support of teachers, each kit provides a guide rich of suggestions
for conducting the activities according to the principles described
above.

Before moving on to the description of the kits we need one last
observation. At a superficial glance one can think that
student-centered instruction has the effect of diminishing the role of
the teacher. In fact a work of this type requires a greater
involvement of the teacher, who has to assume the role of an
unobtrusive director.
Researchers have pointed out ``action'' that the teacher should be
ready to take~\cite[for a review e.g. see][]{schoenfeld92}, in this
context we like to quote the words of~\citet[p.~1]{howtosolve}
which indeed reveal the complexity of the task.
\begin{quote}
  One of the most important tasks of the teacher is to help his
  students. This task is not quite easy; it demands time, practice,
  devotion, and sound principles.

  The student should acquire as much experience of independent work as
  possibile. But if he is left alone with his problem without any help
  or with insufficient help, he may make no progress at all. If the
  teacher helps too much, nothing is left to the student. The teacher
  should help, but not too much and not too little, so that the
  student shall have a \emph{reasonable share of the work}.

  If the student is not able to do much, the teacher should leave him
  at least some illusion of independent work. In order to do so, the
  teacher should help the student discreetly, \emph{unobtrusively}.

  The best is, however, to help the student naturally. The teacher
  should put himself in the student's place, he should see the
  student's case, he should try to understand what is going on in the
  student's mind, and ask a question or indicate a step that
  \emph{could have occurred to the student himself}.
\end{quote}
In particular, the teacher must ``get involved'' personally and be
prepared to deal with unforeseen situations.

\section{The themes}

The laboratory activities proposed in \matemm{}'s kits usually start
with a concrete problem related to a
topic which is interesting from a mathematical point of view
(and as such is usually included in the school mathematics
curriculum). Often the problem is built around an object available
with the kit.

A detailed description (in Italian) of the kits can be found on
\matemm{}'s website
\begin{center}
  \url{http://specchi.mat.unimi.it/matematica/}.
\end{center}
Here we just try to give the relevant information needed for a better
understanding of the teachers' reports.
\begin{table}[tbp]
  \centering
\begin{tabular}{|p{0.13\textwidth}|p{0.80\textwidth}|}
  \hline
  \multicolumn{2}{|c|}{\href{http://specchi.mat.unimi.it/matematica/torri_kit.html}{\textbf{Torri, serpenti e \dots{} geometria}}} \\
\hline
  Level & Grades 1-5 \\ \hline
  Themes &
  Measure perimeters, areas and volumes
  \\ \hline
  Objects &
  `Nice to touch' foam tiles in shape of equilateral triangle for the 2D
  activities, wooden cubes for the 3D
  activities (for pictures see
  p.~\pageref{pic:torri:4} and p.~\pageref{pic:torri:9}).
  \\ \hline
  Activities
  &
  The worksheets describing the activities are available on
  \matemm{}'s website:
  \href{http://specchi.mat.unimi.it/matematica/torri/Schede_I_Elem_240409.pdf}{grade
    1},
  \href{http://specchi.mat.unimi.it/matematica/torri/Schede_II_Elem_240409.pdf}{grade
    2},
  \href{http://specchi.mat.unimi.it/matematica/torri/Schede_III_Elem_240409.pdf}{grade
    3},
  \href{http://specchi.mat.unimi.it/matematica/torri/Schede_IV_Elem_240409.pdf}{grade
    4}, and
  \href{http://specchi.mat.unimi.it/matematica/torri/Schede_V_Elementare_240409.pdf}{grade
    5},
  \noitemsep{}
  \begin{itemize}
  \item \attivita{Per cominciare}: estimate perimeters, areas and
  volumes of given geometric figures, without manipulatives.
\item Main 2D activities: with the tiles build a variety of polygon
  and measure perimeters and areas (discover that perimeter and area
  are independent quantities, solve simple optimization problems).
\item Main 3D activities: with the cubes build a variety of polyhedra
  and measure areas and volume (and deal with the analogous of the 2D
  problems).
\item \attivita{Per concludere}: review and summarize.
\end{itemize}
\\ \hline
\end{tabular}
\caption{Towers, snakes, and \dots{} geometry}
\end{table}
\begin{table}[tbp]
  \centering
\begin{tabular}{|p{0.13\textwidth}|p{0.80\textwidth}|}
  \hline
  \multicolumn{2}{|c|}{\href{http://specchi.mat.unimi.it/matematica/forme1.html}{\textbf{Giocare con le forme}}} \\
\hline
  Level & Grades 1-5 \\ \hline
  Themes & classification of shapes according to their geometric properties
  \\ \hline
  Objects &An ill structured set of `nice to touch' foam shapes
  (for pictures see pp.~\pageref{pic:formemorbide:2},
  \pageref{pic:formemorbide:8}, and~\pageref{pic:formemorbide:9},
  complete lists are available on \matemm{}'s
  website:
  \href{http://specchi.mat.unimi.it/matematica/forme1/tabellone_forme_I_II_attivit1_240409.pdf}{grades
    1-2},
  \href{http://specchi.mat.unimi.it/matematica/forme1/tabellone\%20forme_III_IV_V_attivit1_240409.pdf}{grades
    3-5});
  bingo tickets
  (\href{http://specchi.mat.unimi.it/matematica/forme1/cartelle_tombola_I_II_\%20III_240409.pdf}{grades
    1-3},
  \href{http://specchi.mat.unimi.it/matematica/forme1/cartelle_tombola_IV_V_240409.pdf}{grades
    4-5})
  \\ \hline
  Activities
  & %
  The worksheets describing the
  activities are available on
  \matemm{}'s website:
  \href{http://specchi.mat.unimi.it/matematica/forme1/Scheda\%20A\%20e\%20B\%20Giocare\%20con\%20le\%20forme\%20\%28BASSA\%29_180610.pdf}{grades
    1-2},
  \href{http://specchi.mat.unimi.it/matematica/forme1/Scheda\%20C\%20Giocare\%20con\%20le\%20Forme\%20\%28BASSA\%29_190610.pdf}{grade
    3}, and
  \href{http://specchi.mat.unimi.it/matematica/forme1/scheda_D_12_11_09.pdf}{grades
  4-5}.
\begin{itemize}
\item A pupil picks a shape without showing it, the other pupils have
  to guess which one it is through yes/no question about its geometric
  properties (the use of geometric figures names is not allowed)
\item Given two shapes, list difference and similarities between them
\item Bingo-like game with shapes
\end{itemize}
\\ \hline
\end{tabular}
\caption{Playing with shapes}
\end{table}
\begin{table}[tbp]
  \centering
\begin{tabular}{|p{0.13\textwidth}|p{0.80\textwidth}|}
\hline
\multicolumn{2}{|c|}{\href{http://specchi.mat.unimi.it/matematica/regolarivno.html}{\textbf{Diamo forma alla geometria: Regolari o no?}}} \\
\hline
Level & Grades 6-8 \\ \hline
Themes &Polyhedra, trying to establish the notion of regularity
\\ \hline
Objects &A set of \href{http://www.polydron.co.uk/}{Polydron} construction
\\ \hline
Activities
&
  The worksheets describing the
  activities are available on
  \matemm{}'s website.
  \begin{itemize}
  \item
    \href{http://specchi.mat.unimi.it/matematica/regolarivno/regolarivno_schedaA_15072010.pdf}{\attivita{Scheda~A}}:
    analyze examples of polyhedra
  \item
    \href{http://specchi.mat.unimi.it/matematica/regolarivno/regolarivno_schedaB_15072010.pdf}{\attivita{Scheda~B}}:
    write down your definition of
    regular polyhedron and check it
    against given examples
  \item
    \href{http://specchi.mat.unimi.it/matematica/regolarivno/regolarivno_schedaC_15072010.pdf}{\attivita{Scheda~C}}:
    weaker definition of regularity
  \item
    \href{http://specchi.mat.unimi.it/matematica/regolarivno/regolarivno_schedaD_15072010.pdf}{\attivita{Scheda~D}}:
    use regolar polygons to build
    examples of ``regular'' wallpaper patterns
  \end{itemize}
\\ \hline
\end{tabular}
\caption{Let's shape geometry: Regular or not?}
\end{table}
\begin{table}[tbp]
  \centering
\begin{tabular}{|p{0.13\textwidth}|p{0.80\textwidth}|}
  \hline
  \multicolumn{2}{|c|}{\href{http://specchi.mat.unimi.it/matematica/grandevpiccolo.html}{\textbf{Diamo forma alla geometria: Grande o piccolo?}}} \\
\hline
  Level & Grades 6-8 \\ \hline
  Themes &Polyhedra, trying to solve measure problems
  \\ \hline
  Objects &A set of \href{http://www.polydron.co.uk/}{Polydron}
  construction and a set of foam shapes
  \\ \hline
  Activities
  &
  The worksheets describing the activities are available on \matemm{}'s website.
  \begin{itemize}
  \item
    \href{http://specchi.mat.unimi.it/matematica//matematica/grandevpiccolo/SCHEDAA_19_12_08.pdf}{\attivita{Scheda~A}}:
    analyze examples of polyhedra
  \item
    \href{http://specchi.mat.unimi.it/matematica//matematica/grandevpiccolo/SCHEDAB_19_12_08.pdf}{\attivita{Scheda~B}}:
    build polyhedra and compare volumes
  \item
    \href{http://specchi.mat.unimi.it/matematica//matematica/grandevpiccolo/SCHEDAC_19_12_08.pdf}{\attivita{Scheda~C}}:
    build cubes and compare
    surfaces and volumes
  \item
    \href{http://specchi.mat.unimi.it/matematica//matematica/grandevpiccolo/SCHEDAD_19_12_08.pdf}{\attivita{Scheda~D}}:
    puzzles with polyhedra
  \item
    \href{http://specchi.mat.unimi.it/matematica//matematica/grandevpiccolo/SCHEDAE_19_12_08.pdf}{\attivita{Scheda~E}}: Pythagorean Theorem
  \end{itemize}
\\ \hline
\end{tabular}
\caption{Let's shape geometry: Big or small?}
\end{table}
\begin{table}[tbp]
  \centering
\begin{tabular}{|p{0.13\textwidth}|p{0.80\textwidth}|}
\hline
\multicolumn{2}{|c|}{\href{http://specchi.mat.unimi.it/matematica/viaggio.html}{\textbf{Il viaggio segreto: giochi di aritmetica}}} \\
\hline
Level & Grades 7-8 \\ \hline
Themes &
\begin{minipage}{1.0\linewidth}
  \begin{itemize}
  \item fractions
  \item criptography
  \end{itemize}
\end{minipage}
\\ \hline
Objects & The games
\\ \hline
Activities
&
  The rules of the games are available on \matemm{}'s website.
\begin{itemize}
\item Domino game with fractions: \href{http://specchi.mat.unimi.it/matematica/viaggio/tessere_Domino.PDF}{deck}
\item \attivita{Realizziamo una vetrata}, a bingo like game with
  fractions (see p.~\pageref{pic:vetrate:3} for pictures)
\item \attivita{Viaggio segreto}, decrypt the name of a city
\end{itemize}
\\ \hline
\end{tabular}
\caption{The secret journey}
\end{table}
Each kit provides a set of activities, built around a ``worksheet''
(the ``scheda''), guiding the pupils from the disclosure of the
problem to its solution, through a series of hints and required
steps. Some of the kits clearly name the activities (and such names
can be found in the activities descriptions written by the teachers),
others just number the activities with no names (and this sometimes
makes it difficult to follow the teachers' thread). Tables~1 to~5 sum
up the activities proposed by the kits used in this trial. Please note
that teachers were supposed to discuss with the tutors any changes to
the worksheets or to the activities, but, as it clearly emerges from
their reports, they sometimes did not.

\section{Comments and conclusions}

The aim of this report is to show how the idea of a practice of active
student-centred instruction is welcome in the reality of
teaching. This work is meant to provide a sort of snapshot of the
class implementation of the theory-based teaching recommendations, so
that any readers can draw their conclusions on their own.

Nevertheless I think it is worth pointing out two main aspects clearly
emerging from the teachers' reports.

Even when supported with structured activities, some of the teachers
seem to miss the main point of the activity and tend to focus on
marginal matters (marginal for a true mathematical point of
view). This is especially true with the kit ``Playing with shapes''
where you can see some of the teachers focusing on set terminology,
even if not needed for the activity, and on superficial knowledge,
e.g. having the children memorize the names of the shapes, rather than
have them work on properties of such shapes. There are two types of
possible intervention on this point. One way is to try to better
organize the kit. The other is to try to strengthen the mathematical
skills of teachers, with a suitable training, in order to make them
more aware of the fact that mathematics is not just a matter of
memorizing names and procedures.

Finally I have to point out that sometimes the environmental
constraints are a real barrier to these kinds of activities. I
quote in full the following particularly significant episode
(see p.~\pageref{ex:alessia}):
\begin{quotation}
  \noindent{}\emph{[close to the end of the activity]}

  \noindent{}A colleague comes into the classroom, she seems to show
  no interest at all in what we have been doing. Still holding the
  worksheet I invite the pupils to explain to the other teacher what we
  are doing, as she is joining us for the first time. I ask Rebecca
  (one of the kids with a very good profile) to speak
  \begin{studente}[Rebecca]
    \noindent{}our teacher gave us some coloured tiles and together
    with our classmates we used them to fill in the worksheet and
    build the figures
  \end{studente}
  \noindent{}My colleague notices that Alessia, a girl sometimes a bit
  childish but a willing student capable of great commitment, is
  chuckling with another girl and so asks her to explain what the
  lecture was about. Enthusiastically Alessia says
  \begin{studente}[Alessia]
    \noindent{}I PLAYED with the tiles
  \end{studente}
  \noindent{}My colleague reacts severely
  \begin{tutor}[ ]
    \noindent{}So, Alessia, did you play or did you work?
  \end{tutor}
  \noindent{}Suddenly Alessia becomes sad and falls silent. My
  colleague goes on, speaking to the whole class
  \begin{tutor}[ ]
    \noindent{}because, you know, you do not come here to play!
  \end{tutor}
  \noindent{}The climate has changed. Pupils have become frightened,
  they look down, their voices are trembling, [\dots{}] they are no
  longer able to work.

  \noindent{}My mood has changed too, from rage (with my colleague) to
  frustration, for not being able to cheer up the children.

  \noindent{}[\dots{}] I have many more things to say, but I think
  that the episode speaks for itself.
\end{quotation}

And I agree with her.

\clearpage{}%
\section{Credits}

I wish to thank
\begin{itemize}
\item the \emph{Ufficio Scolastico Regionale della Lombardia} for
  sponsoring the course;
\item the University of Milano-Bicocca, for supporting the activities
  by providing meeting rooms and computer facilities;
\item the tutors who supervised the groups and revised the teachers'
  reports: Carla Abrigo, Laura Bassani, e Mara Costa;
\item the teachers who tested the activities in their classrooms:
  Adele Alberti, Maria Aprigliano, Maria Aprile, Denise Bambara, Paola
  Brezigia, Antonella Bustreo, Patrizia Cappelletti, Federica
  Carluccio, Rosanna Carugo, Maria Chieffa, Antonella Cortina, Maria
  Giuseppa De~Sio, Rossella Fiamingo, Miriam Filippini, Eleonora
  Gatti, Sergio Gentili, Carmen Girola, Vilma Rosa Lampugnani, Franca
  Landi, Sabrina Legnani, Maria Pia Lettieri, Fedra Mandelli,
  Giuseppina Mellina, Gemma Negro, Alessandra Pacchioli, Giovanna
  Pizzo, Laura Sardi, Alessandra Sazio, Erminia Scolari, Claudio
  Sold\`a, Carla Speranza, Patrizia Tonso, Maria Giovanna Turconi,
  Antonietta Ventriglia, Raffaella Visconti.
\end{itemize}

\nocite{brena11:giocare,brena10:regolari,locatelli08:grande,locatelli09:viaggio}
\bibliography{ho}

\begin{thebibliography}{13}

\bibitem[Bertolini et~al.(2003)Bertolini, Cazzola, Dedò, Di~Sieno, Frigerio,
  Luminati, Poldi, Rampichini, Tamanini, Todesco, and Turrini]{matemilano2003}
M.~Bertolini, M.~Cazzola, M.~Dedò, S.~Di~Sieno, E.~Frigerio, D.~Luminati,
  G.~Poldi, M.~Rampichini, I.~Tamanini, G.~M. Todesco, and C.~Turrini.
\newblock \emph{matemilano, percorsi matematici in città}.
\newblock Springer-Verlag Italia, Milano, revised second edition, 2003.

\bibitem[Brena and Locatelli(2010)]{brena10:regolari}
A.~Brena and O.~Locatelli.
\newblock Diamo forma alla geometria: Regolari o no? (per la scuola secondaria
  di {I} grado).
\newblock Quaderni di Laboratorio,
  \url{http://specchi.mat.unimi.it/matematica/}, 2010.

\bibitem[Brena and Locatelli(2011)]{brena11:giocare}
A.~Brena and O.~Locatelli.
\newblock Giocare con le forme (per le classi della scuola primaria).
\newblock Quaderni di Laboratorio,
  \url{http://specchi.mat.unimi.it/matematica/}, 2011.

\bibitem[Cazzola(2005)]{czz2005:matemilano}
M.~Cazzola.
\newblock matemilano, mathematical explorations of the city.
\newblock In \emph{Proceedings Math\&Art conference}, Boulder Colorado USA,
  June 2005.

\bibitem[Cazzola(2007)]{czz2007}
M.~Cazzola.
\newblock Fare esperienza di matematica a scuola.
\newblock In \emph{Conorovesciato. Un esperimento di didattica per problemi
  nella scuola primaria}, Materiale per i Quaderni a Quadretti. Mimesis, 2007.

\bibitem[Cazzola(2008)]{czz2008b:madrid}
M.~Cazzola.
\newblock Problem-based learning and mathematics: Possible synergical actions.
\newblock In L.~G. Chova, D.~M. Belenguer, and I.~C. Torres, editors,
  \emph{ICERI2008 Proceedings}, Valencia (Spain), 2008. IATED (International
  Association of Technology, Education and Development).
\newblock ISBN 978-84-612-5091-2.

\bibitem[Hmelo-Silver(2004)]{hmelosilver04:whatlearn}
C.~E. Hmelo-Silver.
\newblock Problem-based learning: What and how do students learn?
\newblock \emph{Educational Psychology Review}, 16\penalty0 (3):\penalty0
  235--266, September 2004.

\bibitem[Locatelli(2006)]{locatelli06:torri}
O.~Locatelli.
\newblock \emph{Torri, serpenti e\dots{} geometria}.
\newblock Quaderni a quadretti--Quaderni di laboratorio. Mimesis, 2006.

\bibitem[Locatelli(2008)]{locatelli08:grande}
O.~Locatelli.
\newblock Diamo forma alla geometria: Grande o piccolo? (per la scuola
  secondaria di {I} grado).
\newblock Quaderni di Laboratorio,
  \url{http://specchi.mat.unimi.it/matematica/}, 2008.

\bibitem[Locatelli(2009)]{locatelli09:viaggio}
O.~Locatelli.
\newblock Viaggio segreto (per le classi seconde e terze della scuola
  secondaria di {I} grado).
\newblock Quaderni di Laboratorio,
  \url{http://specchi.mat.unimi.it/matematica/}, 2009.

\bibitem[Polya(1945)]{howtosolve}
G.~Polya.
\newblock \emph{How to solve it: a new aspect of mathematical method}.
\newblock Princeton University Press, 1945.

\bibitem[Savery(2006)]{savery06:metabpl}
J.~R. Savery.
\newblock Overview of problem-based learning: Definition and distinctions.
\newblock \emph{The Interdisciplinary Journal of Problem-based Learning},
  1\penalty0 (1):\penalty0 9--20, 2006.

\bibitem[Schoenfeld(1992)]{schoenfeld92}
A.~H. Schoenfeld.
\newblock Learning to think mathematically: Problem solving, metacognition, and
  sense-making in mathematics.
\newblock In D.~Grouws, editor, \emph{Handbook for Research on Mathematics
  Teaching and Learning}, pages 334--370. MacMillan, New York, 1992.

\end{thebibliography}
\bibliographystyle{modplainnat}

\mainmatter{}
\selectlanguage{italian}

\cleardoublepage{}
 \newcounter{saveenum}

\bigskip

{\centering \includegraphics[height=3cm]{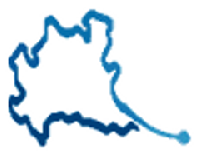}}

\bigskip

\bigskip

\bigskip

{\centering\selectlanguage{italian}\sffamily\bfseries INDICAZIONI PER
  IL CURRICOLO
\par}

{\centering\selectlanguage{italian}\sffamily\bfseries AREA
  MATEMATICO-SCIENTIFICO-TECNOLOGICA
\par}

\bigskip

{\centering\selectlanguage{italian}\sffamily FORMAZIONE PER DOCENTI
  DELLE SCUOLE DELL'INFANZIA, PRIMARIA E SECONDARIA DI
  PRIMO GRADO
\par}

\bigskip

{\centering\selectlanguage{italian}\sffamily\bfseries SECONDA FASE
\par}

{\centering\selectlanguage{italian}\sffamily\bfseries LABORATORI
\par}

\bigskip

\bigskip

{\centering\selectlanguage{italian}\huge\sffamily\bfseries LABORATORIO
  \\``Kit di laboratorio: Tra numeri e
  forme''
\par}

\bigskip

\chapter*{Relazione del docente universitario responsabile}
\addcontentsline{toc}{chapter}{Relazione del docente universitario responsabile}

\section{Valore scientifico del percorso del laboratorio}

Il laboratorio ``Kit di laboratorio: tra numeri e forme'' si proponeva
di condurre i docenti della scuola dell'obbligo a una sperimentazione
in classe dei materiali strutturati (i ``kit'') prodotti dal centro
\matemm{} per la realizzazione di un ``laboratorio di matematica'',
secondo le modalità didattiche suggerite nelle indicazioni nazionali
oggetto del ciclo di seminari (e in realtà suggerite nei programmi
ministeriali da lunga data, ma di cui mai veramente si è vista una
affermazione nella pratica scolastica generale).

La prima parte del corso si è svolta con una serie di incontri in
presenza, in cui si sono riviste le indicazioni e si è cercato di far
sì che i docenti iscritti al seminario vivessero in prima persona una
esperienza di quel ``fare matematica'' tanto auspicata non solo dalle
indicazioni stesse, ma anche dalla ricerca didattica. Si sono quindi
distribuiti i 6 kit di cui l'ufficio scolastico regionale ci ha dotato
per permetterne una esplorazione. Durante gli incontri in presenza si
è infine iniziata la progettazione dei percorsi didattici,
progettazione proseguita con interazione a distanza attraverso
internet. I docenti hanno svolto nelle classi il loro percorso
didattico e, via internet, hanno relazionato, producendo, sotto la
guida dei coordinatori dei gruppi, una ricca documentazione (vedi
capitoli successivi). La documentazione prodotta offre a mio avviso il
punto di partenza per una riflessione sulle azioni che occorre
compiere per una effettiva applicazione in classe dei principi
didattici auspicati dalle indicazioni. Se da un lato i docenti
partecipanti hanno colto l'importanza di un certo tipo di attività e
hanno condiviso lo spirito che ha ispirato il seminario, nondimeno
dalle loro riflessioni contenute nelle rendicontazioni finali (inviate
direttamente all'ufficio) emerge la difficoltà nell'organizzare i
tempi e modi per attività come quelle proposte, per vari motivi:

\begin{itemize}
\item le attività organizzate in modo laboratoriale richiedono più
  tempo di quelle organizzate in modo trasmissivo; una buona parte
  degli iscritti ha scelto di non proseguire la sperimentazione quando
  ha realizzato il tipo di impegno che era richiesto;
\item i docenti si sentono vincolati al ``programma'' e hanno scelto
  di sperimentare i kit immediatamente riconoscibili come
  riconducibili ai ``contenuti'' e ai ``traguardi'', mentre altri temi
  non sono stati presi in considerazione;
\item i docenti che hanno scelto di partecipare attivamente al
  seminario sono in generale tra i più attivi e aperti alle
  innovazioni didattiche e sono stati in generale supportati dalle
  rispettive strutture scolastiche; nondimeno alcuni hanno trovato nei
  colleghi ``indifferenza'' o peggio hanno avuto difficoltà a far
  accettare la sperimentazione di un percorso ``diverso'' nella
  propria scuola.
\end{itemize}

Grazie al seminario si è potuto creare un ``luogo'' di confronto tra
docenti di scuole diverse (e tra docenti e ricercatori universitari),
creando la possibilità di dare vita a progetti comuni. I ``kit'' sono
stati riconosciuti dai docenti come strumenti estremamente validi,
tanto che non pochi tra i docenti partecipanti hanno auspicato che la
loro scuola si potesse dotare di tali materiali.

Dei 65 iscritti iniziali, 35 hanno portato a termine la
sperimentazione.

La sperimentazione ha potuto aver corso sia grazie all'unità
Città~Studi del centro \matemm{}, che ha prodotto i kit, sia grazie
all'università di Milano-Bicocca (in particolare il Dipartimento di
matematica e applicazioni, la Facoltà di scienze della formazione e il
Settore laboratori didattici di ateneo) che ha dato supporto logistico
e informatico (offrendo per i lavori a distanza la piattaforma
\textit{moodle} di ateneo, \url{http://ateneo.elearning.unimib.it},
piattaforma su cui è disponibile la documentazione completa, sia della
fase di progettazione che della fase di realizzazione dei vari
percorsi).

\section{Eventuali spunti di discussione sulle Indicazioni nazionali}

{Oltre agli aspetti già evidenziati al punto precedente, una modalità
  didattica laboratoriale richiede al docente un grosso sforzo per
  ``mettersi in gioco'' e mettere in discussione una prassi didattica
  consolidata. Nei resoconti dei docenti emerge la difficoltà del
  trovarsi davanti a situazioni inaspettate, a reazioni impreviste nei
  bambini e a dubbi sul cammino intrapreso (``\textit{che cosa ho
    combinato?}'', ``\textit{quanta confusione ho creato nei miei
    allievi?}''), dubbi rientrati nei passi successivi del
  percorso. Credo che per una effettiva attuazione delle indicazioni
  nazionali sia necessario dare ai docenti un supporto in questi
  momenti di impasse (come spero che sia stato questo seminario) e non
  lasciare tutto il peso dell'innovazione didattica sulle spalle della
  buona volontà di singoli.}

\section{Valutazione complessiva e prospettive di sviluppo}

Complessivamente ritengo si siano poste le basi per un contatto
scuola-università che porti nella direzione auspicata dalle
indicazioni. Tutti i docenti che hanno portato a termine la
sperimentazione hanno manifestato entusiasmo per i risultati ottenuti
nelle proprie classi, anche solo dal punto di vista della motivazione
degli allievi (spesso anche e soprattutto di quelli ritenuti
``deboli''). Tutti i docenti si sono impegnati ben oltre l'orario
``ufficiale'' del seminario per completare la documentazione e
condividere la propria esperienza.

Per una riproposizione del percorso ci sono alcuni aspetti da
riconsiderare:

\begin{itemize}
\item per quel che riguarda la scelta degli argomenti, ci sono stati
  alcuni ``kit'' che sono stati richiesti da molti docenti e che
  quindi sono stati messi a disposizione delle classi per un tempo
  appena sufficiente a condurre la sperimentazione, mentre forse
  sarebbe stato opportuno condurre una esplorazione più approfondita
  (in una impostazione laboratoriale il tempo è una variabile chiave
  mentre molti docenti si sono trovati costretti a terminare
  tassativamente l'attività il giorno xx perché poi il kit doveva
  essere consegnato a un'altra scuola). Per una riproposizione, se da
  un lato questa sperimentazione porta a concludere che il ``kit''
  costituisca uno strumento didattico valido (e quindi sia opportuno
  che le scuole, nei limiti del possibile, come molti docenti nei loro
  resoconti auspicano, si attrezzino), dall'altra sarebbe auspicabile
  che i docenti si cimentino anche con quegli argomenti che a prima
  vista hanno considerato fuori dalla loro portata;
\item la comunicazione via internet si è rivelata molto difficoltosa e
  credo che sia opportuno prevedere più incontri dal vivo, con ovvie
  difficoltà logistico-organizzative (tra le altre cose, spesso i
  docenti non riescono a avere il pomeriggio libero per partecipare
  agli incontri).
\end{itemize}
Per un futuro proseguimento dell'attività molto docenti auspicano si
possa anche pensare a un percorso per verificare la permanenza nei
loro allievi dei concetti acquisiti con questa modalità didattica,
ponendo una questione di indubbio interesse.

\section{Materiale}

La documentazione completa, sia della fase di progettazione, che della
fase di realizzazione dei percorsi è stata raccolta sul sito
\url{http://ateneo.elearning.unimib.it}. Questo fascicolo contiene la
rendicontazione e rielaborazione dei vari percorsi, estratta da tale
sito.

\setcounter{secnumdepth}{1}
\chapter{Torri, serpenti e \dots{} geometria}


\section[Sperimentazione \#1: quarta primaria]{Sperimentazione \#1:
  classe quarta primaria, novembre--dicembre~2009}

\subsection{Osservazioni generali}

\subsubsection{Presentazione della classe}

Si tratta di una quarta composta da 23 alunni. Gli alunni hanno, in
termini geometrici, solo un concetto intuitivo di perimetro e nessuno
di area; a ciò si aggiunge il fatto che molti sono ottimi esecutori ma
scarsi ``rielaboratori'' concettuali.

\subsubsection{Composizione dei gruppi}
Gruppi omogenei (4 alunni per gruppo, stabiliti
dall'insegnante). Decido di formare gruppi formati ciascuno da 4
alunni (un gruppo da 3, il ``più alto'') più o meno di eguale
livello. Questo perché in attività di gruppo precedenti si era
verificato che i ``bravi'' troppo velocemente eseguissero da soli
tutto il lavoro non dando la possibilità ai compagni in difficoltà di
lavorare e comprendere.

I gruppi sono 6 e di ``livello'', all'interno di ciascun gruppo si è
cercato di mantenere un margine di eterogeneità:
\begin{itemize}
  \item uno di livello alto (A)
  \item uno di livello medio alto (B): un bambino è non italofono,
    cinese con notevoli difficoltà legate alla lingua ma
    cognitivamente e didatticamente MOLTO bravo
  \item due di livello medio (C e D)
  \item due di livello medio basso (E e F)
\end{itemize}
I gruppi sono rimasti invariati per tutto il percorso.

\subsubsection{Insegnanti presenti}

Presente solo l'insegnante di classe.

\subsubsection{Calendarizzazione degli incontri}
\begin{calendario}
  \begin{itemize}
\item 24 novembre
\item {}[\dots{}]\footnote{Non sono disponibili le date dei due
    incontri successivi. Si tratta della prima sperimentazione e
    l'insegnante non aveva ancora a disposizione lo schema di
    presentazione, delle attività che è stato costruito in seguito;
    per questo motivo non sono presenti tutte le informazioni raccolte
    invece per le sperimentazioni successive.}
\end{itemize}
\end{calendario}

\subsection{Primo incontro}

\begin{description}\item[Alunni presenti:]23 (tutti)
\item[Tempo effettivo di lavoro:]dalle 14.30 alle 16.30
\end{description}

\begin{consegna}
  In gruppo ideazione, disegno e individuazione del nome della
  mascotte del gruppo, consegna della scheda 1 \attivita{Per
    cominciare} (per le classi terze, anche se i bambini sono di IV,
  vedi poi)

  \materiali{}
  I materiali previsti dal kit
\end{consegna}

\subsubsection{Osservazioni}
L'attività si divide in più fasi.

Fase~1 (formazione dei gruppi):
Comunico agli alunni che lavoreremo laboratorialmente e scoppia
l'entusiasmo generale (sono sempre MOLTO felici di attività di gruppo
e manipolative) iniziando a tempestarmi di domande per sapere cosa
faremo.

Accenno loro che sarà un laboratorio di geometria in cui giocheranno
con del materiale di cui dovranno avere molta cura.

I bambini protestano venendo a sapere che i gruppi sono stati formati
da me dicendo che vogliono stare con i loro ``amichetti''. Ricordo
loro le difficoltà emerse durante l'ultima attività di gruppo in cui
alcuni si erano lamentati di aver dovuto lavorare da soli perché gli
altri componenti chiacchieravano e, altri si erano lagnati di non
essere riusciti a lavorare perché i compagni avevano eseguito tutto da
soli; spiego loro che la mia decisione, unita al fatto che ogni gruppo
avrà max 4 componenti, altro non è che ``una prova per vedere se il
problema si risolve''. E le proteste si trasformano in assensi e
consensi.

Anche davanti alle ``nuove'' regole di lavoro (cioè quelle suggerite
nella guida docente del kit, ``nuove'' per questa classe) sembrano
timorosi e necessitano di essere rassicurati sul fatto che i ruoli, se
pur ben definiti, possano essere scambiati.

Questa fase di informazione richiede circa 40 minuti! Più tempo del
previsto.

Fase~2 (individuazione nome e mascotte):
Ciascun gruppo è invitato a pensare un nome e disegnare una mascotte
che lo identifichi e lo rappresenti.

Immediatamente emergono le prime difficoltà: i gruppi con i bimbi più
deboli faticano a accordarsi. Anche i 2 gruppi ``medi'' hanno
difficoltà a trovare armonia.

Più volte devo intervenire per suggerire soluzioni e sedare gli
animi\dots{} (non immaginavo tante difficoltà in questa che doveva
essere la fase ``più semplice'').

A questo punto è trascorsa più d un'ora e qualcuno sembra essere già
stanco.

Fase~3 (\attivita{Per cominciare}):
Consegno loro la scheda del kit.

Solo il gruppo B la compila correttamente.

Gli altri contano solo le tessere della punta della stella!!

Fase~4: attività \attivita{Tessere 2}
\begin{description}
\item[attività 2.1] gli alunni sono entusiasti del materiale loro
  dato\dots{} dopo qualche minuto di osservazione e manipolazione
  libera consegno loro la scheda.
  \begin{itemize}
  \item gruppo A: si accorge di aver sbagliato la scheda precedente:
    \bambini{maestra non ho contato le tessere dentro la stella!!}. non
    incontra problemi il gruppo A: collaborano e si divertono.
  \item gruppo B: il bambino non italofono si annoia e sbuffa\dots{}
    (%
    i compagni gli hanno dato il ruolo di lettore, ma non comprende
    ciò che ``legge''), più volte intervengo per cercare di stimolarlo
    ma invano: è frustrato perché non capisce cosa deve fare, anche i
    compagni sono preoccupati per lui e perplessi perché non sanno cm
    aiutarlo, dico loro di stare tranquilli e che sicuramente andando
    avanti troveremo il modo\dots{} l'attività in sé procede
    bene\dots{}
  \item gruppo C: collaborano, provano, riprovano e discutono\dots{}
  \item gruppo D: disegna la figura non nella parte strutturata del
    foglio ma sul margine bianco\dots{} i componenti faticano a
    collaborare.
  \item gruppo E: un bambino non capisce cosa deve fare e non si fida
    dei compagni, che esasperati richiedono il mio intervento\dots{} la
    bambina con maggiori difficoltà sembra non interessata\dots{}
  \item gruppo F: litigano perché non riescono a collaborare, tutti
    vorrebbero usare le tessere e disegnare.
  \end{itemize}
\item[attività 2.2]
  gli alunni fanno fatica, sono stanchi, lo capisco dai loro visi
  e dal fatto che chiedono l'ora, (qualcuno lo dice apertamente\dots{})
  ricercano il mio aiuto e si crea molto rumore. Vista l'ora (sono già
  le 16:20) e la confusione decido di interrompere le attività.
  Anche io sono molto stanca!

  INTEGRAZIONE CONSIDERAZIONI PRIMO LAB: avendo una classe con
  all'interno 2 alunni molto deboli cognitivamente avevo pensato d
  partire con le schede della seconda e terza elementare, per tre
  motivi fondamentali:
  \begin{enumerate}
  \item dare la possibilità a tutti gli alunni di eseguire/seguire il
    più possibile almeno la prima parte del laboratorio con successo,(
    soprattutto ai bimbi con maggiore difficoltà);
  \item far ``rodare'' il gruppo e creare un clima di collaborazione,
    poiché a mio avviso, le attività previste erano abbastanza
    semplici da far sì che le uniche difficoltà emergenti fossero di
    tipo collaborativo/organizzativo;
  \item verificare il livello dell'intero gruppo classe.
\end{enumerate}
\end{description}
Nelle mie previsioni avremmo dovuto eseguire tutte le attività
previste per il primo incontro\dots{}  Così non è stato per diversi
motivi:
\begin{itemize}
\item per costituire i gruppi abbiamo impiegato più tempo del
  previsto;
\item le difficoltà di tipo collaborativo sono state diverse e hanno
  richiesto interruzioni dell'attività nel singolo gruppo per cercare
  insieme il modo di risolverle;
\item nel pomeriggio il rendimento degli alunni cala
  considerevolmente.
\item è stato quasi impossibile seguire e stimolare adeguatamente
  tutti e 6 i gruppi di lavoro,
\item ciascun gruppo richiedeva il mio intervento per risolvere
  controversie o appianare malumori, avere conferme sul lavoro svolto
  o semplicemente per farmi vedere le figure costruite\dots{}
\item ogni gruppo richiedeva insomma un certo margine di tempo da
  dedicargli spesso contemporaneamente e io non ci sono
  riuscita\dots{}

  ad esempio, mentre il gruppo D litigava per chi dovesse disegnare la
  figura e per quale figura disegnare, il gruppo E si arrovellava e
  chiedeva il mio intervento perché non riusciva a disegnare la figura
  (visto che, nonostante le mie spiegazioni, avevano accostato le
  figure non per i lati ma per le ``punte'' e, anche se un bimbo lo
  faceva presente, continuavano a volerla disegnare in quel modo\dots{})
\item la fase di preparazione e formazione dei gruppi con nome e
  mascotte richiede più tempo del previsto\dots{}
\item gli alunni sono entusiasti del materiale consegnato!
\item riusciamo a eseguire solo parte del lavoro programmato\dots{}
\end{itemize}

\subsubsection{Considerazioni finali}

\begin{itemize}
\item le attività hanno richiesto più tempo di quanto immaginassi,
\item è stato quasi impossibile seguire e stimolare adeguatamente
  tutti e 6 i gruppi di lavoro,
\item le ore pomeridiane non sono probabilmente le più adeguate per
  questo tipo di attività, mi riprometto di effettuare il prossimo
  laboratorio di mattina\dots{}
\end{itemize}

\subsection{Secondo incontro}

\begin{description}
\item[Alunni presenti:]23 (tutti)
\item[Tempo effettivo di lavoro:]
Dalle ore 9 alle ore 10:30
\end{description}

\begin{consegna}
  Si affrontano le attività sulle TESSERE previste per la CLASSE II e
  la CLASSE III seguendo le INDICAZIONI DEL MANUALE presente nel kit.
  \materiali{}%
  Il materiale messo a disposizione di ciascun gruppo è costituito da
  16 tessere, 3 cordicelle, 1 penna, 1 gomma, 1 matita e pastelli
  colorati.
\end{consegna}

\subsubsection{Osservazioni}
Il lavoro nei gruppi si è svolto in modo lineare e sereno senza
le difficoltà emerse nell'incontro precedente.

Sono stati confermati i gruppi iniziali (eterogenei tra di loro)
costituiti dall'insegnante per i motivi già elencati.

Rispetto al primo incontro in cui i bambini meno brillanti non
riuscivano a apportare il loro contributo e i bambini più
accentratori rifiutavano di confrontarsi, in questo laboratorio i
primi hanno partecipato più attivamente se pur con contributi non
sempre utili, e i secondi hanno iniziato a aprirsi alla
collaborazione.

I ``bravi'' non hanno
fatto alcuna fatica a accettare spunti e consigli stimolando e
coinvolgendo anche i compagni in difficoltà.

Molti alunni hanno cercato di adottare diverse tecniche di
autoverifica quali confrontare le risposte, verificare concretamente
le diverse soluzioni prima di scegliere quella giusta, rileggere
attentamente il testo.

Altri, i più immaturi, i più deboli e i più accentratori non
hanno ricercato alcun tipo di verifica, anche messi davanti all'errore
hanno faticato a accettare un consiglio o a comprendere il perché
dell'errore.

Tutti hanno dimostrato un impegno e un interesse alto nell'esecuzione
delle attività e un'alta motivazione data dal ``gusto della scoperta
e della riuscita''. Solo alcuni che risultano essere ottimi esecutori
ma poco autonomi nella rielaborazione dei contenuti hanno avuto
momenti di disinteresse, rientrati grazie all'intervento mirato a
ottenere una risposta positiva cui far seguire una gratificazione.

Gli alunni più deboli hanno resistito fino in fondo grazie anche ai
compagni che li incitavano e stimolavano e alla maestra che passando
chiedeva loro qualche chiarimento su quanto stessero facendo.

Non vi è stata competizione tra i diversi gruppi.

In generale i componenti dei gruppi con più difficoltà collaborative
hanno cercato maggiormente di collaborare anche se al loro interno
c'era ancora un clima di competizione.

Nessuna difficoltà di comprensione del testo a opera dei brillanti e
dei ``medi'', qualcuna per i più deboli, superata all'interno del
gruppo o con richiesta di conferma all'insegnante.

\subsection{Terzo incontro}

\begin{description}
\item[Alunni presenti:]23 (tutti)
\item[Tempo effettivo di lavoro:] 10.45-12.30 (tempo effettivo 1 ora)
\end{description}

\begin{consegna}
  \par Le attività \attivita{Per cominciare}, \attivita{attività
    tessere}, \attivita{per concludere} previste per la classe IV.
  \materiali{} %
  I materiali previsti dal kit
\end{consegna}

\subsubsection{Osservazioni}
Chiedo di iniziare facendo mente locale a quanto effettuato
precedentemente. Vedendo il viso perplesso di alcuni alunni consegno a
ciascun gruppo il proprio fascicolo di schede e insieme cerchiamo di
fare il punto della situazione\dots{} chiedo loro se ricordano\dots{}

Anche così gli alunni più deboli sembrano cadere delle nuvole\dots{} un
paio addirittura dicono di non ricordare di aver svolto quelle
attività\dots{}  Solo gli alunni più brillanti dicono di ricordare bene
le attività svolte.

Davide ci racconta che:
\begin{studente}[Davide]%
  lavorando con la clessidra, la lumaca e il gatto mi sono accorto
  che pure se sono uguali in qualcosa sono diversi in
  qualcos'altro\dots{} cioè che pure se sono uguali non è che sono uguali
  in tutto proprio\dots{} cioè non per forza in tutto\dots{}
\end{studente}
io gioisco ma i compagni gli dicono che non hanno capito niente di
quello che ha detto!

Li invito a osservare allora la tabella dell'attività 2 (Classe III) e
di provare a capire, discutendo tra loro, se quello che aveva detto
Davide poteva avere un senso\dots{} io intanto riproduco alla lavagna la
tabella\dots{}

Qualcuno esclama \bambini{è vero!} qualcun altro risponde \bambini{ma
  cosa?}.

Invito tutti a guardare la tabella e i dati che essa contiene e pongo
loro la domanda:
\begin{tutor}[Ins.]
  cosa notate?
\end{tutor}
\begin{studente}[Bambini]%
\begin{itemize}
  \item \bambini{gatto e lumaca valgono uguali}
  \item  \bambini{gatto e lumaca sono con 8 tessere}
  \item  \bambini{però il contorno non vale uguale!}
  \item  \bambini{però clessidra e lumaca hanno il perimetro uguale!}
  \item  \bambini{è vero! una ne ha uguale e una no!}
  \item  \bambini{anche il contorno di clessidra e lumaca è uguale}
  \item  \bambini{però valgono diversi!}
\end{itemize}
\end{studente}
Così, piano piano arriviamo alla conclusione che perimetro e area non
per forza hanno la stessa misura\dots{} PERCHÉ?
\begin{studente}[Bambini]%
  perché il valore ce lo danno il numero di tessere e il contorno non
  centra col numero di tessere
\end{studente}
Sono stra felice\dots{}

Entra in classe la collega che non mostra alcun interesse.

Schede in mano, suggerisco ai bambini di raccontare alla collega che è
presente per la prima volta di cosa si tratta, come hanno lavorato e,
invito Rebecca (una delle bambine con ottime capacità sia concettuali
che linguistiche) a dire la sua.
\begin{studente}[Rebecca]
  la maestra ci ha dato alcune tessere colorate e insieme (coi i
  compagni) le abbiamo usate per completare le schede e costruire
  figure
\end{studente}
La collega si rivolge a Alessia, bimba solare un po' infantile ma
dotata di grande impegno e volontà di fare, che stava ridacchiando con
una compagna e le domanda lei cosa avesse fatto. All'affermazione
entusiasta della bimba:
\begin{studente}[Alessia]
  ho GIOCATO con le tessere
\end{studente}
segue il richiamo saccente della collega\label{ex:alessia} che la incalza:
\begin{tutor}[Ins.]
  \dots{} ah sì Alessia? hai giocato o lavorato?
\end{tutor}
Alessia si intimorisce e ammutolisce e la collega continua
rivolgendosi a tutti
\begin{tutor}[Ins.]
  no, perché qui non si viene per giocare!
\end{tutor}
Intervengo allora spiegando che Alessia aveva ragione perché abbiamo
eseguito dei lavori giocando con alcune tessere.

Il clima della classe è ormai mutato\dots{}  gli alunni sono
intimoriti, i loro sguardi sono bassi e la voce dei bambini ai quali
la collega si rivolge chiedendo che lavori avessero effettuato è
tremolante.

Taglio corto con un ``ok dai, cominciamo con l'attività così la
maestra si fa un'idea di cosa facciamo'' e consegno loro tessere e
scheda \attivita{PER CONCLUDERE 3}\dots{}

Proprio non riescono a collaborare nonostante i miei inviti di
lavorare insieme, di alzarsi e spostarsi dalla sedia per avvicinarsi e
vedere e agire meglio\dots{}

Anche il mio umore è variato passando dalla rabbia (per la collega)
alla frustrazione per l'incapacità di ridare serenità ai bambini.

Dico loro di riconsegnare appena avranno terminato le scheda e, dopo
averle contate, le tessere.

Distraggo la collega proponendole un'attività per il lavoretto di
Natale e le propongo di cominciarla alla riconsegna delle
schede\dots{}

\subsection{Osservazioni finali }
Avrei diverse frasi da riportare e osservazioni in merito all'accaduto
ma ho deciso che la situazione illustrata parla da sola\dots{}

Prego di non dimenticare mai che ho davanti bambini ai quali devo
passare prima il gusto dell'imparare\dots{}  un programma da svolgere
perché deve essere svolto lascia il tempo che trova\dots{}

Mi domando se attività laboratoriali inserite in contesti rigidi forse
non disorientino i bambini\dots{}


\section[Sperimentazione \#2: seconda primaria]{Sperimentazione \#2:
  classe seconda primaria, gennaio/febbraio~2010}

\subsection{Osservazioni generali}

\subsubsection{Presentazione della classe} 21 alunni.

\subsubsection{Composizione dei gruppi}
La classe è stata divisa in 5 gruppi di lavoro, liberamente
costituiti, di 4 elementi ciascuno; i gruppi sono risultati eterogenei
tra loro e sostanzialmente omogenei al loro interno (gruppi di
livello). I gruppi sono rimasti invariati per tutto il percorso.

\subsubsection{Insegnanti presenti}

In alcuni incontri l'insegnante di classe è affiancata dall'insegnante
di sostegno.

\subsubsection{Calendarizzazione degli incontri}
\begin{calendario}
  \begin{itemize}
  \item 21 gennaio
  \item 28 gennaio (presenza insegnante sostegno)
  \item 4 febbraio (presenza insegnante sostegno)
  \item 11 febbraio (presenza insegnante sostegno)
  \end{itemize}
\end{calendario}

\subsection{Primo incontro}
\begin{description}
\item[Alunni presenti:]Sono presenti 20 alunni, un alunno è assente.
\item[Tempo effettivo di lavoro:]Il tempo di lavoro effettivo è stato di 2 ore.
\end{description}

\mbox{}
\begin{consegna}
  Ogni mattina spiego alla classe il programma di lavoro ipotizzato,
  quanto tempo abbiamo a disposizione e come lo useremo. Ho presentato
  il nuovo percorso di lavoro contestualizzandolo nella partecipazione
  a una sperimentazione con l'università, come aiuto che la classe
  darà allo studio e alla riflessione di adulti che vogliono imparare
  come insegnare meglio la matematica anche attraverso l'uso di un
  materiale particolare, il kit; scoprendo e imparando nuovi concetti
  e nuove informazioni aiuteranno gli adulti a imparare. Ho presentato
  in modo dettagliato le modalità di lavoro. Sono state proposte le
  esperienze con le tessere triangolari: \attivita{per cominciare} e
  \attivita{prima esperienza} seguendo le indicazioni riportate nel
  testo Torri, serpenti e \dots{} geometria. Le consegne sono state lette e
  comprese prima a livello collettivo, poi a livello di gruppo come
  verifica e rinforzo.

  \materiali{} %
  Ogni gruppo ha a disposizione sul proprio tavolo 16 tessere e le
  schede di lavoro. Gli alunni non hanno chiesto altro materiale.
\end{consegna}

\subsubsection{Osservazioni}
l ``lavorare in gruppo'' è stato il problema; è stato necessario il mio
intervento talvolta per rilanciare la richiesta di collaborare e di
tenere conto del lavoro di tutti e talvolta per agire sulle dinamiche
relazionali. I bambini si dimostravano interessati, desiderosi di fare
ma un eccesso di individualismo da parte di alcuni ha condizionato la
concentrazione sul lavoro. Tutti i gruppi hanno eseguito il compito
assegnato iniziando a manipolare e a conoscere il
materiale. Qualificante è stato il momento finale di confronto in
grande gruppo.

\subsubsection{Consigli per i colleghi che vogliono proporre le stesse
  attività }
La presenza di (almeno) due insegnanti.

\subsection{Secondo incontro}
\begin{description}
\item[Alunni presenti:]Tutti gli alunni sono presenti.
\item[Tempo effettivo di lavoro:]Il tempo di lavoro preventivato è di
  2 ore. Tempo effettivo tre ore
\end{description}
\mbox{}
\begin{consegna}
  Sono state proposte le esperienze con le tessere triangolari:
  \attivita{seconda e terza esperienza} seguendo le indicazioni
  riportate nel testo Torri, serpenti e \dots{} geometria. Le consegne
  sono lette attentamente e comprese in un momento di lavoro
  collettivo cioè che coinvolge tutto il gruppo classe nel suo
  insieme.

  \materiali{} %
  Ogni gruppo ha a disposizione sul proprio tavolo 16 tessere, le tre
  corde e le schede di lavoro. Gli alunni non hanno chiesto altro
  materiale.
\end{consegna}

\subsubsection{Osservazioni}
Il tempo di lavoro preventivato è di 2 ore ma risulta insufficiente:
il momento finale di confronto in grande gruppo è stato solo
avviato. Gli alunni accolgono l'attività con dichiarata soddisfazione
e motivazione; esplicitano il loro piacere di poter utilizzare il
materiale per costruire forme, anche se non più liberamente ma
finalizzato all'esecuzione di consegne precise. La presenza delle
corde li incuriosisce. Elemento di rilievo è sicuramente la presenza
di due insegnanti, presenza che sembra determinare un'atmosfera di
lavoro più distesa e meno dispersiva. I bambini contano sulla
possibilità di raccontare ciò che fanno, di chiedere conferme e/o
informazioni. La concentrazione e l'impegno appaiono più diffusi e più
adeguati. I livelli di elaborazione e consapevolezza dell'esperienza
effettuata appaiono molto diversificati; per alcuni alunni la
manipolazione, l'osservazione e il parlare con i compagni sollecita
intuizioni, scoperte e curiosità autonome; per altri invece la
manipolazione è quasi fine a se stessa e solo l'intervento guida
dell'insegnante porta a osservare, riconoscere e scoprire.

\subsubsection{Consigli per i colleghi che vogliono proporre le stesse
  attività}
L'esperienza con le corde dovrebbe essere proposta in un incontro successivo.

\subsection{Terzo incontro}
\begin{description}
\item[Alunni presenti:] Tutti gli alunni sono presenti.
\item[Tempo effettivo di lavoro:] Il tempo di lavoro preventivato è di
  2 ore
\end{description}

\begin{consegna}
  Inizialmente è stato dedicato tempo alla ripresa collettiva delle
  esperienze precedenti per ricomporle in un percorso unitario,
  offrire una ulteriore occasione per fare emergere chiarezze e
  incomprensioni e fondare l'attività odierna sul ``già fatto''. Sono
  state proposte le esperienze con le tessere triangolari:
  \attivita{quarta e quinta esperienza} seguendo le indicazioni
  riportate nel testo ``Torri, serpenti e\dots{} geometria''. Le
  consegne sono lette attentamente e comprese in un momento di lavoro
  collettivo, particolare attenzione è stata posta alla domanda: ``Ve
  lo aspettavate?''. Durante il lavoro tutti i gruppi hanno richiesto
  l'intervento dell'insegnante come conferma, chiarificazione o
  ulteriore spiegazione di alcuni termini o del significato globale
  della consegna, in particolare dei quesiti (ciò non garantirà una
  reale comprensione e una corretta risposta ai quesiti)

  \materiali{}%
  Ogni gruppo ha a disposizione sul proprio tavolo 16 tessere, le tre
  corde e le schede di lavoro 4 e 5. Ho offerto la scheda 3 da
  consultare per costruire la piramide.
\end{consegna}

\subsubsection{Osservazioni}
Il tempo di lavoro preventivato è di 2 ore ma risulta insufficiente:
il momento finale di confronto in grande gruppo risulta sempre molto
qualificante e meriterebbe più tempo. I gruppi rimangono invariati e
questo sembra consentire una maggiore collaborazione o una maggiore
disponibilità alla collaborazione; tutti hanno portato il proprio
contributo anche se non sempre adeguato o ``pensato''. La presenza di
due insegnanti permette una gestione più puntuale dei conflitti in
modo che non condizionino troppo l'andamento del lavoro. Un gruppo di
bambini, solitamente brillanti, ha raggiunto una buona capacità di
operare discutendo e confrontandosi in modo pertinente e
funzionale. Realizzare delle figure con le tessere riscuote sempre
successo, i bambini si divertono durante e dopo il lavoro osservando e
mostrando con soddisfazione le proprie opere; in particolare oggi,
hanno creato nuove figure spostando solo quattro tessere! Tuttavia
l'aspetto ludico rimane per alcuni alunni troppo preminente. Il
disegno della nuova figura è difficoltoso per due gruppi di alunni,
solitamente più insicuri; dopo alcuni tentativi falliti ricorrono al
colore delle tessere per individuarne la posizione e riprodurre
correttamente la figura; questi gruppi non riusciranno a realizzare la
quinta esperienza. Tutti i gruppi non hanno subito ``sparato'' delle
risposte e hanno cercato delle soluzioni utilizzando il materiale;
solo uno tuttavia ha tenuto conto del lavoro appena eseguito (abbiamo
solo cambiato posizione a quattro tessere, non le abbiamo tolte, non
le abbiamo aggiunte quindi \dots{}), altri hanno utilizzato strategie
scorrette fondate su ``immagini'' scorrette: la corda delimita lo spazio
occupato dalla figura quindi la sua lunghezza mi dice quanto spazio
occupa oppure indico con tre punti i vertici della piramide e immagino
lo spazio che occupa, vi inserisco la nuova figura e\dots{}, altri si
fermano al dato percettivo superficiale: mi sembra che\dots{} Nel
momento di discussione finale il confronto tra le diverse risposte ha
suscitato l'interesse e la partecipazione della classe; in particolare
il gruppo che ha fornito la risposta esatta ha spiegato il proprio
``ragionamento'' con chiarezza e determinazione, sollecitando interventi
e convincendo i compagni.

\subsubsection{Consigli per i colleghi che vogliono proporre le stesse
  attività}
Non è un consiglio ma una proposta di riflessione: dare tempi di
lavoro più lunghi e la presenza di più insegnanti per consentire una
maggiore elaborazione dell'esperienza.

\subsection{Quarto incontro}
\begin{description}
\item[Alunni presenti:] Tutti gli alunni sono presenti.
\item[Tempo effettivo di lavoro:] Il tempo di lavoro effettivo è di 2
  ore.
\end{description}

\begin{consegna}
  Presento l'esperienza odierna come momento conclusivo della
  sperimentazione in classe che non prevede l'uso del materiale:
  occorre ripensare alle esperienze fatte e utilizzarle solo
  ``immaginando'' le tessere per controllare così ciò che ricordiamo e
  abbiamo incominciato a capire. Gli alunni scelgono di leggere e
  comprendere le consegne in gruppo, da soli. Sollecito il
  controllo/verifica delle risposte date, attraverso strategie che
  dovranno poi essere in grado di spiegare ai compagni per dimostrare
  la correttezza del proprio lavoro.

  \materiali{} %
  Scheda di lavoro: \attivita{Per concludere}.
\end{consegna}

\subsubsection{Osservazioni}
La prima parte del lavoro viene eseguita facilmente da tutti i gruppi;
alcuni mostrano perplessità circa la necessità di precisare
``\dots{} tessere rosse\dots{}'' (a cosa serve indicare il colore?). Per
contare le tessere che compongono l'ochetta e il cane vengono usate
diverse modalità: guardare e immaginare le tessere, ripassare e
completare il contorno delle tessere con il dito, tracciare i lati
mancanti delle tessere con la matita. Difficoltà emergono invece nella
seconda parte del lavoro già a partire dalla comprensione della
consegna. Quattro gruppi su cinque chiedono l'intervento delle
insegnanti per chiarire dubbi, per verificare la comprensione o per
essere guidati nella comprensione, gli alunni hanno incertezze a
individuare i lavori che devono eseguire e a isolare la domanda a cui
rispondere. L'ostacolo sembra essere riconducibile alla formulazione
linguistica. Un solo gruppo non pone richieste di spiegazione,
evidenzia la parola ``spostando'' dalla quale deducono la risposta
esatta, solo successivamente lavoreranno con la figura come strategia
di verifica. I due gruppi più insicuri faticano a lavorare senza poter
manipolare il materiale, chiedono se non sia proprio possibile averlo
a disposizione e ``lo sostituiscono'' con la richiesta della presenza
della maestra. Due gruppi giungono alla risposta corretta dopo aver
lavorato sulla figura. Il momento di discussione finale è stato
partecipato, in questo incontro tuttavia ho rilevato una minore tenuta
dell'ascolto da parte dei bambini più insicuri (avevano già speso
molte energie in un'attività dal livello di astrazione maggiore delle
precedenti). Buona è per alcuni la capacità di argomentare il proprio
operato e le conclusioni che esplicitano: \bambini{cambia la forma,
  non il numero delle tessere che ho usato}, \bambini{leggi la parola
  spostano, ti dice che non togli e non aggiungi\dots{}}

\subsection{Riflessioni di fine percorso}
Con la collega che ha condotto la Sperimentazione \#3 (vedi
sperimentazione successiva) abbiamo scelto di mettere a confronto le
nostre esperienze avendo entrambe lavorato sulle classi seconde e
pervenire così a un documento comune:
\begin{itemize}
\item per l'efficacia della manipolazione (lavorare riflettendo,
  mentre costruisco, mentre disegno\dots{} capisco) è necessaria la
  presenza di più animatori; il manipolare, anche con un materiale
  strutturato, non garantisce da solo né la comprensione né la
  riflessione su ciò che si sta facendo.
\item Agire il doppio ruolo di animatore e osservatore è stato
  parecchio complesso.
\end{itemize}

Abbiamo operato nella costituzione dei gruppi scelte diverse:
\begin{itemize}
\item nei gruppi eterogenei i bambini in difficoltà non sono riusciti
  a trarre vantaggio dalla presenza di compagni ``bravi'', non erano
  pronti a cogliere ciò che veniva loro riferito, lo accettavano ma
  non era una loro conquista: tempi e capacità diverse;
\item i gruppi omogenei si rivelano più funzionali ma necessitano
  della presenza dell'insegnante come supporto a disposizione dei
  gruppi meno pronti.
\end{itemize}
I tempi di esecuzione indicati dal testo non si sono mai dimostrati
adeguati. Abbiamo sempre avuto l'impressione di aver bisogno di più
tempo, ci chiediamo se questo dipende forse da un approccio troppo
``scolastico'' e occorre dare a queste attività più ``leggerezza''.


\section[Sperimentazione \#3: seconda primaria]{Sperimentazione \#3:
  classe seconda primaria, gennaio/febbraio~2010}

\subsection{Osservazioni generali}
L'attività si svolge nelle mie due seconde, ma la documentazione si
riferisce solo a una delle due classi.

\subsubsection{Presentazione della classe}

La classe è composta da 21 alunni, è presente un alunno segnalato.

\subsubsection{Composizione dei gruppi}

Gli alunni sono stati da me divisi in cinque gruppi eterogenei e
ugualmente numerosi. Di solito lavorano a coppie o in gruppi di tre
bambini ciascuno per libera aggregazione o su mia indicazione.

\subsubsection{Insegnanti presenti}

Agli incontri è presente solo l'insegnante di classe.

\subsubsection{Calendarizzazione degli incontri}
\begin{calendario}
  \begin{itemize}
  \item 22 gennaio
  \item 25 gennaio
  \item 1 febbraio
  \item 5 febbraio
  \end{itemize}
\end{calendario}

\subsection{Primo incontro}

\begin{description}
\item[Alunni presenti:] presenti 20 alunni su 21 totali (1 assente:
  l'alunno segnalato)
\item[Tempo effettivo di lavoro:] L'attività è durata 55 minuti.
\end{description}

\begin{consegna}
  Ho avviato l'attività di laboratorio con \attivita{Per cominciare} e
  la prima esperienza in entrambe le classi seconde in cui lavoro come
  docente di matematica in modo da offrire ai miei alunni le medesime
  opportunità. Ogni consegna e/o domanda è stata letta insieme.

  \materiali{}%
  Ho predisposto in aula per ciascun gruppo le ``postazioni di
  lavoro'', il materiale e le schede riguardanti le esperienze
  proposte secondo le indicazioni date.
\end{consegna}

\subsubsection{Osservazioni}
Ho gestito da sola l'attività di laboratorio. Gli alunni hanno subito
dimostrato interesse e curiosità, voglia di fare ma sono stati
distratti e rallentati nel lavoro dall'assegnazione dei ruoli
all'interno del gruppo e, per la prima esperienza, dal materiale:
forte la tendenza a accaparrarsi le tessere sul tavolo per giocare o
per costruire individualmente le figure richieste senza tener conto
delle proteste dei compagni e dei loro suggerimenti.

\attivita{Per cominciare}: gli alunni hanno facilmente individuato il
numero delle tessere (\bambini{basta contare}, \bambini{ci sono le
  divisioni}, \bambini{è come contare i quadrati unità delle torri
  decina}) non così per i mattoni. La rappresentazione a tre
dimensioni li ha disorientati: \bambini{bisogna immaginare},
\bambini{non capisco il disegno}, \bambini{ma ci sono mattoni più
  grandi e mattoni più piccoli?}, \bambini{perché ci sono quattro
  mattoni e sotto è vuoto?}. Un solo gruppo ha risposto correttamente:
un alunno ha usato il righello accostandolo allo spigolo più lungo
della torre ma valutando a occhio. Ci sono stati dei tentativi di
disegnare i mattoni all'interno della torre.

La prima esperienza ha creato un clima di gioco e di competizione tra
i gruppi non solo nel costruire le figure ma anche nel dare loro dei
nomi: \bambini{pesce a bocca aperta}, \bambini{corona da re}. Nella
discussione finale gli alunni hanno spiegato le strategie usate e
manifestato le loro difficoltà; è stato per loro difficile il
confronto tra le diverse risposte, quelle corrette e quelle
sbagliate. Tendevano a rivolgersi direttamente all'insegnante cercando
una risposta o una conferma.

\subsection{Secondo incontro}
\begin{description}
\item[Alunni presenti:] presenti 18 alunni su 21 totali. Ci sono tre
  assenti e, tra questi, l'alunno segnalato.
\item[Tempo effettivo di lavoro:]
\par L'attività è durata due ore
\end{description}
\begin{consegna}
  Agli alunni sono state proposte la seconda e la terza esperienza con
  le tessere triangolari accompagnate, a livello collettivo, da una
  lettura attenta e puntuale delle richieste. È stato loro comunicato
  che l'attività durerà due ore.

  \materiali{}%
  A ciascun gruppo sono state date le 16 tessere triangolari, le tre
  corde e le schede relative alla seconda e terza attività secondo le
  indicazioni date. L'alunna, che mi aiuta a contare le tessere da
  distribuire a ciascun gruppo, mi suggerisce di formare gruppi con
  tessere di diverso colore per essere più veloce perché il
  \bambini{colore delle tessere non serve per il lavoro}.
\end{consegna}

\subsubsection{Osservazioni}
Sono contenti di riprendere l'attività di laboratorio, si ricordano
vicendevolmente i ruoli che si sono dati (all'interno di un gruppo
vengono ridistribuiti in modo più funzionale al lavoro) e le modalità
di utilizzo delle tessere, sono piacevolmente incuriositi dalle forme
da riprodurre. Durante il lavoro i gruppi richiedono spesso
l'intervento dell'insegnante per superare momenti di difficoltà o di
stallo, per comunicare scoperte e osservazioni, per avere conferme nel
procedere. Talvolta necessitano di essere guidati nell'osservazione
delle figure da riprodurre, di essere sostenuti nei ripetuti tentativi
o nelle difficoltà pratiche (nel fare, ad esempio, aderire la corda al
bordo della stella) in modo da mantenere l'attenzione e la
motivazione, nel rileggere le consegne per soddisfare correttamente le
richieste (la parola ``circondare'' per alcuni bambini è così
``pregnante'' da far scomparire ``in modo che aderiscano bene al
bordo''). All'interno dei gruppi i bambini più intuitivi a volte
faticano a esplicitare in modo chiaro e a condividere con i compagni
le loro scoperte (non so che parole usare); altri bambini sembrano
accettare le soluzioni dei compagni che riconoscono come bravi. La
presenza di un solo insegnante è penalizzante. Occorre ``ributtare''
all'interno del gruppo le domande degli alunni ponendo domande
opportune e calibrate in modo da sollecitare le loro riflessioni e e
di far vivere come conquiste personali le loro scoperte. Ciò richiede
tempo e attenzione e crea situazioni di stallo che non aiutano a
mantenere alta l'attenzione. Un gruppo non ha completato la terza
attività; Non c'è stato spazio per la discussione nel gruppo classe:
sarà il punto da cui partire nella prossima esperienza.

\subsection{Terzo incontro}

\begin{description}
\item[Alunni presenti:] Sono presenti tutti gli alunni
\item[Tempo effettivo di lavoro:]2 ore
\end{description}

\begin{consegna}
  Parte del tempo a disposizione è stata utilizzata per riprendere a
  livello collettivo le esperienze precedenti con il supporto delle
  schede sulle quali gli alunni avevano lavorato. È stato utile in
  particolare per gli alunni assenti il 25 gennaio perché hanno avuto
  la possibilità di avere un percorso lineare e di lavorare
  serenamente. Nella lettura delle consegne per alcuni alunni è stato
  ancora necessario puntualizzare la differenza tra spazio occupato
  sul tavolo e bordo da circondare e ricordare le strategie utilizzate
  per misurare le caratteristiche prese in considerazione di volta in
  volta. La domanda ``ve lo aspettavate?'' ha suscitato in alcuni
  bambini delle perplessità perché non avevano chiaro se dovevano
  prima misurare e poi rispondere, perché non era esplicitamente
  richiesto come nelle esperienze precedenti, o rispondere senza
  misurare.

  \materiali{}%
  A ogni gruppo sono state consegnate le 16 tessere triangolari, le
  tre corde e le schede relative alla quarta e alla quinta
  esperienza. Alla fine di ogni attività lascio a ogni gruppo le
  schede di lavoro per renderli più consapevoli e partecipi del
  percorso intrapreso e perché ho osservato che rivedono e mostrano
  agli altri gruppi con piacere e soddisfazione le figure che hanno
  costruito o che le utilizzano per motivare o rinforzare le loro
  scoperte. I gruppi sono invariati; l'alunno segnalato è stato
  inserito in un gruppo da me individuato per la presenza di bambini
  disponibili a ascoltare e a aiutare.
\end{consegna}

\subsubsection{Osservazioni}
Gli alunni costruiscono facilmente più figure, si divertono a
riconoscere nelle forme oggetti, animali ecc. e a dare loro dei nomi,
discutono su quale sia la ``più bella'' da disegnare. L'aspetto ludico
dell'attività fa perdere di vista a alcuni gruppi la richiesta ``solo
cambiando posizione a quattro tessere della piramide'' e li devo
richiamare a rileggere attentamente la consegna. Alla prima domanda i
gruppi rispondono correttamente e velocemente, diverse sono le
strategie usate: due gruppi sentono la necessità di contare le tessere
della piramide o della stella (o costruendola o utilizzando la
rappresentazione della \attivita{scheda 3}) e della nuova figura; un
gruppo risponde correttamente ma poi conta le tessere ``per essere più
sicuri''; due gruppi rispondono subito con sicurezza senza contare le
tessere. Il lavoro con le corde richiede molto tempo; ha per i bambini
delle difficoltà pratiche e li devo sostenere nel mantenere
l'attenzione e l'interesse. Un gruppo richiede il mio intervento
perché un bambino manifesta la sua difficoltà a lavorare in gruppo
distruggendo le figure sul banco. Per facilitare e velocizzare il
lavoro costruisco la stella e la piramide con le tessere del kit non
utilizzate che metto a disposizione dei gruppi perché hanno bisogno di
misurare più volte i bordi e di confrontare i risultati delle
misurazioni. Nonostante l'esperienza precedente (\attivita{Scheda 3})
sono sorpresi dai risultati delle loro misurazioni e dalle scoperte
fatte:
\begin{studente}[ ]
  il bordo non si comporta come lo spazio
\end{studente}
Un alunno è molto interessato allo strumento di misura: mi dice che la
corda non va bene per misurare i bordi, prende la squadra (non è
richiesta ai bambini come materiale, ma viene usata da alcuni bambini
come un righello per costruire tabelle, per tracciare percorsi ecc.)
e la accosta ai lati della piramide:
\begin{studente}[ ]
  non si muove, sta attaccata, ha una salita e una discesa come la
  piramide ma non va bene per la stella
\end{studente}
L'alunno segnalato viene coinvolto dai compagni nel lavoro ma ne
coglie solo l'aspetto manipolativo e ludico. Ai compagni che gli
spiegano pazientemente più volte ciò che stanno facendo risponde
sorridendo con un \bambini{mi arrendo}. Tuttavia rimane nel
gruppo. Osserva e ascolta i compagni senza disturbarli
nell'attività. La presenza dell'insegnante di sostegno gli avrebbe
permesso di fare un percorso a sua misura. Due gruppi non riescono a
completare il lavoro nel tempo prestabilito. Rimando alla fase
\attivita{Per concludere} la discussione nel grande gruppo.

\subsubsection{Consigli per i colleghi che vogliono proporre le stesse
  attività}
Solo alcune riflessioni: l'opportunità della presenza di più animatori
e dell'insegnante di sostegno se vi sono alunni segnalati; una
maggiore efficacia di gruppi omogenei al loro interno in modo da
permettere a ognuno di procedere con il proprio passo senza
rallentamenti o senza anticipazioni inopportune; più tempo per la fase
di manipolazione del materiale.

\subsection{Quarto incontro}

\begin{description}
\item[Alunni presenti:] Sono presenti tutti gli alunni.
\item[Tempo effettivo di lavoro:]1 ora
\end{description}

\begin{consegna}
  Comunico che oggi si conclude la sperimentazione in classe e che si
  richiede loro uno ``sforzo'' maggiore: lavoreranno senza le tessere
  e ciò suscita rammarico perché amano l'aspetto ludico e manipolativo
  dell'attività. Nella lettura collettiva delle consegne mi soffermo
  soprattutto sulla seconda dalla formulazione linguistica complessa:
  la riformulo slivellando le richieste.

  \materiali{}%
  A ciascun gruppo consegno la scheda \attivita{Per concludere}. I
  gruppi sono invariati nella composizione.
\end{consegna}

\subsubsection{Osservazioni}
I gruppi rispondono correttamente e agevolmente alla prima domanda o
tracciando con le dita o disegnando con la matita le tessere
all'interno delle figure. Mi accorgo che alcuni bambini, quando mi
avvicino, cancellano le divisioni tracciate. Li rassicuro dicendo loro
che ogni gruppo può scegliere la strategia che ritiene più
opportuna. Per l'ochetta molti scrivono subito la risposta corretta
ricordando la scheda \attivita{Per cominciare}.

Per la seconda parte della scheda alcuni gruppi richiedono ancora il
mio aiuto per decodificare la consegna per poi procedere speditamente
e correttamente. Un gruppo ricorre ancora alla rappresentazione:
disegna due tessere e copre con le dita le due tessere che immaginano
di spostare ma, ascoltando le loro osservazioni, sembra più un bisogno
di verificare per ``rispondere giusto''. Al ``potete dirlo senza
contare'' due gruppi rispondono ``no'' ma, chiamati a giustificare il
no, scrivono che \bambini{occupano tutte e due lo stesso spazio hanno
  lo stesso numero di tessere} e \bambini{non aumentano o diminuiscono
  cambiano solo il posto}. Nella discussione nel grande gruppo i
bambini sottolineano soprattutto il piacere di lavorare, di capire e
di scoprire insieme: \bambini{i lavori erano alcuni facili alcuni
  difficili ma divertenti}, \bambini{non ci arrendevamo mai
  continuavamo a lottare}, \bambini{mi è piaciuto disegnare e
  capire}. Alcuni bambini invece sottolineano con soddisfazione ciò che
hanno imparato dimostrando una buona conoscenza dei concetti trattati
e consapevolezza del percorso fatto: \bambini{abbiamo lavorato con lo
  spazio e il bordo delle figure}, (riferendosi alla forma)
\bambini{lo spazio è suo e non può essere occupato}, \bambini{il suo
  spazio finisce quando arrivi ai bordi}, \bambini{puoi anche non
  usare la corda, basta contare i bordi} (riferendosi ai lati di
ciascuna tessera), ecc.

\subsection{Riflessioni di fine percorso}
Con la collega che ha condotto la Sperimentazione \#2 (vedi
sperimentazione precedente) abbiamo scelto di mettere a confronto le
nostre esperienze avendo entrambe lavorato sulle classi seconde e
pervenire così a un documento comune:
\begin{itemize}
\item per l'efficacia della manipolazione (lavorare riflettendo,
  mentre costruisco, mentre disegno\dots{} capisco) è necessaria la
  presenza di più animatori; il manipolare, anche con un materiale
  strutturato, non garantisce da solo né la comprensione né la
  riflessione su ciò che si sta facendo. Agire il doppio ruolo di
  animatore e osservatore è stato parecchio complesso.
\item I tempi di esecuzione indicati dal testo non si sono mai
  dimostrati adeguati. Abbiamo sempre avuto l'impressione di aver
  bisogno di più tempo, ci chiediamo se questo dipende forse da un
  approccio troppo ``scolastico'' e occorre dare a queste attività più
  ``leggerezza''
\end{itemize}

Abbiamo operato nella costituzione dei gruppi scelte diverse:
\begin{itemize}
\item nei gruppi eterogenei i bambini in difficoltà non sono riusciti
  a trarre vantaggio dalla presenza di compagni ``bravi'', non erano
  pronti a cogliere ciò che veniva loro riferito, lo accettavano ma
  non era una loro conquista: tempi e capacità diverse;
\item i gruppi omogenei si rivelano più funzionali ma necessitano
  della presenza dell'insegnante come supporto a disposizione dei
  gruppi meno pronti.
\end{itemize}


\section[Sperimentazione \#4: terza primaria]{Sperimentazione \#4:
  classe terza primaria, gennaio/febbraio~2010}

\subsection{Osservazioni generali}

\subsubsection{Presentazione della classe}
La sperimentazione avviene in due classi terze.

Nella prima (``A''): 19 alunni, di cui 1 segnalato per disturbi di
apprendimento e 1 con problemi di apprendimento legati a un forte
disagio psico-affettivo.

Nella seconda (``B''): 19 alunni, di cui 1 con gravi difficoltà di
relazione con i compagni e 1 con problemi di apprendimento in via di
accertamento.

I bambini di entrambe le classi lavorano spesso in piccoli gruppi.

\subsubsection{Composizione dei gruppi}
Abbiamo chiesto ai bambini di formare liberamente i gruppi, con
l'unica condizione che all'interno di ognuno ci fossero maschi e
femmine. In ogni classe si sono formati 4 gruppi da 4 bambini e 1 da
3, tutti abbastanza eterogenei al loro interno. I gruppi sono gli
stessi per le 4 esperienze. Ognuna delle insegnanti conduce l'attività
nella propria classe, senza compresenza.

\subsubsection{Insegnanti presenti}
In ognuna delle due classi è presente solo l'insegnante di classe.

\subsubsection{Calendarizzazione degli incontri}
\begin{calendario}
  \begin{itemize}
  \item 15 gennaio
  \item 22 gennaio
  \item 27 gennaio
  \item 3 febbraio
  \end{itemize}
\end{calendario}

\subsection{Primo incontro}
\begin{description}\item[Alunni presenti:] Tutti (19 in ogni classe)
\item[Tempo effettivo di lavoro:] meno di 2 ore
\end{description}

\begin{consegna}
  Consegna delle schede \attivita{Per cominciare} di classe terza,
  prima scheda \attivita{Tessere} di classe seconda, prima scheda
  \attivita{Tessere} di classe terza.

  Abbiamo comunicato ai bambini che avrebbero cominciato un
  laboratorio di geometria utilizzando del nuovo materiale e che
  avrebbero ``fatto delle scoperte''. Abbiamo letto insieme i comandi,
  senza però dare particolari indicazioni.

  \materiali{} %
  12 tessere triangolari per ogni gruppo. Abbiamo specificato che il
  numero delle tessere era quello necessario per svolgere le consegne.
\end{consegna}

\subsubsection{Osservazioni}
Fase 1: Abbiamo chiesto ai bambini di formare liberamente i gruppi,
imponendo l'unica condizione che all'interno di ognuno ci fossero
maschi e femmine. In ogni classe si sono formati 4 gruppi da 4 bambini
e 1 da 3, tutti abbastanza eterogenei al loro interno. Per questa
operazione ci sono voluti pochi minuti. Abbiamo chiesto poi di
scegliere un nome per il gruppo e poi dividersi i ruoli di:
\begin{itemize}
  \item leggere i quesiti e scrivere le risposte concordate
  \item disegnare/rappresentare le figure
  \item porre eventuali domande all'insegnante
  \item riferire alla classe nel momento collettivo.
\end{itemize}
Queste operazioni si sono svolte in una decina di minuti senza
particolari problemi, anche perché i bambini lavorano spesso in
piccolo gruppo.

Fase 2: Consegna delle schede \attivita{Per cominciare} di classe
terza. Tempo di esecuzione 10 minuti. I bambini hanno partecipato
tutti attivamente.

Fase 3: Consegna delle 12 tessere triangolari per ogni gruppo. Abbiamo
invitato i bambini a costruire diverse forme e a dare loro un nome: in
questa fase di manipolazione abbiamo notato partecipazione da parte di
tutti. Abbiamo consegnato la prima scheda \attivita{Tessere} di classe
seconda, sulla quale i bambini dovevano rappresentare graficamente due
delle forme realizzate. In questa fase dell'esperienza alcuni (tra
quelli più ``deboli'') hanno manifestato calo di interesse (chi si
avvicinava a compagni di altri gruppi, chi giocherellava con il
materiale, chi appoggiava la testa sul banco\dots{}). La durata di
questa fase è stata di circa 30 minuti.

Fase 4: Consegna della prima scheda \attivita{Tessere} di classe
terza: in questa fase i bambini dovevano costruire due figure indicate
e confrontarle sia dal punto di vista dello spazio occupato, sia da
quello del contorno (dovevano usare delle corde predisposte). Durante
questa attività (durata circa 20 minuti) gli alunni ``distratti'' hanno
ritrovato motivazione e interesse.

Fase 5: Confronto tra i lavori di gruppo. Ogni gruppo riferisce le
risposte date ai quesiti. I bambini si rendono conto di aver dato
risposte diverse. Si mettono in discussione, si pongono domande:
\bambini{abbiamo contato male}, \bambini{abbiamo contato solo il
  contorno (le tessere del contorno) e non all'interno},\dots{} La
verifica delle risposte viene rimandata a una fase successiva.

\subsection{Secondo incontro}
\begin{description}\item[Alunni presenti:]
\par Tutti (19 alunni in ogni classe)
\item[Tempo effettivo di lavoro:]
\par Circa 2 ore
\end{description}

\begin{consegna}
  seconda scheda \attivita{Tessere} di classe terza

  terza scheda \attivita{Tessere} (punto 4)

  quarta scheda \attivita{Tessere} (punto 5)

  Abbiamo letto insieme le consegne, senza dare ulteriori indicazioni.

  \materiali{} %
  Abbiamo consegnato ai gruppi le 12 tessere triangolari.
\end{consegna}

\subsubsection{Osservazioni}
Fase 1: Abbiamo consegnato ai gruppi le 12 tessere triangolari e la
seconda scheda \attivita{Tessere} di classe terza. In entrambe le
classi tutti i bambini hanno cominciato a lavorare con interesse e
partecipazione. All'interno dei gruppi sono nate discussioni intorno
alle domande ``quale figura occupa più o meno spazio'', ``quale ha il
contorno più o meno lungo''. In ogni gruppo, in un caso anche con
l'intervento dell'insegnante, si è giunti a un accordo sul
significato da attribuire a ogni termine geometrico (in classe non
abbiamo affrontato con sistematicità l'argomento). I tempi di
completamento dell'attività cominciano a differenziarsi
significativamente tra i vari gruppi. Gli alunni ``più deboli'' hanno
ancora difficoltà a differenziare ``lo spazio'' occupato dalle forme dal
loro ``contorno''. Tempo destinato a questa fase: 30 minuti.

Fase 2: Consegna della terza scheda \attivita{Tessere} (punto 4). I
vari gruppi hanno avuto difficoltà a soddisfare la richiesta
(costruire una figura che abbia il contorno che misura tanto quanto
quello del serpente). In una classe un gruppo, dopo aver provato
ripetutamente, ha chiesto aiuto all'insegnante; un altro gruppo
nell'altra classe ha rinunciato dopo pochi tentativi. Tra due gruppi
sono nati conflitti e discussioni (\bambini{avete copiato la nostra
  forma}\dots{} \bambini{non è bella}\dots{}). L'attività si protrae
per 30 minuti circa. Notiamo che nel momento della rappresentazione i
bambini più ``fragili'' si distraggono.

Fase 3: Consegna della quarta scheda \attivita{Tessere} (punto
5). Alcuni gruppi cominciano a provare e riprovare senza concordare
strategie. Nei gruppi dove ci sono gli alunni ``più alti'', scoprono che
per avere il contorno più corto devono compattare le tessere (metterle
il più vicino possibile), per il contorno più lungo invece devono
\bambini{fare la figura lunga}, \bambini{con più insenature}. 30/35
minuti il tempo di esecuzione per questa fase. Ancora momenti di
stanchezza da parte di alcuni alunni nel momento della
rappresentazione.

Fase 4: Condivisione del lavoro di gruppo (10/15 minuti). Si
confrontano le costruzioni dell'ultima scheda. I bambini riferiscono
le misure delle figure rappresentate e alcuni dei criteri seguiti. Dal
confronto qualcuno critica la propria costruzione e chiede di poterla
modificare per fare meglio. Qualcuno sottolinea di aver ``imparato'' a
misurare senza usare il righello, usando il lato delle tessere.

\subsection{Terzo incontro}

\begin{description}
\item[Alunni presenti:]tutti presenti
\item[Tempo effettivo di lavoro:] 1 ora e 30 minuti
\end{description}

\begin{consegna}
  le prime 2 schede \attivita{Cubi}. (Gli alunni dovevano costruire a
  loro scelta 3 figure e rappresentarne almeno una) Prima di avviare
  il lavoro abbiamo letto insieme la consegna.

  \materiali{}%
  cubi in legno (8 cubi di legno per ogni gruppo). I bambini non hanno
  chiesto altro materiale.
\end{consegna}

\subsubsection{Osservazioni}
Fase 1: Consegna di 8 cubi in legno per ogni gruppo e le prime 2
schede \attivita{Cubi}. (Gli alunni dovevano costruire a loro scelta 3
figure e rappresentarne almeno una). Prima di avviare il lavoro
abbiamo letto insieme la consegna.

Durante il lavoro di costruzione, in entrambe le classi si è notato
entusiasmo e interesse, inizialmente anche nei bambini più
``deboli''. In un gruppo nella classe ``B'' sono nate discussioni per
accordarsi su quale figura costruire (\bambini{Facciamo un ponte, no
  una pistola}). Nella classe ``A'' due gruppi si sono confrontati e
scontrati sulle figure costruite (\bambini{la nostra è più bella!},
\bambini{non assomiglia a \dots{}}). Tutti i gruppi hanno provato a
rappresentare le figure tridimensionali costruite. Qualcuno ha deciso
di rappresentare in due dimensioni. Qualcuno ha chiesto aiuto
all'insegnante. Qualcuno voleva rinunciare, ma con aiuto ha
rappresentato almeno una costruzione. Nella fase rappresentativa i
bambini problematici si sono distratti: chi giocava con i cubi, chi si
distraeva con oggetti personali\dots{}

Fase 2: consegna della terza e quarta scheda \attivita{Cubi}. Abbiamo
invitato gli alunni a contare le facce delle loro costruzioni. Tutti
hanno incontrato delle difficoltà a contare, perché non riuscivano a
trovare delle strategie. Dopo circa 10 minuti di tentativi alcuni
alunni hanno proposto \bambini{di contare le facce di `una parte' e
  poi moltiplicare per il numero delle parti uguali}; \bambini{\dots{}
  ma poi dobbiamo aggiungere le facce sotto e sopra}. Qualcuno fa
notare che non bisogna \bambini{contare le parti che combaciano}.
\begin{studente}[ ]
  Proviamo a contare una parte, poi quella uguale e poi contiamo le
  altre che restano fuori
\end{studente}
(all'esterno).

Fase 3: costruzione del dado. Non abbiamo notato particolari
difficoltà, forse perché avevano già lavorato lo scorso anno con
alcuni solidi. Poche difficoltà anche nella rappresentazione e nel
conteggio delle facce.

Fase 4: Costruzione della figura con 24 facce. Inizialmente hanno
avuto tutti difficoltà, anche se la motivazione era ancora alta (a
eccezione dei bambini più ``deboli'' che davano segni di stanchezza e di
insofferenza). Nonostante la richiesta delle insegnanti di utilizzare
il numero di cubi che volevano, tutti hanno provato con gli 8 cubi in
loro possesso. Un gruppo ha pensato poi di usarne meno, riuscendo così
a soddisfare la richiesta.

Fase 5: Nel confronto collettivo sono emerse le difficoltà soprattutto
nel contare le facce e nella rappresentazione. Alcuni hanno preso
consapevolezza degli errori.

\subsection{Quarto incontro}

\begin{description}
\item[Alunni presenti:]tutti presenti
\item[Tempo effettivo di lavoro:] 1 ora e 30 minuti
\end{description}

\begin{consegna}
  conversazione sul percorso precedente, schede \attivita{Per concludere},
  verifica collettiva sull'esperienza del laboratorio

  \materiali{}%
  12 tessere a ogni gruppo.
\end{consegna}

\subsubsection{Osservazioni}
Fase 1: Conversazione sul percorso fatto nei 3 precedenti
incontri. Riportiamo alcune affermazioni dei bambini:
\begin{studente}[Bambini]
  \begin{itemize}
  \item \bambini{Abbiamo usato delle tessere per fare forme e
      confrontarle}
  \item \bambini{Abbiamo costruito con i cubi}
  \item \bambini{Abbiamo contato le facce}
  \item \bambini{Abbiamo visto che le forme con i cubi occupano uno
      spazio `alto' e non `piano' come le forme costruite con le
      tessere}
  \end{itemize}
\end{studente}
Noi siamo intervenute nella conversazione solo come
moderatori. Durata: 10/15 minuti.

Fase 2: Distribuzione delle schede \attivita{Per concludere}. Abbiamo
invitato gli alunni a leggere attentamente le richieste, a
confrontarsi e a completare. I bambini ricordano e assumono i ruoli
che si erano assegnati negli incontri precedenti (nella classe ``B''
un'alunna non ha voluto assumere il ruolo concordato, è nata una
discussione, terminata con l'intervento dell'insegnante). Inizialmente
l'interesse all'interno dei vari gruppi è alto. I gruppi con gli
alunni ``più alti''
hanno memoria delle esperienze dei laboratori precedenti. I bambini
``più deboli'' non partecipano, si adeguano e qualcuno si distrae. Le
difficoltà maggiori sono emerse nel completamento della scheda
\attivita{Cubi}. Alcuni hanno provato a cercare delle strategie per
``contare i piani delle torri'' usando il righello oppure provando a
disegnare i cubi su ogni piano. 2 gruppi hanno contato i cubi del
piano in alto e quelli appoggiati al piano. Durata: 15/20 minuti.

Fase 3: Abbiamo pensato a questo punto di fare una verifica collettiva
di quanto scritto nelle schede \attivita{Per concludere}. Abbiamo
consegnato a ogni gruppo 12 tessere e abbiamo chiesto di provare a
rispondere nuovamente usando il materiale e di fare un confronto su
quanto scritto precedentemente. Ogni gruppo ha poi letto alla classe
le sue risposte. Qualcuno ha trovato gli errori, altri hanno
confermato. Abbiamo notato un riaccendersi dell'interesse alla
consegna del materiale. Abbiamo fatto collettivamente la verifica
delle costruzioni con i cubi (non erano sufficienti per tutti i
gruppi). È la fase in cui sono emerse maggiori difficoltà. Solo 2
gruppi per classe sono riusciti a contare esattamente i cubi delle due
torri. Pochi hanno contato esattamente il numero delle facce
esterne. Tutti hanno avuto difficoltà nel motivare le loro risposte
(``PERCHÉ?''). Durata 30 minuti circa.

Fase 4: Abbiamo chiesto ai bambini di esprimere le loro considerazioni
sull'esperienza di laboratorio. Elenchiamo le più significative:
\begin{studente}[Bambini]
  \begin{itemize}
  \item \bambini{Abbiamo scoperto e capito nuove cose divertendoci}
  \item \bambini{Abbiamo imparato a scoprire da soli}
  \item \bambini{Parlando con i compagni è più facile capire}
  \item \bambini{È stato un po' difficile dopo aver discusso mettersi
      d'accordo con i compagni su una risposta comune}
  \item \bambini{Abbiamo avuto difficoltà a contare le facce dei cubi
      e a misurare le figure piane usando le corde}
  \item \bambini{Ci è piaciuto lavorare con i compagni}
  \end{itemize}
\end{studente}
Durata 10 minuti.

\subsection{Osservazioni di fine percorso}
La sperimentazione effettuata nelle nostre due classi terze è stata
complessivamente positiva. Il gradimento e il coinvolgimento degli
alunni durante le attività sono stati alti. L'approccio ai nuovi
argomenti in modo laboratoriale ha suscitato curiosità, interesse e
motivazione all'apprendere. Riteniamo positiva la possibilità che i
bambini hanno avuto di confrontarsi nel piccolo gruppo per concordare
delle soluzioni comuni ai quesiti proposti.

\begin{figure}[pht]
  \centering
  \label{pic:torri:4}
  \begin{tabular}{cc}
    \includegraphics[width=0.47\textwidth]{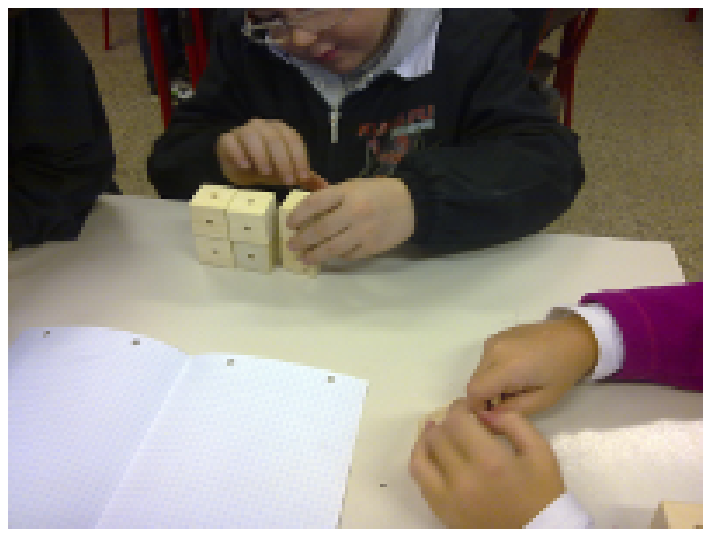} &
    \includegraphics[width=0.47\textwidth]{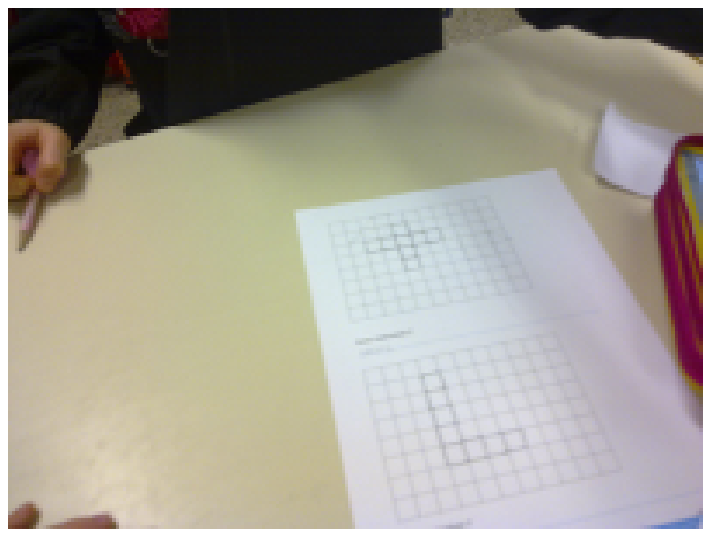} \\
    \includegraphics[width=0.47\textwidth]{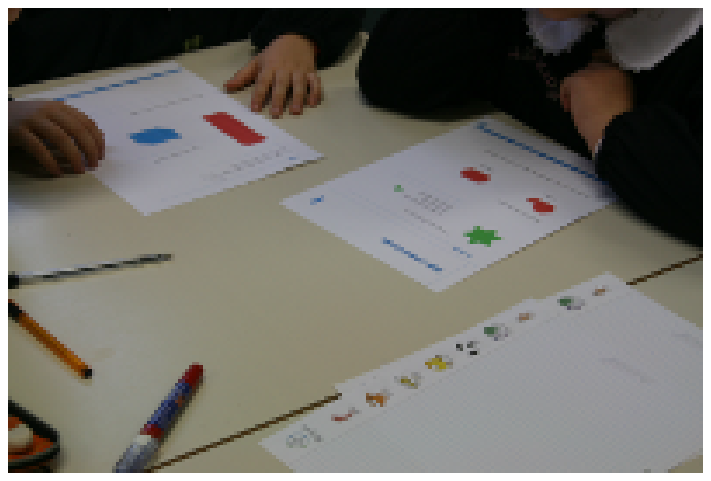} &
    \includegraphics[width=0.47\textwidth]{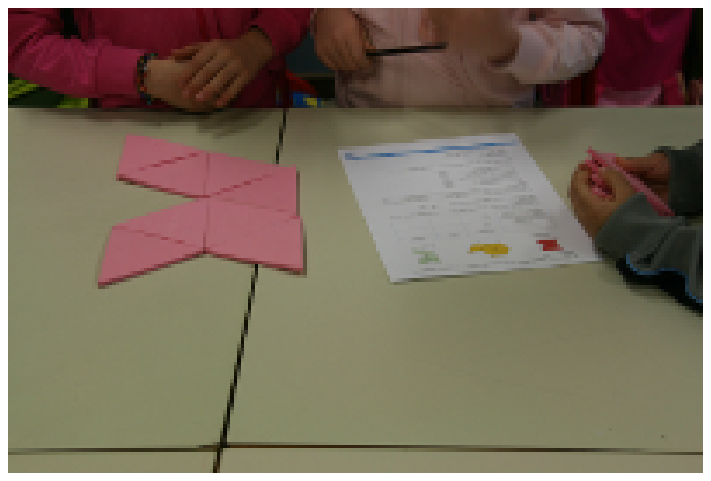} \\
    \includegraphics[width=0.47\textwidth]{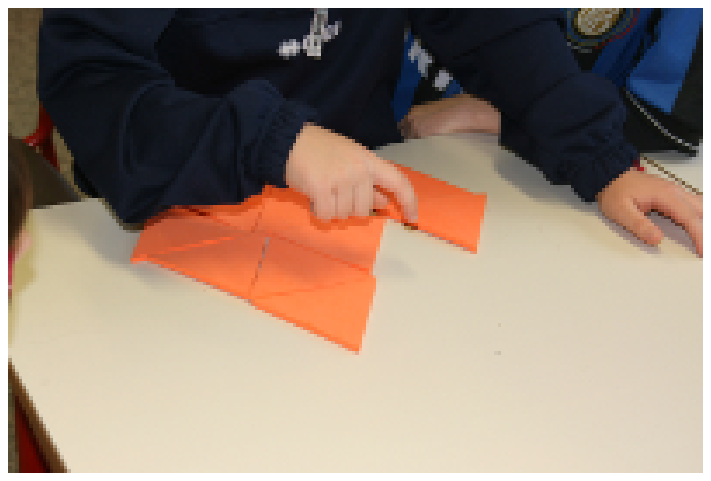} &
    \includegraphics[width=0.47\textwidth]{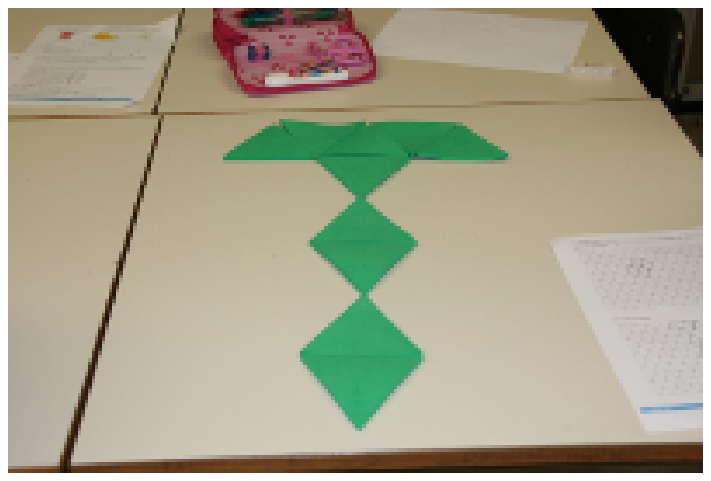} \\
    \includegraphics[width=0.47\textwidth]{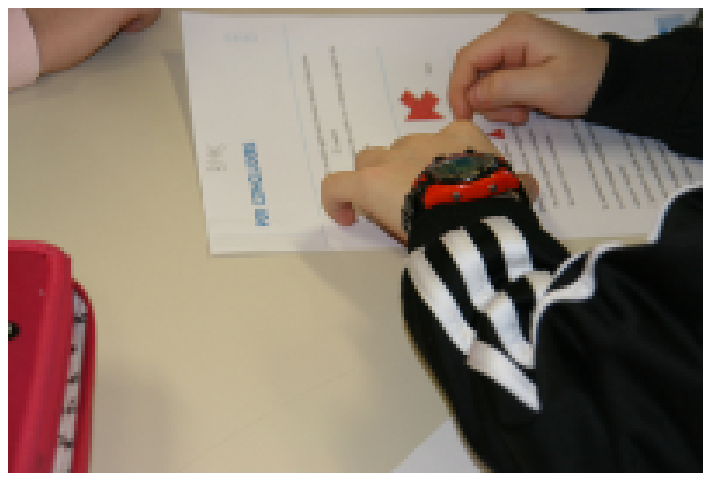} &
    \includegraphics[width=0.47\textwidth]{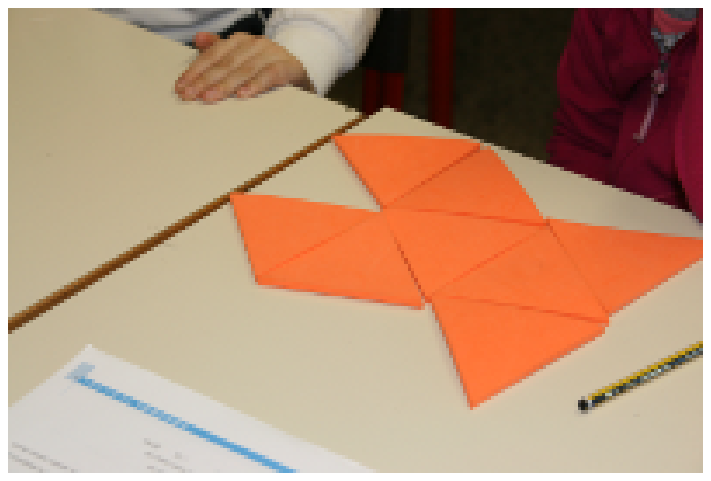} \\
  \end{tabular}
\end{figure}
Non avendo potuto proporre le attività alla classe durante le ore di
compresenza, abbiamo avuto alcune difficoltà nel doppio ruolo di
conduttore e osservatore. Ci siamo accorte di ``aver perso'' alcune fasi
di lavoro di alcuni gruppi.

Pensiamo che il materiale proposto dal kit sia riuscito a catturare
l'interesse dei ragazzi, aiutandoli a conseguire quegli obiettivi che
ci eravamo proposte all'inizio del percorso. Riteniamo però che il
poco tempo a disposizione (1 giorno alla settimana per 4 settimane)
abbia un po' limitato le potenzialità del kit stesso, soprattutto per
il consolidamento dei concetti acquisiti.

Per quanto riguarda la nostra formazione riteniamo che siano stati
molto positivi gli incontri ``laboratoriali'' effettuati in
Bicocca%
; la presenza del tutor e di altri colleghi ci ha offerto
l'opportunità di un confronto diretto e di un approccio alla
progettazione.

Per quanto riguarda l'uso della piattaforma, è stato positivo il poter
riflettere sul percorso effettuato nelle nostre classi, anche grazie
ai suggerimenti dei tutor; abbiamo invece trovato difficoltà a
confrontare la nostra esperienza con quella dei colleghi (anche perché
effettuate in tempi diversi)\dots{} pensiamo inoltre che sia più
efficace un confronto immediato e\dots{} ``verbale''.


\section[Sperimentazione \#5: quinta primaria]{Sperimentazione \#5:
  classe quinta primaria, gennaio/febbraio~2010}

\subsection{Osservazioni generali}

\subsubsection{Presentazione della classe}
Si tratta di una classe quarta composta da 24 alunni di cui 2
segnalati per disturbi dell'apprendimento e 2 alunni stranieri.

\subsubsection{Composizione dei gruppi}
I gruppi, eterogenei formati da me, sono 5. I gruppi saranno sempre
gli stessi per i 4 incontri ma sarà possibile scambiarsi i ruoli.

\subsubsection{Insegnanti presenti}

Durante le attività è presente l'insegnante di sostegno

\subsubsection{Calendarizzazione degli incontri}
\begin{calendario}
  \begin{itemize}
  \item 27 gennaio
  \item 29 gennaio
  \item 3 febbraio
  \item 4 febbraio
  \end{itemize}
\end{calendario}

\subsection{Primo incontro}

\begin{description}
\item[Alunni presenti:] 24 alunni presenti (tutti)
\item[Tempo effettivo di lavoro:] Tempo 2 ore circa
\end{description}

\begin{consegna}
Scheda \attivita{Per cominciare} di classe quarta.

\materiali{}%
Consegna delle 16 tessere a ogni gruppo e della prima parte delle
schede di lavoro
\end{consegna}

\subsubsection{Osservazioni}
Fase 1: comunico ai bambini che parteciperanno a un laboratorio di
geometria. Utilizzeranno materiale fornito dall'insegnante che
permetterà loro di avvicinarsi a alcuni concetti geometrici. Saranno
divisi in gruppi già stabiliti. Ogni gruppo dovrà darsi un nome e
scegliere le persone che avranno il ruolo di scrivere, disegnare,
riferire alla classe. Preciso che il laboratorio avrà la durata di 4
incontri di 2 ore ciascuno; che i gruppi saranno sempre gli stessi, ma
sarà possibile scambiarsi i ruoli.

Fase 2: consegna schede \attivita{Per cominciare} di classe quarta
relativa al lavoro con le tessere. Tempo di esecuzione 10 minuti.

Alcuni bambini riferiscono che è difficile rispondere senza avere a
disposizione il materiale. Li tranquillizzo dicendo che hanno ragione,
ma di cercare di dare una risposta.

Fase 3: consegna delle 16 tessere e della prima parte delle schede di
lavoro (pag. 3-4). Subito si mettono al lavoro. I primi dubbi arrivano
quando devono indicare ``Quant'è lungo il contorno di\dots{}''  avendo a
disposizione solo il materiale dato: tessere, schede, matita, gomma,
pastelli. Chiedono di poter usare il righello (quest'anno abbiamo
rivisto le misure di lunghezza). Rispondo di trovare un modo diverso
per calcolarlo. Dopo varie discussioni e in tempi diversi i gruppi
arrivano alla stessa conclusione:
\begin{studente}[ ]
  Contiamo i lati esterni dei triangoli, tanto abbiamo verificato che
  sono tutti uguali
\end{studente}
Il resto del lavoro procede più veloce.

Durante il lavoro siamo dovute intervenire per fare in modo che
leggessero con attenzione quanto richiesto e che tutti
partecipassero. Gli alunni più deboli cognitivamente hanno avuto
difficoltà a seguire i ragionamenti dei compagni i quali sono stati
invitati a spiegare con calma le loro conclusioni. I 2 bambini
segnalati non sono riusciti a dare il loro contributo sebbene
stimolati da noi e dai compagni. Il confronto delle risposte viene
rimandato alla data successiva (venerdì).

\subsubsection{Consigli per i colleghi che vogliono proporre le stesse
  attività}
Utile lavorare in compresenza se si ha una classe numerosa e con
problemi. Ciò permette di intervenire in tempi brevi per aiutare a
risolvere dubbi o incertezze, evitando che gli alunni perdano la
motivazione o l'attenzione.

\subsection{Secondo incontro}

\begin{description}
\item[Alunni presenti:] 24 alunni (tutti)
\item[Tempo effettivo di lavoro:] 2 ore circa
\end{description}

\begin{consegna}
  Schede del kit (pag. 3-4-5-6)

  \materiali{}%
  Consegno il materiale (24 Tessere per ogni gruppo) e le schede
  (pag.3-4-5-6). I bambini non hanno chiesto altro materiale
\end{consegna}

\subsubsection{Osservazioni}
I bambini, divisi in gruppi, si ridistribuiscono i ruoli. Consegno il
materiale e le schede (pag.3-4-5-6). Li invito a rivedere il lavoro
precedente per poi continuare. Molti si dimostrano entusiasti di
continuare l'attività; alcuni, inizialmente passivi, si lasciano poi
coinvolgere. I due bambini segnalati sono in grande difficoltà: per
loro sarebbe stato opportuno svolgere il lavoro iniziando dalle schede
di seconda. Le richieste impegnano molto gli alunni. Tendono a
realizzare figure allungate, poco compatte; ricordano la soluzione
dell'esercizio 3 di pag. 4 e fanno fatica a sganciarsi da quello
schema. Interveniamo suggerendo di provare a fare figure meno estese
. Li abbiamo invitati a osservare bene come era stata formata la
lumaca. Dopo vari tentativi ci riescono con loro grande
soddisfazione. Capito come fare, le altre proposte richieste dal n. 5
vengono svolte con più facilità. Dicono:
\begin{studente}[ ]
  \begin{itemize}
\item \bambini{Abbiamo capito che tenendo sempre la stessa forma, una
    specie di palla, potevamo fare il lavoro che ci chiedeva la
    scheda}
\item \bambini{Le figure hanno la stessa forma, ma cambiando numero di
    tessere diventano sempre più grandi}
\item \bambini{Tutte hanno una forma circolare perché le tessere devono
    essere il più vicino possibile}
\end{itemize}
\end{studente}
Gli ultimi 20 minuti vengono utilizzati per lo scambio di informazioni
Dicono
\begin{studente}[ ]
  abbiamo costruito figure diverse ma con lo stesso numero di tessere
\end{studente}
\begin{studente}[ ]
  Anche noi e abbiamo visto che allora occupano lo stesso spazio
\end{studente}
Dico loro che è corretto e che quelle figure si dicono ``equiestese''
\begin{studente}[ ]
  Ma siamo riusciti a fare anche figure diverse, ma con lo stesso
  contorno
\end{studente}
Spiego che il contorno si dice perimetro e che quelle figure si dicono
``isoperimetriche'', hanno cioè uguale perimetro.

\subsubsection{Consigli per i colleghi che vogliono proporre le stesse
  attività}
\`E importante seguire con attenzione i lavori dei gruppi così da
risolvere subito eventuali dubbi invitandoli a continuare.

\subsection{Terzo incontro}

\begin{description}
\item[Alunni presenti:] Presenti 24 alunni (tutti)
\item[Tempo effettivo di lavoro:] Tempo 2 ore circa
\end{description}

\begin{consegna}
  Consegna scheda \attivita{Per cominciare} relativa ai cubi

  \materiali{}%
  Consegna materiale: 9 mattoni e alcuni pioli. Schede
  (pag. 7-8-9-10). I bambini non hanno chiesto altro materiale
\end{consegna}

\subsubsection{Osservazioni}
Fase 1: Divisione in gruppi e ridefinizione dei ruoli. Consegna scheda
\attivita{Per cominciare} relativa ai cubi. Tempo 10 minuti

I bambini sono in grande difficoltà, non sanno come
rispondere. Intervengo dicendo di stare tranquilli, che è sicuramente
un lavoro difficile non avendo il materiale a disposizione, ma, come
per le tessere, avrebbero risolto i dubbi in seguito. Solo 2 gruppi
danno risposte ai quesiti.

Fase 2: Consegna materiale: 9 mattoni e alcuni pioli. Schede
(pag. 7-8-9-10) e li informo che non devono terminare tutto nella
giornata,

Costruiscono facilmente i 2 oggetti (aereo e cane). Faticano a contare
le facce esterne. Li invitiamo a trovare una strategia che gli
permetta di contare senza dimenticare nulla. Dicono:
\begin{studente}[ ]%
  \begin{itemize}
\item \bambini{Proviamo a contare prima quelli davanti a me, poi
    quelli davanti a te, poi sopra e infine sotto}
\item \bambini{Per il cane contiamo da una parte poi raddoppiamo; poi
    quelli di contorno compresi i 2 che appoggiano}
\end{itemize}
\end{studente}
Due gruppi arrivano a intuire in poco tempo che le due figure
occupano lo stesso spazio
\begin{studente}[ ]
  È come per le tessere!
\end{studente}
gli altri ci arrivano dopo aver discusso molto. La richiesta di pag. 8
viene risolta dai gruppi abbastanza facilmente, mentre la n. 3 di
pag. 9 è per loro impegnativa. Ci dicono che è difficile. Ci mostrano
figure troppo allargate e ogni volta diciamo loro di provare a
metterle più unite
\begin{studente}[ ]
  Dobbiamo quindi metterle il più unite possibile come per le tessere
\end{studente}
Finalmente, con grande soddisfazione formano il cubo. È passata 1 ora
e 30 minuti. Sono stanchi e alcuni mostrano segni di nervosismo. I 2
bambini segnalati intervengono solo per eseguire i suggerimenti dei
compagni. Decidiamo di interrompere l'attività. Spiego che la stessa
verrà ripresa il giorno successivo

\subsubsection{Consigli per i colleghi che vogliono proporre le stesse attività }

Un numero maggiore di cubi a disposizione li avrebbe aiutati a
comparare le due figure iniziali

\subsection{Quarto incontro}

\begin{description}
\item[Alunni presenti:] 24 (tutti)
\item[Tempo effettivo di lavoro:] 2 ore circa
\end{description}

\begin{consegna}
  Seconda parte lavoro cubi (da pag.~9) e schede \attivita{Per
    concludere}

\materiali{}%
mattoni 10
\end{consegna}

\subsubsection{Osservazioni}
Fase 1:  Distribuzione materiale (mattoni) e relative schede.

I bambini si mettono subito al lavoro: rivedono ciò che avevano
eseguito e cominciano a costruire riprendendo da pag.~9. La proposta 4
li mette a dura prova: nonostante tutti i loro tentativi non trovano
soluzione. A questo punto interveniamo per non rischiare lo
scoraggiamento dicendo che in effetti non vi è soluzione possibile al
quesito e li facciamo riflettere sul perché. Affrontano quindi
l'ultima proposta. Qualsiasi figura realizzata dà sempre come
risultato un numero pari di facce. Qualcuno comincia a guardare
attentamente i mattoni e dice
\begin{studente}[ ]
  Il cubo ha 6 facce. Se provo a metterne vicini 2 ottengo 10 facce
  cioè un numero pari
\end{studente}
Invito a fare la prova con 3, con 4 cubi. Finalmente arrivano alla
soluzione
\begin{studente}[ ]
  Se il cubo ha 6 facce quando li unisco nascondo 2 facce alla
  volta. Quindi non ci può essere una figura con un numero dispari di
  facce
\end{studente}

Fase 2: Ritiro il materiale e distribuisco la scheda \attivita{Per
  concludere}, pag. 11-12. Tempo 10 minuti

I quesiti riguardanti le tessere vengono svolti esattamente e
velocemente. Quelli riguardanti i cubi ancora con qualche difficoltà

\subsection{Riflessioni a fine percorso}
Presenti sempre 24 alunni. I gruppi sono rimasti invariati, mentre i
ruoli venivano ridistribuiti di volta in volta. Entrambe le insegnanti
(di classe e di sostegno) sono intervenute su richiesta degli alunni
per chiarire, stimolare le riflessioni, sollecitare la
collaborazione. Attenzione particolare è stata rivolta ai 2 gruppi con
inseriti i bambini segnalati. Abbiamo dovuto chiedere al gruppo di
spiegare sempre in modo chiaro i ragionamenti e i diversi passaggi
così che tutti seguissero il lavoro. Abbiamo più volte sollecitato
l'intervento dei bambini segnalati, ma per loro le richieste erano
troppo alte. Sono solo riusciti a contare il numero delle tessere
usate per la formazione delle varie figure.

Alcune riflessioni: la manipolazione di materiale strutturato ha
permesso agli alunni di sperimentare e scoprire nuovi concetti
mettendo alla prova le loro capacità logiche e intuitive. Sapere che
non vi sarebbe stata alcuna valutazione, li ha aiutati a esprimersi
spontaneamente all'interno del gruppo. I bambini segnalati avrebbero
la necessità di iniziare il percorso partendo dalle attività di classe
seconda e di lavorare individualmente con l'insegnante. I bambini
solitamente più bravi hanno faticato a adeguare i loro tempi di
lavoro (più svelti) con quelli di chi era in difficoltà (più
lenti). \`E stato faticoso gestire le dinamiche relazionali e
riportare sempre il gruppo a una effettiva collaborazione.

In conclusione: È stata una faticaccia, ma ce l'abbiamo
fatta!


\section[Sperimentazione \#6: quinta primaria]{Sperimentazione \#6:
  classe quinta primaria, gennaio/febbraio~2010}

\subsection{Osservazioni generali}

\subsubsection{Presentazione della classe}
La classe è composta da 20 alunni, di cui 4 sono extracomunitari (una
di loro, da poco inserita, ancora non parla italiano). La
comunicazione di questa attività è stata accolta con entusiasmo dagli
alunni che amano lavorare in gruppo e manipolare materiali strutturati
e non.

\subsubsection{Composizione dei gruppi}
L'insegnante ha deciso di formare lei stessa gruppi eterogenei di
lavoro perché vi sono alcuni alunni con una forte personalità e altri
che si adeguano passivamente alle loro decisioni. I gruppi sono
rimasti invariati per tutto il percorso, sono invece cambiati i ruoli
al loro interno. La classe è formata da 20 alunni; l'insegnante dunque
ha formato 4 gruppi da 5 per poter osservare meglio le dinamiche al
loro interno (ha condotto le attività da sola).

\subsubsection{Insegnanti presenti}
Agli incontri è presente solo l'insegnante di classe.

\subsubsection{Calendarizzazione degli incontri}
\begin{calendario}
  \begin{itemize}
  \item 18 gennaio
  \item 21 gennaio
  \item 25 gennaio
  \item 28 gennaio
  \end{itemize}
\end{calendario}

\subsection{Primo incontro}

\begin{description}
\item[Alunni presenti:]20 alunni (tutti)
\item[Tempo effettivo di lavoro:]dalle ore 14.35 alle ore 16.35
\end{description}

\begin{consegna}
  È stata utilizzata la scheda \attivita{Per cominciare} del
  kit. Essendo una classe V ho sempre dato le schede senza alcuna
  indicazione, anche perché volevo verificare il loro grado di
  comprensione dei comandi.

  \materiali{}%
  Ho consegnato ai gruppi tutti i triangoli disponibili
\end{consegna}

\subsubsection{Osservazioni}
Dopo aver comunicato loro le modalità di lavoro, ho lasciato 10 minuti
per organizzarsi, scegliere il nome del gruppo e i ruoli dei diversi
componenti; in ogni gruppo il compito di scrivere sulle schede è stato
dato a alunne femmine in genere molto ordinate, mentre quello di
rappresentare le figure con il disegno era ambito da molti. Su loro
richiesta si è deciso di poter cambiare i ruoli all'interno dei gruppi
le volte successive.

Ho quindi consegnato la scheda \attivita{Per cominciare} relativa al
lavoro con le tessere triangolari lasciando 10 minuti per rispondere
ai quesiti. Tutti gli alunni hanno dato il loro contributo per la
soluzione tranne la bimba da poco inserita che ancora non capisce bene
l'italiano. Completata la scheda ho consegnato il materiale a ogni
gruppo per poter verificare le loro ipotesi risolutive.

Ho poi distribuito la scheda della prima esperienza. Mi è stata
richiesta la spiegazione del termine ``area'' in quanto ancora non era
stato presentato agli alunni (avevo solo introdotto il concetto di
superficie). Ho notato che alcuni alunni iniziavano a costruire figure
con il materiale senza aspettare la completa lettura della richiesta,
per questo non operavano in modo esatto. Sono comunque stati ``ripresi''
dagli altri componenti del gruppo. Nei momenti operativi anche
l'alunna extracomunitaria cercava di dare il suo contributo, formando
figure con caratteristiche simili a quelle dei suoi compagni. Un
alunno ha chiesto il mio intervento perché i compagni non stavano
operando in modo corretto, ma non volevano ascoltare la sua
proposta. Completata la scheda ho lasciato 15 minuti per costruire e
disegnare altre figure isoperimetriche e equiestese. Ogni gruppo ha
poi relazionato e confrontato le risposte date ai quesiti.

\subsection{Secondo incontro}
\begin{description}
\item[Alunni presenti:]20 alunni (tutti)
\item[Tempo effettivo di lavoro:]dalle ore 8.30 alle ore 10
\end{description}

\begin{consegna}
  Sono state utilizzate le schede \attivita{Per cominciare} relativa
  ai cubi e la scheda \attivita{per la prima attività}
  (finestra/podio) del kit. Essendo una classe V ho sempre dato le
  schede senza alcuna indicazione, anche perché volevo verificare il
  loro grado di comprensione dei comandi.

  \materiali{}%
  ho dato solo i cubi necessari per la singola scheda (essendo
  concetti nuovi per loro e più complicati non volevo troppa
  confusione e distrazione).
\end{consegna}

\subsubsection{Osservazioni}
Gli alunni si sono organizzati in gruppi in pochi minuti, ridefinendo
i ruoli al loro interno. Nel momento operativo di costruzione di
figure è nata una ``gara'' spontanea tra i gruppi per chi riusciva a
costruire figure con il contorno di lunghezza minore e maggiore
possibile, dato un determinato numero di tessere. Ho lasciato fare
perché ciò era uno stimolo positivo a trovare nuove soluzioni. Nella
verbalizzazione comune dell'attività i gruppi hanno cercato di
spiegare le strategie adottate: \bambini{per il perimetro minore
  bisogna costruire figure compatte, con poche punte che escono fuori
  perché ogni punta ha due lati\dots{}}, \bambini{per il perimetro
  maggiore bisogna costruire figure molto estese e con più punte che
  escono\dots{}}. Riposte le tessere triangolari ho distribuito la
scheda \attivita{Per cominciare} relativa ai cubi; gli alunni hanno
avuto maggiore difficoltà a rispondere ai quesiti. Per verificare ho
distribuito i cubi e 3 gruppi su 4 si sono accorti di non aver
risposto correttamente (non avevano calcolato la giusta profondità dei
solidi).

Ho quindi consegnato la scheda con la prima attività
(finestra/podio). Il podio è stato costruito da tutti i gruppi con 9
cubi, di conseguenza hanno risposto in modo diverso dalla previsione.

Al termine del loro lavoro ho sollecitato la costruzione del podio con
solo 8 cubi (ci sono voluti diversi tentativi) per poter confrontare
figure solide con uguale volume e area esterna diversa.

Il calcolo dell'area esterna del podio è stato complicato per gli
alunni; i gruppi hanno ripetuto il calcolo più volte prima di arrivare
a una strategia corretta. Un gruppo si è fidato del calcolo eseguito
da un solo compagno (quello che in genere non sbaglia mai), purtroppo
questa volta è stato troppo frettoloso!

Nelle verbalizzazioni al termine delle attività gli alunni hanno
incontrato maggiori difficoltà a usare termini corretti, perché erano
legati a concetti nuovi per loro.

Nonostante questo tutti gli alunni hanno riposto il materiale a
malincuore chiedendo di riprendere presto le attività.

\subsection{Terzo incontro}

\begin{description}
\item[Alunni presenti:]20  (tutti)
\item[Tempo effettivo di lavoro:]dalle ore 14.30 alle ore 16.30
\end{description}

\begin{consegna}
  comandi 2-3-4 delle schede. Essendo una classe V ho sempre dato le
  schede senza alcuna indicazione, anche perché volevo verificare il
  loro grado di comprensione dei comandi.

  \materiali{}%
  mattoncini.
\end{consegna}

\subsubsection{Osservazioni}
Gli alunni sono stati invitati a dividersi in gruppo e a predisporre
da soli i banchi per le attività. Si mostrano contenti di riprendere
il lavoro. Hanno lavorato con i mattoncini per rispondere ai comandi
2-3-4 delle schede. La costruzione di figure solide richiede più
tentativi delle figure piane con i triangoli e quindi più tempo. Più
mani che spostano i mattoni creano confusione, per questo due gruppi
hanno deciso di costruire a turno le figure e di confrontarsi al loro
interno. In un gruppo invece si tende a delegare i due elementi più
``forti'' a trovare la soluzione. Tutti hanno scoperto facilmente che
il cubo è il solido con l'area esterna minore dato un determinato
numero di mattoni. Un alunno ha confermato a voce alta che si applica
anche nei solidi la regola che una figura compatta occupa meno
spazio. Il comando 5 chiedeva di costruire una figura con 17 facce. Due
alunni, ancor prima di iniziare la manipolazione, hanno intuito la non
possibilità, ma in genere, prima di arrivare alla soluzione, hanno
fatto molti tentativi; la spiegazione delle loro scoperte però era
poco chiara, quindi abbiamo cercato insieme la regola aritmetica per
confermare il numero pari di facce esterne di mattoncini.

\subsection{Quarto incontro}

\begin{description}
\item[Alunni presenti:]20 alunni (tutti)
\item[Tempo effettivo di lavoro:]dalle ore 9.00 alle ore 10.30.
\end{description}

\begin{consegna}
  Scheda \attivita{Per concludere}. Essendo una classe V ho sempre
  dato le schede senza alcuna indicazione, anche perché volevo
  verificare il loro grado di comprensione dei comandi.

  \materiali{}%
  Il materiale previsto dal kit
\end{consegna}

\subsubsection{Osservazioni}
Ho deciso di dare la scheda \attivita{Per concludere} a ogni alunno
per verificare quanto abbia inciso questa esperienza in ognuno di
loro. I bambini più deboli leggono la consegna più volte prima di
rispondere, in generale tutti faticano a completare la parte relativa
ai mattoni. Non è semplice capire da quanti mattoni è composta una
figura; qualcuno usa il righello e la divide. Riuniti in gruppo poi
ogni alunno ha confrontato le proprie risposte con quelle dei compagni
e verificato con il materiale.

\subsection{Riflessioni a fine percorso}

Apprendere attraverso una metodologia ``attiva e collaborativa'' rende
gli alunni più motivati e coinvolti: i bambini con più difficoltà
hanno potuto imparare ``costruendo'', a volte anche solo per imitazione;
quelli con buone capacità, invece, hanno potuto consolidare,
approfondire e essere più consapevoli del proprio sapere. All'interno
dei gruppi c'è stata collaborazione, anche se i più ``bravi'' e veloci
cercavano di imporre comunque la propria idea. In genere i maschi
hanno manipolato maggiormente il materiale, mentre le femmine
guardavano, suggerivano e poi, dopo di loro, ripetevano l'attività. I
maschi, in questa classe hanno quasi tutti una forte personalità e
sono competitivi tra loro.

L'alunna extracomunitaria ha lavorato per imitazione.

Nei momenti di difficoltà in ogni gruppo è stato chiesto il mio
intervento.

Avere già il materiale sui banchi prima della scheda portava qualcuno
a distrarsi nel momento di lettura dei comandi.

Ho dato importanza al momento di confronto finale tra i gruppi; i
bambini più deboli faticavano a mantenere l'attenzione; per questo
chiedevo a più elementi di ripetere la riflessione comune.

Ora fisseremo i concetti acquisiti con esempi e verbalizzazioni su
cartelloni e sui loro quaderni.


\section[Sperimentazione \#7: prima primaria]{Sperimentazione \#7:
  classe prima primaria, febbraio/marzo~2010}

\subsection{Osservazioni generali}

\subsubsection{Presentazione della classe}
 26 alunni

\subsubsection{Composizione dei gruppi}
Gruppi eterogenei scelti dall'insegnante: 5 gruppi, di questi quattro
formati da 5 alunni e uno da 6.

Su indicazione dell'insegnante, durante ogni incontro, ogni gruppo ha
scelto un responsabile.

Durante questa scelta alcuni alunni si sono imposti, altri sono stati
scelti in accordo tra i bambini secondo i seguenti criteri ``il
responsabile deve essere attento, sveglio, deve saper
scrivere\dots{}''.

Gli alunni hanno sperimentato il lavoro di gruppo per attività manuali
solo alla scuola materna.

\subsubsection{Insegnanti presenti}

La docente sperimentatrice è sempre stata affiancata dalla collega di
classe.

\subsubsection{Calendarizzazione degli incontri}
\begin{calendario}
  \begin{itemize}
  \item 25 febbraio
  \item 1 marzo
  \item 8 marzo
  \item 9 marzo
  \end{itemize}
\end{calendario}

\subsection{Primo incontro}

\begin{description}
\item[Alunni presenti:] 25 presenti,1 assente
\item[Tempo effettivo di lavoro:] Circa 2 ore: dalle 10,45 alle 12,30
\end{description}

\begin{consegna}
  Ho iniziato la sperimentazione presentandola come un nuovo modo di
  imparare, giocando con del materiale nuovo.

  Ho proposto l'attività \attivita{Per cominciare} (solo la parte
  relativa alle tessere), in seguito ho fatto manipolare liberamente
  il materiale e creare pavimentazioni.

  Ho proposto la prima esperienza con le tessere triangolari
  presentata dal kit e alla fine ho invitato gli alunni di ogni gruppo
  a scegliere una forma costruita, a cui dare un nome e a riprodurla
  su un cartellone, dove in precedenza era stata predisposta una
  griglia avente come unità di misura il triangolo del kit.

  \materiali{}%
  Ogni gruppo ha 12 tessere e la scheda di lavoro. Gli alunni non
  hanno chiesto altro materiale.
\end{consegna}

\subsubsection{Osservazioni}
Gli alunni si sono dimostrati subito entusiasti delle attività
proposte: hanno partecipato attivamente, lavorato con più impegno
rispetto al solito e si sono divertiti.

Tutti hanno dimostrato interesse. Ognuno voleva portare il proprio
contributo, anche quelli che solitamente durante le lezioni non
partecipano mai.

L'attività manuale li ha stimolati positivamente.

Le difficoltà si sono rivelate nel lavoro di gruppo, poiché alcuni
alunni hanno faticato a collaborare e a ascoltare i suggerimenti dei
compagni, a causa di un eccesso di individualismo.

L'attività \attivita{Per Cominciare} è stata proposta secondo le
indicazioni del kit.

Ogni consegna e/o domanda è stata letta e spiegata dall'insegnante.

Dopo una prima osservazione della scheda gli alunni hanno commentato:
\begin{studente}[ ]
  Non sembra un'ochetta\dots{} è fatta con dei triangoli\dots{} è
  fatta con delle forme geometriche\dots{} è tutta rossa
\end{studente}
Alla domanda dell'insegnante:
\begin{tutor}[ ]
  Cosa significa ricoprire l'ochetta?
\end{tutor}
un alunno risponde:
\begin{studente}[ ]
  È come quando hai un panino, metti sopra la marmellata e hai
  ricoperto il panino
\end{studente}
Tre gruppi su cinque hanno individuato facilmente il numero di tessere
\begin{studente}[ ]
  Sono 6\dots{} bisogna contarle\dots{} è stato facile\dots{} basta
  pensare a tre più tre
\end{studente}
Un gruppo ha valutato a occhio
\begin{studente}[ ]
  Sono 7\dots{} sono 8\dots{} no, bisogna mettere il dito su ogni
  tessera
\end{studente}
mentre un altro ha disegnato sulla scheda altre tesserine vicino
all'ochetta.

L'esperienza delle tessere ha creato nella classe un clima di gioco e
una forte tendenza di ogni alunno a impossessarsi del maggior numero
di tessere. In un gruppo un alunno ha diviso le tessere per colore e
le ha distribuite ai compagni equamente.

Solo un gruppo ha mostrato subito un atteggiamento collaborativo
costruendo insieme le forme.

All'interno di alcuni gruppi si sono create situazioni di stallo,
dovute a difficoltà di collaborazione e a eccessi di individualismo
da parte di alcuni alunni desiderosi di ``fare tutto''.

Faticosa per alcuni gruppi è stata la scelta di una forma comune,
dovuta a difficoltà di adeguarsi alle decisioni della maggioranza.

Le forme decise sono state: la stella, il razzo, il serpente delle
montagne russe, il gatto selvaggio e Nemo.

Più collaborativi si sono invece dimostrati durante la riproduzione
sui cartelloni: ognuno colorava il triangolo per formare la figura
scelta.

Durante questa attività i gruppi hanno posizionato le tessere della
loro figura sul cartellone, poi, togliendole una alla volta, hanno
colorato lo spazio interno. Alcuni hanno ripassato il contorno della
figura.

Il momento del confronto è stato molto qualificante: ogni gruppo ha
mostrato alla classe il cartellone preparato con la forma scelta dal
gruppo, spiegando le difficoltà incontrate durante questa scelta e
durante la sua realizzazione:
\begin{studente}[ ]
  All'inizio non sapevamo come fare, poi ci è venuta l'idea, abbiamo
  messo insieme le tessere, ma ogni tanto le tessere si spostavano,
  allora abbiamo fatto il contorno e colorato dentro i triangoli e la
  figura si è formata
\end{studente}

\subsubsection{{Consigli per i colleghi che vogliono proporre le
    stesse attività }}
La presenza di due insegnanti.

\subsection{Secondo incontro}
\begin{description}\item[Alunni presenti:]
 Tutti gli alunni sono presenti
\item[Tempo effettivo di lavoro:]
 Circa 2 ore: dalle 10,45 alle 12,30
\end{description}

\begin{consegna}
  Agli alunni sono state proposte la seconda e la terza esperienza con
  le tessere triangolari secondo le indicazioni riportate nel testo
  ``Torri, serpenti \dots{} e geometria''.

  Ogni esperienza è stata preceduta da una lettura da parte
  dell'insegnante e da una puntualizzazione delle richieste della
  scheda.

  \materiali{}%
  Ogni gruppo ha 12 tessere e la scheda di lavoro.

  Gli alunni non hanno chiesto altro materiale.
\end{consegna}

\subsubsection{Osservazioni}
Tutti gli alunni sono contenti di riprendere l'attività.

Ricordano le modalità di utilizzo delle tessere e iniziano a costruire
le figure richieste.

Sono molto divertiti da questa attività: inizialmente qualcuno
costruisce le figure senza osservare i modelli, ma poi viene
richiamato dagli altri componenti del gruppo.

Dopo aver costruito le figure gli alunni discutono e si confrontano
tra loro sul numero di tessere usate per costruirle.

Durante la fase iniziale di questa attività alcuni gruppi necessitano
di essere guidati nella riproduzione delle figure, altri chiedono solo
conferme.

Tutti i gruppi, anche se in tempi diversi, riescono a ricostruire le
figure.

Durante la presentazione della terza attività gli alunni notano subito
che le figure presenti sulla scheda
\begin{studente}[ ]
  non hanno più le righe per far vedere quante sono le tessere
\end{studente}
Un gruppo che nella seconda attività aveva avuto qualche problema
nella ricostruzione delle figure, ora è il primo a finire.

Alla fine di questa attività alcuni bambini riferiscono:
\begin{studente}[ ]
  Il lavoro è stato un po' difficoltoso perché alcuni costruiscono,
  mentre altri disfano il lavoro fatto
\end{studente}
Durante il momento del confronto, i livelli di consapevolezza delle
esperienze a cui sono giunti gli alunni sono stai diversi. Alcuni
hanno avuto bisogno della guida dell'insegnante per arrivare a
scoperte, altri sono riusciti a verbalizzare ciò che hanno fatto e le
difficoltà incontrate
\begin{studente}[ ]
  \begin{itemize}
  \item \bambini{abbiamo avuto difficoltà a costruire il serpente
      perché c'erano posizioni che non riuscivamo a capire}
  \item \bambini{noi abbiamo fatto fatica a ricostruire l'aquilone,
      non riuscivamo a capire se le punte andavano in su o in giù }
  \item \bambini{il gatto è stato difficile!}
  \item \bambini{le orecchie abbiamo capito come farle, soltanto che
      quando ci giravamo qualcuno le spostava}
  \end{itemize}
\end{studente}
Durante questo momento i gruppi si sono confrontati e, vedendo che
alcune risposte erano diverse tra loro, hanno provato, su suggerimento
di alcuni compagni, a ricostruire le figure ricontando le tessere;
infine hanno abbandonato le loro posizioni iniziali.

\subsubsection{Consigli per i colleghi che vogliono proporre le
    stesse attività}
 Lasciare più spazio al confronto tra gli alunni nel piccolo gruppo.

\subsection{Terzo incontro}

\begin{description}
\item[Alunni presenti:] Tutti gli alunni sono presenti
\item[Tempo effettivo di lavoro:] Circa 2 ore: dalle 10,45 alle 12,30
\end{description}

\begin{consegna}
  Agli alunni è stata proposta la quarta esperienza con le tessere
  triangolari secondo le indicazioni riportate nel testo ``Torri,
  serpenti \dots{} e geometria''.

  L'esperienza è stata preceduta da una lettura da parte
  dell'insegnante e da una puntualizzazione delle richieste della
  scheda.

  \materiali{}%
  Ogni gruppo ha 12 tessere e la scheda di lavoro.

  Gli alunni non hanno chiesto altro materiale.
\end{consegna}

\subsubsection{Osservazioni}
Gli alunni ricordano ciò che è stato fatto nelle lezioni precedenti e
sono entusiasti di iniziare un'altra attività.

Durante questa esperienza i vari gruppi non hanno difficoltà a
ricostruire l'immagine dell'ochetta, ma alcuni devono essere guidati
nella ricostruzione del pinguino.

Alcuni gruppi richiedono la presenza dell'insegnane per mostrare le
loro costruzioni e per avere conferme sul lavoro svolto, a volte anche
per intervenire nelle dinamiche relazionali.

Ampio spazio viene dato al momento del confronto, per cercare di
individuare le caratteristiche delle due immagini. Due gruppi
riferiscono di aver fatto fatica nella ricostruzione del pinguino
perché
\begin{studente}[ ]
  \dots{} Non riuscivamo a capire da quanti triangoli era
  fatto,\dots{} perché non l'avevamo mai fatto, \dots{} era più
  difficile, \dots{} siamo partiti dalla testa, poi facevamo fatica a
  fare il corpo
\end{studente}
Tutti i gruppi individuano che il pinguino occupa più spazio sul
tavolo, ma non tutti sanno spiegare il perché. Alcuni riferiscono
\begin{studente}[ ]
  \dots{} Occupa più spazio perché ha più tessere, \dots{} è più
  grande
\end{studente}
Osservando le due figure alcuni alunni riferiscono che le immagini
hanno anche forme diverse. Un'alunna osserva
\begin{studente}[ ]
  \dots{} Possiamo guardare anche il bordo e vedere da quale misura è
  fatto
\end{studente}
Un alunno specifica
\begin{studente}[ ]
  Il bordo possiamo chiamarlo contorno o confine della figura \dots{}
  e per contare il bordo della figura contiamo i suoi lati
\end{studente}
Alcuni bambini hanno fatto fatica a partecipare a questo momento,
lasciavano parlare gli altri e si distraevano, mentre altri erano
molto interessati e continuavano a fare osservazioni.

\subsubsection{Consigli per i colleghi che vogliono proporre le stesse attività}

Coinvolgere tutti gli alunni durante il momento del confronto, anche i
più restii.

\subsection{Quarto incontro}

\begin{description}
\item[Alunni presenti:] Tutti gli alunni sono presenti
\item[Tempo effettivo di lavoro:] 2 ore: dalle 8,30 alle 10,30
\end{description}

\begin{consegna}
  Agli alunni è stata proposta la quinta esperienza con le tessere
  triangolari secondo le indicazioni riportate nel testo ``Torri,
  serpenti\dots{} e geometria''.

  L'esperienza è stata preceduta da una lettura da parte
  dell'insegnante, da una puntualizzazione delle richieste della
  scheda e successivamente è stata riprodotta su un cartellone, dove
  in precedenza era stata predisposta una griglia avente come unità di
  misura il triangolo del kit.

  Consegna della scheda \attivita{Per Concludere} (solo parte relativa
  alle tessere).

  \materiali{}%
  Ogni gruppo ha 12 tessere, la scheda di lavoro e il cartellone. Gli
  alunni non hanno chiesto altro materiale.
\end{consegna}

\subsubsection{Osservazioni}
I gruppi, rimasti invariati, riescono a collaborare tra loro e a
adeguarsi alle decisioni della maggioranza. In questo modo non hanno
difficoltà a costruire, a scegliere una figura diversa dal pinguino e
a riprodurla sia sulla scheda che sul cartellone.

Alla domanda
\begin{tutor}[ ]
  Cosa è successo ai vari elementi che formavano il pinguino?
\end{tutor}
gli alunni rispondono
\begin{studente}[ ]
  \dots{} Li abbiamo spostati\dots{} ma solo tre\dots{} così è venuta
  fuori una nuova figura,\dots{} è come se si è trasformata, ma è
  rimasto uguale il numero di tessere,\dots{} lo spazio occupato sul
  banco dalla nuova figura è uguale a quello del pinguino
\end{studente}
A queste osservazioni, guidate dall'insegnante, gli alunni sono giunti
dopo un'attenta osservazione delle figure.

Durante il lavoro conclusivo non emergono particolari difficoltà: i
gruppi accettano di non utilizzare il materiale. Dicono
\begin{studente}[ ]
  È come un mettersi alla prova
\end{studente}

Tutti i gruppi individuano subito il numero di tessere dell'ochetta,
ricordando la scheda \attivita{Per Cominciare}, mentre solo tre gruppi
quelle del fantasma.

Le modalità per contare le tessere sono le seguenti: segnare ogni
tessera con le dita, disegnare con la matita le tessere, immaginare le
tessere.

Durante il momento di confronto finale, avendo dato risposte diverse,
i gruppi discutono tra loro e chiedono di provare a ricostruire le
figure con le tessere, per verificare praticamente le loro risposte.

\subsubsection{Consigli per i colleghi che vogliono proporre le stesse
  attività}

In classi così numerose è utile la formazione del minor numero di
gruppi, per consentire agli animatori un'osservazione più puntuale e
dettagliata.


\section[Sperimentazione \#8: terza primaria]{Sperimentazione \#8:
  classe terza primaria, febbraio/marzo~2010}

\subsection{Osservazioni generali}

\subsubsection{Presentazione della classe}
La classe è composta da 21 alunni, 9 maschi e 12 femmine.

Nella classe sono da segnalare: un bambino dislessico, un alunno
egiziano con difficoltà di inserimento nel gruppo classe e gravi
difficoltà di apprendimento.

\subsubsection{Composizione dei gruppi}
Gruppi eterogenei. 3 gruppi da 5 bambini, 1 gruppo da 4
bambini. L'insegnante ha dato a 4 bambini con difficoltà, il compito
di formare i gruppi, scegliendo a turno un compagno/a

\subsubsection{Insegnanti presenti}

In alcuni incontri l'insegnante sperimentatrice è affiancata
dall'altra insegnante di classe (compresenza).

\subsubsection{Calendarizzazione degli incontri}
\begin{calendario}
  \begin{itemize}
  \item 26 febbraio
  \item 2 marzo (compresenza)
  \item 3 marzo (compresenza)
  \item 5 marzo
  \end{itemize}
\end{calendario}

\subsection{Primo incontro}

\begin{description}
\item[Alunni presenti:] Bambini presenti 19 - assenti 2
\item[Tempo effettivo di lavoro:] 2 ore
\end{description}

\begin{consegna}
  Ogni gruppo ha avuto il materiale/schede del kit nella modalità e
  nei tempi stabiliti dal testo

  \materiali{}%
  I materiali del kit vengono distribuiti dopo la presentazione
  dell'attività e la formazione dei gruppi.
\end{consegna}

\subsubsection{Osservazioni}

Fase~1:
Ai bambini viene comunicato che ``parteciperanno a una
sperimentazione matematica'' viene detto loro che dovranno collaborare
per provare a imparare in un modo nuovo con nuovi materiali.

Tutti intraprendono con entusiasmo l'attività, in ogni gruppo i
bambini provano a suddividersi i compiti, nascono le prime discussioni
perché gli alunni con difficoltà, scelti per formare i gruppi,
vogliono essere "i capi".

È necessario l'intervento dell'insegnante che ribadisce la necessità
di collaborare senza la presenza di un capo in questa prima fase del
lavoro.

Fase~2:
Vengono distribuite le schede \attivita{Per iniziare} e tra i gruppi
non c'è rivalità, ma in 3 di essi si evidenziano difficoltà nel
trovare strategie e soluzioni condivise.

Solo un gruppo prova a completare la scheda confrontandosi e
verificando le diverse idee.

Fase~3:
Vengono date 12 tessere a ogni gruppo, con la richiesta di formare
figure a loro piacere, 2 gruppi si suddividono le tessere equamente
tra maschi e femmine, 2 gruppi costruiscono figure insieme usando
tutte le tessere a disposizione.

Alcuni bambini provano a usare le tessere per costruire piramidi.

Da più bambini vengono richieste un maggior numero di tessere.

Le tessere vengono contese, i bambini con difficoltà tendono a
utilizzarle come giochi dimenticandosi dell'attività da svolgere.

Quando viene chiesto di scegliere alcune figure piane e di riprodurle
(pag.20), senza discussioni il compito è affidato dal gruppo ai
bambini più abili nel disegno.

Fase~4:
Vengono confrontati i disegni, ogni gruppo racconta quello che ha
costruito, la classe ascolta e partecipa con osservazioni pertinenti
alla discussione. I bambini con difficoltà non intervengono, si
soffermano sui nomi scelti per le figure costruite.

Emergono osservazioni in merito alle tessere:
\begin{studente}[ ]
  \begin{itemize}
  \item \bambini{sono morbide}
  \item \bambini{sono triangoli uguali}
  \item \bambini{sono piatte}
  \end{itemize}
\end{studente}
riguardo alle figure dicono:
\begin{studente}[ ]
  \begin{itemize}
  \item \bambini{sono figure simili costruite con numeri di tessere
      differenti}
  \item \bambini{sono figure diverse fatte con un uguale numero di
      triangoli}
  \end{itemize}
\end{studente}
(osservazioni rilevanti).

\subsubsection{Consigli per colleghi che vogliono proporre le stesse
  attività}
Per eseguire l'attività nel migliore dei modi sarebbe stata necessaria
la presenza della collega di classe purché adeguatamente informata
circa l'attività da svolgere.

\subsection{Secondo incontro}

\begin{description}
\item[Alunni presenti:] presenti 19, assenti 2
\item[Tempo effettivo di lavoro:]2 ore
\end{description}

\begin{consegna}
  Schede kit per terza elementare

  \materiali{}%
  Vengono distribuite ai gruppi le schede previste dal testo per il
  confronto tra serpente, stella\dots{}, 12 tessere e le corde.

  I bambini non hanno chiesto altro materiale.
\end{consegna}

\subsubsection{Osservazioni}

Fase~1:
Vengono distribuiti il materiale e le schede. Si nota maggior
collaborazione all'interno di ogni gruppo ma minor coinvolgimento dei
bambini con difficoltà.

Durante la compilazione delle schede a un gruppo sembra ambigua la
richiesta numero 3.  Non sa se osservare tutte le 5 figure o solo le 3
della pagina in questione.

Fase~2:
Tutti i gruppi incontrano difficoltà nel far aderire le corde al
contorno delle figure, alcuni bambini sottolineano la praticità di
contare i lati per misurare il contorno.

Tra i gruppi nasce il desiderio di terminare prima degli altri la
compilazione delle schede, i bambini più competitivi prendono il
sopravvento, viene meno la condivisione delle scelte.

Fase~3:
La collaborazione riprende quando i gruppi incontrano difficoltà nel
costruire figure con 10 tessere dal contorno più lungo e più corto
possibile.

Tutti i bambini si mettono in gioco, provano diverse soluzioni e
chiedono più tempo per assicurarsi che le loro figure siano le più
adeguate alla richiesta.

Alcune osservazioni fatte a voce alta dai bambini più veloci, riguardo
ai lati che si toccano e non si toccano aiutano quasi tutti a
arrivare a una soluzione adeguata.

Fase~4:
Quando tutti i gruppi consegnano le schede compilate, si passa al
confronto.

A turno un bambino scelto dal gruppo relaziona quanto fatto e legge le
risposte date.

Due gruppi hanno trovato con 10 tessere figure dal contorno di 8 e 30
pezzi, due gruppi hanno trovato contorni diversi.

L'insegnante invita tutti i bambini al confronto per trovare la
soluzione corretta, alcune bambine provano a convincere i compagni che
8 e 30 sono le lunghezze esatte perché:
\begin{studente}[ ]
  nel primo caso le tessere sono TUTTE VICINE mentre nel secondo
  NESSUN TRATTINO È APPOGGIATO SU UN ALTRO
\end{studente}
la classe all'unanimità accetta la spiegazione delle compagne. La
classe sembra contenta del lavoro svolto, i gruppi che hanno dato
risposte diverse non manifestano delusione%
.

Alcuni bambini non sembrano comprendere l'importanza del confronto,
non ascoltano le soluzioni dei compagni e tendono a disturbare il
momento di rielaborazione. L'insegnante più volte deve intervenire per
riportare silenzio e concentrazione.

\subsection{Terzo incontro}

\begin{description}
\item[Alunni presenti:] presenti 18,  assenti 3
\item[Tempo effettivo di lavoro:]2 ore
\end{description}

\begin{consegna}
  schede kit \attivita{Per iniziare} (attività con i cubi)

  \materiali{}%
  Le schede sono date subito dopo la formazione dei gruppi; i cubi, 8
  per ciascun gruppo, dopo la compilazione delle schede per
  iniziare. Tutti i gruppi chiedono più materiale.
\end{consegna}

\subsubsection{Osservazioni}

Fase~1:
Ai gruppi vengono date le schede \attivita{Per iniziare} (attività con
i cubi), tutti i gruppi contano velocemente i cubetti della torre
rossa ma incontrano difficoltà a contare i mattoni della torre blu.

Le insegnanti incoraggiano i bambini a provare a contare usando varie
strategie, nessuno trova difficile rispondere alla domanda relativa
allo spazio occupato.

Tutti ricordano: ``più tessere più superficie'' così senza pensare
affermano
\begin{studente}[ ]
  più mattoni più superficie
\end{studente}
L'insegnante ne approfitta per dire a tutti i gruppi che le tessere
che erano piatte occupavano una superficie, mentre i mattoni che sono
spessi occupano uno spazio.

Vengono distribuiti 8 mattoni a ciascun gruppo, il materiale suscita
interesse e entusiasmo, si nota maggior coinvolgimento e
concentrazione nel completare le schede del testo.

Il lavoro svolto con le tessere ha lasciato la consapevolezza che è
necessario intraprendere con metodo e collaborazione il lavoro.

Fase~2:
Ai bambini viene chiesto di scegliere tre figure costruite con otto
mattoni e di disegnarle, tutti i gruppi fanno un ottimo lavoro,
cercando di dare il senso della profondità. Il lavoro prosegue con un
momento di confronto, ogni portavoce dei gruppi descrive le figure
disegnate.

Sono interessanti le osservazioni fatte da alcuni bambini che
riprendono l'esperienza delle tessere e sottolineano che anche in
questo caso ci sono figure diverse costruite con lo stesso numero di
mattoni.

Una bambina ripete quanto detto in merito a superficie e spazio, altri
bambini ricordano un'esperienza fatta in classe prima, con oggetti
aventi forma di solidi che hanno lasciato la loro impronta. Quasi
tutti i bambini ricordano quale oggetto avevano portato e quale
impronta aveva lasciato.

L'insegnante riporta la discussione al lavoro svolto, un bambino nota
che tre gruppi hanno costruito e disegnato \bambini{un cubo fatto con
  8 cubi} dice che sono figure \bambini{proprio uguali}. I bambini
sembrano stanchi e solo pochi sono ancora attenti.

La discussione viene fermata e si avvia la terza fase dell'attività.

Fase~3:
I gruppi devono completare le schede del libro contando e confrontando
il numero delle facce esterne delle figure disegnate. Tutti incontrano
difficoltà, ogni bambino cerca di contare a voce alta generando
confusione. La presenza di due insegnanti in classe rende più semplice
il controllo dei gruppi che sono guidati nel contare con metodo le
facce.

Risulta immediato per tutti i gruppi risolvere l'esercizio 4, poiché
durante il momento di confronto la classe si era soffermata a
osservare il GRANDE DADO.

I bambini sembrano sicuri ormai nel confrontare le figure osservando
lo spazio occupato e il numero delle facce esterne, solo pochi bambini
manifestano dubbi.

Fase~4:
Ogni gruppo legge le risposte date completando le schede, tutto
avviene molto velocemente, le risposte coincidono, non è necessario
discutere ulteriormente.

I bambini sono stanchi, il lavoro sarebbe stato da dividere in due
giornate.

\subsubsection{Consigli per i colleghi che vogliono proporre le stesse
  attività }
Più tempo da dedicare alla manipolazione e più cubi a disposizione

\subsection{Quarto incontro}

\begin{description}
\item[Alunni presenti:] 19 bambini presenti, 2 assenti
\item[Tempo effettivo di lavoro:]1 ora
\end{description}

\begin{consegna}
  schede kit \attivita{Per concludere}

  \materiali{}%
  A ciascun gruppo vendono date le schede \attivita{Per concludere}
  prese dal testo, i bambini avrebbero voluto ancora i materiali,
  tessere e cubi, ormai restituiti.
\end{consegna}

\subsubsection{Osservazioni}

Fase~1:
Vengono consegnate le schede \attivita{Per concludere}, tutti i gruppi
eseguono senza difficoltà la parte relativa alle tessere.

Nasce qualche dubbio nell'esercizio di confronto delle torri. Nei
gruppi nasce tensione perché i bambini non riescono a accordarsi. Un
solo gruppo decide di provare a disegnare i mattoni per facilitare il
conteggio delle facce esterne. Gli alunni con difficoltà si
distraggono e non partecipano alla compilazione delle schede, il
lavoro di gruppo sembra aver perso per loro il fascino suscitato nei
primi incontri. Tutti i gruppi compilano la parte finale della scheda
esprimendo le loro difficoltà:
\begin{studente}[ ]
  \begin{itemize}
  \item \bambini{Abbiamo dovuto contate perché le torri hanno forme
      diverse}
  \item \bambini{Abbiamo dovuto contare perché non riuscivamo a
      immaginare}
  \item \bambini{Abbiamo dovuto contare perché non sapevamo come poter
      scoprire il numero delle facce}
  \item \bambini{Abbiamo dovuto contare e ci è servito il righello}
  \end{itemize}
\end{studente}

Fase~2:
Nel momento del confronto emerge che tutti i gruppi hanno eseguito il
lavoro con le tessere allo stesso modo, trovando risposte
adeguate. Nella compilazione della parte riguardante i cubi, solo un
gruppo, quello che ha disegnato i mattoni, ha risposto che le torri
hanno lo stesso numero di facce esterne.

Tre gruppi hanno detto che LA TORRE VERDE ha più facce esterne e
sostengono di aver contato con attenzione. Sarebbe stata necessaria a
questo punto una verifica con il materiale ma il kit era stato
consegnato%
. L'insegnante decide
di usare i cubetti bianchi dei regoli per ricostruire le torri.

Contando insieme la facce esterne i bambini si convincono degli errori
e qualche bambino sottolinea che
\begin{studente}[ ]
  CON GLI OGGETTI È PIÙ FACILE
\end{studente}


\section[Sperimentazione \#9: terza primaria]{Sperimentazione \#9:
  classe terza primaria, febbraio/marzo~2010}

\subsection{Osservazioni generali}

\subsubsection{Presentazione della classe}
La classe è composta da 20 alunni:
\begin{itemize}
\item uno con il sostegno per deficit legati alla difficoltà di
  comprensione, ha seguito terapia logopedia e psicomotoria;
\item una bimba dislessica seguita dalla logopedista;
\item è attualmente in via di accertamento la probabile dislessia in
  un altro alunno;
\item è stata inserita nella classe a gennaio dello scorso anno una
  bimba proveniente da un orfanotrofio del Brasile e, pur mostrando
  notevoli progressi, questa alunna evidenzia ancora una certa
  difficoltà di comprensione e capacità di attenzione limitate, che si
  traducono, a volte, in atteggiamenti di disturbo;
\item una bimba di origini rumene, attenta e con desiderio di
  partecipazione, ma con difficoltà nel linguaggio sia orale che
  scritto, è stata inserita nella nostra classe in prima come
  ripetente, è stata seguita dalla logopedista, ma non è ancora stato
  determinato se le sue difficoltà siano dovute a una lingua madre
  diversa o a un problema specifico del linguaggio;
  \item da pochi mesi è stato inserito un bambino di origini cinesi,
    ma nato a Milano, è attento e con buone capacità di apprendimento.
\end{itemize}

Gli alunni hanno un concetto intuitivo di perimetro che chiamiamo
confine, legato ai concetti di regione interna e esterna. Hanno
giocato in modo pratico con figure solide in legno che ho messo loro a
disposizione, giocando a fare i timbri con le facce dei solidi abbiamo
ottenuto le figure piane, questa attività è stata svolta in modo
giocoso e i bambini non conoscono ancora i termini geometrici
specifici.

\subsubsection{Composizione dei gruppi}
Abbiamo chiesto ai bambini di formare liberamente quattro gruppi con
la presenza sia di maschi che di femmine. Si sono formati in pochi
minuti tre gruppi composti da cinque alunni e uno da quattro (c'era un
alunno assente). Uno di questi gruppi era costituito completamente da
bambini con buone capacità, gli altri tre comprendevano livelli
diversi di apprendimento.

\subsubsection{Insegnanti presenti}
In alcuni incontri l'insegnante sperimentatrice è affiancata
dall'insegnante di sostegno e/o da alcune tirocinanti.

\subsubsection{Calendarizzazione degli incontri}
\begin{calendario}
  \begin{itemize}
  \item 1 marzo (presenza insegnante sostegno e tirocinante per il
    sostegno)
  \item 2 marzo (presenza tirocinante per il sostegno), con
    integrazione il 3 marzo
  \item 4 marzo (presenza di due tirocinanti)
  \item 8 marzo (presenza tirocinante per il sostegno)
  \end{itemize}
\end{calendario}

\subsection{Primo incontro}

\begin{description}
\item[Alunni presenti:] 19 (un alunno è assente)
\item[Tempo effettivo di lavoro:] Dalle 8.45 alle 10.25
\end{description}

\begin{consegna}
  La classe già sapeva che stavo frequentando questo corso presso
  l'Università e era informata che avremmo sperimentato un materiale
  nuovo e avremmo lavorato in gruppo. All'inizio della lezione ho
  cercato di creare motivazione e interesse, spiegando loro che
  giocando con il materiale avremmo dovuto verificare cosa era
  possibile capire e imparare, ho detto che non c'era nessun problema
  in caso di errore, perché a volte dagli errori si possono scoprire
  tante cose, in fondo anche l'America è stata scoperta per errore:
  Cristoforo Colombo pensava di trovarsi in India!

  Essendo presenti in tre adulti abbiamo ritenuto di non spiegare
  collettivamente le istruzioni, ma abbiamo invitato i bambini a
  leggere con molta attenzione le richieste e di richiedere
  spiegazioni solo in caso di difficoltà.

  \materiali{}%
  Abbiamo distribuito le schede \attivita{Per cominciare} di classe
  terza e a ogni gruppo, nel momento in cui tutti hanno terminato le
  schede, sono state consegnate 12 TESSERE triangolari e fogli di
  carta da pacco bianca.
\end{consegna}

\subsubsection{Osservazioni}
I gruppi si sono formati abbastanza velocemente, solo due bambini sono
rimasti seduti ``in attesa'': il bambino con sostegno, che è stato
subito chiamato da un compagno e da una compagna per invitarlo a
unirsi al loro gruppo e lui ha immediatamente aderito; una bambina che
spesso all'inizio ha questo tipo di comportamento non collaborativo e
sul cui atteggiamento stiamo cercando di intervenire, pertanto l'ho
invitata a accordarsi con i compagni e invitata da loro si è unita a
un gruppo. Una volta formati i gruppi, ho chiesto di organizzare i
banchi in modo adeguato al lavoro, abbiamo deciso insieme di metterli
in gruppi di quattro. Ho chiesto poi di dare un nome al gruppo e di
dividersi i ruoli per
\begin{itemize}
  \item leggere i quesiti
  \item scrivere le risposte
  \item disegnare le figure
  \item fare eventualmente domande all'insegnante
  \item riferire alla classe nel momento collettivo
\end{itemize}

Tutto si è svolto in pochi minuti senza particolari problemi. È stato
osservato in modo particolare l'alunno con diritto al sostegno, si è
notato che non ha mai fatto domande, ma ha risposto quando gli
venivano poste, non si è distratto e ha accettato l'incarico datogli
dal gruppo di riferire nel momento collettivo,
\begin{figure}[thp]
  \centering
  \label{pic:torri:9}
  \begin{tabular}{cc}
    \includegraphics[width=0.47\textwidth]{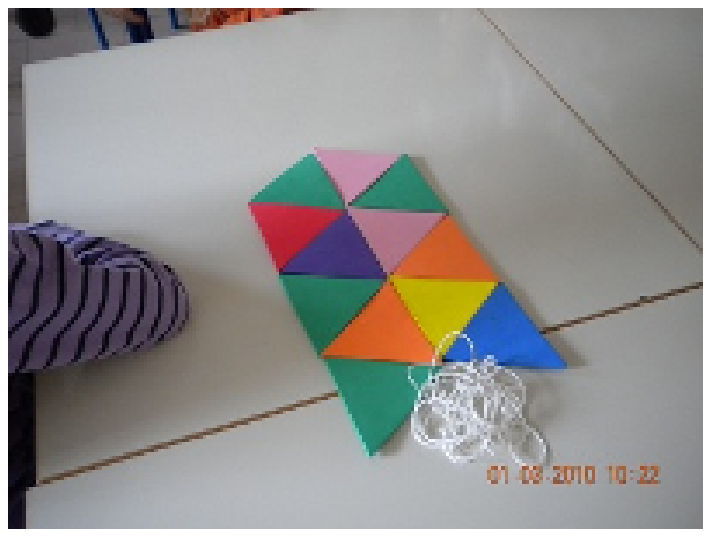} &
    \includegraphics[width=0.47\textwidth]{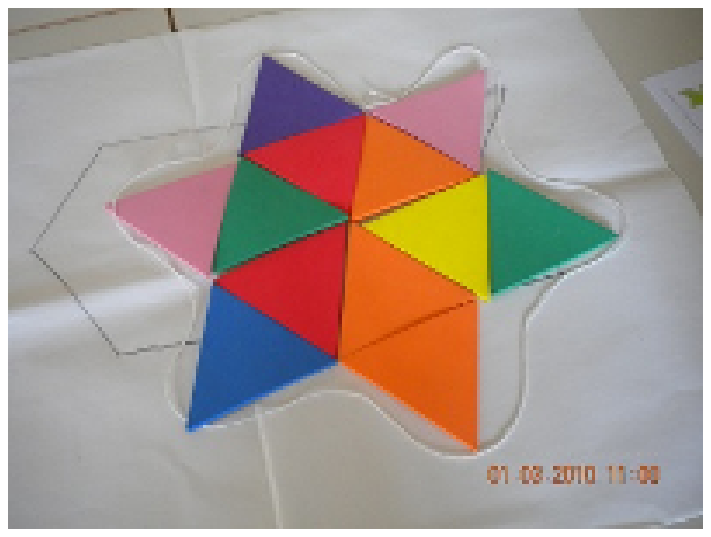} \\
  \end{tabular}
\end{figure}

\begin{wrapfigure}{R}{0pt}
  \includegraphics[width=0.47\textwidth]{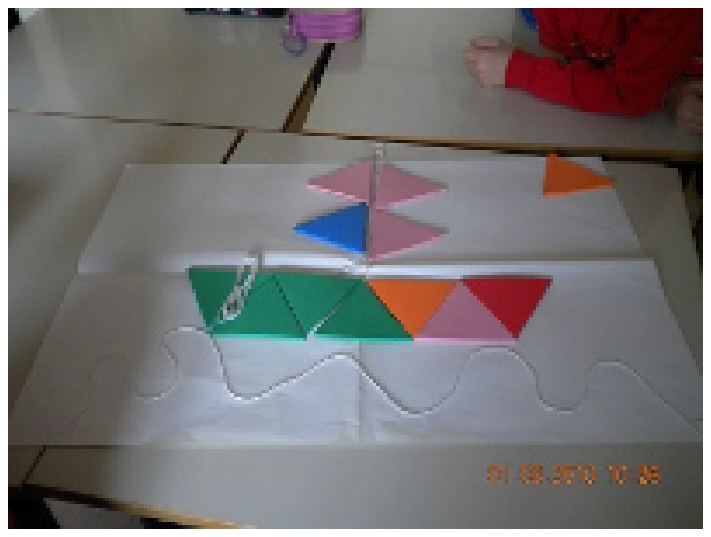}
\end{wrapfigure}
Il lavoro nei vari gruppi è iniziato con la scheda \attivita{Per
  cominciare}. Tutti i bambini hanno partecipato, anche se qualcuno
solo come ascolto e non esprimendosi attivamente. In un gruppo una
prima lettura superficiale aveva causato un conteggio errato delle
tessere utili per ricoprire la stella, ma una bambina ha portato
all'attenzione dei compagni del gruppo il significato della parola
``ricoprire'' e si sono pertanto corretti. Solo un gruppo si è ``aiutato''
disegnando le tessere all'interno della stella e, nell'esercizio
successivo, spostando con il disegno due lati e facendo diventare la
figura della lumaca uguale a quella del diamante. Un gruppo ha contato
le tessere necessarie solo per ricoprire le ``punte'' della stella.

Successivamente sono state distribuite le tessere e i fogli di carta
da pacco bianca, dove disegnare le figure costruite. Gli alunni hanno
subito iniziato a sperimentare, qualcuno ha sovrapposto le tessere, ma
poi ha visto che non riusciva a riportare la figura sul foglio e
allora ha modificato. In questa fase di lavoro si è verificato qualche
conflitto all'interno dei gruppi, conflitti risolti senza l'intervento
dell'adulto. Tre gruppi hanno disegnato sul foglio bianco le figure
con i triangoli che le componevano, uno ha disegnato solo il contorno
delle figure.

Abbiamo poi distribuito la scheda \attivita{Tessere} dove i bambini
hanno riportato una forma e eseguito come richiesto. Tutti hanno
cercato di realizzare forme con significato e alcuni, nel riportare la
figura sulla scheda, hanno completato in modo creativo con onde, sole
``triangolare''\dots{}

A questo punto la classe mi è sembrata stanca e ho ritenuto opportuno
non procedere con il momento di confronto collettivo.

\subsection{Secondo incontro}

\begin{description}
\item[Alunni presenti:] 19 (un alunno è assente)
\item[Tempo effettivo di lavoro:] Dalle 8.45 alle 10.30
\end{description}

\begin{consegna}
  Chiedo se l'organizzazione dei banchi il giorno precedente fosse
  funzionale, non si evidenziano particolari problemi, ma si decide di
  utilizzare solo due banchi per stare più vicini. Al termine della
  sistemazione chiedo ai bambini se non notano nulla che non va, un
  alunno si accorge che un banco è stato lasciato davanti alla porta
  d'ingresso e questo non si può fare per motivi di sicurezza, subito
  il banco viene sistemato in modo opportuno. Chiedo agli alunni di
  trovare la giusta concentrazione e iniziamo il lavoro senza
  ulteriori spiegazioni.

 \materiali{}%
 Si consegnano 12 tessere e le corde, le schede dove confrontare le
 figure del serpente, della stella e poi la clessidra, la lumaca e il
 gatto. Al termine di questo lavoro si consegnano le schede dove, dopo
 aver costruito, riportare le figure.
\end{consegna}

\subsubsection{Osservazioni}
Le figure proposte vengono realizzate abbastanza facilmente. Poco
funzionale per i bambini si rivela l'uso delle corde, che non
aderiscono bene e non si differenziano molto, decidono che è più
comodo contare i trattini. Ogni tanto ci accorgiamo che la loro
lettura è un po' superficiale: chiediamo di chiarire bene quando si
tratta di spazio e di tessere utilizzate e quando si tratta di
contorno; in qualche gruppo dobbiamo richiamare l'attenzione per
rispondere in modo corretto, ma comunque questa parte del lavoro viene
svolta con la partecipazione di tutti. In modo particolare, qualche
alunno, che durante il lavoro di classe non ci sembra attivo e
raramente interviene, l'abbiamo visto esprimere il proprio punto di
vista nel gruppo e partecipare: questo ci ha fatto molto
piacere%
.

Alcuni conflitti si sono evidenziati quando hanno dovuto costruire una
figura con il contorno che misurasse tanto quanto quello del serpente,
soprattutto perché volevano rappresentare qualcosa di reale e non
volevano assomigliasse a quella realizzata dagli altri gruppi. Allora
si sono sentite frasi del tipo \bambini{Non sei mai contenta di
  nulla}, \bambini{Non ti va mai bene\dots{}} Abbiamo dovuto invitarli
a riflettere bene e a confrontarsi per rispondere alla domanda:
\begin{tutor}[ ]
  Avete dato la risposta giusta prima di contare le tessere?
\end{tutor}
Durante l'attività svolta per realizzare le figure con il contorno più
corto e più lungo possibile, alcuni bambini si sono proprio
distratti. Ad esempio la bimba di origini brasiliane a questo punto è
stata motivo di disturbo coinvolgendo anche un altro elemento del
gruppo, un'altra bimba si è posta come leader e ha continuato il
lavoro con l'aiuto di un compagno.

Il gruppo formato da bambini con buone capacità ha realizzato una
figura con 8 trattini come contorno e una con 30, ma si è rivelato
conflittuale ``perdendo'' così molto tempo. Negli altri gruppi hanno
realizzato figure con contorno più corto e più lungo tra di loro.

In un gruppo due alunni hanno cercato di imporre le loro risposte
senza fornire spiegazioni accettabili dagli altri, per cui anche se
alcune conclusioni erano esatte, non sono state fatte proprie dal
gruppo.

Al termine del lavoro la bimba brasiliana ha dimostrato tutto il suo
impegno nel colorare le figure con molta attenzione, trovando così il
suo spazio positivo all'interno del gruppo.

I bambini mi sembrano stanchi e termino la lezione.

\subsubsection{Discussione finale}
\begin{wrapfigure}{R}{0pt}
  \includegraphics[width=0.48\textwidth]{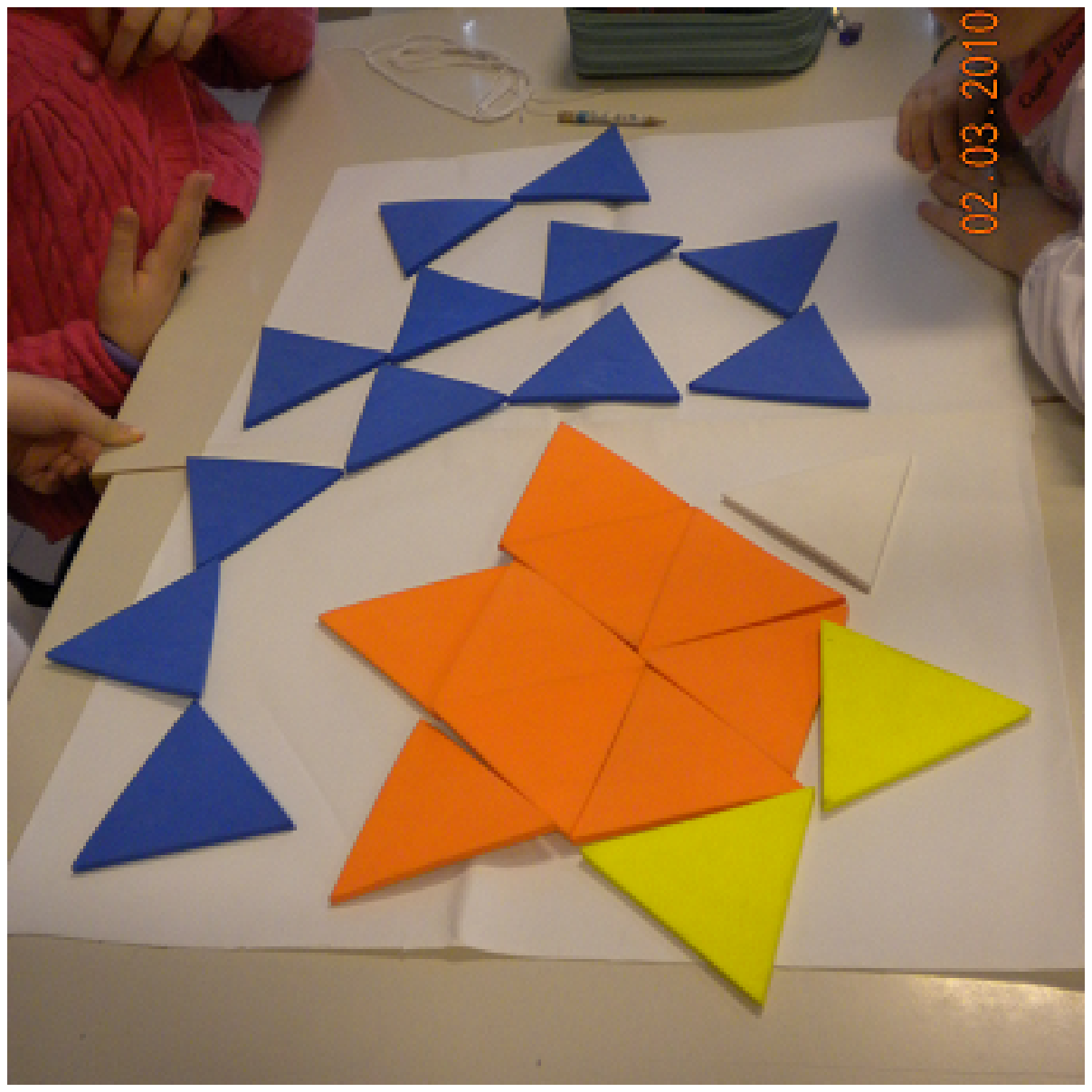}
\end{wrapfigure}
Nella giornata successiva, 3 marzo, durante un'ora esterna alle ore
dedicate alla sperimentazione, avviene l'esposizione collettiva del
lavoro di gruppo riguardante i due incontri, sono presenti entrambe le
insegnanti di classe.

Gli alunni incaricati hanno in mano le loro schede e illustrano i
risultati conseguiti. Tutti i compagni possono intervenire alzando la
mano e questo avviene.

L'incaricato del primo gruppo racconta come, grazie a una compagna
che ha evidenziato il significato di una parola (ricoprire), siano
riusciti a correggere un errore, per cui io evidenzio come
l'attenzione e la comprensione di ciò che si legge sia fondamentale.

A mano, a mano che il confronto prosegue si rendono conto di aver dato
risposte diverse, li invito a illustrare i ragionamenti, mi sembra che
ognuno prenda coscienza dei concetti proposti da questa attività.

Quando è la volta che deve illustrare l'attività del suo gruppo il
bambino con diritto al sostegno, racconta ai compagni come nel suo
gruppo abbiano sbagliato a contare le tessere della stella, ma abbiano
anche capito il perché:
\begin{studente}[ ]
  Abbiamo contato le tessere intorno e non quelle della pancia, ci
  siamo accorti facendo le costruzioni
\end{studente}
Comunque, nonostante il timore di noi adulti (quando ho visto il tipo
di incarico che gli era stato dato, per un attimo sono stata tentata
di intervenire, poi non l'ho fatto per una questione di rispetto,
avendogli i compagni dato l'incarico e avendolo lui accettato,
riservandomi così eventualmente di guidarlo)%
, questo alunno è riuscito a illustrare il lavoro con qualche
intervento dei compagni e soddisfazione delle insegnanti, che si sono
ripromesse di riferire alla mamma.

Quando ha illustrato il proprio lavoro il gruppo che ha realizzato le
figure con il contorno di 8 e 30 lati, i compagni sono stati
particolarmente attenti, hanno condiviso ciò che era stato scoperto:
per avere il contorno più corto si devono mettere il più possibile
vicino le tessere e per ottenere il contorno più lungo le tessere
devono toccarsi il meno possibile.

Alcuni bambini vorrebbero poter ``ritoccare'' le loro schede o i
disegni, ma spiego che non è possibile e che comunque la cosa
importante è che abbiano capito.

Gli alunni hanno detto di aver trovato il primo incontro divertente,
il secondo incontro meno, secondo me è solo perché si sono distratti
più facilmente di fronte a una richiesta per loro più difficile.

Mi è sembrato importante questo momento di rielaborazione del lavoro
svolto, sia per la condivisione dei ragionamenti e delle strategie,
sia per aver esplicitato in modo condiviso le regole di comportamento
per un lavoro di gruppo proficuo.

\subsection{Terzo incontro}

\begin{description}
\item[Alunni presenti:] 18 (due assenti)
\item[Tempo effettivo di lavoro:] Dalle 8.50 alle 10.30
\end{description}

\begin{consegna}
  Ho invitato gli alunni a ripensare alle riflessioni fatte e a
  svolgere questa parte del lavoro come meglio potevano. Non abbiamo
  letto insieme le consegne, ma ho ricordato che eravamo eventualmente
  disponibili per spiegare ciò che c'era scritto sulle schede.

  \materiali{}%
  Abbiamo consegnato 8 cubi e una quantità non definita di bastoncini
  dimostrando come potevano essere utilizzati, poi abbiamo consegnato
  le schede sui cubi.
\end{consegna}

\subsubsection{Osservazioni}
\begin{wrapfigure}{R}{0pt}
  \includegraphics[width=0.48\textwidth]{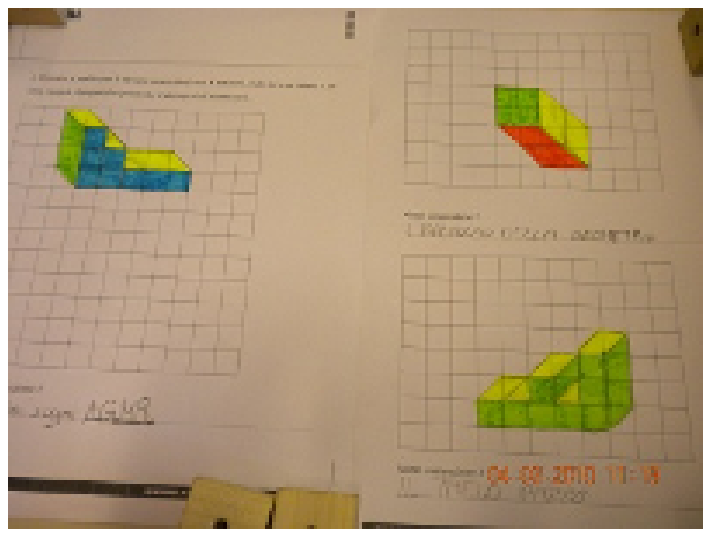}
\end{wrapfigure}
Dopo mia richiesta, una bambina di un gruppo formato da cinque alunni
si è resa disponibile a spostarsi dal proprio gruppo a quello dove
erano rimasti in tre a causa delle assenze.

Tutti i gruppi hanno costruito le figure, concordato, cercato la
costruzione più bella\dots{}

Tutti i bambini hanno partecipato a questo momento di
costruzione. Qualche difficoltà l'hanno incontrata quando hanno dovuto
disegnare la costruzione sulle schede. Due gruppi hanno rappresentato
solo due dimensioni e poi hanno chiesto agli adulti come rappresentare
lo ``spessore''. Un altro gruppo ha fatto il disegno in due dimensioni e
poi ha raddoppiato il contorno per realizzare lo spessore. Gli adulti
hanno suggerito praticamente a tutti i gruppi come disegnare le tre
dimensioni. Tre gruppi hanno colorato le figure dando un colore
diverso allo spessore e evidenziando così la tridimensionalità.

Quando nella terza scheda è stato richiesto di contare le facce, i
bambini hanno ricostruito le figure che ormai avevano distrutto, hanno
incontrato difficoltà, hanno cercato strategie, contato prima tutte le
facce da una parte e poi girato la figura\dots{} qualcuno cercava di
tenere il dito sulla faccia già contata\dots{} Poi un gruppo ha
iniziato a inserire un bastoncino in ogni faccia contata e la
strategia si è rivelata efficace e subito copiata dagli altri
gruppi. Durante il lavoro di conteggio delle facce alcuni alunni si
sono distratti e hanno iniziato a giocare con il materiale.

Nella costruzione della figura con 24 facce, inizialmente hanno avuto
tutti difficoltà, anche perché, quando ho detto che potevano decidere
loro quanti mattoni usare, hanno cercato di utilizzarne tanti, poi un
po' alla volta hanno diminuito la quantità di mattoni e ci sono
riusciti.

Nel confronto collettivo è emerso che avrebbero voluto non disfare
ogni volta le costruzioni e che era stato difficile contare le
facce. Mi sembra comunque che ci si stata una presa di consapevolezza
degli errori commessi.

\subsection{Quarto incontro}

\begin{description}
\item[Alunni presenti:] 20 (tutti presenti)
\item[Tempo effettivo di lavoro:] 1 ora
\end{description}

\begin{consegna}
  Ho chiesto ai bambini di riflettere sul lavoro svolto
  precedentemente, di leggere con attenzione le richieste, di
  rispondere alle domande senza l'uso del materiale.

  \materiali{}%
  Ho distribuito le schede \attivita{Per concludere}, a lavoro
  ultimato ho distribuito i cubi per verificare le risposte.
\end{consegna}

\subsubsection{Osservazioni}
\begin{wrapfigure}{R}{0pt}
  \includegraphics[width=0.48\textwidth]{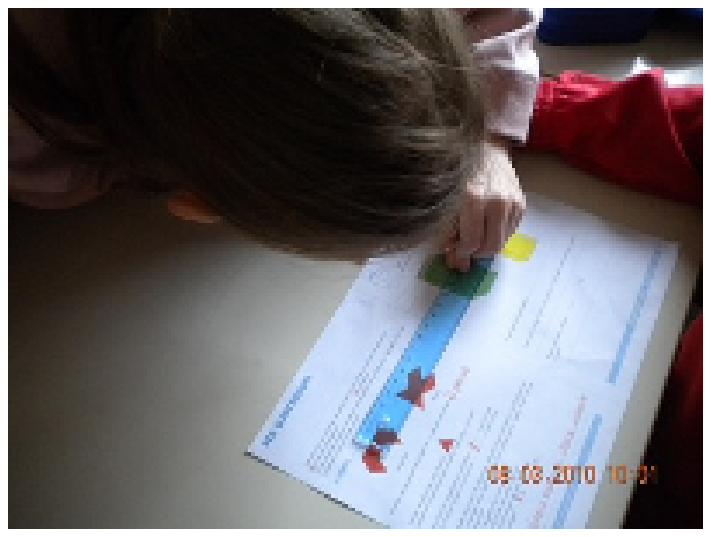}
\end{wrapfigure}
Si sono riformati i gruppi come nel primo incontro e, nell'unico
gruppo da 4, si è spontaneamente aggiunto il bambino assente negli
altri incontri. Autonomamente si sono ridistribuiti i ruoli.

La parte dell'attività riguardante le tessere è stata svolta
abbastanza facilmente, un gruppo ha provato a ``cancellare'' con una
riga 2 tessere e a ridisegnarle in un'altra posizione, altri hanno
discusso, ma dato la risposta esatta senza nessuna prova pratica. Più
problematico si è dimostrato il rispondere alle domande riguardanti i
mattoni, alcuni hanno disegnato, altri hanno misurato in modo
arbitrario, altri hanno misurato con il righello per quantificare i
mattoni. Si è dimostrata una difficoltà il conteggio delle facce
esterne. Solo due gruppi hanno dato le risposte corrette.

A questo punto abbiamo distribuito i mattoni, i bambini hanno
verificato praticamente il lavoro e dove avevano sbagliato si sono
autocorretti.

Poi il relatore di ogni gruppo ha esposto il lavoro svolto. Tutti
hanno espresso soddisfazione per aver lavorato in gruppo e per
l'attività, hanno dimostrato di aver capito i concetti proposti.

Questi sono stati alcuni commenti:
\begin{studente}[ ]
  \begin{itemize}
\item \bambini{Ho imparato a contare da furbo senza facce disegnate}
\item \bambini{Ho imparato che per collaborare non dobbiamo litigare}
\item \bambini{Abbiamo capito che anche figure diverse possono occupare lo stesso spazio}
\item \bambini{Ho capito bene cos'è lo spazio e ho imparato a
    disegnare lo spessore per lo spazio alto}
\item \bambini{Ho capito che con lo stesso spazio ho contorni più
    lunghi o più corti}
\item \bambini{Ho capito che mettendo le tessere vicine faccio la
    figura con il contorno più corto, se avvicino poco le tessere il
    contorno diventa più lungo}
\end{itemize}
\end{studente}


\section[Sperimentazione \#10: terza primaria]{Sperimentazione \#10:
  classe terza primaria, febbraio/marzo~2010}

\subsubsection{Presentazione classe}
La classe è formata da 21 alunni, di cui 1 disabile, 3 dislessici
certificati, 1 con disturbi nell'apprendimento, 1 con disturbo
nell'abilità scolastiche, 2 stranieri.

\subsubsection{Composizione gruppi}
Gruppi omogenei fra loro e eterogenei al loro interno, stabiliti
dall'insegnante per dare la possibilità agli alunni in difficoltà di
ricevere stimoli durante tutto il percorso. Nei gruppi c'è mescolanza
di maschi e femmine.Tre gruppi sono formati da 5 alunni e 1 gruppo da
6. I gruppi resteranno gli stessi per le cinque esperienze. All'interno
di ciascun gruppo è stato nominato un capogruppo per ogni singola
esperienza con compiti di coordinare, suddividere i seguenti ruoli:
\begin{enumerate}
  \item leggere i quesiti;
  \item compilare i quesiti;
  \item disegnare le figure;
  \item fare domande all'insegnante;
  \item riferire alla classe nel momento collettivo.
\end{enumerate}

\subsubsection{Insegnanti presenti}
In quasi tutti gli incontri l'insegnante sperimentatrice è affiancata
dall'altra insegnante di classe (compresenza).

\subsubsection{Calendarizzazione degli incontri}
\begin{calendario}
  \begin{itemize}
  \item 15 febbraio (mattino, compresenza)
  \item 15 febbraio (pomeriggio, compresenza)
  \item 17 febbraio (mattino, compresenza)
  \item 22 febbraio (mattino, compresenza)
  \item 23 febbraio (pomeriggio)
  \item 24 febbraio (mattino, compresenza)
  \end{itemize}
\end{calendario}

\subsection{Primo incontro}

\begin{description}
\item[Alunni presenti:] 19 alunni (2 assenti)
\item[Tempo effettivo di lavoro:] 1 h e 30 minuti
\end{description}

\begin{consegna}
  Ho avviato le attività presentando la sperimentazione alla
  classe. Ho spiegato che avrebbero partecipato a un laboratorio di
  matematica con l'opportunità di acquisire alcuni concetti di
  geometria attraverso l'utilizzo di due nuovi materiali. Ho
  presentato in modo dettagliato le modalità di lavoro. Ho precisato
  che le attività sarebbero state cinque per la durata di circa due
  ore. Ho predisposto in aula per ciascun gruppo le postazioni di
  lavoro, il materiale e le schede riguardanti le esperienze
  predisposte precedentemente seguendo le indicazioni riportate nel
  testo ``Torri, serpenti e\dots{} geometria''. Dopo di che li ho
  suddivisi in gruppi e ho consegnato a ciascun gruppo la scheda
  \attivita{Per cominciare} che hanno letto e completato
  autonomamente. Solo due gruppi, però, l'hanno compilata
  correttamente. Subito dopo ho distribuito le tessere triangolari (12
  per ciascun gruppo). Gli alunni hanno familiarizzato con il nuovo
  materiale divertendosi a costruire liberamente figure. Poi, seguendo
  le mie istruzioni, hanno costruito almeno quattro figure e a ciascuna
  di esse hanno dato un nome. Una fra queste l'hanno scelta come
  mascotte al fine di identificarli e rappresentarli. Infine hanno
  ripassato il contorno delle figure ottenute su carta bianca e le
  hanno riprodotte su scheda predisposta (pag. 20).
\end{consegna}

\subsubsection{Osservazioni}
Durante la presentazione del lavoro, gli alunni hanno dimostrato
subito interesse e curiosità. Di norma sono bambini che esternano
sempre molto entusiasmo nello sperimentare nuove attività e a lavorare
in gruppo. Nella fase di manipolazione del materiale si sono divertiti
molto, ma hanno subito richiesto di avere a disposizione più
tessere. Due gruppi su quattro hanno collaborato molto fra di loro,
hanno provato, discusso e si sono divertiti. Negli altri due è stato
richiesto spesso l'intervento dell'insegnante per cercare di stimolare
chi non partecipava o per risolvere controversie. La bambina disabile
è stata stimolata a partecipare più volte sia dai compagni, sia dalle
insegnanti. Alla fine dell'attività, invitati dall'insegnante, ciascun
gruppo ha preso visione delle figure costruite dagli altri compagni e
sono emerse diverse osservazioni:
\begin{enumerate}
\item confronto numero tessere utilizzate per ciascuna figura
\item numero lati di ciascuna figura e strategie per contarli
\end{enumerate}

\subsubsection{Consigli per i colleghi che vogliono proporre le stesse
  attività}
È stato vantaggioso formare pochi gruppi anche se più numerosi al
loro interno; è stata comunque indispensabile la presenza di due
insegnanti durante lo svolgimento dell'attività.

\subsection{Secondo incontro}
\begin{description}
\item[Alunni presenti:] 19 alunni (due assenti)
\item[Tempo effettivo di lavoro:] 2 ore circa
\end{description}

\begin{consegna}
  Ho avviato l'attività distribuendo 16 tessere triangolari a ciascun
  gruppo. Ho spiegato che dovevano formare figure con minor/maggior
  contorno ma lo spazio occupato sul piano doveva essere lo stesso per
  ciascuna figura. Hanno disegnato la sagoma delle figure ottenute e
  le hanno riprodotte su scheda predisposta (pag. 20).
\end{consegna}

\subsubsection{Osservazioni}
I bambini divisi in gruppi si sono ridistribuiti i ruoli e hanno
cominciato a lavorare con interesse e impegno. È stato interessante
vedere come si sono organizzati subito nel lavoro. Due gruppi su
quattro hanno diviso a metà le tessere e ha costruito due figure con
stessa superficie e diverso contorno. Un gruppo non ha utilizzato
tutte le tessere ma ha rispettato comunque la consegna data. Solo un
gruppo ha costruito con le tessere date due figure con stessa
superficie e contorno. Gratificante e stimolante è stato il momento
collettivo: quando ciascun gruppo ha preso visione dei lavori
effettuati dagli altri compagni e ha condiviso e confrontato la
costruzioni ottenute. Dal confronto sono emersi i criteri utilizzati
per la misurazione delle figure, diversi hanno sottolineato che si può
misurare il contorno usando il lato delle tessere. A questo punto il
gruppo che aveva costruito due figure con uguale contorno ha chiesto
di poterne modificare una. Dopo vari tentativi è riuscito a costruire
la figura giusta spiegando che bisognava creare solo più
\bambini{insenature}. La bambina disabile ha faticato a seguire i
ragionamenti dei compagni e, nonostante l'aiuto sia dei compagni, sia
dell'insegnante, non è riuscita a apportare il proprio contributo.

\subsection{Terzo incontro}

\begin{description}
\item[Alunni presenti:] 19 alunni (2 assenti)
\end{description}

\begin{consegna}
  Ho avviato l'attività facendola precedere da un momento di
  riflessione collettiva per rivedere tutti insieme i risultati
  raggiunti negli incontri precedenti. Ho poi consegnato a ciascun
  gruppo le due schede di lavoro (pagg. 54 - 55) e ho distribuito le
  tessere necessarie per la costruzione di entrambe le figure
  rappresentate sulla prima scheda (pag. 54). Infine ho messo a
  disposizione a ciascun gruppo due corde della misura giusta. Ogni
  gruppo ha letto e completato autonomamente le due schede.
\end{consegna}

\subsubsection{Osservazioni}
Parte del tempo a disposizione è stato utilizzato per riprendere a
livello collettivo le esperienze precedenti. È stato utile soprattutto
per gli alunni assenti il 15 febbraio perché sono riusciti attraverso
anche il supporto delle schede su cui avevano lavorato i propri
compagni a recuperare il percorso fin qui svolto. Tutti i gruppi hanno
completato le schede senza incontrare particolari difficoltà anche se
è stato necessario puntualizzare durante l'attività la differenza tra
spazio occupato sul tavolo e bordo da circondare. Importantissimo è si
rivelato il momento conclusivo durante il quale è avvenuto il
confronto delle risposte date e la verifica dei risultati ottenuti.

\subsubsection{Consigli per i colleghi che vogliono proporre le stesse
  attività}
Il tempo a disposizione non è stato sufficiente
per effettuare le esperienze 4 e 5 (pagg.~56-57).

\subsection{Quarto incontro}
\begin{description}
\item[Alunni presenti:] 20 alunni (1 assente)
\item[Tempo effettivo di lavoro:] due ore
\end{description}

\begin{consegna}
  Ho avviato il lavoro consegnando a ciascun gruppo la scheda
  \attivita{Per Cominciare} (pag. 53). Gli alunni l'hanno letta e
  compilata autonomamente. Solo due gruppi l'hanno completata
  correttamente. Poi ho distribuito il nuovo materiale (8 mattoni a
  ciascun gruppo) invitando tutti a costruire figure
  liberamente. Hanno poi realizzato e disegnato 3 costruzioni diverse
  sulle schede predisposte (pagg. 58-59)
\end{consegna}

\subsubsection{Osservazioni}
Nella fase di manipolazione del materiale gli alunni hanno dimostrato
interesse e entusiasmo. Quando li ho invitati a contare le facce
esterne delle loro costruzioni, 3 gruppi su quattro non hanno
incontrato difficoltà nel riconoscere le facce di ogni costruzione e
guidati sono riusciti a trovare le strategie per contarle. Durante il
lavoro, un gruppo ha fatto notare che non bisognava contare le parti
che combaciavano. Prima di eseguire le esperienze 3 e 4 delle schede
(pagg. 60-61) è stato necessario creare un momento collettivo utile
soprattutto per il gruppo che non era riuscito a riconoscere e a
contare le facce delle proprie costruzioni, rendendosi così
consapevole dei propri errori.

\subsubsection{Consigli per i colleghi che vogliono proporre le stesse
  attività}
Il tempo preventivato per effettuare l'intera consegna è risultato
insufficiente. Durante il momento collettivo è stato utile rivedere il
lavoro sui solidi presentato a novembre%
.

\subsection{Quinto incontro}
\begin{description}
\item[Alunni presenti:] 20 alunni (1 assente)
\item[Tempo effettivo di lavoro:] 2 ore
\end{description}

\begin{consegna}
  Dopo aver realizzato e disegnato con una quantità definita di
  materiali (8 mattoni) le tre costruzioni rispettando la consegna
  delle schede pagg. 58-59 (vedi consegna quarto incontro, prima
  parte), gli alunni hanno costruito il dado e lo hanno confrontato
  con le costruzioni precedenti considerando il numero delle facce
  esterne e lo spazio occupato da ciascuna costruzione. Infine hanno
  provato a realizzare una costruzione diversa dal dado ma che avesse
  lo stesso numero di facce (schede pagg. 60-61).
\end{consegna}

\subsubsection{Osservazioni}
Questa seconda parte del lavoro li ha entusiasmati
tantissimo. L'impegno, l'attenzione, la voglia di riuscire a
realizzare quanto richiesto li ha sicuramente aiutati a raggiungere il
risultato finale. Tre gruppi su quattro sono riusciti senza alcuna
difficoltà e autonomamente a effettuare le esperienze proposte nelle
pagine 60 e 61. Solo in un gruppo le dinamiche relazionali e la
lentezza nell'esecuzione dei lavori hanno inciso negativamente sul
percorso di lavoro. È prevalsa comunque la voglia di riuscire e, anche
se in tempi di esecuzione diversi, anche i bambini di questo gruppo,
il mattino seguente, prima di iniziare l'ultimo incontro, hanno
chiesto di poter provare ancora a realizzare la costruzione su cui
avevano faticato il pomeriggio precedente ottenendo ottimi risultati.

\subsection{Sesto incontro}
\begin{description}
\item[Alunni presenti:] 21 alunni
\item[Tempo effettivo di lavoro:] 1 ora e trenta minuti
\end{description}

\begin{consegna}
  Schede \attivita{Per concludere}.
\end{consegna}

\subsubsection{Osservazioni}
Gli alunni hanno letto e completato autonomamente le due schede. Tutti
e quattro i gruppi sono riusciti a compilare correttamente la scheda
per concludere sulle tessere; mentre solo due gruppi su quattro sono
riusciti a trovare da soli le strategie per contare i mattoni sia
della torre verde, sia della torre gialla e quindi a eseguire
correttamente anche la seconda scheda (\attivita{Per concludere sui
  mattoni}). A mio parere, in effetti, non era semplice calcolare il
numero di mattoni contenuti in ciascuna torre anche a causa di
problemi legati alla stampa delle due rappresentazioni.


\section[Sperimentazione \#11: terza primaria]{Sperimentazione \#11:
  classe terza primaria, marzo~2010}

\subsection{Osservazioni generali}

\subsubsection{Presentazione della classe}
La classe è formata da 25 alunni.

\subsubsection{Composizione dei gruppi}
All'inizio del primo incontro ho formato 5 gruppi eterogenei (5
bambini in ogni gruppo). Ho formato personalmente i gruppi, perché
temevo che venisse impiegato troppo tempo e poi perché ultimamente 3
bambini risultano ``poco graditi'' al gruppo classe e sovente vengono
``scelti'' per ultimi. Questi gruppi sono rimasti invariati nell'arco
dei 4 incontri, modificandosi solo in base alle assenze.

\subsubsection{Insegnanti presenti}

Agli incontri è presente solo l'insegnante sperimentatrice.

\subsubsection{Calendarizzazione degli incontri}
\begin{calendario}
  \begin{itemize}
  \item 11 marzo
  \item 17 marzo
  \item 22 marzo
  \item 31 marzo
  \end{itemize}
\end{calendario}

\subsection{Primo incontro}

\begin{description}
\item[Alunni presenti:]22 bambini presenti, 3 assenti
\item[Tempo effettivo di lavoro:] 1 ora e 1/2 (11/13.30)
\end{description}

\begin{consegna}
  scheda \attivita{Per cominciare} di seconda e di terza.

  \materiali{}%
  Scheda \attivita{Per cominciare} di seconda, messa a disposizione
  sui tavoli.
\end{consegna}

\subsubsection{Osservazioni}
Siccome abbiamo lavorato poco per piccoli gruppi, prima di dare il via
all'attività ho cercato, con l'aiuto dei bambini, di riportare alla
luce le modalità del lavoro di gruppo.

Dal momento che c'era una scheda da compilare i ruoli individuati sono
stati:
\begin{itemize}
  \item lettore
  \item scrittore
  \item relatore-comunicatore
  \item postino
\end{itemize}
La difficoltà dell'interiorizzare i ruoli è emersa fin da subito:
tutti i componenti del gruppo hanno preso la matita per scrivere; al
momento della consegna tutto il gruppo si alzava, per recarsi da me e
lo stesso avveniva durante il prelievo e il ritiro del materiale.

L'attribuzione dei ruoli è avvenuta in merito alle attitudini dei
singoli, riconosciute dal gruppo, e solo in pochi casi è stata
affidata alla casualità. In un gruppo, nel quale è presente un bambino
``caratteriale'' che vuole essere leader, ma non è riconosciuto come
tale dai compagni, hanno fatto ``pari e dispari'' per decidere
l'assegnazione degli incarichi. In questo gruppo un bambino,
introverso e poco attivo, ha fatto sentire la sua voce nel tentativo
di mettere a tacere il presunto leader.

La condivisione vera delle risposte ai quesiti, prima delle scrittura
delle stesse sulla scheda non è avvenuta:
\begin{studente}[ ]
  \begin{itemize}
\item \bambini{Non è colpa mia; ho scritto quello che mi hai detto tu}
\item \bambini{Te l'avevo detto io!}
\item \bambini{Vuoi sempre avere ragione e hai sbagliato}
\item \bambini{È colpa sua se abbiamo sbagliato}
\end{itemize}
\end{studente}
I ruoli del ``bravo'', piuttosto che del ``tontolone'' o del
``prepotente'' hanno influenzato l'andamento dei lavori, nel senso che
i bambini mi sono sembrati ben disposti a accettare le risposte dei
``bravi'', (che il più delle volte non venivano nemmeno messe in
discussione), mentre mi è parso che non ascoltassero i consigli di
coloro che sono considerati poco capaci.

In un gruppo un bambino ha fatto osservare ai compagni che l'immagine
della stella che stavano riproducendo era sbagliata, perché
\bambini{la riga che fanno i triangolini sulla scheda è dritta e la
  vostra è storta}. L'osservazione corretta, sebbene articolata poco
chiaramente, è stata completamente ignorata dai compagni, due dei
quali, in ambito matematico, hanno prestazioni migliori rispetto al
bambino che era intervenuto.

Fase~1 (scheda \attivita{Per cominciare} di seconda):
Consegna: leggere la scheda, compilarla in gruppo e di riconsegnarla.

La compilazione ha richiesto pochi minuti: giusto il tempo di leggere
e di scrivere i 2 numeri richiesti.

Quando ho visto che per 4 gruppi la torre era formata da 4 cubi
(anziché 16) e per un gruppo era formata da 5 cubi, ho pensato di
restituire la scheda ai bimbi e di chiedere che verificassero le loro
risposte, utilizzando il materiale (che a quel punto ho mostrato e
distribuito in base alla quantità da loro richiesta).

Inutile dire che c'è stata la corsa al banco distribuzione.

Con il materiale alla mano, i bambini hanno provato a riprodurre
l'immagine e successivamente hanno modificato le loro risposte, ma
solo due gruppi sono arrivati alla soluzione esatta.

Uno di questi, usando il righello, ha misurato l'altezza di un cubo
sulla scheda e poi l'ha riprodotta lungo l'altezza della torre.

Prima di farmi restituire le schede ho chiesto ai bambini di spiegare
ai compagni le tecniche utilizzate per verificare l'ipotesi. Quattro
gruppi hanno detto di aver guardato la scheda e rifatto l'immagine
uguale uguale, mentre ``il gruppo del righello'' ha illustrato la
strategia utilizzata.

Dato che questa era l'unica proposta individuata, tutti i gruppi hanno
provato a applicarla trovando così il numero corretto dei mattoni
impiegati. Da qui la condivisione del risultato.

Fase~2 (scheda \attivita{Per cominciare} di terza):
Consegna:
\begin{enumerate}
\item leggere e compilare la scheda in ogni sua parte;
\item verificare le risposte attraverso l'utilizzo del materiale a
  disposizione (cubi, triangoli e corde)
\end{enumerate}
Ecco le risposte scritte sulle schede secondo la prima richiesta
\begin{center}
  \begin{tabular}{|l|c|l|c|c|l|}
    \hline
    &
    \multicolumn{1}{p{2.3cm}|}{\raggedright{}Quante tessere nella
      stella \par (Risposta: 12)}
    &
    \multicolumn{1}{p{2.3cm}|}{\raggedright{}Contorno più lungo\par
      (Risposta: uguali)}
    &
    \multicolumn{1}{p{2.3cm}|}{\raggedright{}Mattoni nella torre rossa\par
      (Risposta: 14)}
    &
    \multicolumn{1}{p{2.3cm}|}{\raggedright{}Mattoni nella torre blu\par
      (Risposta: 12)}
    &
    \multicolumn{1}{p{2.3cm}|}{\raggedright{}Torre che occupa più spazio\par
      (Risposta: rossa)}
    \\ \hline
    Gruppo 1&
    12&
    Diamante&
    14&
    16&
    La blu\\ \hline
    Gruppo 2&
    6&
    La lumaca&
    21&
    10&
    La rossa\\ \hline
    Gruppo 3&
    6&
    Diamante&
    21&
    12&
    La blu\\ \hline
    Gruppo 4&
    6&
    La lumaca&
    15&
    12&
    La rossa\\ \hline
    Gruppo 5&
    6&
    La lumaca&
    14&
    9&
    La rossa\\ \hline
  \end{tabular}
\end{center}

Ecco le risposte dopo aver effettuato la verifica con il materiale
(sono evidenziate le risposte che sono state modificate rispetto alla
tabella precedente)
\begin{center}
  \begin{tabular}{|l|c|l|c|c|l|}
    \hline
    &
    \multicolumn{1}{p{2.3cm}|}{\raggedright{}Quante tessere nella stella \par
      (Risposta: 12)}
    &
    \multicolumn{1}{p{2.3cm}|}{\raggedright{}Contorno più lungo \par
      (Risposta: uguali)}
    &
    \multicolumn{1}{p{2.3cm}|}{\raggedright{}Mattoni nella torre rossa \par
      (Risposta: 14)}
    &
    \multicolumn{1}{p{2.3cm}|}{\raggedright{}Mattoni nella torre blu
      (Risposta: 12)}
    &
    \multicolumn{1}{p{2.3cm}|}{\raggedright{}Torre che occupa più spazio \par
      (Risposta: rossa)}
    \\ \hline
    Gruppo 1&
    12&
    Diamante&
    14&
    \hl{12}&
    \hl{La rossa}\\ \hline
    Gruppo 2&
    \hl{12}&
    \hl{Diamante}&
    \hl{14}&
    \hl{8}&
    La rossa\\ \hline
    Gruppo 3&
    \hl{12}&
    Diamante&
    \hl{14}&
    12&
    La blu\\ \hline
    Gruppo 4&
    \hl{12}&
    Lumaca&
    \hl{14}&
    12&
    La rossa\\ \hline
    Gruppo 5&
    \hl{12}&
    Lumaca&
    14&
    \hl{13}&
    \hl{La blu}\\ \hline
  \end{tabular}
\end{center}
Nel cercare di riprodurre la stella ho notato che i bambini creavano
solo il contorno (6 triangoli) senza riempirne l'interno (per formare
il quale sarebbero serviti altri 6 triangoli).

Sono rimasta senza parole nell'osservare come è risultato difficile
conteggiare il numero dei cubi sulla scheda (anche da parte dei
``bravi'') e confrontare lo spazio occupato dalle due
torri%
.

I gruppi convinti che la torre che occupa più spazio è quella blu
hanno giustificato la loro affermazione dicendo che è più grande
sotto.
Un'altra frase, sentita da lontano, che mi ha interrogata è:
\bambini{Il contorno è quello che occupa più spazio}.
Per effettuare la verifica del contorno più lungo, avevo messo a
disposizione per ogni gruppo due corde, ma\dots{} non sono state
utilizzate. \underline{\textbf{Dove suggerire io il possibile
    utilizzo?}}

Il calcolo della lunghezza del contorno è stato tradotto in conteggio
dei triangoli che formano la figura.

\medskip{}
In entrambe le fasi di lavoro,
quasi tutti i presenti hanno partecipato attivamente al lavoro: solo 3
bambini hanno vissuto l'esperienza da osservatori, rispetto
all'attività didattica vera e propria, ma hanno comunque rivestito e
svolto il ruolo che era stato loro assegnato nel gruppo.

Mi è sembrato che tutti si siano divertiti. Quando ho chiesto di
ritirare molti mi hanno chiesto come mai, non essendosi resi conto che
era già giunta l'ora di andare a pranzo!

Rispetto alla gratificazione\dots{} non ho avuto modo di compiere
questa osservazione. Le verifiche effettuate hanno smentito talvolta
le ipotesi formulate in partenza, quindi penso che i bambini (più
abituati a confrontarsi sui risultati che non sui processi) non si
siano sentiti gratificati.

Il clima nei gruppi è stato complessivamente positivo.

In un gruppo, il desiderio di affermarsi di un componente ha creato
tensione, al punto che questo ha deciso di mettere la testa sul banco
per una decina di minuti, nel tentativo di essere richiamato in gioco
dai compagni, che al termine dell'attività mostravano una crescente
intolleranza.

In un altro gruppo, invece, il desiderio di leadership di un bimbo è
stato messo a tacere dall'unione degli altri 3 componenti che lo hanno
più volte sollecitato a ascoltare anche la loro versione e a
rispettare i ruoli.

Presa visione delle risposte, mi auguro che le abbiano ``sparate''
subito, nel desiderio di arrivare primi alla soluzione.

La strategia, adottata da diversi bambini, di disegnare sulla scheda i
contorni delle figure per facilitarne il conteggio non si è rivelata
per tutti efficace, perché i triangoli (o i cubi) che alcuni bambini
hanno disegnato non corrispondevano a quelli del modello, quindi il
conteggio è risultato comunque sbagliato.

Un gruppo, quello che dopo la prima scheda aveva proposto l'utilizzo
del righello, si è accanito nell'usare il suddetto strumento non
cercando strategie alternative, forse più idonee.

Io\dots{} mi sono letteralmente messa le mani nei capelli: ma quanta
confusione hanno in testa o meglio, quanta confusione sono riuscita a
creare nelle loro teste?!?

\subsection{Secondo incontro}

\begin{description}
\item[Alunni presenti:] 20 bambini presenti, 5 assenti
\item[Tempo effettivo di lavoro:] circa 2 ore 13.40/14.30 -
  15.30/16.30
\end{description}

\begin{consegna}
  5 gruppi eterogenei (4 bambini in ogni gruppo). I gruppi non sono
  stati modificati rispetto al primo incontro, anche perché all'invito
  ``Mettiamoci in gruppo per cominciare l'attività'' i bambini si sono
  raggruppati dando per scontata la formazione proposta la settimana
  precedente.

  \materiali{}%
  Scheda \attivita{Tessere} -- 2 corde di lunghezze diverse a
  disposizione sui tavoli
\end{consegna}

\subsubsection{Osservazioni}
In tutti i gruppi i bambini hanno rimesso subito ``le mani in pasta''
con entusiasmo; in alcuni gruppi i ruoli assegnati durante il primo
incontro sono stati modificati.

Mi è parso che i bambini si ``accusassero'' meno e qualcuno si è
sforzato di non imporre il proprio pensiero.

Inizialmente i due gruppi, dove sono stati inseriti i bambini che
erano assenti durante il primo incontro, hanno lavorato con maggiore
lentezza; uno dei due poi ha assunto ritmi ``regolari'', mentre
l'altro ha eseguito l'attività con estrema rapidità, ma con scarsa
precisione.

Il ``nuovo'' mi ha comunicato la sua gioia nel fare \bambini{geometria
  in questo modo divertente}, ma più volte ho notato che giocava con
le tessere da solo, anche quando i compagni stavano conteggiando,
piuttosto che misurando.

Fase~1:
Osservazione e riproduzione di serpente e stella; conteggio delle aree
(n. tessere) e confronto dei perimetri (con utilizzo di corde), ovvero
attività di pag.~3.

La costruzione del serpente e della stella e il successivo conteggio
delle tessere ha richiesto pochi minuti, mentre difficoltosa è
risultata la misurazione del perimetro con l'utilizzo delle corde che
erano state consegnate.

Le difficoltà emerse erano di due tipi:
\begin{itemize}
\item le tessere si spostavano molto facilmente e le 8 mani che i
  bambini avevano a disposizione risultavano insufficienti a
  bloccarle. Un bambino ha pensato di fissare le tessere al tavolo
  mettendovi sopra il suo astuccio, ma essendo troppo ingombrante gli
  ostacolava l'attività di misurazione, quindi l'ha tolto non
  riuscendo però comunque a raggiungere il suo scopo;
\item la corda slittava dalla parete delle tessere (bassa e non
  sufficientemente rigida) e finiva sotto le stesse, rendendo
  scorretta la misurazione.
\end{itemize}
Le punte della stella hanno messa a dura prova la tenacia dei bambini,
alcuni dei quali si sono comodamente seduti, aspettando che i compagni
fornissero una risposta anche a nome loro.

Un gruppo ha risposto che i contorni del serpente e della stella erano
uguali e quando hanno dovuto motivare la loro affermazione hanno
detto:
\begin{studente}[ ]
  La volta scorsa abbiamo fatto un esercizio quasi uguale e la
  risposta giusta era questa!
\end{studente}
Ognuno dei cinque ``relatori'', al termine della compilazione della
prima pagina dell'attività, ha spiegato ai compagni della classe le
tecniche adottate dal gruppo per costruire le figure proposte, per
effettuare il conteggio delle tessere e per misurare il perimetro.

Fase~2:
Costruzione e osservazione di 3 figure (clessidra, lumaca, gatto),
misurazione di aree (n. tessere), di perimetri (n. lati delle tessere
contenuti nel contorno) e confronto (attività di pag.~4).

\underline{Sulle attività, a pagina 4, manca l'immagine del
  gatto}

Gli interrogativi emersi durante la compilazione di questa scheda sono
stati:
\begin{itemize}
\item \bambini{Cosa vuol dire numero dei lati di tessera contenuti nel
    contorno?}
\item \bambini{Per figure che abbiamo costruito intendi anche quelle
    che abbiamo fatto nella prima pagina?}
\end{itemize}
Tutti i gruppi hanno compilato la tabella relativa al perimetro,
conteggiando i lati della figura e non il numero dei lati delle
tessere contenuti nel contorno.

Ho sentito che una bambina si è interrogata sul senso della consegna,
ma non mi ha coinvolto perché è stata ``convinta'' dai compagni che
\begin{studente}[]
  dire lati della figura e lati delle tessere era la stessa cosa
\end{studente}

Nella fase del confronto collettivo, a seguito di alcune mie domande,
che avevano lo scopo di far riflettere i bambini sulla consegna, mi è
stata chiesta la possibilità di rivedere le risposte scritte sulla
scheda \bambini{altrimenti è tutto sbagliato!}

Chiaramente la risposta è stata affermativa.

Un gruppo ha persistito nell'errore (prima aveva conteggiato tutti i
lati delle tessere, non solo quelli del contorno, poi ha conteggiato i
lati della figura); gli altri hanno rettificato correttamente le loro
risposte nella tabella, ma due non hanno effettuato la correzione
nell'esercizio sottostante che richiedeva lo stesso tipo di conteggio.

Qualche imprecisione è emersa anche nel confronto degli spazi occupati
(``quale figura occupa più spazio sul tavolo?''): 2 gruppi hanno
risposto correttamente; 2 gruppi hanno dato una risposta parziale
(hanno citato solo il serpente); 1 gruppo ha scritto ``gatto''.

Nel confronto tra i contorni delle figure, 2 gruppi hanno risposto
correttamente (per uno di questi non si sa bene come perché una
risposta non era coerente con i dati raccolti nell'esercizio
precedente) mentre gli altri tre hanno risposto in maniera errata, uno
perché ha considerato solo le figure presentate nella pagina; gli
altri due perché hanno confuso ``tessere'' con ``lati delle tessere''.

Le strategie condivise al termine di questa seconda esperienza sono
state:
\begin{itemize}
\item per contare i lati delle tessere senza dimenticarne alcuno è
  bene segnarli;
\item prima di procedere nel confronto è necessario accertarsi di cosa
  si deve confrontare;
\item i dati raccolti, possono essere utilizzati (per evitare di fare
  più volte lo stesso lavoro).
\end{itemize}
Conclusioni a cui si è giunti:
\begin{itemize}
\item contorno e spazio sono due ``cose'' diverse;
\item con lo stesso numero di tessere si possono fare figure diverse.
\end{itemize}

La conclusione che ci sono figure con il perimetro della stessa
lunghezza che sono formate da un numero di tessere diverso non è
giunta spontanea.

\medskip{}
In questo secondo incontro,
come prevedevo, l'attività è stata ripresa con entusiasmo, da parte di
tutti.

Prima di proporre il secondo step è emersa la necessità di raccontare
il vissuto del primo incontro, in quanto alcuni alunni non erano
presenti.

Sono stati ricordati e riassegnati i ruoli tra i vari componenti del
gruppo (il fatto che ci fossero 4 alunni per gruppo è stata una
coincidenza ottimale) e si è cercato di riportare alla memoria le
strategie individuate durante il primo incontro, perché diventassero
patrimonio comune.

Con mia grande sorpresa la strategia del misurare (che aveva
condizionato il lavoro di un gruppo per tutto il primo incontro) non è
stata ricordata nemmeno dal bambino che l'aveva proposta; tutti
avevano ben presente di avere sbagliato a individuare il numero dei
cubi utilizzati per costruire la torre, ma nessuno riusciva a spiegare
la tecnica utilizzata per verificare le ipotesi numeriche che erano
state scritte.

Un bimbo ha raccontato che nel suo gruppo il guardare bene la
posizione delle linee dei triangoli prima di ``disegnare'' la figura
richiesta è stato utile per riprodurre correttamente l'immagine e un
altro ha invitato a non toccare tutti insieme le tessere, altrimenti
uno costruisce e l'altro distrugge e non si conclude niente! Il clima
nei gruppi è stato positivo, anche nel gruppo di L (chiamiamo così il
bambino che desidererebbe essere riconosciuto come leader).

La presenza costante del materiale sul tavolo favorisce
l'``isolamento'' dei bambini che faticano a inserirsi nel gruppo o di
quelli che hanno un tempo di concentrazione limitato.

I bambini avrebbero voluto costruire contemporaneamente tutte le
figure ``per vederle meglio'', ma la quantità di tessere a
disposizione non ha consentito tale verifica. La presenza di
un'insegnante sola non consente un'osservazione puntale delle
dinamiche e tanto meno un intervento preciso e rispondente ai bisogni
dei singoli gruppi.

\subsection{Terzo incontro}

\begin{description}
\item[Alunni presenti:] 21 bambini presenti, 4 assenti
\item[Tempo effettivo di lavoro:] circa 2 ore 13.45/15.30
\end{description}

\begin{consegna}
  5 gruppi eterogenei (2 gruppi da 5; 2 gruppi da 4; 1 da 3)

  Non ho ritenuto opportuno spostare un bambino del gruppo dei 5 nel
  gruppo dei 3, perché non è emersa come necessità

  \materiali{}%
  Scheda \attivita{Tessere} a disposizione sui tavoli.
\end{consegna}

\subsubsection{Osservazioni}

I bambini si sono suddivisi rapidamente nei rispettivi gruppi, uguali
per formazione e anche per collocazione all'interno dello
spazio-aula.

Fase~1:
Realizzare una costruzione con contorno uguale al serpente (14 linee),
ovvero attività di pag.~5.

Prima di cominciare l'attività ho ricordato ai bambini che dovevano
provvedere a leggere la consegna e a eseguire il compito con
autonomia, ma che non dovevano farsi remore a chiedere il mio
intervento, qualora all'interno del gruppo non fossero riusciti a
trovare un accordo sul da farsi.

Quattro gruppi su 5 hanno tradotto il termine ``realizzare'' con
costruire, quindi dopo aver letto la consegna si sono diretti verso la
cattedra per prendere il numero di tessere concordato all'interno del
gruppo.

Un gruppo ha deciso di conteggiare il numero dei lati delle tessere
direttamente sulla griglia stampata sul foglio (utilizzando la
simbologia adottata per conteggiare i contorni del serpente e della
stella) e solo successivamente ha recuperato il numero esatto delle
tessere di cui aveva bisogno per realizzare la costruzione inventata.

Un bambino avendo osservato che il proprio gruppo aveva realizzato una
costruzione piuttosto grande ha detto:
\begin{studente}[ ]
  Abbiamo fatto una figura grandissima che ha il contorno come il
  serpente. Ora disegniamola!
\end{studente}
ma subito dopo ha aggiunto
\begin{studente}[ ]
  Proviamo a controllare prima
\end{studente}
(Il bambino apparteneva al gruppo che nell'incontro precedente aveva
lavorato in modo molto impreciso).

Nello stesso gruppo, lo ``scrittore'', disegnando la costruzione
realizzata collettivamente, ha raddoppiato la grandezza delle tessere
triangolari; i membri del gruppo hanno scoperto l'errore quando,
ricontando i lati delle tessere, si sono resi conto che era maggiore
di 14 e hanno prontamente provveduto a rettificare lo scritto.

Fase~2:
Costruzione di una figura con il contorno più corto e più lungo
possibile, usando uno stesso numero di tessere (attività di pag.~6).

Non sono state date regole prima di procedere alla costruzione delle
figure, ma è nata la necessità di stabilirle durante la fase di
confronto, in quanto alcuni bambini non riconoscevano come figure
valide quelle costruite da altri gruppi.

I bambini hanno così concordato il significato della parola ``figura''
(le tessere che la costituivano avrebbero dovuto avere almeno un lato
in comune con quella vicina), dopo di che hanno proposto di rifare
l'esercizio alla luce della regola stabilita.

Al termine dell'attività i gruppi hanno comunicato agli altri la
lunghezza del contorno delle proprie figure (a cui avevano assegnato
un nome) e chi si è accorto di non aver trovato la figura avente le
caratteristiche richieste ha provveduto a rettificare la propria
costruzione.

Un bambino ha cercato (volontariamente) di svelare il trucco per
trovare il contorno più corto e quello più lungo e ha detto:
\begin{studente}[ ]
  Per costruire una figura con il contorno corto bisogna circondare
  una tessera con le altre e per costruire quella con il contorno
  lungo si deve mettere una tessera di fianco all'altra
\end{studente}
Al termine di questa terza esperienza la quasi totalità dei bambini si
è detta propensa a manipolare il materiale per cercare di risolvere il
problema, piuttosto che a immaginare e a verificare la soluzione con
l'uso di carta e penna.

È stato interessante osservare che di fronte alla presa di coscienza
dell'errore i bambini non manifestavano il desiderio di vedere la
soluzione trovata dai compagni, ma si davano da fare per scovarne
autonomamente una valida%
.

Non si è giunti a conclusioni nuove, ma si sono ribadite quelle
intuite nel precedente incontro:
\begin{itemize}
\item con lo stesso numero di tessere si possono fare figure diverse;
\item con numero di tessere diverse si possono costruire figure aventi
  lo stesso contorno.
\end{itemize}

Per consentire la sperimentazione di queste affermazioni ho proposto
all'interno dei gruppi di giocare con le tessere formando figure equi
estese e/o figure che avessero lo stesso perimetro (sfida a coppie).

\medskip{}
Anche in questo terzo incontro,
come nella ``puntata'' precedente, l'aggiornare i numerosi bambini
assenti ha permesso alla classe di ricordare le ``scoperte'' fatte
nell'incontro numero due.

Stavolta non è stato necessario che io invitassi i bambini a
assegnare/riassegnare i ruoli all'interno del gruppo, in quanto ho
sentito che si organizzavano con autonomia.

Forse non è un caso che alla richiesta di relazionare l'attività
svolta mercoledì 17, si sono fatti avanti due bambini ``relatori'',
cioè due bambini a cui era stato assegnato tale compito.

Si è ricordato:
\begin{itemize}
\item la strategia di segnare i lati, per non dimenticarne alcuno nel
  conteggio;
\item la necessità di comprendere bene la consegna per poterla
  eseguire correttamente;
\item le trappole nelle quali erano incappati i vari gruppi nel
  calcolare il perimetro;
\item la diversità tra il concetto di perimetro e quello di
  superficie;
\item la possibilità di realizzare figure diverse aventi la stessa
  superficie;
\item la possibilità di realizzare figure con uguale perimetro e
  superfici diverse.
\end{itemize}
Il gruppo che nel precedente incontro aveva lavorato frettolosamente,
commettendo numerosi errori, ha lavorato con maggiore
attenzione. Probabilmente il cambiamento è da ricondurre anche al
fatto che uno degli assenti è un bambino piuttosto esuberante, che
influenza l'andamento dell'attività.

Il clima nei gruppi è stato come sempre positivo.

\subsection{Quarto incontro}

\begin{description}
\item[Alunni presenti:]19 bambini presenti, 6 assenti
\item[Tempo effettivo di lavoro:]circa 2 ore 13.45/14.30 - 15.30/16.30
\end{description}

\begin{consegna}
 scheda \attivita{Cubi}

 sempre 5 gruppi eterogenei (2 gruppi da 5; 1 gruppi da 4; 1 da 3; 1 da
 2). Come nell'incontro precedente i bambini non hanno chiesto di
 spostarsi, anche il gruppo di 2 unità non si è posto il problema di
 lavorare in coppia, quindi ho preferito non proporre cambi, anche per
 il fatto che potevo disporre di due ore spezzate (organizzazione un
 po' dispersiva) e i bambini avrebbero dovuto impiegare del tempo per
 cercare l'equilibrio nella nuova situazione.

 \materiali{}%
 scheda \attivita{Cubi} -- 8 cubi per ogni gruppo.
\end{consegna}

\subsubsection{Osservazioni}
Nella distribuzione del materiale, fatta da un paio di bambini nella
pausa gioco, i gruppi sono stati collocati come negli incontri
precedenti.

Nei gruppi di 2, di 3 e di 4 non ho sentito commenti rispetto ai nulla
facenti, mentre in uno dei gruppi da 5 ho osservato che due bambine
tendevano a isolarsi e a parlottare tra loro; nell'altro mi è stato
chiesto di intervenire perché un bimbo non faceva nulla (preciso che
la collega, in servizio la mattina, mi aveva segnalato che quel
bimbetto le sembrava un po' frastornato e ``perso'').

L'attività del disegnare i solidi, per la quale era espressamente
citata sulla scheda la possibilità di rivolgersi all'animatore, ha
richiesto la mia presenza all'interno di tutti i gruppi (durante
quest'anno scolastico abbiamo disegnato solo un paio di volte le
figure solide e il lavoro di arte e immagine svolto dalla collega non
prevedeva di affrontare questo argomento).

Fase~1:
Realizzare 3 diverse costruzioni con 8 mattoni, dare loro un nome e
disegnarle (attività di pag.~7).

Prima di cominciare l'attività, memori dell'esperienza vissuta nel
precedente incontro, i bambini hanno domandato \bambini{Quali
  costruzioni possiamo fare?}, cioè hanno chiesto che fossero
esplicitati i criteri che avrebbero rese accettabili le costruzioni.

Abbiamo concordato che i cubi dovevano essere uniti tra loro tramite
l'affiancamento di una faccia (e non solo dello spigolo o del
vertice).

Un gruppo ha disegnato la costruzione e poi l'ha costruita (come aveva
fatto con le tessere triangolari).

Tutti i gruppi hanno realizzato la prima costruzione sviluppandola
orizzontalmente (non avevo fornito i pioli di legno, riservandomi di
farlo solo su esplicita richiesta) e l'hanno disegnata ``vedendola''
dall'alto.

Ho allora precisato che i cubi potevano anche essere sovrapposti e non
solo affiancati così la fantasia ha spaziato maggiormente.

La rappresentazione su foglio della profondità è stata un po'
problematica, per due ragioni:
\begin{itemize}
\item non riuscivano a capire dove orientarla (i più hanno deciso di
  indicare la profondità verso destra, sul lato destro della
  costruzione, e verso sinistra su quello sinistro);
\item il disegno risultava diverso dalla costruzione inventata (il
  cubo riprodotto sulla griglia quadrata sembrava un
  parallelepipedo). Non mi è sembrato il caso di spiegare come si
  dovesse procedere per effettuare una corretta trasposizione grafica,
  perché ho ritenuto che non fosse un'informazione utile al
  raggiungimento del fine che mi ero proposta: riflettere sul volume e
  sull'estensione della superficie delle figure solide.
\end{itemize}

Fase~2:
Conteggio delle facce esterne dei mattoni che formano ciascuna
costruzione (attività di pag.~9).

Prima di procedere al conteggio delle facce esterne, un gruppo mi ha
chiesto conferma rispetto alla loro interpretazione della consegna.
\begin{studente}[G]
  Le facce sono i lati del cubo?
\end{studente}
mi ha chiesto G.
\begin{tutor}[ ]
  Cosa intendi per i lati del cubo?
\end{tutor}
ho domandato a mia volta
\begin{studente}[G]
  Questo
\end{studente}
e mi ha indicato la faccia.
\begin{tutor}[ ]
  Sì, per facce intende proprio quello che tu mi hai indicato che non
  si chiama lato, perché il lato\dots{}
\end{tutor}
\begin{studente}[R]
  Il lato è del quadrato, non del cubo
\end{studente}
ha concluso R, un bambino dello stesso gruppo
\begin{tutor}[ ]
  Perfetto, allora potete procedere!
\end{tutor}
Passando per i vari gruppi mi sono accorta che il conteggio delle
facce avveniva in modo molto impreciso, allora ho chiesto che ogni
gruppo comunicasse agli altri la strategia adottata.

In verità i bambini non stavano applicando nessun metodo, ma contavano
a ruota libera: qualcuno saltava dall'alto al basso, qualche altro
contava per alcuni cubi tutte le facce e per altri solo quelle
laterali riservandosi di conteggiare in un secondo tempo quelle
superiori e quelle inferiori\dots{}

Avendo verificato che gli errori di conteggio erano molto numerosi, ho
proposto di inventare una scaletta da seguire per ``vincere'' la
distrazione o la fretta. Questo è l'iter che si è deciso di usare per
effettuare la verifica:
\begin{itemize}
\item contare le facce laterali;
\item procedere con le superiori;
\item e concludere con le inferiori.
\end{itemize}
I bambini si sono accorti che in molte costruzioni le facce superiori
e inferiori coincidevano come numero, quindi hanno proposto di
osservare la costruzione, prima di procedere nel conteggio, perché
\begin{studente}[ ]
  Se le facce sopra e sotto sono uguali non serve contarle due
  volte. Basta farlo una volta e poi si fa per 2
\end{studente}
Procedendo secondo la suddetta scaletta i bambini hanno rettificato i
loro conteggi giungendo a risultati decisamente più corretti dei
precedenti.

Alla domanda, ``Quale costruzione occupa più spazio?''
\begin{itemize}
\item due gruppi hanno scritto subito \bambini{Nessuna} o
  \bambini{Tutte}, perché abbiamo usato sempre 8 cubi;
\item un gruppo ha scritto \bambini{Sono uguali, perché hanno lo
    stesso numero di cubi e di facce} (in questo gruppo il numero
  delle facce esterne di tutte e tre le costruzioni era effettivamente
  34)
\item due gruppi hanno citato una delle costruzioni inventate.
\end{itemize}
Qualche bambino ha collegato quest'attività con quella del conteggio
delle tessere (equiestensione) effettuata in uno dei precedenti
laboratori e ha manifestato meraviglia nell'osservare che alcuni
gruppi avevano sbagliato la risposta.

La costruzione del cubo più grande possibile con 24 dadi ha richiesto
un po' di tempo.

Tutti i gruppi, a rotazione, costruivano parallelepipedi e poi, quando
si rendevano conto che non erano cubi li disfavano.

Per limitare i tempi di attesa, ho proposto che i gruppi da 2 e da 3
si unissero agli altri, per lavorare in contemporanea.

Ho osservato con attenzione un gruppo e ho annotato le frasi più
significative che hanno portato poi a individuare quale fosse il cubo
più grande possibile che si potesse realizzare.
\begin{studente}[ ]
  \begin{itemize}
\item \bambini{No, non va bene. Così è un parallelepipedo}
\item \bambini{Il cubo deve avere la stessa quantità (di mattoncini)
    sotto, sopra e di lato}
\item \bambini{Ce ne mancano tre. Come facciamo a fare il cubo?}
\item \bambini{Ma non è detto che dobbiamo utilizzarli tutti}
\item \bambini{Guardate che è così con due cubi su ogni lato}
\item \bambini{Se da tre non si può, non si può nemmeno con 4, con 5,
    con 6, quindi si può solo con 2 come ha detto R}
\end{itemize}
\end{studente}
E così hanno concordato che il cubo più grande possibile fosse quello,
già costruito da qualche gruppo, con 8 mattoncini.

Anche in quest'ultimo incontro,
come ormai di consueto prima di affrontare la nuova proposta si è
fatto il riassunto delle puntate precedenti.  Si è ricordato che:
\begin{itemize}
\item è bene chiarirsi sui termini, prima di procedere in un'attività;
\item è possibile realizzare figure diverse aventi la stessa
  superficie e perimetro diverso;
\item è possibile realizzare figure con uguale perimetro e superfici
  diverse.
\end{itemize}
I bambini vivono l'attività come un gioco e continuano a accostarvisi
con piacere e entusiasmo.

Mi sembra che il fatto di non essere valutati singolarmente, li porti
a affrontare l'esperienza con serenità. Ho notato che un bimbo, che
non si espone mai frenato dall'eccessivo timore di sbagliare, durante
il lavoro di gruppo ha invece fatto delle osservazioni
interessanti%
.

Le conclusioni a cui si è giunti al termine di questo quarto incontro
sono state:
\begin{itemize}
\item le figure formate da uno stesso numero di cubi occupano lo
  stesso spazio;
\item non è detto che le figure che occupano uno stesso spazio hanno
  uguale superficie esterna.
\end{itemize}
Il clima nei gruppi è stato come sempre positivo.

\chapter{Giocare con le forme}


\section[Sperimentazione \#1: seconda primaria]{Sperimentazione \#1: classe seconda primaria, gennaio~2010}

\subsection{Osservazioni generali }

\subsubsection{Presentazione della classe }

Si tratta di una seconda composta da 19 alunni.

\subsubsection{Composizione dei gruppi }
4 gruppi eterogenei formati da 4-5 bambini

\subsubsection{Insegnanti presenti}

Alcuni incontri si sono svolti in periodi di compresenza.

\subsubsection{Calendarizzazione degli incontri}
\begin{calendario}
  \begin{itemize}
  \item 19 gennaio (compresenza, due docenti)
  \item 21 gennaio (compresenza, due docenti)
  \item 25 gennaio (solo l'insegnante di classe)
  \item 1 febbraio (solo l'insegnante di classe)
  \end{itemize}
\end{calendario}

\subsection{Primo incontro }

\begin{description}
\item[Alunni presenti:]18 presenti, 1 alunna assente
\item[Tempo effettivo di lavoro:]1 ora e mezza, dalle 14.40 alle 16.15
  circa: 15 minuti per l'invenzione delle figure esemplificative di
  ogni gruppo, 45 minuti per la classificazione delle figure negli
  insiemi, 40 minuti per il disegno degli insiemi trovati
\end{description}

\begin{consegna}Inizialmente le insegnanti dicono ai bambini di
  osservare le forme e di manipolarle; quindi invitano ogni gruppo a
  costruire una figura con le forme messe a disposizione, disegnarla e
  ripassarne i contorni. Il nome della figura costruita darà il nome
  anche al gruppo.

  In un secondo momento viene detto ai bambini di ``mettere in ordine
  le figure per formare degli insiemi'' e disegnare tutti gli insiemi
  formati seguendo il contorno delle forme.

  \materiali{} %
  A ogni gruppo è stato consegnato il materiale previsto per la
  classe seconda, con l'aggiunta del rettangolo e del parallelogramma.
\end{consegna}

\subsubsection{Osservazioni }
Nella prima parte dell'attività
(invenzione della figura esemplificativa di ogni gruppo) inizialmente
ogni bambino cerca di appropriarsi di alcune forme e di costruire
singolarmente una figura. Dopo un'ulteriore spiegazione
dell'insegnante i bambini provano a costruire insieme più figure fino
a sceglierne una. Le insegnanti hanno dovuto stabilire un tempo
massimo (15 minuti circa) entro cui scegliere altrimenti questa prima
fase del lavoro si sarebbe protratta più a lungo.

Nella seconda parte dell'attività (ordinamento delle figure negli
insiemi) si nota che i bambini fanno fatica a concentrarsi e
inizialmente qualcuno sembra che ordini a caso, senza saper dare una
spiegazione logica. Le insegnanti intervengono ricordando cosa
significa ``formare insiemi''%
, facendo molti
esempi e ricordando le attività svolte l'anno precedente. Una volta
chiaritisi le idee i bambini cominciano a formare insiemi.

Ogni gruppo ha classificato per colore e per il numero dei lati (che
solo un bimbo chiama lati mentre gli altri dicono che le forme hanno
gli ``spigoli'' o le ``punte'', quindi formano insiemi di forme con 3
punte, 4 punte, 5 punte\dots{}).

Un gruppo ha notato che si possono mettere insieme tutte le forme che
hanno le punte: rimane fuori il ``cerchio'' perché non ha le punte.

\subsubsection{Consigli per i colleghi che vogliono proporre le stesse
  attività}
Sarebbe meglio effettuare questa attività in compresenza come ho fatto
io per poter seguire meglio i vari gruppi, soprattutto se c'è qualche
bimbo che ha poca concentrazione e gli risulta difficile lavorare in
gruppo.

\subsection{Secondo incontro }

\begin{description}
\item[Alunni presenti:]18 alunni presenti, 1 assente
\item[Tempo effettivo di lavoro:]1 ora, dalle 11.00 alle 12.00 circa
\end{description}
\begin{consegna}
  Ai bambini è stato chiesto di osservare di nuovo gli insiemi formati
  nell'incontro precedente e individuare ulteriori caratteristiche,
  per arrivare a classificare le forme secondo il numero dei lati.

  \materiali{}%
  A ogni gruppo è stato consegnato il materiale previsto per la
  classe seconda con l'aggiunta del rettangolo e del parallelogramma.
\end{consegna}

\subsubsection{Osservazioni }
I vari gruppi già nell'incontro precedente avevano effettuato la
classificazione richiesta, che aveva come criterio il numero dei lati;
si è trattato quindi di nominare gli insiemi scoprendo i termini
geometrici: triangoli, quadrilateri, pentagoni, esagoni. Tranne che
per esagono, ottagono e romboide
 alcuni bambini hanno dimostrato di
conoscere i nomi delle altre figure piane: pentagono, quadrato,
trapezio, rettangolo, rombo, triangolo.

È stato interessante annotare i termini che i vari gruppi hanno usato
per definire l'insieme di figure con quattro punte:
\bambini{quadrupedo} (perché ricordava loro gli animali a 4 zampe),
\bambini{quadrupedia}, \bambini{cubo}, \bambini{quarantina},
\bambini{quattrina}, \bambini{quattrangolo},
\bambini{quadrallelepipedo}, \bambini{quarto}.

Per l'insieme delle figure con 5 punte i bambini hanno inventato i
nomi: \bambini{cinquelati}, \bambini{pentacoli}, \bambini{pentali},
\bambini{quinto}, \bambini{pentagoli}, \bambini{pentagologici},
\bambini{pentacinque}, \bambini{pentalolo}, \bambini{cinquelini},
\bambini{cinqueldi}, \bambini{pentagolinque}, \bambini{pentagonqui},
\bambini{cinquina}, \bambini{penquini}. Alla fine sono stati forniti i
termini esatti.

\subsubsection{Consigli per i colleghi che vogliono proporre le stesse
  attività}
Questa attività non ha richiesto molto tempo perché era già stata
effettuata una classificazione; meglio così perché i tempi di
concentrazione e di attenzione nei bambini di quest'età sono brevi.

\subsection{Terzo incontro }

\begin{description}
\item[Alunni presenti:]19 presenti (tutti)
\item[Tempo effettivo di lavoro:]1 ora, dalle 14.30 alle 15.30
\end{description}
\begin{consegna}
  I bambini giocano a \attivita{Indovina chi}: un bambino sceglie una
  forma che non deve rivelare, poi un compagno deve porre delle
  domande scegliendo ogni volta una sola caratteristica della forma
  (per es.: è rossa, ha tre lati\dots{}) a cui il primo bambino
  risponderà con SÌ o NO. Per totalizzare 1 punto si deve indicare la
  forma, per totalizzarne 2 si deve dire anche il nome corretto
  (es. pentagono). Alla fine vince il gruppo che ha realizzato più
  punti. L'attività si è svolta solo a livello orale. L'attività
  riprende e modifica quella proposta dal kit.

  \materiali{}%
  A ogni gruppo è stato consegnato il materiale previsto per la
  classe seconda con l'aggiunta del rettangolo e del parallelogramma.
\end{consegna}

\subsubsection{Osservazioni}
Prima
di iniziare il gioco, nella classe sono stati ripetuti più volte i
nomi delle figure. Nel gioco un paio di gruppi, che si sono accordati
velocemente su chi doveva iniziare per primo e hanno organizzato dei
turni, sono riusciti a fare il giro più volte realizzando così tanti
punti.

Un gruppo, pur indicando le forme esatte, si è penalizzato perché non
ricordava i nomi delle forme.

L'insegnante è passata fra i gruppi, senza intervenire, per
controllare che i bambini non ingannassero sui punti e per osservare
che la maggior parte di loro conosceva i nomi delle figure.

\subsubsection{Consigli per i colleghi che vogliono proporre le stesse
  attività}
Anche per questa attività sarebbe meglio essere in due insegnanti per
controllare meglio lo svolgimento.

\subsection{Quarto incontro }

\begin{description}
\item[Alunni presenti:]19 presenti (tutti)
\item[Tempo effettivo di lavoro:]1 ora, dalle 14.30 alle 15.30
\end{description}
\begin{consegna}
  Gioco: \attivita{Tombola delle forme}. L'insegnante spiega che
  quando estrarrà un cartoncino con una forma i bambini dovranno
  osservarlo bene e, se nella loro cartella riconosceranno la forma,
  dovranno colorarla. Vince chi per primo colora tutta la
  cartella. L'attività riprende quella prevista nel kit.

  \materiali{} %
  I materiali previsti dal kit per l'attività: una cartella per ogni
  coppia di bambini e le forme fatte di cartoncino.
\end{consegna}

\subsubsection{Osservazioni }
Inizialmente l'insegnante pensava di dare una cartella per ogni
gruppo, facendo colorare le forme a rotazione; poi, valutando
l'eccitazione e l'impazienza di partecipare dei bambini, ha consegnato
una cartella ogni 2 bambini.

Quando inizia l'estrazione i bambini sono molto attenti; a ogni
estrazione osservano la forma e la colorano se la ritrovano in
cartella, mentre per le forme che conoscono meno (pentagono,
esagono\dots{}) contano i lati. Dopo ogni estrazione la classe dice il
nome del poligono; qualcuno ha difficoltà, non ha ancora memorizzato i
nomi.

A un certo punto due bambini si accorgono che qualcosa non quadra
nella loro tabella: infatti hanno colorato una forma sbagliata
(pentagono anziché esagono), per cui all'estrazione del pentagono lo
ritrovano già colorato. Il gioco procede anche quando una coppia ha
colorato tutta la cartella per dare la possibilità a tutti di
completare.

Il gioco è stato molto gradito e i bambini si sono
divertiti. L'insegnante pensa di riproporre il gioco di tanto in
tanto, non mostrando più le forme ma dicendo i loro nomi, per
rinforzare nei bambini la conoscenza dei nomi delle figure.

\subsubsection{Consigli per i colleghi che vogliono proporre le stesse
  attività}
Giocare almeno due volte: io ritengo sia meglio a coppie perché i
bambini si compensano. Riproporre la tombola a distanza di tempo.


\section[Sperimentazione \#2: seconda primaria]{Sperimentazione \#2:
  classe seconda primaria, gennaio/febbraio~2010}

\subsection{Osservazioni generali}

\subsubsection{Presentazione della classe}
21 alunni; di questi 14 frequentano il tempo pieno (40 ore
settimanali), mentre 7 frequentano 27 ore settimanali e nelle giornate
di martedì e venerdì la loro attività didattica termina alle
12.30. \`E presente un'alunna, L., con un grave handicap cognitivo.

\subsubsection{Composizione dei gruppi}
Gruppi eterogenei scelti dall'insegnante; 3 gruppi da 5 alunni e un
gruppo da 6 con la presenza di L. Gli alunni hanno già sperimentato il
lavoro di gruppo per attività pratiche e manuali.

\subsubsection{Insegnanti presenti}

Alcuni incontri si sono svolti in compresenza con l'ausilio
dell'insegnante di sostegno di L.

\subsubsection{Calendarizzazione degli incontri}
\begin{calendario}
  \begin{itemize}
  \item 27 gennaio  (compresenza)
  \item 3 febbraio  (compresenza) con integrazione il 9 febbraio
  \item 10 febbraio  (compresenza)
  \item 11 febbraio  (solo l'insegnante di classe)
  \end{itemize}
\end{calendario}

\subsection{Primo incontro}

\begin{description}
\item[Alunni presenti:]19 presenti, 2 assenti
\item[Tempo effettivo di lavoro:]2 ore, dalle 14.30 alle 16.30: un'ora
  per la prima parte dell'attività e una per la seconda parte.
\end{description}
\begin{consegna}
  Prima parte dell'attività: iniziale manipolazione libera del
  materiale; poi invito a costruire una forma condivisa dal gruppo, a
  cui dare un nome; quindi riproduzione della forma su di un
  cartellone, con il ripasso dei contorni.

  Seconda parte dell'attività: ogni gruppo con le forme che ha a
  disposizione prova a costruire insiemi e a nominarli; quindi disegna
  gli insiemi trovati.  L'attività riprende quella prevista dal kit
  per le classi seconde.

  \materiali{}%
  Forme morbide previste dal kit per l'attività
\end{consegna}

\subsubsection{Osservazioni}
Inizialmente i bambini hanno manipolato le forme considerandone la
morbidezza al tatto e, spontaneamente ma singolarmente, hanno
cominciato a utilizzarle per costruire disegni e oggetti piani. I
bambini si sono divertiti e sono riusciti, anche i più restii a a
farlo, a apportare il proprio contributo. Si è rilevata grande
serietà e partecipazione.

Nella seconda parte dell'attività gli alunni immediatamente hanno
diviso le forme in relazione al colore. Hanno poi individuato
l'insieme ``triangoli''. Dopo riflessioni che hanno richiesto grande
messa in comune di idee un gruppo (tartaruga) ha focalizzato
l'attenzione sui lati; le forme sono state quindi classificate in base
al numero dei lati: 3 - 4 - 5. Non sono state considerate, invece, le
forme a 6 lati e a 8 perché elementi singoli.

Un alunno, dopo lunga riflessione e con grande difficoltà di
verbalizzazione del pensiero, ha associato esagono e ottagono perché
\bambini{hanno i lati che coincidono}. Successivamente, con l'aiuto
dei compagni ha specificato il suo \bambini{coincidono} con
\bambini{lati lunghi uguali}. Non ha saputo però vedere altre forme da
poter inserire nell'insieme.

Gli altri gruppi hanno, in momenti diversi, individuato insiemi con 3
- 4 - 5 - 6 \bambini{punte}. Tutti i gruppi hanno inserito il cerchio
nell'insieme \bambini{senza punte e senza lati}.

In tutto questo tempo l'alunna L. ha osservato e manipolato le forme.

Al termine della discussione finale la classe ha concordato un elenco
di classificazioni, che sono diventate il punto di partenza per il
secondo incontro.

\subsection{Secondo incontro}

\begin{description}
\item[Alunni presenti:]20 presenti, 1 assente
\item[Tempo effettivo di lavoro:]2 ore, dalle 14.30 alle 16.30
\end{description}
\begin{consegna}
  Si invitano gli alunni a formare nuovamente gli insiemi già creati,
  per affinare l'osservazione e favorire una classificazione più
  specifica delle forme, riprodotta anche graficamente.

  \materiali{}%
  Forme morbide previste dal kit per l'attività
\end{consegna}

\subsubsection{Osservazioni}

In classe sono stati appesi cartelloni con le rappresentazioni degli
insiemi formati nel primo incontro. Ciascun gruppo nomina i propri
insiemi e tutti classificano le forme in relazione al numero dei lati.

I bambini osservano insieme all'insegnante le ``punte'' delle forme e
alcuni suggeriscono di non chiamarle ``punte'', ma che \bambini{dove
  si incontrano due lati è come l'angolo delle strade}. Si scopre che
se la forma ha 4 lati ha anche 4 angoli; la classificazione si fa più
specifica per ogni insieme.

La classe passa poi a osservare le forme di ciascun insieme e prova a
trovare somiglianze e differenze. Nell'insieme dei 4 lati e 4 angoli
gli alunni confrontano il quadrato (di cui già conoscono il nome in
quanto forma dei blocchi logici) con il trapezio (che però non sanno
nominare) e scoprono che il trapezio ha un lato \bambini{storto};
\bambini{si dice obliquo} precisa un alunno.

Un altro alunno sottolinea allora che il rombo (forma conosciuta nei
diagrammi di flusso) \bambini{ha tutti e 4 i lati obliqui}. Le
osservazioni continuano in modo incalzante e con grande
entusiasmo. Emerge anche che ci sono forme che hanno tutti i lati
\bambini{dritti o obliqui} e tutti gli angoli \bambini{larghi o
  stretti} uguali.

Quello che manca agli alunni è la nomenclatura corretta, sanno fare
osservazioni interessanti ma usano un linguaggio semplice e concreto:
\begin{studente}[]
  Se le forme con tre lati e tre angoli si chiamano triangoli allora i
  4 lati e i 4 angoli si chiamano `quadriangoli'\dots{}no, si possono
  chiamare anche `quadrilateri'. Allora i 5 lati si chiamano
  `cinquelateri'\dots{}
\end{studente}
L'insegnante ritiene opportuno a questo punto fornire la corretta
nomenclatura, della quale i bambini prendono atto come una nuova
conoscenza.

Si realizza quindi un grande insieme dove entrano tutte le figure,
tranne il cerchio che gli alunni hanno escluso dall'inizio del lavoro
perché senza lati e senza angoli; quindi a ogni insieme si associa il
nome corretto: triangoli, quadrilateri, pentagoni, esagoni,
ottagono. Poi i bambini decidono che le figure che hanno lati e angoli
uguali possono stare in un insieme intersezione:
\begin{studente}[]
  hanno tutto uguale\dots{} sono regolari!!
\end{studente}
L'insegnante suggerisce il nome dell'insieme che contiene tutte le
forme, l'insieme dei ``poligoni'', spiegando il significato del
termine: molti angoli. Tutti gli elementi dell'insieme rispondono a
quella caratteristica.

Insieme si decide di realizzare un unico grande cartellone con gli
insiemi trovati e sistemati correttamente. Qualche alunno suggerisce
di creare una tabella a doppia entrata dove collocare le forme. In
quest'occasione manca il tempo materiale per farlo, ma la classe la
realizzerà nell'incontro successivo.

Gli alunni in questo secondo incontro hanno partecipato con tanto
interesse, nonostante l'attività fosse meno stimolante dal punto di
vista manuale. La curiosità era dettata, più che dal fare, dal trovare
le ``giuste'' parole.

L'alunna L., durante questo secondo incontro, ha manipolato le forme
contando le ``punte''.

Grande soddisfazione nel vedere i bambini organizzarsi, confrontarsi,
valutare le riflessioni e concordare la risposta.

\subsection{Attività integrativa al secondo incontro}

\begin{description}
\item[Alunni presenti:]L'attività è stata realizzata di martedì
  pomeriggio, ovvero con la presenza solo degli alunni frequentanti il
  tempo pieno. Di questi 14, vi erano 11 presenti e 3 assenti.
\item[Tempo effettivo di lavoro:]2 ore, dalle 14.30 alle 16.30.
\end{description}
\begin{consegna}
  Si ricordano le riflessioni emerse dalla discussione al termine del
  secondo incontro; si decide di realizzare insieme la tabella a
  doppia entrata che era stata proposta in quel momento. Si tratta di
  un'attività nata dalle esigenze della classe e non proposta dal kit.

  \materiali{}%
  Carta, pennarelli e colori.
\end{consegna}

\subsubsection{Osservazioni}
Si decide di impostare la tabella in questo modo:
\begin{center}
  \begin{tabular}{|l|c|c|c|c|c|}
    \hline
    lati/angoli & 3 & 4 & 5 & 6 & 8 \\ \hline
    3 &&&&& \\ \hline
    4 &&&&& \\ \hline
    5 &&&&& \\ \hline
    6 &&&&& \\ \hline
    8 &&&&& \\ \hline
  \end{tabular}
~~
\raisebox{-0.5\totalheight}{\includegraphics[width=0.48\textwidth]{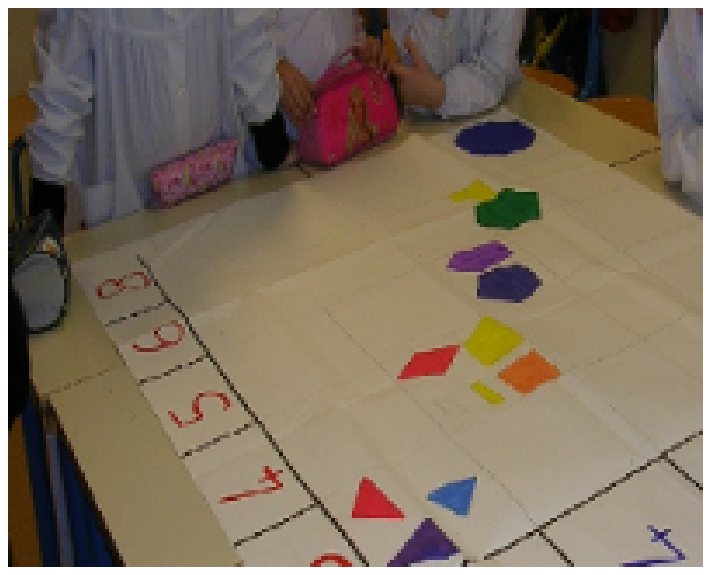}}
\end{center}

All'interno della tabella sono state poi inserite le forme
ripassandone i contorni; con grande facilità i bambini hanno
organizzato le forme negli spazi corretti%
.

Un alunno osservando il cartellone degli insiemi si è reso conto che
nella tabella a doppia entrata non è possibile evidenziare i poligoni
regolari perché
\begin{studente}[]
  non c'è la cella adatta
\end{studente}
Un altro alunno invece è rimasto con il cerchio in mano, perché non ha
trovato l'\bambini{incrocio} dove poterlo inserire. Qualcuno, con
l'aria di chi adesso ha imparato tante cose gli ha detto:
\begin{studente}[]
  Ma quello non ce lo puoi mettere\dots{} non vedi che questa è la
  tabella dei poligoni?\dots{} e il cerchio mica ha gli angoli!!!!!
\end{studente}
Risposta del
primo bambino:
\begin{studente}[]
  Hai ragione\dots{} non ci ho pensato!
\end{studente}
rimettendo il cerchio nel sacchetto.

L'insegnante osserva che è questo il bello del lavorare con la
spontaneità e la curiosità dei bambini: quando hanno capito si sentono
importanti e grandi.

\begin{figure}[pht]
  \centering
  \begin{tabular}{cc}
    \includegraphics[width=0.47\textwidth]{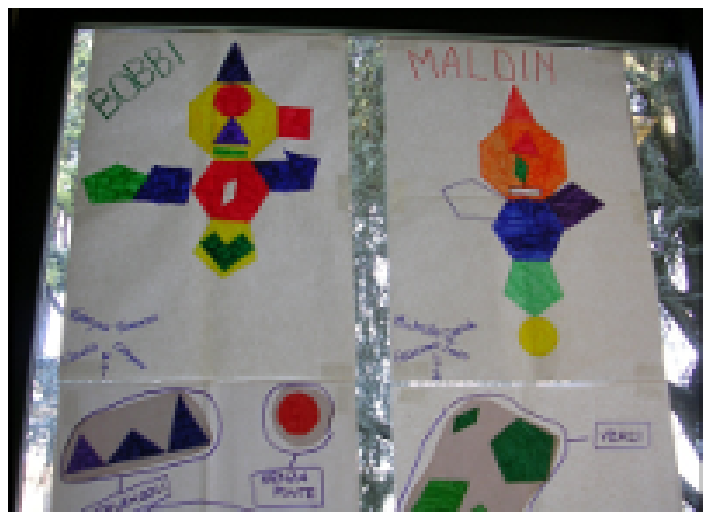} &
    \includegraphics[width=0.47\textwidth]{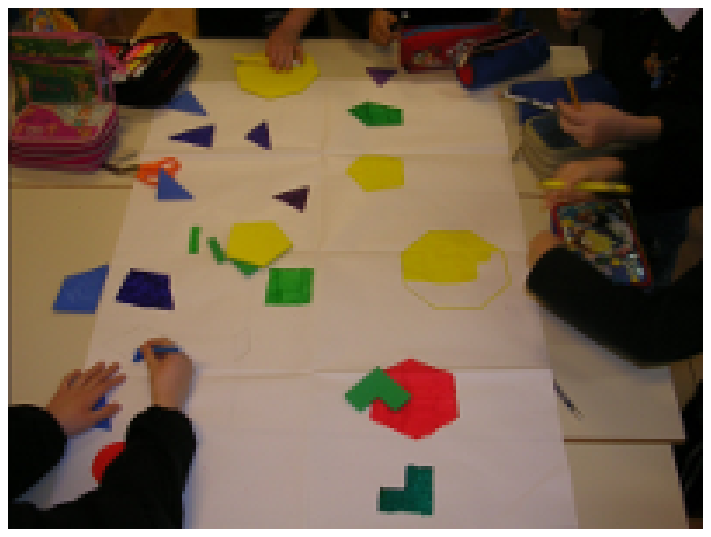} \\
    \includegraphics[width=0.47\textwidth]{\imgdir/forme-2-3} &
    \includegraphics[width=0.47\textwidth]{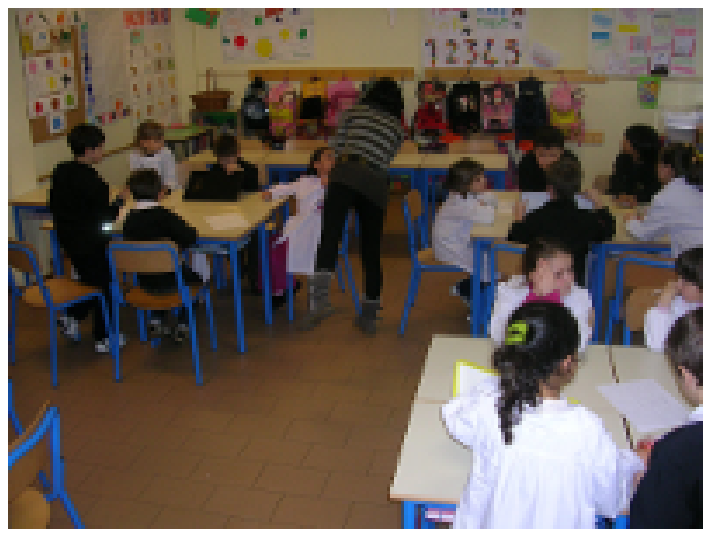} \\
  \end{tabular}
  \label{pic:formemorbide:2}
\end{figure}

\subsection{Terzo incontro}

\begin{description}
\item[Alunni presenti:]19 presenti, 2 assenti.
\item[Tempo effettivo di lavoro:]2 ore, dalle 14.30 alle 16.30.
\end{description}
\begin{consegna}
  \`E presentata l'attività numero 3 prevista dal kit: \attivita{Indovina
    la forma}.

  \materiali{}
  I materiali previsti dal kit per l'attività.
\end{consegna}

\subsubsection{Osservazioni}
Prima di iniziare l'attività alcuni alunni presenti il giorno
precedente hanno raccontato ai compagni ciò che è stato realizzato nel
pomeriggio: la tabella a doppia entrata. Hanno spiegato come è stata
costruita e il perché di alcune omissioni (i poligoni regolari non
sono evidenziati, il cerchio non è presente). Da questa conversazione
l'insegnante si è resa conto che tanti sono i messaggi che sono stati
interiorizzati e molte le nuove conoscenze. Ad esempio, quando hanno
spiegato il motivo per cui non erano stati evidenziati i poligoni
regolari, i bambini hanno mostrato di aver interiorizzato l'uso e la
funzione della tabella a doppia entrata. \`E emersa anche una chiarezza
di idee sui poligoni regolari, i quali pur avendo una caratteristica
specifica in più non potevano essere ``incasellati'' nella tabulazione
costruita.

Quindi, una volta spiegate le regole del gioco, gli alunni hanno
cominciato a ``giocare'' e da subito si è notato l'impegno con il
quale cercavano di scoprire la forma nascosta. Le domande poste erano
rispondenti alla richiesta del gioco e i bambini si sono divisi bene i
ruoli. Si riporta un esempio di domanda pertinente: se un componente
del gruppo poneva la domanda
\begin{studente}[]
  ha 4 lati?
\end{studente}
il compagno faceva precedere la propria domanda successiva con un
\begin{studente}[]
  se non ha 4 lati non avrà nemmeno 4 angoli quindi ti chiedo: ha 5
  angoli?
\end{studente}

Nella seconda parte del gioco, quando i bambini sono stati chiamati
alla riflessione sulle differenze tra forme, si sono confrontati molto
e sono arrivati a ``concordare'' le risposte da dare.

Durante il momento di gioco, il ruolo delle insegnanti è stato quello
di osservatori; deliberatamente sono intervenute il meno possibile per
lasciare ai gruppi la possibilità di esprimersi liberamente e
spontaneamente. Grande soddisfazione nel vedere i bambini
organizzarsi, confrontarsi, valutare le riflessioni e concordare la
risposta.

\subsection{Quarto incontro}

\begin{description}
\item[Alunni presenti:]20 presenti; è assente l'alunna L.
\item[Tempo effettivo di lavoro:]1 ora e mezza, dalle 11 alle 12.30
\end{description}
\begin{consegna}
  Si propone agli alunni il gioco della tombola, previsto come
  attività del kit. Gli alunni non sono però divisi a gruppi ma
  giocano singolarmente.

  \materiali{}
  Forme e cartelle della tombola previste dal kit per l'attività.
\end{consegna}

\subsubsection{Osservazioni}
Non sono necessarie ampie spiegazioni perché i bambini conoscono il
gioco, anche se qualcuno preferisce chiamarlo ``bingo''. L'insegnante
estrae le figure e, senza alcun commento, le mostra alla classe; gli
alunni riconoscono la figura sulla propria tabella e provvedono a
colorarla.

Alla quinta estrazione alcuni alunni fanno tombola, ma il gioco
continua sino a quando tutti riescono a colorare la tabella. Qualcuno
chiede di poter scrivere il nome delle figure sulla propria tabella e
successivamente chiedono di ripetere il gioco. Vengono quindi
distribuite altre cartelle e si ricomincia l'estrazione.

Il gioco, in quanto tale, li ha divertiti molto. Nessuno ha mostrato
difficoltà né nel riconoscimento delle forme, né nell'attribuzione
del giusto nome.

\subsection{Considerazioni conclusive al termine della
  sperimentazione}
La sperimentazione, a mio avviso, si è rilevata molto interessante sia
nelle modalità di attuazione, sia per quanto riguarda i risultati
ottenuti. Gli alunni hanno partecipato con costante e attivo
interesse; ho potuto constatare che gli alunni più coinvolti sono
stati coloro i quali vivono la matematica come qualcosa di ostico e
troppo difficile per loro.

Intitolare queste quattro attività come ``Imparare sperimentando'' ha
fatto sì che l'approccio risultasse svincolato da qualsiasi timore. Un
alunno particolarmente emotivo ha più volte voluto sottolineare che
questi lavori gli piacevano tanto con frasi come:
\begin{studente}[]
  Mi sento importante perché il mio pensiero servirà ai signori di una
  università di matematica di Milano
\end{studente}

Tante sono state le conoscenze acquisite per scoperta e non per
trasmissione. Un consiglio che vorrei dare a chi ancora deve
``sperimentare'' è di non porsi assolutamente un obiettivo finale, non
domandarsi mai ``ma dove voglio arrivare?''. Il percorso di lavoro
deve nascere dagli alunni%
, l'evoluzione
deve seguire i loro ritmi, ciò che verrà acquisito, tanto o poco, sarà
frutto delle loro scoperte comuni. Il lavoro di gruppo come messa in
comune di pensieri e riflessioni li porterà a interiorizzare in modo
profondo.

Sottolineo che solo alcuni dei miei alunni hanno 8 anni, la maggior
parte di loro ne ha solo 7 ma con questa esperienza hanno dimostrato
una grande ``maturità''.


\section[Sperimentazione \#3: seconda primaria]{Sperimentazione \#3:
  classe seconda primaria, gennaio~2010}

\subsection{Osservazioni generali}

\subsubsection{Presentazione della classe}
La classe è formata da 23 alunni.

\subsubsection{Composizione dei gruppi}
Gruppi eterogenei, ciascuno con 4, 5 o 6 componenti, non costanti nei
vari incontri a causa degli alunni assenti.

\subsubsection{Insegnanti presenti}
In quasi tutti gli incontri la docente sperimentatrice è affiancata
dall'altra docente di classe (compresenza).

\subsubsection{Calendarizzazione degli incontri}

\begin{calendario}
  \begin{itemize}
  \item 27 gennaio (compresenza)
  \item 28 gennaio (compresenza)
  \item 3 febbraio (compresenza)
  \item 8 febbraio
  \end{itemize}
\end{calendario}

\subsection{Primo incontro}

\begin{description}
\item[Alunni presenti:]21 alunni presenti, 2 assenti
\item[Tempo effettivo di lavoro:]2 ore, dalle 10.30 alle 12.30. Trenta
  minuti sono stati dedicati alla discussione e un'ora e mezza al
  lavoro di gruppo.
\end{description}

\begin{consegna}
  Inizialmente si è data la possibilità di manipolare liberamente il
  materiale, poi dopo una decina di minuti l'insegnante ha chiesto ai
  gruppi di ``creare attraverso l'accostamento, la sovrapposizione di
  forme una figura che li identificasse''. A ogni figura è stato dato
  un nome e i bambini ne hanno ripassato i contorni su un foglio. Poi
  l'insegnante, dopo aver ribadito il concetto di insieme, ha chiesto
  ai bambini di provare a formarne qualcuno e a dare loro dei
  nomi. Alla fine gli alunni hanno disegnato gli insiemi formati.

  L'attività proposta riprendeva quanto previsto dal kit per le classi
  prima e seconda.

  \materiali{}
  A ogni gruppo è stato consegnato un sacchetto con le forme dentro
  (sacchetti già preparati prima). Le forme erano quelle previste dal
  kit per le attività delle classi prima e seconda; l'insegnante ha
  aggiunto il rettangolo e il parallelogramma previsti per le attività
  delle classi terza, quarta e quinta.
\end{consegna}

\subsubsection{Osservazioni}
All'inizio
i bambini hanno tirato fuori le forme e hanno cominciato a
manipolarle, a giocare e anche a divertirsi. Quindi ogni gruppo ha
creato la sua mascotte con l'accostamento di forme e le ha dato un
nome. Quando è stato chiesto di formare gli insiemi con le figure, i
bambini:
\begin{itemize}
\item hanno dapprima diviso le forme in relazione al colore.
\item hanno riconosciuto l'insieme dei triangoli e, dopo l'invito a
  un'osservazione più attenta, li hanno classificati come figure a 3
  lati. Quindi l'insegnante ha spostato l'attenzione sui lati delle
  altre figure e sono state classificate quelle con 4 e 5 lati.
\item un alunno ha messo insieme il pentagono, l'esagono e l'ottagono
  perché ogni figura, ha affermato, ha i lati uguali tra loro.
\item un bambino del gruppo Pippo, osservando alcune forme, le ha
  messe insieme perché sosteneva che avessero i lati paralleli;
  l'insegnante e la sua collega hanno verificato che il bambino
  conosceva il significato del termine. Il termine è stato poi ripreso
  e spiegato dall'insegnante alla classe, usando l'immagine dei binari
  del treno.
\end{itemize}

\subsubsection{Consigli per i colleghi che vogliono proporre le stesse
  attività}
Lasciarli giocare perché è proprio attraverso il gioco che si facilita
l'apprendimento.

\subsection{Secondo incontro}

\begin{description}
\item[Alunni presenti:] 23 presenti (tutti)
\item[Tempo effettivo di lavoro:]2 ore, dalle 10.30 alle 12.30
\end{description}

\begin{consegna}
  Dopo aver appeso alla parete della classe i cartelloni che
  identificano i gruppi e i disegni degli insiemi trovati la lezione
  precedente, l'insegnante ha chiesto ai bambini di riguardarli e
  cercare di fare un'ulteriore classificazione. Gli insiemi sono stati
  poi disegnati su carta da pacco.

  \materiali{}
  Forme previste dal kit e carta da pacco.
\end{consegna}

\subsubsection{Osservazioni}
Dopo un momento di smarrimento i bambini hanno capito che avrebbero
dovuto formare dei sottoinsiemi. La classe ha ripassato la definizione
di sottoinsieme e i bambini si sono messi al lavoro.

A partire dall'insieme dei triangoli gli alunni hanno formato il
sottoinsieme dei triangoli con lati diversi (lasciando fuori quello
equilatero). Partendo dalle forme con quattro lati un gruppo ha creato
il sottoinsieme delle forme con i lati ``storti'' (rombo, trapezio) e un
altro quello con i lati uguali tra loro (quadrato, rombo). Il gruppo
dell'Arancia è partito dall'insieme delle forme con cinque punte e ha
formato il sottoinsieme che conteneva la figura con tutti i lati
uguali. Il gruppo Pippo, guardando l'insieme con le figure con cinque,
sei e otto lati, ha formato il sottoinsieme contenente l'esagono e
l'ottagono; i bambini hanno notato che tutti i lati di queste due
figure hanno un altro lato ``di fronte'' (paralleli), mentre il
pentagono non soddisfa questo criterio.

A questo punto l'insegnante, dopo aver messo sui banchi un foglio di
carta da pacco, vi ha disposto sopra tutte le forme, ha disegnato con
un pennarello un grande insieme e man mano ha sistemato le figure nei
sottoinsiemi, che poi sono stati nominati. La classe ha escluso
dall'insieme il cerchio.

L'insegnante ha invitato i bambini a osservare ancora tutti i
sottoinsiemi; così è stato notato che ogni gruppo di forme ne aveva
una con tutti i lati e le punte uguali (poligoni regolari). Insieme
hanno racchiuso queste ultime in un ulteriore gruppo, individuando
così l'intersezione tra gli insiemi. \`E allegata alla quarta attività
una fotografia che riassume tutto il lavoro svolto.

\subsubsection{Consigli per i colleghi che vogliono proporre le stesse
  attività}
Lasciare i bambini liberi di esprimersi.

\subsection{Terzo incontro}

\begin{description}
\item[Alunni presenti:]23 presenti (tutti)
\item[Tempo effettivo di lavoro:]1 ora circa, dalle 14.30 alle 15.30
\end{description}

\begin{consegna}
  \`E stata proposta l'attività \attivita{Indovina la forma} prevista
  dal kit.

  \materiali{}
  A ogni gruppo sono state consegnate la scheda fotocopiata con tutte
  le forme presentate durante le attività precedenti e la scheda da
  compilare durante il gioco stesso. Sulla cattedra sono state
  predisposte le forme colorate e non quelle presenti nel sacchetto
  nero.
\end{consegna}

\subsubsection{Osservazioni}
\begin{wrapfigure}{R}{0pt}
  \includegraphics[width=0.48\textwidth]{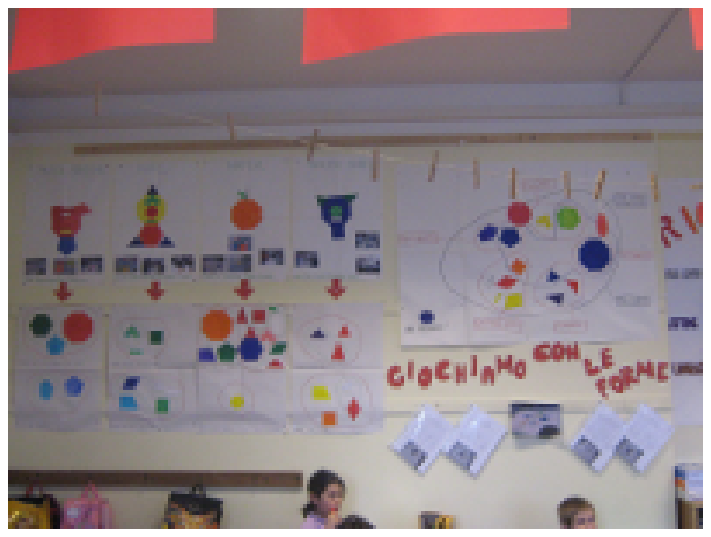}
\end{wrapfigure}
Dopo aver predisposto sulla cattedra le forme, l'insegnante ha
spiegato le regole e il funzionamento del gioco. Un bambino alla volta
poteva andare alla cattedra e scegliere
tra le tante forme una da far indovinare ai compagni; questi dovevano
arrivare alla forma ponendo delle domande secche, alle quali si poteva
rispondere solo con un sì o con un no. Le domande, generalmente, si
riferivano alle conoscenze acquisite durante le prime due attività, ad
esempio:
\begin{itemize}
\item Ha tre punte?
\item Ha quattro lati?
\item \`E rotondo?
\item Ha i lati storti?
\end{itemize}
Ogni partecipante poteva porre una sola domanda; vinceva chi
indovinava la forma e la nominava correttamente. I bambini non hanno
avuto difficoltà a riconoscere le forme e nemmeno a nominarle anche
se, per qualcuno, l'esagono è diventato \bambini{esagolo} e l'ottagono
\bambini{ottangolo}.

I bambini si sono divertiti molto e l'attività è servita loro per
memorizzare meglio i nomi delle figure e coglierne, ancora una volta,
le caratteristiche.

\subsection{Quarto incontro}

\begin{description}
\item[Alunni presenti:]22 presenti, 1 assente
\item[Tempo effettivo di lavoro:]1 ora, dalle 14.30 alle 15.30
\end{description}

\begin{consegna}
  \`E stata proposta l'attività \attivita{La tombola delle forme}
  prevista dal kit; l'insegnante però non ha diviso la classe a
  gruppi, come previsto nel kit, ma ha consegnato una cartella a ogni
  bambino. In questo modo si intendeva fornire a ogni bambino la
  possibilità di essere maggiormente partecipe nell'attività.

  \materiali{}
  Cartelle della tombola e forme per l'estrazione.
\end{consegna}

\subsubsection{Osservazioni}
L'insegnante ha iniziato a estrarre le forme, che i bambini man mano
riconoscevano sulla propria cartella e coloravano. Il gioco è
continuato finché tutti i bambini hanno completato la propria tabella.

A questo punto, per avere un ulteriore riscontro dell'apprendimento
dei bambini, l'insegnante ha chiesto loro di scrivere il nome
geometrico di ciascuna figura. Escluso qualche errore di ortografia,
tutti hanno riconosciuto le forme e sono stati in grado di
denominarle.


\section[Sperimentazione \#4: seconda primaria]{Sperimentazione \#4:
  classe seconda primaria, gennaio/febbraio~2010}

\subsection{Osservazioni generali}

\subsubsection{Presentazione della classe}
20 alunni, di cui uno seguito da un'insegnante di sostegno

\subsubsection{Composizione dei gruppi}
4 gruppi eterogenei, costanti negli incontri: 3 gruppi da 5 bambini e
1 gruppo da 4.

\subsubsection{Insegnanti presenti}

Negli incontri la docente sperimentatrice è affiancata dall'altra
docente di classe (compresenza) o dalla collega di sostegno.

\subsubsection{Calendarizzazione degli incontri}
\begin{calendario}
  \begin{itemize}
  \item 28 gennaio (compresenza)
  \item 29 gennaio (con insegnante sostegno)
  \item 4 febbraio (compresenza)
  \item 5 febbraio (con insegnante sostegno)
  \end{itemize}
\end{calendario}

\subsection{Primo incontro}

\begin{description}
\item[Alunni presenti:]Tutti, a eccezione del bambino seguito
  dall'insegnante di sostegno
\item[Tempo effettivo di lavoro:]2 ore, dalle 10:30 alle 12:30;
  all'incirca 1 ora dedicata al lavoro di gruppo e 1 ora alla
  discussione
\end{description}

\begin{consegna}
  Iniziale manipolazione e osservazione libera delle figure
  consegnate. In seguito l'insegnante ha chiesto a ogni gruppo di
  realizzare, accostando queste forme, una mascotte. I bambini hanno
  riportato il contorno della figura su un foglio, l'hanno colorata e
  hanno scelto un nome. In seguito, dopo aver ribadito il concetto
  d'insieme, l'insegnante ha invitato ogni gruppo a formare degli
  insiemi e a nominarli. Come attività conclusiva i bambini hanno
  ricalcato su un foglio le forme ricostruendo così gli insiemi
  trovati. L'attività segue quanto proposto dal kit.

  \materiali{}
  A ogni gruppo sono state consegnate le forme geometriche previste
  dal kit per le attività delle classi prima e seconda.
\end{consegna}

\subsubsection{Osservazioni}
Gli alunni erano molto soddisfatti e contenti dell'attività proposta e
si sono anche molto divertiti.

Nell'attività di classificazione delle figure e formazione degli insiemi:

\begin{itemize}
\item Tutti i gruppi hanno inizialmente classificato per colore.
\item Invitati a trovare altri criteri hanno quindi classificato
  contando le ``punte'' 3, 4, 5 (senza attribuire il termine
  geometrico).
\item Tutti i gruppi hanno individuato l'insieme dei triangoli e
  l'hanno nominato.
\item Un gruppo ha formato l'insieme delle figure con ``punte'' e
  ``spigoli'' lasciando fuori il cerchio.
\item Tre gruppi hanno messo insieme anche figure con ``spigoli''
  della stessa lunghezza.
\end{itemize}
I bambini hanno usato spontaneamente i termini ``punte'' e
``spigoli''; qualcuno ha introdotto anche il termine ``lati''.

\begin{itemize}
\item L'insegnante ha sollecitato i bambini a osservare anche gli
  angoli, a confrontarli per sovrapposizione, ma è risultato
  difficile.
\item Qualcuno ha notato che il rombo, inserito correttamente
  nell'insieme delle figure con tutti i lati uguali, aveva però
  qualcosa di strano\dots{} Solo alcuni si sono accorti che mentre
  ogni figura aveva tutti gli angoli uguali, il rombo li aveva uguali
  a due a due. Per questo ultimo passaggio l'insegnante ha dovuto
  invitare più volte i bambini a un'osservazione attenta e alla
  chiara verbalizzazione del loro pensiero.
\end{itemize}
Un esempio di scambio verbale tra insegnante e alunni è stato il
seguente:
\begin{tutor}[Ins.]
  Le punte sono tutte uguali?
\end{tutor}
\begin{studente}[]
  No, alcune sono più strette, altre più larghe.
\end{studente}
\begin{tutor}[Ins.]
  Quante strette e quante larghe?
\end{tutor}
\begin{studente}[]
  2 e 2.
\end{studente}
\begin{tutor}[Ins.]
  Le punte larghe sono uguali tra loro?
\end{tutor}
\begin{studente}[]
  Si, e anche quelle strette sono uguali tra loro.
\end{studente}

\subsubsection{Consigli per i colleghi che vogliono proporre le stesse
  attività}
Incoraggiare gli alunni a esprimersi per tentativi
successivi. L'esercizio risulta stimolante e le osservazioni sono
pertinenti.

\subsection{Secondo incontro}

\begin{description}
\item[Alunni presenti:]Tutti, a eccezione del bambino seguito
  dall'insegnante di sostegno
\item[Tempo effettivo di lavoro:] 2 ore, dalle 10:30 alle 12:30. Il
  tempo è stato dedicato a fasi alterne al lavoro di gruppo e alla
  discussione, complessivamente circa 1 ora per il lavoro di gruppo e 1
  ora per la discussione
\end{description}

\begin{consegna}
  L'insegnante ha iniziato l'attività mostrando ai bambini i disegni
  degli insiemi trovati il giorno precedente e invitandoli a ricordare
  i criteri di classificazione adottati. In seguito ha consegnato a
  ogni gruppo le forme e un foglio di carta da pacco, chiedendo di
  sistemare sul foglio le forme per rappresentare gli insiemi
  trovati. L'insegnante ha esplicitato che dovevano essere utilizzate
  tutte le forme, anche il cerchio.

  Le richieste dell'attività erano aggiuntive a quanto previsto dal
  kit.

  \materiali{}
  Carta da pacco e forme, previste dal kit per la prima e la seconda
  classe.
\end{consegna}

\subsubsection{Osservazioni}
All'inizio senza problemi ogni gruppo ha formato l'insieme
\bambini{forme con spigoli e punte}, escludendo il cerchio.

Quindi l'insegnante ha domandato:
\begin{tutor}[]
  All'interno di questo insieme è possibile trovarne altri?
\end{tutor}
Senza difficoltà i bambini hanno ricordato il concetto di sottoinsieme
e hanno formato quelli con le figure con 3, 4, 5, 6 punte.

Poi l'insegnante ha chiesto di costruire anche l'insieme con punte e
spigoli uguali trovato il giorno precedente. Per formare quest'ultimo
insieme i bambini hanno smontato quelli costruiti in precedenza,
perché una forma (es. quadrato) poteva stare o nell'insieme 4 punte
oppure nell'insieme punte e spigoli uguali.

A questo punto l'insegnante ha chiesto di trovare un modo di
rappresentare tutti gli insiemi. Ogni gruppo ha fatto diversi
tentativi: i bambini hanno discusso tra loro, provato, smontato e così
via.

L'insegnante ha ricordato loro un lavoro svolto all'inizio dell'anno
sui luoghi di villeggiatura, che aveva utilizzato per introdurre
l'insieme intersezione%
.

Il gruppo del PINGUINO ha subito capito e una bambina ha suggerito di
formare prima il sottoinsieme delle forme spigoli uguali e poi di
``incrociare'' gli altri. In un attimo anche dal punto di vista
dell'organizzazione spaziale i bambini hanno disegnato più insiemi
intersezione (pur ricordandone solo alla fine il nome).

Il gruppo COLORATO ha avuto bisogno di un po' più di tempo ma è
riuscito senza problemi.

Il gruppo PINCHINO ha capito cosa avrebbe dovuto fare ma non riusciva
a organizzarsi dal punto di vista spaziale e
grafico. Nell'eccitazione di aver capito i bambini facevano
confusione. Finalmente una bambina con molta calma ha preso in mano la
situazione e ha organizzato il lavoro, formando quattro intersezioni.

L'ultimo gruppo è riuscito a formare i sottoinsiemi con 3, 4, 5, 6
punte, ma le intersezioni le ha fatte solo dopo aver visto gli altri
gruppi.

Tutti hanno ricalcato, colorato le forme, ripassato con colori diversi
le linee degli insiemi, scritto i cartellini anche agli insiemi
intersezione. Li ho visti molto soddisfatti e entusiasti del lavoro
svolto.

\subsection{Terzo incontro}

\begin{description}
\item[Alunni presenti:]19 presenti, 1 assente
\item[Tempo effettivo di lavoro:]1 ora circa, dalle 11:30 alle 12:30
\end{description}

\begin{consegna}
  L'insegnante ha presentato il gioco \attivita{Indovina la forma},
  come suggerito dal kit. L'unica variante che è stata introdotta:
  ogni componente del gruppo ha scelto una forma tra quelle colorate
  (erano più tipi di forme rispetto a quelle plastificate) da far
  indovinare ai compagni del gruppo stesso.

  \materiali{}
  A ogni gruppo sono state distribuite la scheda predisposta e una
  fotocopia di tutte le forme presentate nelle lezioni precedenti,
  analogamente a quanto previsto nel kit.
\end{consegna}

\subsubsection{Osservazioni}
I bambini
hanno rispettato le regole del gioco; le domande che hanno posto si
riferivano però quasi sempre solo al numero dei lati della figura. I
bambini non hanno avuto difficoltà a indovinare le forme ma qualche
difficoltà a nominarle, in particolare per l'esagono e il pentagono,
che indicavano sulla fotocopia.

I gruppi hanno funzionato: i bambini hanno stabilito chi doveva
scrivere sulla scheda e hanno concordato le risposte alle domande. Gli
alunni si sono molto divertiti.

\subsection{Quarto incontro}

\begin{description}
\item[Alunni presenti:]19 presenti, 1 assente
\item[Tempo effettivo di lavoro:]1 ora circa, dalle 11:30 alle 12:30
\end{description}

\mbox{}
\begin{consegna}
  L'insegnante ha presentato il gioco la \attivita{Tombola delle
    forme} come suggerito dal kit. Unica variante introdotta: questo
  gioco è stato considerato come una verifica, pertanto i bambini non
  sono stati divisi in gruppo, ma a ogni alunno è stata consegnata
  una cartella. L'insegnante ha chiesto di scrivere sotto alla forma
  riconosciuta e colorata anche il suo nome geometrico.

  \materiali{}
  Cartelle e forme previste dal kit per quest'attività.
\end{consegna}

\subsubsection{Osservazioni}
I bambini non hanno avuto nessuna difficoltà a riconoscere le forme;
solo 3 bambini, tra cui l'alunno segnalato, hanno invertito i termini
triangolo-rettangolo.

A conclusione di questo lavoro sono stati terminati e appesi in classe
i cartelloni realizzati: le MASCOTTE, gli INSIEMI TROVATI e le FORME
GEOMETRICHE ricalcate con i rispettivi nomi%
.


\section[Sperimentazione \#5: seconda primaria]{Sperimentazione \#5:
  classe seconda primaria, gennaio/febbraio~2010}

\subsection{Osservazioni generali}

\subsubsection{Presentazione della classe}
19 bambini

\subsubsection{Composizione dei gruppi}
Quattro gruppi eterogenei, decisi dall'insegnante: 3 da 5 bambini e 1
da 4.

\subsubsection{Insegnanti presenti}

La docente sperimentatrice è affiancata dall'altra docente di classe
(compresenza).

\subsubsection{Calendarizzazione degli incontri}
\begin{calendario}
  \begin{itemize}
  \item 18 gennaio  (compresenza)
  \item 19 gennaio  (compresenza)
  \item 25 gennaio  (compresenza)
  \item 26 gennaio  (compresenza)
  \end{itemize}
\end{calendario}

\subsection{Primo incontro}
\begin{description}
\item[Alunni presenti:] 19 presenti
\item[Tempo effettivo di lavoro:] 2 ore, dalle 10,30 alle 12,30.
\end{description}

\begin{consegna}
  All'inizio i bambini erano liberi di toccare, manipolare, comporre e
  scomporre figure. Poi l'insegnante ha chiesto loro di comporre con
  più forme una figura che sarebbe diventata la mascotte del loro
  gruppo.

  Successivamente si è chiesto di formare degli insiemi e di spiegare
  quale fosse il criterio per il quale quelle forme potessero stare
  insieme.

  \materiali{} %
  Sono state consegnate a ogni gruppo le figure geometriche previste
  dal kit
\end{consegna}

\subsubsection{Osservazioni}
All'inizio i bambini hanno faticato a organizzarsi, poi hanno
interagito e trovato tutti una figura che li rappresentasse.

Nell'attività di formazione degli insiemi, tutti hanno creato subito
gruppi per colore; dopo la richiesta dell'insegnante di cercare nuovi
criteri, i bambini hanno preso in considerazione il numero dei lati.

L. ad esempio ha formato un gruppo di triangoli e un secondo gruppo
con tutte le altre figure; alla richiesta di esplicitazione del
criterio ha risposto: \bambini{Triangoli e non
  triangoli}. L'insegnante ha chiesto allora di formare più
sottoinsiemi nel gruppo dei non triangoli. Dopo molti tentativi i
bambini si sono accordati, individuando nuovi gruppi per il numero dei
lati: 0, 4, 5, 6, 8.

Gli alunni sono apparsi molto soddisfatti e contenti del risultato
ottenuto; inoltre hanno dimostrato di conoscere il nome di alcune
figure: \bambini{triangolo}, \bambini{rotondo}, \bambini{rettangolo},
\bambini{quadrato}, \bambini{rombo}.

\subsection{Secondo incontro}
\begin{description}
\item[Alunni presenti:] 18 presenti, 1 assente
\item[Tempo effettivo di lavoro:]2 ore, dalle ore 10,30 alle 12,30.
\end{description}

\begin{consegna}
  L'insegnante ha chiesto di formare ancora gli insiemi in base al
  numero dei lati della figura; successivamente ha chiesto ai vari
  gruppi se conoscevano il nome di alcune figure. Ha cercato di farli
  riflettere sul significato della parola triangolo e, sulla base
  delle loro riflessioni, ha chiesto di provare a dare dei nomi anche
  alle figure che i bambini non conoscevano.

  \materiali{} %
  Sono state consegnate a ogni gruppo le figure geometriche previste
  dal kit
\end{consegna}

\subsubsection{Osservazioni}
Sono riportate di seguito alcune ipotesi dei bambini sui nomi delle
figure piane: \bambini{seiangoli}, \bambini{seilateri},
\bambini{ottoangoli}, \bambini{cinqueangoli}\dots{} Dalla discussione
all'interno dei gruppi è scaturita in seguito l'importanza di un nome
che fosse condiviso e riconoscibile da tutti; solo in quel momento
l'insegnante ha comunicato loro i nomi di queste nuove figure.

All'interno del grande gruppo dei quadrilateri i bambini hanno
imparato i nomi del trapezio e del romboide%
. Hanno inoltre intuito, osservando i due pentagoni e i due esagoni,
il quadrato e il rombo%
, il concetto di poligono regolare.

\subsection{Terzo incontro}
\begin{description}
\item[Alunni presenti:] 19 presenti
\item[Tempo effettivo di lavoro:] 2 ore, dalle 10,30 alle 12,30
\end{description}

\begin{consegna}
  Prima parte dell'attività: l'insegnante ha distribuito a ciascun
  gruppo le forme e ha chiesto di lasciarle al centro del tavolo di
  lavoro. Quindi diceva a alta voce il nome di una figura che i
  bambini dovevano individuare fra tutte le altre forme. Il gioco è
  stato ripetuto più volte.

  Seconda parte dell'attività: è stato proposto il gioco
  dell'\attivita{Indovina chi?}. In ogni gruppo un bambino per volta
  doveva scegliere una forma senza però comunicarla ai suoi compagni;
  i compagni per individuarla dovevano porre delle domande a cui poter
  rispondere solo sì o no, del tipo: è rosso? \`E grande? Ha quattro
  lati? Ogni partecipante poteva cercare di individuare la forma solo
  dopo aver posto una domanda; vinceva chi individuava e chiamava la
  forma con il giusto nome.

  \materiali{} %
  Sono state consegnate a ogni gruppo le figure geometriche previste
  dal kit
\end{consegna}

\subsubsection{Osservazioni}
Nella prima parte dell'attività gli alunni si sono divertiti molto;
dopo varie prove si è data loro la possibilità di alternarsi al
comando del gioco. I bambini hanno memorizzato e riconosciuto con
facilità le figure.

Anche nella seconda parte dell'attività i bambini si sono molto
divertiti, hanno memorizzato ulteriormente i nomi delle figure, hanno
colto alcuni aspetti peculiari di ciascuna figura (il quadrato ha
tutti i lati uguali, il rettangolo ha i lati uguali a due a
due\dots{}) e hanno allenato le loro capacità logiche nello sforzo di
porre domande, sempre più mirate, che li aiutassero a individuare la
figura giusta.

\subsection{Quarto incontro}
\begin{description}
\item[Alunni presenti:] 19 presenti
\item[Tempo effettivo di lavoro:] 2 ore, dalle 10,30 alle 12,30.
\end{description}

\begin{consegna}
  \`E stata proposta la tombola delle forme prevista dal kit; a
  differenza del kit però il gioco non è stato proposto al piccolo
  gruppo ma a tutti i bambini della classe. L'insegnante ha spiegato
  le regole: quando il bambino vedeva la figura estratta riprodotta
  sulla sua cartella, la doveva colorare; avrebbe vinto il bambino che
  per primo avesse completato la coloritura di tutte le figure della
  cartella.

  \materiali{} %
  \`E stata distribuita a ciascun bambino una cartella della tombola
  delle forme, prevista dal kit.
\end{consegna}

\subsubsection{Osservazioni}
Dopo che il primo bambino ha fatto tombola, l'insegnante ha deciso di
continuare l'estrazione delle figure per consentire a tutti di
completare la coloritura della cartella.

Il gioco ha talmente coinvolto e divertito i bambini che si è deciso
di fare un'altra partita: l'insegnante ha distribuito altre cartelle
e ha proceduto alla nuova estrazione. Questa volta, però, non
mostrava subito le figure ma ne diceva prima il loro nome e
controllava che i bambini colorassero la figura giusta;
successivamente mostrava anche la figura e consentiva loro,
eventualmente, di correggersi.

L'insegnante ha notato con piacere che i bambini non hanno avuto
bisogno di correggersi. Questo gioco le ha permesso di verificare
quanto i bambini avessero appreso nel corso delle tre attività
precedenti. Il risultato è stato indubbiamente positivo e ha
confermato che il ``divertirsi facendo'' è un'ottima motivazione a
apprendere.


\section[Sperimentazione \#6: prima primaria]{Sperimentazione \#6:
  classe prima primaria, febbraio~2010}

\subsection{Osservazioni generali}

\subsubsection{Presentazione della classe}
23 alunni (un'alunna è assente durante tutta la sperimentazione)

\subsubsection{Composizione dei gruppi}
Nel primo, secondo e terzo incontro i gruppi erano 4, formati da 5 o 6
alunni, eterogenei, scelti dall'insegnante. Uno di questi è stato
composto volutamente da bambini un po' più ``tranquilli'', con tempi di
interiorizzazione e di lavoro leggermente più lenti degli
altri. L'intento era di permettere loro di agire senza fretta e senza
l'esuberanza di qualche compagno che avrebbe potuto inibire la loro
creatività. L'insegnante ha avuto dei dubbi sull'effettuare o meno
questa scelta, ma a posteriori la riproporrebbe. È stata la prima
occasione che la classe ha avuto per lavorare in
gruppi%
.

Durante la quarta attività i bambini hanno lavorato a coppie,
mantenendo quelle già presenti con la disposizione dei banchi
utilizzata durante le consuete lezioni scolastiche.

\subsubsection{Insegnanti presenti}

Solo l'insegnante di classe. All'ultimo incontro era prevista la
partecipazione di una collega, ma poi quest'ultima è stata destinata
alla sostituzione di un'insegnante assente.

\subsubsection{Calendarizzazione degli incontri}
\begin{calendario}
  \begin{itemize}
  \item 22 febbraio
  \item 23 febbraio
  \item 25 febbraio
  \item 26 febbraio
  \end{itemize}
\end{calendario}

\subsection{Primo incontro}
\begin{description}
\item[Alunni presenti:] 22 presenti, 1 assente
\item[Tempo effettivo di lavoro:] Circa un'ora e mezza, dalle 14.30
  circa alle 16.15 circa
\end{description}

\begin{consegna}
  I bambini non hanno mai lavorato con dei poligoni, per cui
  l'insegnante intende proporre la sperimentazione come primo
  approccio alla geometria piana. Inoltre decide di prendere,
  soprattutto in questo primo incontro, il tempo necessario per
  lavorare con calma e creare l'atmosfera più proficua. Prima di
  cominciare l'insegnante spiega che si lavorerà con la geometria,
  spiegando a grandi linee che cosa tratta questa disciplina. A
  posteriori si rende conto di non aver fatto un collegamento con il
  lavoro sulle linee chiuse e aperte, regioni interne e regioni
  esterne proposto qualche settimana prima. A questo punto
  l'insegnante spiega che i bambini dovranno fare ordine tra le
  figure, che sono tante e diverse ma se guardate bene possono essere
  messe insieme, perché alcune hanno qualcosa in comune. Esprime
  subito la regola che non si può raggruppare per colore; a posteriori
  pensa che forse non avrebbe dovuto anticiparlo e che sarebbe stato
  meglio far emergere questa esigenza dalla
  sperimentazione. L'insegnante richiama un'attività di
  classificazione di animali svolta di recente con l'insegnante di
  scienze e poi dà inizio al lavoro.

  L'attività ha ripreso quella del kit, con la sola esclusione della
  figura del cerchio.

  \materiali{}%
  Al centro del tavolo di lavoro di ogni gruppo l'insegnante pone alla
  rinfusa i poligoni del kit, cercando di darne a tutti lo stesso
  numero con la stessa varietà di forme. Volutamente lascia da parte
  il cerchio (con il quale lavoreranno nel terzo incontro), sia perché
  desidera che gli alunni dirigano la loro attenzione sugli angoli e
  sui lati, sia perché ritiene che la consegna sia già alquanto
  impegnativa e immagina che il cerchio potrebbe porre dei quesiti che
  devierebbero la riflessione.
\end{consegna}

\subsubsection{Osservazioni}
I bambini appaiono fin da subito molto curiosi: ancora prima di
iniziare l'attività hanno notato la scatola arancione sul banco vicino
alla cattedra e vogliono scoprire cosa c'è di nuovo. Il lavoro dei
gruppi parte a stento, non si sa bene da dove cominciare\dots{}
qualcuno si impossessa di quante più figure può perché vuole creare
una figura composta. Questo fa prendere coscienza all'insegnante che
sarebbe stato meglio partire con la manipolazione dei
poligoni. Addirittura una bimba non resiste a assaggiare qualche
piccolo pezzetto di un quadrato giallo (grrr!!).

L'insegnante interviene passando di gruppo in gruppo e chiarendo la
consegna. Dopo un po' di tempo un gruppo fa la prima proposta che
porterà a una bella riflessione: sembra che abbiano deciso di
raggruppare poligoni piccoli, grandi e medi%
.

Mentre gli altri si scervellano ancora un po' spaesati, l'insegnante
chiede a questo gruppo di far vedere quali poligoni apparterrebbero ai
raggruppamenti da loro individuati; a questo punto nasce una
discussione tra i componenti che non riescono a mettersi
d'accordo. L'insegnante fa notare allora che il criterio scelto non
deve essere arbitrario ma riconoscibile da tutti; gli alunni sembrano
convinti e vengono lasciati lavorare.

Passa qualche minuto e una bimba del gruppo arriva raggiante dicendo
che hanno deciso quali siano i piccoli, i medi e i grandi! A questo
punto diventa importante approfondire e l'insegnante chiede alla
classe di interrompere quello che sta facendo, scegliere tra i
poligoni a disposizione quelli piccoli e sollevarli in alto perché
siano ben evidenti. Non tutti hanno optato per la stessa soluzione (il
gruppo promotore ne disapprova un altro dicendo: \bambini{Ma quelle
  sono le forme medie!}).

L'insegnante inizia un breve dibattito ponendo delle domande:
\begin{tutor}[Ins]
  avete scelto i poligoni piccoli? Perché non avete individuato gli
  stessi? In questo modo abbiamo fatto ordine? Si capisce bene quali
  siano i poligoni piccoli o ci sono dei dubbi? Può essere utile
  trovare un metodo che ci dia la sicurezza di raggruppare senza
  incertezze?
\end{tutor}
Gli interventi sono appropriati; le risposte alle domande
dell'insegnante sono centrate ma non espandono il concetto. I bambini
sono pensierosi, in fase di interiorizzazione e ricerca. Non
sembravano molto convinti sull'utilità di mettere ordine, forse solo a
causa della loro età; hanno però portato avanti la consegna con
impegno e entusiasmo.

Dopo aver ascoltato bene il confronto avvenuto un bimbo propone di
contare le ``punte'' per mettere insieme le figure con 3 punte, 4
punte, ecc. I gruppi ricominciano a lavorare, sono concentrati e
arrivano presto alla realizzazione dei raggruppamenti. L'insegnante
chiede di sollevare le forme con 3 punte (successivamente le altre):
tutti hanno le stesse! Il criterio risulta valido!

Poi lo stesso bambino propone di contare anche le ``linee'' (lati) e si
procede. Qualcun altro decide di mettere insieme quelle forme che si
sovrappongono perfettamente.

Alla fine la classe ribadisce le regole scoperte e tutti insieme danno
il nome a ogni poligono (in alcuni casi necessita una votazione):
sembra incoraggiante che propongano dei nomi che riguardano le
caratteristiche geometriche osservate.

I bambini, sono entusiasti e motivati: hanno fatto scoperte attraverso
loro iniziative e esperienze! Sperano di continuare a lavorare con la
geometria!

\subsubsection{Consigli per i colleghi che vogliono proporre le stesse
  attività}
Nei resoconti delle attività l'insegnante ha inserito commenti
personali che pensa possano già costituire una sorta di
idea/consiglio. Aggiunge poi che crede che da tutto questo si debba in
generale cogliere lo stimolo a creare in tutte le discipline una
didattica più centrata sulla ricerca, sulla sperimentazione, sulla
proposta diretta da parte del bambino. Dà ottimi frutti e entusiasmo.

\subsection{Secondo incontro}
\begin{description}
\item[Alunni presenti:] 22 presenti, 1 assente
\item[Tempo effettivo di lavoro:] Circa un'ora e mezza, dalle 9 alle
  10.30.
\end{description}
\mbox{}
\begin{consegna}
  Per il secondo incontro l'insegnante ha pensato di proporre ai
  bambini la terza attività presente nel kit, invertendo pertanto la
  seconda e la terza (la quarta quindi resterà la \attivita{Tombola
    delle forme}). L'insegnante
  ha infatti pensato che ai bambini avrebbe giovato osservare e
  maneggiare ancora le forme con quest'attività ludica, prima di
  procedere all'approfondimento delle differenze tra i poligoni e alla
  loro rappresentazione.

  L'insegnante spiega le regole del gioco: si estrae la figura a turno
  passando il sacchetto e il ``separé'' al compagno di destra; anche
  le domande vanno poste a turno; chi pensa di aver indovinato può
  dire il nome della forma o indicarla sul tabellone. %
  \materiali{}%
  A ogni gruppo sono stati consegnati:
  \begin{itemize}
  \item \textbf{un foglio diviso in tre colonne} (nome del bambino che
    estrae la figura, bambino che indovina e forma individuata
    \dots{} l'insegnante ha riportato sulla lavagna i poligoni con i
    rispettivi nomi);
  \item il \textbf{tabellone} plastificato con i disegni delle figure;
  \item un \textbf{raccoglitore} aperto in piedi sul tavolo per non
    far vedere ai compagni la forma estratta;
  \item della \textbf{pasta} per segnare i poligoni già indovinati;
  \item infine il \textbf{sacchetto} contenente \underline{una figura
      per tipo}; l'insegnante non ha inserito tante forme come era
    previsto dal kit perché le sembrava che potessero dare troppe
    variabili per bambini di sei anni.
  \end{itemize}
  Riflettendo a posteriori, se l'insegnante dovesse ripetere questa
  attività lascerebbe da parte il foglio con le tre colonne, perché i
  bambini essendo di prima hanno fatto fatica a leggere e scrivere. Il
  rischio è stato quello di distogliere l'attenzione dall'obiettivo
  principale, considerando inoltre che questo foglio non era un
  elemento indispensabile all'attività.
\end{consegna}

\subsubsection{Osservazioni}
I gruppi
cominciano spediti e non si fanno distrarre dalla novità della
telecamera che inquadra il loro gioco. L'insegnante nota però che
molti bambini scelgono di non porre la domanda, ma di indicare la
forma direttamente sul tabellone, chiedendo ``è questa?''.

Decide quindi di passare in tutti i gruppi, ricordando di chiedere ciò
che può essere utile per individuare più velocemente possibile la
figura estratta. L'insegnante ha poi notato che i bambini avevano
difficoltà a gestire da soli l'attività, dovendosi concentrare su una
pluralità di aspetti: il sacchetto, il raccoglitore, il foglio per
segnare i nomi, la pasta, il tabellone\dots{} Li ha quindi aiutati
moderando una parte del gioco, senza suggerire gli elementi essenziali
dei quesiti; si è avvicinata a un gruppo alla volta, intervenendo con
sollecitazioni del tipo
\begin{tutor}[Ins]
  a chi tocca? cosa bisogna fare ora?
\end{tutor}
Per i bambini era importante concentrarsi sulla consegna
principale. L'insegnante li ha aiutati con interventi tipo:
\begin{tutor}[Ins]
  \begin{itemize}
  \item Silvia poni una domanda a Jacopo che ti permetta di
    individuare la forma estratta
  \item bene, ora tocca a Lorenzo
  \item Nicol, ti sembra di aver capito? Che nome abbiamo dato al
    poligono che hai in mente?
  \item Indovinato! ORA potete trovarlo sul tabellone e segnarlo con
    un pezzo di pasta
  \end{itemize}
\end{tutor}
L'intervento dell'insegnante è stato utile per rendere più autonomo il
proseguimento del lavoro.

Spesso si ripete una situazione-tipo: il primo bambino che deve
parlare chiede \bambini{Quante punte ha?}
quindi il secondo: \bambini{Quante linee ha?}  (qualcuno usa la parola
righe). Succede, sempre o quasi, che il terzo indovini.

L'attività si conclude con una discussione in grande gruppo, in cui si
provano a analizzare le novità apprese sulle forme geometriche:
intervengono molti bambini. Si propongono soprattutto osservazioni
confrontando il numero di punte presenti sulle varie forme, si usa
l'espressione \bambini{assomiglia} o \bambini{è uguale perché} o
\bambini{sarebbe uguale se}\dots{} c'è attenzione e interesse. Si
``aggiusta'' qualche nome dato alle forme il giorno precedente. Viene
riportata una parte della discussione ripresa dalle telecamere:
\begin{tutor}[Ins]
  Avete scoperto qualcosa di interessante su questi poligoni, su queste forme?
\end{tutor}
\begin{studente}[Stefano]
  Delle forme possono essere anche un po' uguali
\end{studente}
\begin{tutor}[Ins]
  Prova a dirci in che senso\dots{}
\end{tutor}
\begin{studente}[S]
  Il rotondo (ottagono) e il cerchio sono quasi uguali, l'unica differenza è che il cerchio è tutto rotondo e l'ottagono è a strisce
\end{studente}
\begin{studente}[Lorenzo]
  Il PI 6 (esagono) e il PI 5 (pentagono) hanno la stessa forma perché
  hanno la punta in su
\end{studente}
\begin{tutor}[Ins]
  Forse al PI 5 manca però la punta in giù?
\end{tutor}
\begin{studente}[L]
  Sì!
\end{studente}
\begin{studente}[Edoardo]
  PI 5 e matitone (pentagono storto) hanno tutto uguale ma il matitone
  è un po' più alto
\end{studente}
\begin{studente}[Giovanna]
  Il rombo (triangolo isoscele) è solo un po' più alto del triangolo
  (triangolo equilatero)
\end{studente}
\begin{tutor}[Ins]
  E il PI 3 (triangolo rettangolo) ha qualcosa di simile agli altri
  due?
\end{tutor}
\begin{studente}[G]
  Ha le punte come loro: 3 punte, sembra messo al contrario
\end{studente}
\begin{studente}[Matteo]
  PI 4 (parallelogramma rettangolo) e PI 3 (triangolo rettangolo) sono un po' uguali perché se fai un po' alzare, allungare il PI 4 \dots{} sono uguali \dots{}
\end{studente}
(non riesce a proseguire)(\dots{} provo io, senza aver capito bene e
guardo le figure\dots{})
\begin{tutor}[Ins]
  In che modo alziamo\dots{} allunghiamo?
\end{tutor}
(ho capito!)
\begin{tutor}[Ins]
  allunghiamo il lato storto e quello diritto fino a farli
  toccare\dots{}
\end{tutor}
\begin{studente}[M]
  Sì, sì! Così!
\end{studente}
\begin{studente}[Anita]
  Se al PI 6 raddrizziamo le linee che formano la punta in basso
  diventa come il PI 5
\end{studente}
Da una valutazione finale emerge il buon esito dell'attività, seppur
svolta in un momento diverso rispetto a quanto proposto dal kit.

\subsection{Terzo incontro}
\begin{description}
\item[Alunni presenti:]21 presenti, 2 assenti
\item[Tempo effettivo di lavoro:]Due ore, negli intervalli di tempo
  8.45-10.15 e 11.00-11.30
\end{description}

\begin{consegna}
  La terza proposta corrisponde alla seconda del kit. L'insegnante
  riprende con la classe i discorsi dei giorni precedenti ricordando
  le regole che hanno già scoperto: esistono poligoni con tre,
  quattro, cinque\dots{} punte e anche con tre, quattro, cinque\dots{}
  linee. Propone allora di fare bene ordine e raccogliere in ciascun
  cartellone solo le figure geometriche che hanno lo stesso numero di
  punte, poi quelle che hanno lo stesso numero di linee.

  \materiali{}%
  Prima di esplicitare la consegna l'insegnante prepara l'aula con i
  bambini, liberandola dagli arredi per permettere il lavoro sul
  pavimento. C'è entusiasmo per la novità mentre tutti i banchi e le
  sedie vengono spinti verso il perimetro dell'aula. L'insegnante ha
  ritenuto più opportuno far lavorare i bambini sul pavimento, per
  dare loro più spazio per i cartelloni, evitando gli intralci ai
  movimenti che possono dare banchi e sedie.

  Ogni bambino ha una matita e una gomma e a ogni gruppo è consegnato
  un grande foglio bianco.
\end{consegna}

\subsubsection{Osservazioni}
Inizia la classificazione: l'insegnante chiama il primo gruppo, che
deve prendere dal mucchio di tutti i poligoni riuniti alla rinfusa su
un banco solo quelli con 3 punte. In quest'occasione l'insegnante
inserisce per la prima volta il cerchio, non esplicitando ai bambini
questa novità per osservare le loro reazioni.

I bambini sono molto attenti e contano bene le punte, anche facendosi
aiutare da un compagno per essere certi di non sbagliare. L'insegnante
appare stupita, perché immaginava che dopo le attività svolte in
precedenza i bambini procedessero con maggior sicurezza.

Raccolte con precisione tutte le forme i bambini tornano a
sedersi. Via via sfilano gli altri tre gruppi (uno per le forme con
quattro punte, uno per le cinque punte, uno ha il compito di prendere
le forme dalle sei punte in su) e la scelta dei poligoni si svolge con
le stesse modalità descritte per il primo gruppo.

Il momento più stimolante arriva quando un bambino che deve prendere i
poligoni con quattro punte si rivolge all'insegnante dicendo che anche
il cerchio ha quattro punte e toccando quattro punti (più o meno alla
stessa distanza tra loro) della circonferenza di un cerchio trovato
nel mucchio.

Il gruppo si ferma interessato.
L'insegnante propone a questo bambino, e alla compagna che intanto si
è unita a lui, di provare a contare se il cerchio possa avere anche
sei punte. Contano attenti e rispondono di sì incuriositi. Chiede
allora di vedere se possa averne dieci; dopo aver contato questi
rispondono nuovamente in modo affermativo. L'insegnante chiede ancora
se secondo loro il cerchio possa arrivare a avere addirittura venti o
trenta punte e i bambini rispondono \bambini{sì!}. Senza grandi
commenti ma soddisfatti tornano al posto lasciando i cerchi con gli
altri poligoni%
.

Dopo che tutti i gruppi hanno scelto le forme restano sul banco solo i
cerchi; l'insegnante chiede ai bambini il motivo e se secondo loro si
potrebbero inserire in uno dei cartelloni; coinvolge il bambino citato
in precedenza dicendogli che può essere di molto aiuto la sua scoperta
sul cerchio. A questo punto si apre uno scambio breve ma chiaro dove
intervengono alcuni bambini e il cerchio viene assegnato al gruppo che
deve riunire sul cartellone le forme con più di sei punte.

Un altro momento interessante riguarda il momento in cui il gruppo
delle 5 punte non sa se portare via il poligono a ``L'': ha 5 o 6 punte?
L'insegnante prova a utilizzare l'intuito dei bambini e chiede alla
classe di votare: la maggior parte ritiene che abbia 6 punte, così il
poligono viene assegnato al gruppo corrispondente senza ulteriori
spiegazioni in quanto poco dopo la classe ne avrà la conferma
occupandosi dei lati.

Vengono riportate sui quattro cartelloni le forme: ogni gruppo ricalca
e colora quelle a cui è stato assegnato.

Il passo successivo consiste nel preparare un cartellone anche per la
classificazione per numero di lati, come accordato con i bambini
all'inizio del lavoro. L'insegnante propone allora di riportare tutti
i poligoni sul banco dove erano posti all'inizio e chiama il primo
gruppo, che sceglie quelli con tre ``linee'' ; prendendo i primi
triangoli senza indugiare i bambini esclamano con sorpresa:
\begin{studente}[]
  Sono le
stesse figure che abbiamo preso prima!
\end{studente}
La stessa magia si verifica anche per gli altri gruppi e l'insegnante
fa loro i complimenti per aver scoperto delle importanti regole
geometriche. Applauso. Si può fare merenda!

Dopo ricreazione resta circa mezz'ora e l'insegnante propone ai
bambini di sedersi in cerchio a terra per provare a scoprire qualcosa
sulle superfici mettendo a confronto i poligoni: quale occupa una
superficie maggiore? Quale la superficie minore? La scelta di stare
sul pavimento non è casuale e è stata pensata per facilitare una
riflessione sulle superfici. I bambini sono stanchi ma la maggior
parte partecipa con voglia: si prendono le forme che occupano ``più
spazio'' e sovrapponendone altre in vari tentativi si trovano quelle
che occupano ``un po' meno spazio'' e poi quelle che occupano ``meno
spazio'' delle altre. L'insegnante fa poi notare che la mattonella ha
la forma di uno dei poligoni e i bambini individuano giustamente il
``dado'' (quadrato). Chiede allora quale dei due dadi (piastrella o
quadrato del kit) occupa una superficie maggiore\dots{} si conclude
che tutti i poligoni che la classe ha conosciuto possono assumere
dimensioni diverse ma mantenere le stesse caratteristiche e quindi
continuare a appartenere alla classificazione che abbiamo riportato
sui cartelloni.

I bambini sono molto soddisfatti!

All'insegnante resta ancora il dubbio se insegnare ai bambini i nomi
geometrici delle figure e delle parti che li compongono.

\subsection{Quarto incontro}
\begin{description}
\item[Alunni presenti:]21 presenti, 2 assenti
\item[Tempo effettivo di lavoro:] 1 ora, dalle ore 14.30 circa alle
  ore 15.30 circa
\end{description}

\begin{consegna}
  L'insegnante annuncia ai bambini che oggi giocheranno a un gioco
  che è come la tombola: ha infatti le stesse regole di quella che già
  molti conoscono ma non si estrarranno dei numeri, bensì dei
  poligoni! Un bambino distribuisce le cartelle (una per
  coppia). Erano già stati preparati sui tavoli i pennarelli a punta
  grossa e insieme capiscono che verranno utilizzati per colorare le
  forme che a mano a mano verranno estratte dal
  sacchetto. L'insegnante invita i bambini a cambiare colore per ogni
  poligono, mettendosi d'accordo con il vicino.

  \materiali{}%
  Cartelle, pennarelli, forme. Rispetto al materiale fornito dal kit
  l'insegnante preferisce inserire nel sacchetto una figura per ogni
  tipo (un solo quadrato, un solo rettangolo ecc.) sia per velocizzare
  il gioco, sia per evitare che i bambini andassero un po' in
  confusione, vedendo uscire una figura già colorata. Scelta poi
  rivelatasi positiva in quanto l'evidenziarsi di vari dubbi,
  soprattutto prima che iniziasse il gioco, ha rallentato l'attività e
  dimostrato che non per tutti fosse immediato capire.
\end{consegna}

\subsubsection{Osservazioni}
L'insegnante chiede%
, prima di iniziare a giocare, se qualche coppia
abbia notato sulla propria scheda poligoni a cui non avevano ancora
dato un nome, perché non incontrati nei giorni scorsi. Qualcuno alza
la mano per indicare il rettangolo che a maggioranza viene chiamato
col nome di ``salame'' (a questo proposito l'insegnante fa notare: è
un errore del kit o un fatto voluto aver inserito il rettangolo solo
nell'attività della tombola e non nelle
precedenti?%
). Una bimba
sollecitata dall'insegnante e aiutata da altri compagni classifica il
rettangolo considerando punte e linee. Viene anche evidenziata qualche
somiglianza con altri poligoni.

Finalmente il gioco inizia e procede senza particolari situazioni da
commentare. L'unico momento in cui si è dovuto interrompere l'attività
e ragionare insieme è stato quello in cui, durante la prima partita, è
stato estratto l'ottagono: un paio di coppie ha colorato l'esagono e
una terza coppia si è fermata molto incerta sul da farsi. Con calma e
con la partecipazione della classe sono state contate bene le punte e
le linee e si è sciolto con facilità il dubbio. L'insegnante ha poi
tranquillizzato tutti, esplicitando che nel caso in cui fosse stato
estratto l'esagono i gruppi che già lo avevano colorato si sarebbero
trovati il lavoro fatto, in caso contrario non sarebbe valso per
un'eventuale tombola\dots{} la figura estratta subito dopo è proprio
l'esagono (\dots{} e senza barare)!

La classe è riuscita a completare due partite, ha avuto un buon numero
di coppie felici di aver completato la tombola e una coppia vincitrice
di entrambe le manche. Gli alunni sono stati contenti e molto
partecipi anche in questo lavoro, hanno dimostrato dispiacere nel
constatare che il kit stava partendo per andare in altre scuole, ma
l'insegnante ha promesso loro che giocheranno ancora con la geometria.

\subsubsection{Consigli per i colleghi che vogliono proporre le stesse
  attività}
L'insegnante ha trovato bello e utile conservare le schede della
tombola e la settimana successiva, rielaborando e riassumendo con i
bambini il percorso impegnativo affrontato insieme, lasciare sul
quaderno di matematica e geometria una traccia concreta dell'attività
svolta. La classe ha quindi scritto due o tre semplici frasi per
spiegare (anche ai genitori che le avrebbero lette), ricordare e
interiorizzare ancora i concetti trattati e incollato la scheda della
tombola. I bambini si sono divertiti molto anche a tracciare a mano
libera su una pagina di quaderno i poligoni divisi a seconda della
nostra classificazione. L'hanno fatto seguendo le indicazioni orali
dell'insegnante e le linee che contemporaneamente mostrava alla
lavagna (andiamo verso l'alto di due quadretti, poi a sinistra di tre
quadretti ecc.): un ottimo esercizio anche per l'orientamento
spaziale!


\section[Sperimentazione \#7: terza primaria]{Sperimentazione \#7:
  classe terza primaria, marzo~2010}
\subsection{Osservazioni generali}

\subsubsection{Presentazione della classe}
La classe è formata da 24 alunni. Nella classe c'è una bambina cinese
che ancora non conosce la lingua italiana: l'insegnante afferma che è
stato un po' difficoltoso riuscire a coinvolgerla in queste attività.

\subsubsection{Composizione dei gruppi}
Si sono formati 4 gruppi da sei alunni. I gruppi che l'insegnante ha
formato sono eterogenei: i bambini più vivaci sono stati distribuiti
nei diversi gruppi. L'ultima attività è stata svolta a coppie.

\subsubsection{Insegnanti presenti}
Nel primo incontro la docente sperimentatrice si è fatta affiancare
dall'altra docente di classe (compresenza) perché ha ritenuto in
questo modo di poter organizzare meglio i gruppi e sistemare l'aula in
modo adeguato.

\subsubsection{Calendarizzazione degli incontri}
\begin{calendario}
  \begin{itemize}
  \item 3 marzo  (compresenza)
  \item 4 marzo
  \item 5 marzo
  \item 8 marzo
  \end{itemize}
\end{calendario}

\subsection{Primo incontro}
\begin{description}
\item[Alunni presenti:]24 presenti (tutti)
\item[Tempo effettivo di lavoro:]2 ore, dalle 10,30 alle 12,30
\end{description}

\begin{consegna}
  Prima di distribuire il materiale l'insegnante ha detto ai bambini
  che dovevano lavorare nel rispetto del compagno e del materiale che
  avevano in ``cura''. Quindi ha spiegato loro che, se osserviamo con
  attenzione la realtà, ci accorgiamo che esistono forme geometriche
  tra le più svariate e che l'uomo ha deciso di classificarle in base
  a delle caratteristiche, per dare un ordine a tutto ciò che si
  osserva. Le insegnanti hanno poi consegnato un certo numero di
  forme, che gli alunni avrebbero dovuto osservare, analizzare e
  classificare ``come fanno gli adulti''. Dato un po' di tempo per la
  classificazione libera, sono stati consegnati a ogni gruppo un
  foglio e una penna, chiedendo di rappresentare la classificazione
  del materiale e dare un nome anche alle forme che i bambini ancora
  non conoscevano. Successivamente le insegnanti hanno chiesto di
  mettere insieme diverse forme e formare qualche oggetto da
  utilizzare come mascotte del gruppo.

  La prima attività non ha ripreso quella del kit perché l'insegnante
  ha preferito svolgere attività di manipolazione, classificazione e
  conoscenza delle diverse figure geometriche date.

  \materiali{} %
  Le forme geometriche previste nel kit, fogli e penne
\end{consegna}

\subsubsection{Osservazioni}
Quest'anno è la prima volta che gli alunni lavorano in gruppo e la
proposta di lavorare con un nuovo materiale e in gruppo è stata
accettata con entusiasmo. I bambini già conoscevano ``i blocchi logici''
che comprendono quattro figure geometriche. La novità per loro è stata
la diversità del materiale e la varietà delle figure geometriche che
esistono nella realtà; osservando quella grossa scatola non vedevano
l'ora di iniziare.

L'idea di agire come gli adulti entusiasma abbastanza i bambini. La
prima cosa che fanno i diversi gruppi è la manipolazione e
osservazione di diverse figure; sono maggiormente attratti dalle
figure più grandi.

Poi contano quante figure hanno e iniziano a classificarle in base al
colore, alla dimensione e al numero di lati. C'è sempre il bambino che
vuole ``primeggiare'' sugli altri ma, grazie alla discussione avuta in
precedenza su come si lavora in gruppo, gli alunni sono riusciti a
rispettare i tempi e i pareri di tutti i compagni.

All'inizio c'erano anche bambini che usavano un tono di voce più alto;
l'insegnante è intervenuta promettendo agli alunni che alla fine del
lavoro avrebbe premiato il gruppo che avrebbe lavorato meglio e quasi
tutti si sono impegnati per essere ``quel gruppo''.

Nella seconda parte dell'attività, i bambini hanno inventato dei nomi
per le figure che ancora non conoscevano: il trapezio rettangolo lo
hanno chiamato \bambini{scarpa}, l'esagono \bambini{diamante}, il
trapezio isoscele \bambini{tetto}. Alla fine si sono divertiti a
formare ``robot'' dando dei nome di fantasia.

\`E stato un lavoro per loro divertente e hanno scoperto come nella
realtà parecchi oggetti non sono lo specchio di una ``forma geometrica
semplice'', ma l'insieme di più forme geometriche. Ogni gruppo ha poi
deposto il materiale manipolato in un sacchetto e lo ha conservato per
il giorno successivo.

\subsubsection{Consigli per i colleghi che vogliono proporre le stesse
  attività}
L'insegnante ha ritenuto quest'attività molto utile per aiutare il
bambino a familiarizzare con la geometria, in particolare con le
figure piane, attraverso il gioco\dots{} Consiglia di lasciare che i
bambini utilizzino il materiale anche costruendo forme di
fantasia%
; quando descriveranno ciò che hanno realizzato lo
faranno con soddisfazione e maggiore conoscenza.

\subsection{Secondo incontro}
\begin{description}
\item[Alunni presenti:]24 alunni (tutti)
\item[Tempo effettivo di lavoro:]2 ore circa (14:30-16:15)
\end{description}

\begin{consegna}
  L'attività prevista è
  \attivita{Indovina chi?}
  A turno un bambino nel gruppo prende in mano il sacchetto e si
  nasconde dietro una cartelletta, così può osservare con attenzione
  la forma che ha pescato. Gli altri bambini a turno devono fare una
  domanda per indovinare la forma; chi ha l'oggetto in mano deve
  rispondere con un SÌ o NO. Non si può dire il nome della forma nelle
  domande, ma la si può nominare solo quando si tratta dare la
  risposta. Un bambino è incaricato di compilare il modulo scrivendo
  le domande che vengono poste. L'attività ha ripreso in toto quella
  del kit

  \materiali{} %
  A ogni gruppo vengono consegnati: un sacchetto con le forme
  geometriche che hanno manipolato in precedenza, una tabella
  plastificata con le forme disegnate e un foglio con un tabulato da
  compilare man mano che ``pescano'' le forme da far indovinare, così
  come previsto dal kit per l'attività
\end{consegna}

\subsubsection{Osservazioni}
La maggior parte dei bambini ha domandato il colore della figura e poi
ha tirato a indovinare sulla tabella plastificata. Successivamente si
è discusso sulla forma delle figure geometriche, cercando di capire se
il colore è rilevante o se ci sono altre caratteristiche più
importanti per identificare una figura, come il numero di lati. Questo
lavoro potrà essere svolto in maniera migliore dopo aver studiato
l'angolo e le sue caratteristiche%
. L'insegnante ha notato che i bambini hanno mantenuto il nome di
alcune forme che era stato dato il giorno precedente.

\subsubsection{Consigli per i colleghi che vogliono proporre le stesse
  attività}
Prima
di iniziare quest'attività sarebbe utile discutere sulle
caratteristiche che differenziano le diverse forme geometriche

\subsection{Terzo incontro}
\begin{description}
\item[Alunni presenti:]24 presenti (tutti)
\item[Tempo effettivo di lavoro:]2 ore circa (14:30-16:15)
\end{description}

\begin{consegna}
  Si consegnano i sacchetti e l'insegnante fa ripetere il gioco del
  giorno precedente; i bambini dovranno però cercare questa volta di
  fare domande più precise per indovinare la forma. Dopo più di
  un'ora, l'insegnante distribuisce ai bambini coppie di forme come
  quelle che sono segnate sul foglio consegnato all'inizio del gioco
  (rombo-pentagono, ottagono-cerchio); bisogna osservare con
  attenzione e compilare la scheda, rispondendo alle domande. (si
  vedano gli allegati)

  \materiali{}%
  Gli stessi del terzo incontro
\end{consegna}

\begin{figure}[pht]
  \label{ex:scheda:forme:1}
  \centering
  \fbox{\includegraphics[width=0.90\textwidth]{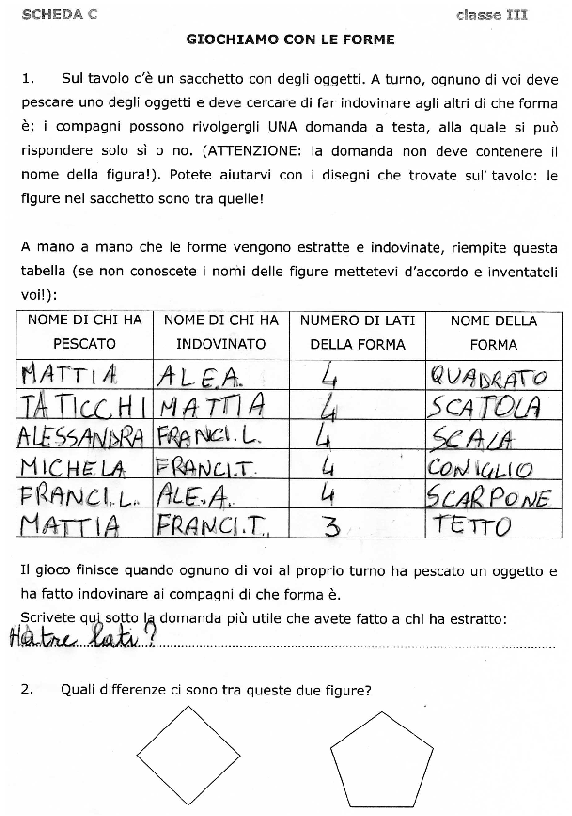}}
\end{figure}
\begin{figure}[pht]
  \label{ex:scheda:forme:2}
  \centering
  \fbox{\includegraphics[width=0.90\textwidth]{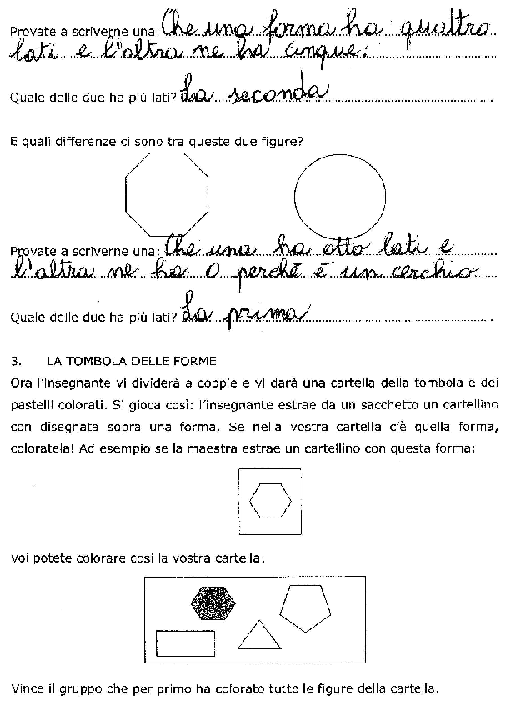}}
\end{figure}

\subsubsection{Osservazioni}
Rispetto all'incontro precedente, i bambini hanno fatto domande più
dettagliate per indovinare le forme, del tipo:
\begin{studente}[ ]
  \begin{itemize}
  \item \bambini{assomiglia a un marciapiede?}
  \item \bambini{assomiglia a una montagna?}
  \item \bambini{ha quattro vertici?}
  \item \bambini{sembra uno scudo?}
  \end{itemize}
\end{studente}
Si nota come le forme geometriche vengono paragonate agli oggetti
presenti nella realtà%
.

Nella seconda parte dell'attività i bambini si sono impegnati con
attenzione e interesse per rispondere correttamente al questionario,
analizzando la figura che avevano in mano. A un gruppo è stato dato
un esagono al posto di un ottagono\dots{} dopo circa dieci minuti i
bambini sono andati dall'insegnante per comunicarle che la forma era
di sei lati al posto di otto e perciò non uguale a quella disegnata
sul questionario.

\subsubsection{Consigli per i colleghi che vogliono proporre le stesse attività}
Prima di far compilare il questionario, discuterne in classe

\subsection{Quarto incontro}
\begin{description}
\item[Alunni presenti:]24 presenti (tutti)
\item[Tempo effettivo di lavoro:]2 ore circa (10:15-12:15)
\end{description}

\begin{consegna}
  Il gioco di oggi è la tombola delle forme, prevista anche nel
  kit. L'insegnante consegna alle coppie di bambini delle fotocopie di
  cartelle con forme geometriche; man mano l'insegnante pescherà da un
  sacchetto una forma e i bambini dovranno colorare la forma pescata
  se è disegnata sulla loro cartelletta. Vince chi completa per primo
  la sua cartella.

  \materiali{} %
  Cartelle della tombola e forme previste dal kit.
\end{consegna}

\subsubsection{Osservazioni}
Quest'attività ha lo scopo di far acquisire ai bambini maggior
confidenza con le forme piane e di consolidare la capacita di
riconoscere le forme; viene inoltre richiesto di saperle
disegnare. L'insegnante ha richiesto ai due bambini di alternarsi nel
colorare la forma e di colorarla correttamente, senza uscire dai
margini o pasticciare.

Il gioco è proceduto bene; l'entusiasmo è aumentato quando a qualche
coppia mancava una forma per completare la cartella. Ogni volta la
coppia che vinceva doveva dire il nome delle forme che aveva colorato;
l'insegnante consegnava loro un piccolo premio e il gioco
ricominciava. \`E stato un gioco divertente, stimolante. I bambini
hanno riconosciuto e denominato tutti le figure senza difficoltà.

\subsubsection{Consigli per i colleghi che vogliono proporre le stesse
  attività}
Molti bambini colorano in modo frettoloso; aggiungere maggiore
precisione nel colorare è da stimolo per una maggiore
concentrazione. Poi man mano che si scoprono le altre figure le forme
possono aumentare e si possono costruire così nuove tabelle.

\subsection{Osservazioni finali}
Il lavoro di gruppo è stato importante per favorire sia lo scambio
delle idee tra i bambini che l'autocontrollo di ogni bambino. L'alunno
doveva, per rendere il lavoro proficuo, abituarsi a controllare sia la
propria voce che quella del compagno in modo che non diventasse troppo
alta; doveva rispettare il proprio turno di parola negli
interventi. \`E stata perciò questa un'attività che ha visto la classe
impegnata nell'apprendimento cooperativo e metacognitivo della
matematica, che potrebbe essere riassunta come ``Impariamo insieme a
conoscere e distinguere il mondo delle forme intorno a noi''. Obiettivo
principale di queste attività è stata sicuramente per ogni bambino la
costruzione del proprio sapere.

\section[Sperimentazione \#8: infanzia]{Sperimentazione \#8: scuola
  dell'infanzia, marzo~2010}

Questa sperimentazione, pur utilizzando il kit ``Giochiamo con le
forme'', introduce in realtà anche attività proposte dal kit ``Torri,
serpenti e\dots{} geometria''.

\subsection{Osservazioni generali }

La sperimentazione si è svolta in una sezione di una scuola
dell'infanzia.

\subsubsection{Presentazione della classe }
Un gruppo eterogeneo composto da 7 bambini di 5 anni e 5 di 4 anni

\subsubsection{Composizione dei gruppi }
Diversi a seconda delle attività: nella prima attività i bambini hanno
lavorato individualmente e poi formato coppie con compagni liberamente
scelti; nella seconda e nella terza attività hanno lavorato in 2
gruppi eterogenei. L'insegnante ha invitato i bambini a organizzarsi
in modo che in ogni gruppo fossero presenti bambini di 4 e 5 anni e
con differenti capacità in modo da bilanciare le risorse interne ai
gruppi. Nell'ultima attività i bambini hanno lavorato individualmente.

\subsubsection{Insegnanti presenti}

Solo l'insegnante di classe

\subsubsection{Calendarizzazione degli incontri}
\begin{calendario}
  \begin{itemize}
  \item 16 marzo
  \item 23 marzo
  \item 31 marzo
  \item 16 aprile
  \end{itemize}
\end{calendario}

\subsection{Primo incontro }

\begin{description}
\item[Alunni presenti:] 11 presenti, 1 assente
\item[Tempo effettivo di lavoro:] 1 ora, dalle 11 alle 12, nel
  laboratorio di linguaggio, più mezz'ora di rielaborazione in classe
\end{description}

\begin{consegna}
  Presentazione dei triangoli equilateri del kit e esplorazione
  libera delle forme, con la possibilità di costruire liberamente
  delle figure. Riproduzione delle figure proposte in ``Torri, serpenti
  e geometria'': stella, vaso, gatto, aquilone, clessidra, serpente,
  della scheda per i bambini di prima elementare.

  \materiali{} %
  I triangoli equilateri del kit ``Giocare con le forme'', triangoli
  equilateri in cartoncino costruiti dalle insegnanti, triangoli
  equilateri in carta colorata per le rielaborazioni successive.
\end{consegna}

\subsubsection{Osservazioni}
Inizialmente si pensava di proporre l'attività solo ai bambini di
cinque anni, ma alcuni di quattro hanno insistito per essere
coinvolti. Il materiale è stato accolto con entusiasmo, soprattutto
quello del kit per i suoi colori e la consistenza
tattile. Inizialmente ogni bambino ha tentato di costruire qualcosa da
solo, anche sottraendo triangoli ai compagni, poi pian piano i bambini
hanno cominciato a lavorare a coppie, alcuni spontaneamente, altri su
indicazione dell'insegnante.

Dopo la libera costruzione sono state presentate le immagini contenute
nel kit: ``Torri, serpenti e geometria''. Alcuni bambini hanno
prodotto figure somiglianti a quelle campione, mentre altri le hanno
riprodotte fedelmente. Una bambina ha supportato un compagno più
piccolo ruotando la struttura costruita e permettendogli di completare
con successo il lavoro; uno dei bambini di cinque anni ha avuto
bisogno di essere stimolato a partecipare.

I bambini di 5 anni, su richiesta dell'insegnante, hanno contato i
triangoli contenuti nelle figure e gli spigoli (punte) esterne; non
hanno incontrato particolari difficoltà con il numero dei triangoli,
ma per gli spigoli qualcuno li ha dovuti contare più volte. In questa
situazione, dopo aver osservato la strategia utilizzata, l'insegnante
è intervenuta fornendo un punto di riferimento all'interno della
figura, da cui partire con il conteggio.

Successivamente in classe è stato proposto ai partecipanti di
riprodurre le figure su un foglio con triangoli di carta colorata
preritagliati, con la possibilità di usufruire o meno di modelli
costruiti precedentemente.

\subsubsection{Consigli per i colleghi che vogliono proporre le stesse
  attività}
L'insegnante ha trovato importante per la realizzazione di
quest'attività le seguenti condizioni:
\begin{itemize}
\item allestimento di un setting idoneo
\item piccoli gruppi
\item precedente acquisizione a livello percettivo della forma proposta.
\end{itemize}

\subsection{Secondo incontro }

\begin{description}
\item[Alunni presenti:] 11 presenti, 1 assente
\item[Tempo effettivo di lavoro:] 1 ora, dalle 11 alle 12, nel
  laboratorio di linguaggio
\end{description}

\begin{consegna}
  Offerta di triangoli equilateri e riproduzione delle figure proposte
  in: ``Torri, serpenti e geometria'': stella, vaso, gatto, aquilone,
  clessidra, serpente, della scheda per i bambini di prima elementare.

  Riproposizione della stessa attività e delle stesse figure da
  realizzarsi con triangoli isosceli e rettangoli. Analisi delle
  differenze fra i triangoli equilateri e quelli proposti in fase
  successiva. Scoperta delle differenze fra i diversi triangoli e dei
  possibili modi in cui si possono usare.

  \materiali{} %
  I triangoli equilateri, isosceli e rettangoli del kit ``Giocare con
  le forme'' e triangoli equilateri in cartoncino costruiti dalle
  insegnanti.
\end{consegna}

\subsubsection{Osservazioni}
I bambini hanno lavorato divisi in due sottogruppi di età eterogenea,
partecipando attivamente. Dapprima hanno riprodotto le figure con i
triangoli equilateri, come nell'attività precedente, poi sono stati
proposti a un gruppo i triangoli rettangoli e all'altro quelli
isosceli.

Con i triangoli isosceli i bambini hanno tentato di riprodurre la
parte centrale della figura del sole (vedi ``Torri, serpenti e
geometria''), notando che il ``cristallo'' era più grande e più lungo
di quello fatto con i triangoli equilateri.

Confrontando i due tipi di triangoli hanno subito osservato che i
triangoli isosceli erano più lunghi di quelli equilateri e rettangoli.

Una bambina di quattro anni si muoveva fra i due gruppi osservando e
riportando commenti ai compagni. Alcuni bambini, in coppia, sono
riusciti a costruire il gatto con i triangoli rettangoli, notando che
aggiungendo o togliendo altri triangoli questo si trasformava in una
casa.

A quest'attività non ha avuto seguito alcuna rielaborazione grafica
successiva.

\subsubsection{Consigli per i colleghi che vogliono proporre le stesse
  attività}
È necessario lavorare in piccoli gruppi, la fascia di età più idonea è
quella dei cinque anni, mentre il coinvolgimento dei bambini più
piccoli è più complesso e richiede un'attenzione e un supporto
maggiore. Nel gruppo dei bambini di quattro anni è più facile
incontrare difficoltà di ordine collaborativo e di organizzazione
spaziale.

\subsection{Terzo incontro }

\begin{description}
\item[Alunni presenti:] 10 presenti, 2 assenti
\item[Tempo effettivo di lavoro:] 1 ora, dalle 11 alle 12, nel laboratorio di linguaggio
\end{description}

\begin{consegna}
  Costruzione di figure con triangoli differenti (isosceli,
  rettangoli, equilateri) e scoperta delle forme che si possono
  realizzare. A entrambi i gruppi è stato chiesto di costruire alcune
  delle figure prodotte precedentemente, dando a chi lo volesse copia
  del disegno.

  Prove di pavimentazione.

  \materiali{} %
  Triangoli equilateri, isosceli e rettangoli del kit e triangoli
  costruiti con il cartoncino. Sono stati consegnati a un gruppo i
  triangoli del kit e all'altro quelli di cartoncino.

\end{consegna}

\subsubsection{Osservazione}
Prima sono stati consegnati i triangoli equilateri e molti bambini
hanno costruito delle stelle. Successivamente sono stati proposti i
triangoli isosceli, con i quali i bambini hanno costruito la stella,
notando che la parte centrale della figura era più grande e i raggi
erano più lunghi. Un bambino ha subito notato che la stella fatta con
triangoli isosceli aveva più pezzi di quella fatta con triangoli
equilateri%
. L'insegnante ha fatto notare agli altri bambini che il ``cristallo''
centrale aveva più pezzi.

Un bambino ha tentato di costruire una stella con triangoli isosceli
mantenendo inalterato il numero dei pezzi usati per costruire quello
precedente, ma è stato costretto a lasciare tra una tessera e l'altra
dello spazio.

Un bambino ha notato che la clessidra la si poteva costruire sia con i
triangoli equilateri che con i triangoli rettangoli%
.

Successivamente sono stati consegnati ai bambini dei fogli da riempire
completamente in modo da formare un ``pavimento''. I bambini sono
stati sollecitati a utilizzare un solo tipo di triangoli e hanno
visto che con quelli isosceli e equilateri non riuscivano a evitare
di lasciare ``buchi''.

In un caso l'insegnante ha aiutato un bambino a orientare nello
spazio due triangoli rettangoli per far loro occupare una posizione
idonea alla pavimentazione.

All'esperienza ha fatto seguito un'attività di pavimentazione a
collage in classe.

\subsubsection{Consigli per i colleghi che vogliono proporre le stesse
  attività}
Utilizzare molto materiale, in modo che tutti possano disporre di un
quantitativo adeguato; far giocare liberamente i bambini, in modo che
siano loro a cogliere le varie possibilità dei triangoli.

\subsection{Quarto incontro }

\begin{description}
\item[Alunni presenti:] 7 bambini di cinque anni. Per quest'attività
  il gruppo di lavoro è stato limitato ai soli bambini di cinque anni,
  sia per motivi di tempo, che per la necessità di soffermare
  l'attenzione su un particolare aspetto: la relazione fra il numero
  dei pezzi e la loro disposizione nello spazio.
\item[Tempo effettivo di lavoro:] 30 minuti, dalle 13,30 alle 14, in classe
\end{description}

\begin{consegna}
  Sono stati consegnati ai bambini i triangoli equilateri di
  cartoncino e è stato loro chiesto di realizzare singolarmente e
  senza avere davanti un modello la forma della
  stella. Successivamente è stato chiesto a ogni bambino di contare
  il numero di pezzi della stella e di costruire con gli stessi pezzi
  un serpente.

  \materiali{} %
  I triangoli equilateri in cartoncino costruiti dalle insegnanti
  simili a quelli del kit ``Giocare con le forme''. Triangoli
  equilateri in carta colorata per le rielaborazioni successive.
\end{consegna}

\subsubsection{Osservazioni}
\begin{figure}[bp]
  \centering
  \includegraphics[width=0.90\textwidth]{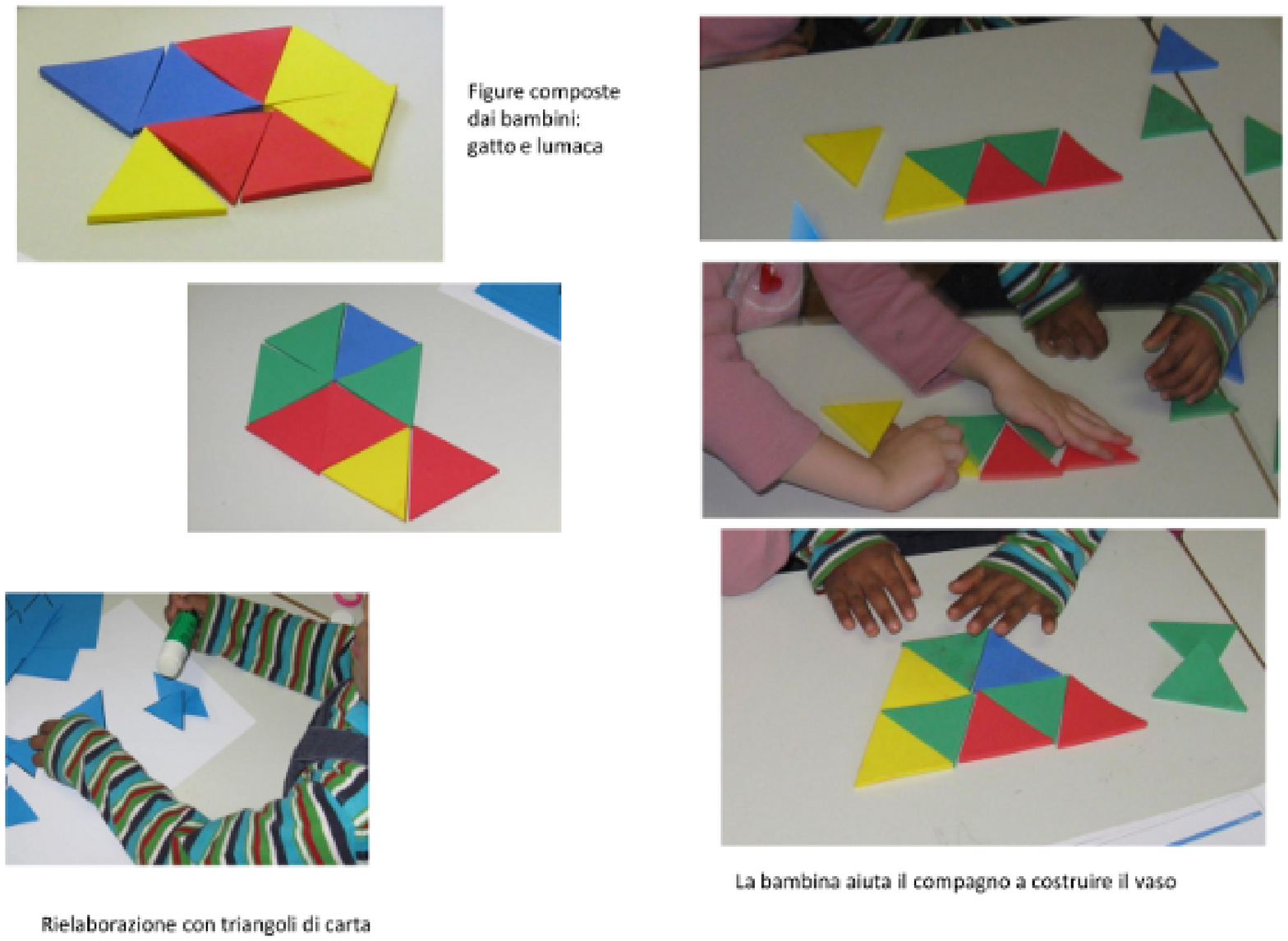}
  \label{pic:formemorbide:8}
\end{figure}
I bambini hanno costruito la stella in autonomia, senza bisogno di
immagini campione. Una bambina che durante l'anno ha avuto una scarsa
frequenza l'ha realizzata parzialmente e è stata aiutata da una
compagna. Un bambino ha iniziato a costruire la stella a partire dalle
punte (contrariamente a quanto aveva fatto quando aveva a disposizione
davanti a sé lo schema), per poi dover allargare i triangoli per
inserire l'esagono centrale; altri bambini hanno usato i triangoli
rimanenti per costruire piccole stelline con due triangoli sovrapposti
aggiungendo:
\begin{studente}[ ]
  Questa è la mamma e queste sono le figlie
\end{studente}

Ogni bambino ha poi iniziato a contare i propri pezzi e tutti sono
stati concordi nell'affermare che la stella era fatta da 12
triangoli. Un bambino ha avuto difficoltà a contare i pezzi perché non
riusciva a individuare il punto da cui era partito per calcolarli.

\`E stato poi chiesto loro di prendere un numero uguale di pezzi e
realizzare un serpente facendo coincidere i lati. La maggior parte dei
bambini ha realizzato il serpente senza difficoltà; solo una bambina
non ha fatto coincidere i lati mettendo in fila i
triangoli. Costatando la differenza fra il suo serpente e quello degli
altri ha provato a rifarlo facendo combaciare i lati.

Alla fine di quest'attività è stato chiesto ai bambini di riportare su
di un foglio le figure. I bambini hanno realizzato prima la stella e
poi il serpente, prendendo lo stesso numero di triangoli di carta
colorata già precedentemente ritagliati e incollandoli sul foglio. \`E
stato poi chiesto loro quale fosse l'oggetto più lungo e tutti hanno
indicato il serpente, mentre l'oggetto più grande è stato individuato
nella stella. Un bambino ha detto che la stella conteneva più pezzi
perché era più grande mentre un altro ha individuato la mancanza di un
pezzo nel serpente di una compagna.

\subsubsection{Consigli per i colleghi che vogliono proporre le stesse
  attività}

Quest'attività è legata non solo alle forme, ma anche al rapporto fra
spazio e quantità. \`E perciò necessario lavorare con bambini di
cinque anni. La capacità di riconoscere quantità è legata anche allo
sviluppo psicologico del bambino, questo però non deve scoraggiare
proposte sfidanti.


\section[Sperimentazione \#9: infanzia]{Sperimentazione \#9: scuola
  dell'infanzia, marzo~2010}

Analogamente alla precedente, questa sperimentazione, pur utilizzando
il kit ``Giochiamo con le forme'', introduce in realtà anche attività
proposte dal kit ``Torri, serpenti e\dots{} geometria''.

\subsection{Osservazioni generali}

La sperimentazione si è svolta di una sezione della stessa scuola
dell'infanzia
in cui si è svolta la sperimentazione precedente.

\subsubsection{Presentazione della classe}
Un gruppo omogeneo di 10 bambini di 5 anni

\subsubsection{Composizione dei gruppi}
I bambini sono stati divisi dall'insegnante in due piccoli gruppi;
seguendo i seguenti criteri:
\begin{itemize}
\item Abbinamento di bambini in base a una precedente buona relazione
\item presenza di leader positivi
\item presenza in ogni singolo gruppo di bambini con differenti
  capacità relazionali, attenzione e abilità cognitive.
\end{itemize}
Si è puntato a avere in ogni gruppo condizioni simili in termini di
dinamicità e di apprendimento.

\subsubsection{Insegnanti presenti}

Solo l'insegnante di classe

\subsubsection{Calendarizzazione degli incontri}
\begin{calendario}
  \begin{itemize}
  \item 18 marzo
  \item 22 marzo
  \item 24 marzo
  \item 29 marzo
  \end{itemize}
\end{calendario}

\subsection{Primo incontro}

\begin{description}
\item[Alunni presenti:] 9 presenti, 1 assente
\item[Tempo effettivo di lavoro:] 1 ora e mezza, dalle 10,15 alle 11,45
\end{description}
\begin{consegna}
  Presentazione dei triangoli equilateri del kit e esplorazione
  libera delle forme, con la possibilità di costruire liberamente
  delle figure. Riproduzione delle figure proposte in ``Torri,
  serpenti e geometria'' Schede per la classe prima, pag. 4 e
  5. (Lumaca, vaso, aquilone, papera)

  \materiali{} %
  I triangoli equilateri del kit ``Giocare con le forme'' e triangoli
  equilateri in cartoncino costruiti dalle insegnanti. Triangoli in
  carta colorata per le rielaborazioni successive.
\end{consegna}

\subsubsection{Osservazioni}
Il materiale è stato accolto con entusiasmo. Inizialmente ciascun
bambino ha cercato di appropriarsi dei triangoli e costruire
singolarmente delle figure tridimensionali: la tenda degli indiani, la
montagna, la piramide.

Quando sono state presentate le immagini contenute nel kit i bambini
hanno trovato modalità collaborative, accoppiandosi spontaneamente,
correggendosi tra di loro, o dando indicazioni verbali o spostando i
pezzi, in modo che le figure risultassero corrispondenti al modello;
inoltre in questo contesto bambini normalmente poco coinvolti hanno
avuto un ruolo attivo all'interno del piccolo gruppo e della coppia.

Per alcuni bambini è risultato molto facile riprodurre le figure,
contare e sommare i triangoli e successivamente sperimentare modalità
più semplici per ricomporre la stessa figura (la stella).

Inizialmente i bambini cominciavano a costruire la figura guardando il
modello e posizionando prima le punte, per poi costruire il corpo
centrale. Poi alcuni bambini hanno individuato che il corpo centrale
di alcune figure (stella, chiocciola, aquilone, vaso) era un esagono e
lo hanno verbalizzato; sono quindi partiti da questa scoperta per
comporre più facilmente le stesse figure, senza che ci fosse bisogno
del modello. Un bambino ha colto il numero di triangoli necessari per
comporre la stella, affermando che ne occorrevano 12 perché
all'interno c'era un esagono.

Successivamente i bambini hanno realizzato figure con triangoli di
carta ritagliati e incollati sul foglio.

\subsubsection{Consigli per i colleghi che vogliono proporre le stesse
  attività}
Valorizzare il materiale, indicandone l'importanza (la provenienza, in
questo caso) e responsabilizzando i bambini circa la tenuta e la
conservazione. Lavorare in piccoli gruppi e se possibile in uno spazio
raccolto e dedicato a esperienze mirate. Costruire rituali, per
favorire la concentrazione e la collaborazione.

\subsection{Secondo incontro}

\begin{description}
\item[Alunni presenti:] 8 presenti, 2 assenti
\item[Tempo effettivo di lavoro:] 1 ora e mezza, dalle 14,15 alle 15,40
\end{description}
\begin{consegna}
  Riproduzione delle figure proposte in: ``Torri, serpenti e
  geometria'' Schede per la classe prima, pag. 4 e 5 (Lumaca, vaso,
  aquilone, papera), con i triangoli equilateri. Riproposizione della
  stessa attività con i triangoli isosceli e rettangoli.

  Analisi delle differenze fra i triangoli equilateri e quelli
  proposti in fase successiva. Scoperta delle differenze fra i diversi
  triangoli e dei possibili modi in cui si possono usare.

  \materiali{} %
  I triangoli equilateri, isosceli e rettangoli del kit ``Giocare con
  le forme'' e triangoli equilateri in cartoncino costruiti dalle
  insegnanti.
\end{consegna}

\subsubsection{Osservazioni}
I bambini sono stati suddivisi in due piccoli gruppi e sono stati
consegnati loro i triangoli isosceli e equilateri. I bambini su
sollecitazione dell'insegnante hanno osservato i pezzi; alla domanda
se anche i nuovi pezzi fossero dei triangoli hanno risposto
affermativamente:
\begin{studente}[ ]
  Soltanto che questi sono più grandi di quelli dell'altra volta
\end{studente}
Alla domanda
\begin{tutor}[ ]
  In che cosa sono diversi?
\end{tutor}
una bambina ha risposto
\begin{studente}[ ]
  I due lati sono uguali più alti
\end{studente}
Successivamente sono stati distribuiti i triangoli rettangoli. I
bambini hanno sovrapposto e confrontato i due tipi di triangoli,
trovando la differenza: nel triangolo rettangolo un lato è più lungo
degli altri; alcuni bambini lo hanno indicato con il dito. Partendo
dalla scoperta fatta da alcuni bambini si è discusso collettivamente,
condividendo le conclusioni.

\subsubsection{Consigli per i colleghi che vogliono proporre le stesse
  attività}
Può essere importante per facilitare il lavoro di bambini piccoli fare
prima delle esperienze con le forme di tipo percettivo utilizzando
altri materiali, ad esempio bastoni per psicomotricità, spiedini senza
punta, bastoncini, ecc.

\subsection{Terzo incontro}

\begin{description}
\item[Alunni presenti:] 8 presenti, 2 assenti
\item[Tempo effettivo di lavoro:] 1 ora e mezza, dalle ore 14 alle
  15,30
\end{description}
\begin{consegna}
  Viene proposto ai bambini di realizzare le figure proposte negli
  incontri precedenti con i tre tipi di triangoli: equilatero,
  isoscele e rettangolo.

  \materiali{} %
  Triangoli equilateri, rettangoli e isosceli del kit.
\end{consegna}

\subsubsection{Osservazioni}
I primi a essere proposti sono i triangoli equilateri e l'insegnante
invita i bambini a comporre una figura a piacere: vengono composti
degli esagoni. L'insegnante invita i bambini a contare i triangoli di
cui è composto l'esagono; alcuni bambini rispondono che l'esagono è
formato da sei triangoli.

Vengono poi offerti sei triangoli rettangoli e viene chiesto ai
bambini di comporre la stessa figura. I bambini ci provano e qualcuno
dice: \bambini{Non viene}; \bambini{No, è un rombo}, risponde un
altro. E un altro ancora: \bambini{È un esagono strano}. L'insegnante
chiede come mai non sia venuta l'immagine ``giusta''. I bambini fanno
notare la differenza tra questi tipi di triangoli, cioè
\bambini{questi triangoli hanno un lato più lungo}.

Viene poi proposto di provare a costruire l'esagono con sei triangoli
isosceli. I bambini provano, ma dicono che non bastano
\begin{studente}[ ]
  Ne mancano ancora, Ce ne vogliono altri tre
\end{studente}
L'insegnante chiede il motivo per cui ce ne vogliono di più e i
bambini rispondono:
\begin{studente}[ ]
  Questi non sono come gli altri (riferendosi agli equilateri) hanno
  le punte molto alte
\end{studente}
Alla fine affermano che la prima figura (esagono) è giusta, mentre la
seconda con i triangoli rettangoli e l'ultima con i triangoli isosceli
sono diverse. I bambini del secondo gruppo hanno costruito la lumaca
seguendo lo stesso procedimento e hanno notato che con otto triangoli
rettangoli si riesce a costruire un quadrato.

\subsubsection{Consigli per i colleghi che vogliono proporre le stesse attività}
Lavorare in piccoli gruppi con le modalità indicate precedentemente,
fornendo materiale a sufficienza e altro materiale in carta, in modo
che i bambini possano lavorare e rielaborare le esperienze.

\subsection{Quarto incontro}

\begin{description}
\item[Alunni presenti:] 8 presenti, 2 assenti
\item[Tempo effettivo di lavoro:] 1 ora e mezza, dalle ore 10,15 alle
  11,45
\end{description}
\begin{consegna}
  In quest'incontro l'insegnante ha sentito la necessità di
  approfondire le attività svolte nell'incontro precedente, in modo da
  consolidare gli obiettivi raggiunti.

  \attivita{Prova di pavimentazione con i tre tipi di triangoli}:
  Riproduzione di pavimentazione su fogli con triangoli di carta
  colorata.  Composizione libera di immagini con triangoli equilateri
  del kit.  Riproduzione libera di immagini su fogli con triangoli
  equilateri in carta colorata.

  \materiali{} %
  Triangoli equilateri, isosceli e rettangoli del kit.  Triangoli in
  carta colorata per le rielaborazioni.
\end{consegna}

\subsubsection{Osservazioni}
\begin{figure}[htp]
  \centering
  \includegraphics[width=0.90\textwidth]{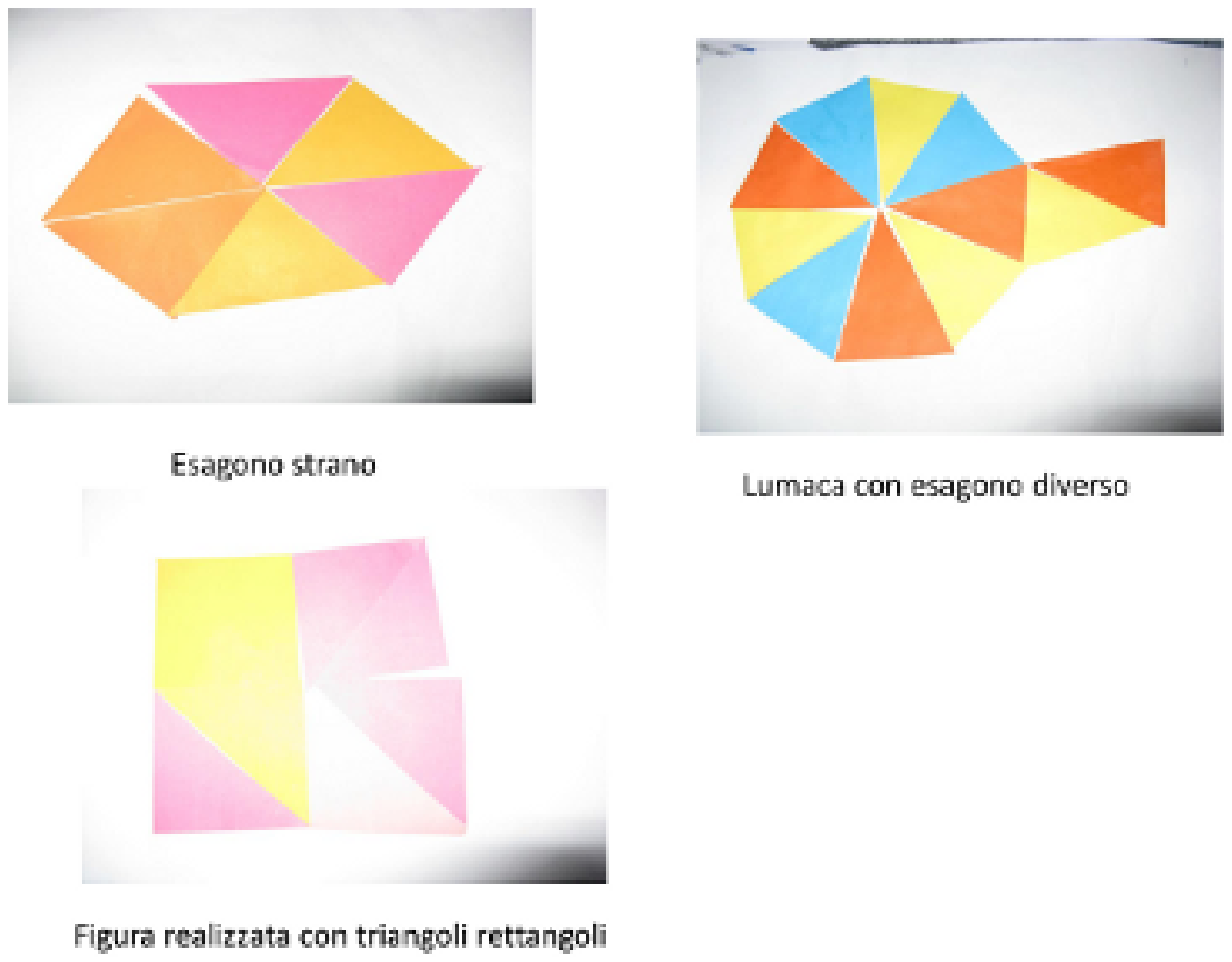}
  \label{pic:formemorbide:9}
\end{figure}
I bambini, suddivisi in due piccoli gruppi come nelle precedenti
attività, accolgono l'attività con sempre maggiore interesse. Dopo
avere distribuito i tre tipi di triangolo a ciascun gruppo,
l'insegnante invita i bambini a osservarli e pone la domanda:
\begin{tutor}[ ]
  Sono tutti e tre dei triangoli?
\end{tutor}
I bambini rispondono affermativamente. L'insegnante chiede di chiarire
la motivazione e gli alunni rispondono che
\begin{studente}[ ]
  Hanno tutti tre punte
\end{studente}

Successivamente li mettono vicini e fanno una classificazione dei
triangoli per altezza, denominandoli come: \bambini{quello alto,
  quello medio e quello basso}.  Poi l'insegnante delimita uno spazio
dei tavolini e suggerisce ai bambini di far finta che quello spazio
sia il pavimento di una stanza; chiede quindi di provare a
piastrellarla. I bambini iniziano prima a utilizzare i triangoli
equilateri, li sistemano sui tavolini e cercano di riempire lo
spazio. Alla domanda
\begin{tutor}[ ]
  Con questi triangoli si può piastrellare il pavimento?
\end{tutor}
i bambini rispondono
affermativamente, ma notano che:
\begin{studente}[ ]
  Ne manca un pochino
\end{studente}

In seguito provano con i triangoli rettangoli. Inizialmente li mettono
un po' a caso, ma una bambina si rende conto che la piastrellatura non
va bene e dice: \bambini{Così non riesce}; prova unendo le due basi
dei triangoli e ottiene un quadrato. Gli altri bambini disfano il
pezzo piastrellato male e continuano su quello proposto dalla
compagna, affermando: \bambini{È venuto bene, si può}.

Per ultimo usano i triangoli isosceli; provano a fare la pavimentazione
in diversi modi e verbalizzano: \bambini{Con questi non si
  riesce!}. Concludono dicendo:
\begin{studente}[ ]
  con i medi, per fare il pavimento ne manca un pochettino, con quelli
  bassi si può e con quelli alti non si riesce
\end{studente}
Con triangoli di carta, precedentemente preparati, i bambini
riprendono l'attività svolta realizzando su dei fogli A3 le
pavimentazioni prima eseguite.

\subsubsection{Consigli per i colleghi che vogliono proporre le stesse
  attività}
Lavorare in piccoli gruppi, in spazi raccolti, con abbondante
materiale, in modo che tutti bambini possano disporre di tessere a
sufficienza per lavorare seguendo i loro bisogni.

\chapter[Diamo forma alla geometria: Regolari o no]{Diamo forma alla geometria: Regolari~o~no?}



\section[Sperimentazione \#1: prima secondaria di primo
grado]{Sperimentazione \#1: classe prima secondaria di primo grado,
  novembre~2009/febbraio~2010}

\subsection{Osservazioni generali}

\subsubsection{Presentazione della classe}
La sperimentazione si è svolta in due classi prime di scuola
secondaria di primo grado.

La prima delle due classi (classe ``A'') è composta da 24 alunni, di
cui una portatrice di handicap e un dislessico certificato; i ragazzi
hanno già provato modalità di lavoro in gruppo.

La seconda delle due classi (classe ``B'') e composta da 24 alunni e
ha una situazione identica con handicap molto grave (autonoma ma
completamente analfabeta, età scolare 1 elementare). Oltre a un
dislessico certificato è presente un alunno molto intelligente ma con
disturbi comportamentali.

\subsubsection{Composizione dei gruppi}
Gruppi creati liberamente dai ragazzi, si formano 5 gruppi di 5 o 4
componenti in entrambe le classi.

Nella classe ``A'' i gruppi sono rimasti invariati; nella classe ``B''
nell'ultima attività (\attivita{scheda D}) ci sono state variazioni
spontanee nella formazione dei gruppi.

Le classi sono abituate a lavorare in gruppo (laboratorio di scienze,
giochi matematica senza frontiere junior) e in genere si lascia che i
gruppi si formino autonomamente. Di solito [gli insegnanti] non
intervengo, tranne se ci sono tensioni o emarginazioni che non si sono
presentate; le variazioni sono scelte libere dei ragazzi.

\subsubsection{Insegnanti presenti}
Durante le attività in ogni classe è presente solo l'insegnante di
classe.

\subsubsection{Calendarizzazione degli incontri}
\begin{calendario}
  \begin{itemize}
  \item 16 novembre (per entrambe le classi)
  \item 17 novembre (``A'') e 19 novembre (``B'')
  \item 23 novembre (``A'') e 26 novembre (``B'')
  \item 30 novembre (``A'') e 3 dicembre (``B'')
  \item 11 gennaio (per entrambe le classi)
  \item 12 gennaio (``A'') e 14 gennaio (``B'')
  \item 25 gennaio (``A'') e 28 gennaio (``B'')
  \item 12 febbraio (``A'') e 15 febbraio (``B'')
  \item 22 febbraio (solo classe ``B'')
  \item 25 febbraio (solo classe ``B'')
  \end{itemize}
\end{calendario}
Per la classe ``A'': si è, quasi sempre, lavorato nella prima e
seconda ora del lunedì mattina; per la classe ``B'' nella prima e
seconda ora del giovedì mattina. In entrambe le classi quando non
necessitavano tempi così lunghi si è utilizzata la quarta ora del
lunedì per la classe ``A'' o la sesta ora del giovedì per la classe
``B''.

\subsection{Primo incontro}

\begin{description}\item[Alunni presenti:] Classe ``A'': 24 alunni,
  classe ``B'': 23 alunni
\item[Tempo effettivo di lavoro:] Tempo previsto 1 ora.
\end{description}
\begin{consegna}
  Nella prima sperimentazione abbiamo lasciato liberi i ragazzi di
  manipolare il materiale contenuto nel kit.
  Abbiamo dato pochissime indicazioni: osservate le forme; guardate
  come poterle unire; guardate i colori, fateci quello che volete,
  senza romperle!

  \materiali{}%
  Abbiamo consegnato a ogni gruppo un certo numero (10 circa) di
  poligoni uguali; lasciati liberi di costruire dopo poco (10 minuti
  circa) i ragazzi si sono scambiati le forme. Potevano anche venire
  un delegato per gruppo dal docente a richiedere ulteriori
  forme. Così abbiamo stimolato gli alunni a usare i nomi corretti
  per i poligoni.
\end{consegna}

\subsubsection{Osservazioni}
Rispetto al percorso proposto dal kit, abbiamo scelto di aggiungere al
percorso del kit un'ora iniziale di conoscenza del
materiale
perché due anni fa, nel lavorare in un progetto di robotica, abbiamo
visto che è necessario dar modo ai ragazzi di sfogare curiosità e
entusiasmo verso il materiale nuovo prima di indirizzarli a un fine
didatticamente più preciso.

Abbiamo lasciato liberi i ragazzi di manipolare il materiale del kit,
per seguire il lavoro dei ragazzi abbiamo compilato una scheda,
ricalcando il protocollo osservativo proposto da questo corso.

\begin{center}
\makebox[0pt]{%
  \begin{tabular}{|*{9}{l|}}
    \hline
    \multicolumn{3}{|l|}{All'inizio del lavoro} &
    \multicolumn{3}{l|}{Dopo mezz'ora} &
    \multicolumn{3}{l|}{Al termine} \\ \hline
    \mbox{}
    &
    \multicolumn{2}{l|}{n. alunni} &
    \mbox{}&
    \multicolumn{2}{l|}{n. alunni} &
    \mbox{}&
    \multicolumn{2}{l|}{n. alunni} \\ \hline
    stupore & 5 & 2 &
    interesse al lavoro & 21 & 19 &
    interesse al lavoro & 15 & 21 \\ \hline
    curiosità&11&2&curiosità&&&curiosità&4& \\ \hline
    desiderio di giocare&2&3&desiderio di giocare&1&3&desiderio di
    giocare&3&1 \\ \hline
    entusiasmo&6&15&entusiasmo&2&&entusiasmo&2& \\ \hline
    indifferenza&&1&indifferenza&&1&indifferenza&&1  \\ \hline
  \end{tabular}%
  }
\end{center}
E in una tabella abbiamo evidenziato alcune questioni toccate dal
lavoro dei ragazzi:

\begin{center}
  \begin{tabular}{|l|p{4cm}|p{4cm}|p{4cm}|}
    \hline
    \mbox{} & \mbox{} &
    \multicolumn{1}{c|}{SÌ} &
    \multicolumn{1}{c|}{NO} \\ \hline
    COLORE & si chiedono se ha significato? &
    Solo due ma decidono che il colore non ha importanza &
    Tutti gli altri non si pongono il problema
    \\ \hline
    &
    si chiedono perché queste forme e non altre? &
    In 5 notano la mancanza di rombi e rettangoli, ma presi dal gioco
    non si chiedono il perché.  &
    \\ \hhline{|~---|}
    &
    ne conoscono il nome?
    &
    Sì per le forme piane. Le forme solide per la classe ``B'' non
    rappresentano solidi geometrici (castello, gatto, astronave,
    topino\dots{}).
    &
    classe ``A'' non danno nomi
    \\ \hhline{|~---|}
    &
    provano a mescolare forme diverse?
    &
    Dopo pochi minuti. (alcuni dicono agli altri che è
    vietato!!). Quando non trovano le tessere giuste lasciano buchi
    nella costruzione. &
    \\ \hhline{|~---|}
    FORME &
    costruiscono spontaneamente pavimentazioni?
    &
    Sì, nella classe ``A'' in 12 iniziano dal piano, nella classe
    ``B'' solo 3 sul piano&
    \\ \hhline{|~---|}
    &
    costruiscono spontaneamente solidi?
    &
    nella classe ``A'' 12 iniziano dai solidi poi seguiti dai
    compagni. nella classe ``B'' partono subito dai solidi in 21&
    \\ \hhline{|~---|}
    &
    le costruzioni sono casuali o cercano regolarità nella costruzione?
    &
    classe ``B'' quasi sempre oggetti di fantasia, ma con un progetto
    in mente. classe ``A'' metà casuali e metà cercava simmetrie
    &
    \\ \hline
  \end{tabular}
\end{center}
Cosa dicono gli alunni: \bambini{Forte}, \bambini{mitico},
\bambini{creativo e collaborativo}, \bambini{serve a esprimere le
  idee e a lavorare insieme}, \bambini{fa venire in mente cose},
\bambini{scatena le idee}, \bambini{qualcosa di diverso}, \dots{}

Alcuni hanno notato che il solido è \bambini{flessibile}, cioè
\bambini{si muove} (ha uno o più gradi di libertà%
) e che il movimento verrebbe bloccato dall'inserimento di una
tessera.

Cosa dicono i professori: molto entusiasmo e collaborazione; difficile
far osservare gli aspetti geometrici perché il desiderio di giocare è
troppo forte. Sarebbe meglio dopo 40 minuti ritirare il materiale e
farli riflettere sui solidi costruiti.

Un alunno si comporta in maniera indisponente: gioca a mettersi i
triangoli come orecchini; corone di esagoni in testa ecc.

Le classi hanno fatto alcune fotografie [che non sono disponibili].

\subsubsection{Consigli per i colleghi che vogliono proporre le stesse
  attività}I consigli sono le solite attenzioni didattiche: a volte,
se il clima era esageratamente effervescente, abbiamo trovato utile
ritirare il materiale, lasciare sulla cattedra i solidi o le
tassellazioni costruite, e ragionare insieme appuntando alla lavagna
le osservazioni fatte.

\subsection{Secondo incontro}

\begin{description}\item[Alunni presenti:]Tutti gli alunni
\item[Tempo effettivo di lavoro:]1 ora e mezza.
\end{description}

\begin{consegna}
  Abbiamo consegnato una scheda per ciascun alunno, più precisamente
  la scheda A \attivita{Per cominciare} del kit nella parte che
  riguarda ``vertici/spigoli/facce''. La consegna non è stata
  commentata a voce.

  \materiali{}%
  Sono stati consegnati esattamente i ``pezzi'' previsti dalle
  istruzioni del kit. Dopo poco i ragazzi hanno cominciato a
  scambiarsi e mescolare il materiale; poi ci hanno richiesto altre
  forme in aggiunta anche per il desiderio di costruire tanti solidi
  diversi.
\end{consegna}

\subsubsection{Osservazioni}

Dall'osservazione della scheda i gruppi hanno compreso la differenza
tra faccia, spigolo, vertice. Hanno poi cominciato a costruire i due
solidi e poi hanno completato la scheda come da richiesta.
Intanto che lavoravano abbiamo colto queste frasi:
\begin{studente}[ ]
  \begin{itemize}
\item \bambini{Non abbiamo fatto niente!!} \makebox{(quasi al pianto)}.
  \item \bambini{Abbiamo provato in tutti i modi ma non ci siamo
      riusciti} (richiesta di aiuto alla prof).
  \item \bambini{Noi abbiamo fatto il più difficile}; \bambini{noi
      abbiamo finito per primi};
  \item \bambini{Olè palla fatta siamo soddisfatti!}
\end{itemize}
\end{studente}
Alla fine di ogni scheda ci siamo sempre ritrovate per discutere e
confrontare i risultati raggiunti e le nostre impressioni, le
osservazioni che seguono riguardano in generale entrambe le classi.

In entrambe le classi, tutti gli alunni sono molto interessati,
compreso il ragazzo (classe ``B'') che nell'incontro precedente era
indisponente. Dopo qualche battibecco hanno deciso di dividersi in due
sottogruppi per costruire separatamente i due solidi.

C'è stata molta differenza nel tempo impiegato dai vari gruppi: i più
veloci in 30 minuti hanno fatto anche le risposte della scheda, altri
hanno utilizzato tutto il tempo e la discussione finale è stata
ripresa successivamente. Alcuni gruppi si sono divisi il lavoro a metà
dando per scontato una simmetria, ciascuno ha costruito metà del
poliedro ``palla'' e poi non sono riusciti a unire i due emisferi. Un
numero significativo di ragazzi (9 su 24) hanno richiesto più volte
che venisse mostrata la foto a colori trovando difficoltà a osservare
la foto in bianco e nero: in particolare era molto in crisi un ragazzo
con discalculia grave e certificata. Pensavamo che vedere le immagini
in cui i poliedri hanno le facce colorate e lavorare con pezzi di
polydron che hanno solo gli spigoli (cioè quelli con le facce vuote)
potesse creare problemi, invece hanno notato la differenza ma non ha
creato problemi.

Nella costruzione della palla un gruppo ha costruito un poliedro di
poco diverso da quello richiesto, simmetrico e solo dal confronto con
quello dei compagni ha colto la differenza, che infatti risulta poco
visibile ma molto utile per le successive osservazioni sulla
simmetria. Per contare le facce qualcuno le ha contate realmente,
qualcuno ha notato la simmetria, qualcuno sapendo il numero iniziale
delle tessere ha lavorato sulle differenze. Per conteggiare quanto
spigoli in totale aveva ogni poliedro, qualcuno ha moltiplicati i
vertici per il numero degli spigoli uscenti e poi si è accorto che
erano troppi e ha diviso per due.

Essendo necessario stare nei tempo previsti, abbiamo trovato utile,
come per ogni altra lezione, riprendere il filo della lezione
precedente prima di continuare il cammino.

\subsection{Terzo incontro}

\begin{description}
\item[Alunni presenti:] nella classe ``B'' era assente l'alunna con
  handicap
\item[Tempo effettivo di lavoro:] L'attività inizialmente prevista per
  1 ora, si è protratta per 2 perché la discussione e la costruzione
  della tabella è stata laboriosa.
\end{description}

\begin{consegna}
  Visto l'interesse dei ragazzi e considerati che volevamo
  approfondire l'argomento prima di passare a \attivita{Mosca cieca},
  abbiamo inserito nel percorso una attività non prevista dal kit in
  cui abbiamo chiesto agli allievi di costruire i 10 poliedri della
  scheda plastificata.

  \materiali{}%
  Non è stato distribuito materiale, i vari capogruppi venivano alla
  cattedra a richiedere le forme desiderate per la costruzione.
\end{consegna}

\subsubsection{Osservazioni}
L'esperienza è stata più difficoltosa del previsto, per una sola
docente è complicato seguire i 5 gruppi che tendono a lavorare
autonomamente e con tempi e modalità diverse.

\begin{itemize}
\item c'è stata viva partecipazione da parte di tutti
\item per alcuni alunni ci sono state difficoltà nel montaggio e
  incastro dei pezzi
\item i tempi sono stati molto diversi anche a seconda della
  complessità dei poliedri
\item molti gruppi si sono divisi il lavoro tenendo conto della
  simmetria del poliedro (i due emisferi) ma poi non tutti hanno
  potuto montare insieme le due parti per un problema di
  ``ganci''
\item i due solidi che nella \attivita{scheda A} sono colorati con
  tonalità diverse ma usando lo stesso colore non sono risultati
  ``chiari'' a vari alunni che hanno trovato difficoltà nella
  ``lettura'' del solido. in particolare un alunno dislessico ha
  dichiarato chiaramente che senza il colore era \bambini{impossibile
    vedere il solido}.
\end{itemize}
Frasi degli alunni:
\begin{studente}[ ]
  \begin{itemize}
  \item \bambini{quelli tutti rossi non si vedono bene}
  \item \bambini{perché le nostre due metà non si incastrano?}
  \item \bambini{ma come è possibile agli altri si incastrano\dots{}?}
  \item \bambini{non sta agganciato\dots{}}
  \item \bambini{uffa appena aggiungo un pezzo si rompe dall'altro
      lato\dots{}}
  \end{itemize}
\end{studente}
Comunque alla fine tutti i gruppi sono riusciti a costruire i 10 poliedri.

A questo punto abbiamo chiesto loro di contare di ognuno di essi il
numero di facce, spigoli, vertici esattamente come avevano fatto nelle
lezione precedente con i due poligoni campione.

E qui per i poliedri più semplici il lavoro è stato veloce ma per i
più complessi le difficoltà sono state insormontabili perché il
conteggio delle facce risultava ancora possibile ma per vertici e
spigoli gli alunni si confondevano.

Senza dare un nome ai poliedri, abbiamo preparato una tabella dei
poliedri della scheda plastificata; li abbiamo semplicemente numerati
e poi i ragazzi hanno contato trovando un accordo sul numero finale da
indicare in tabella:

\begin{center}
  \begin{tabular}{|l|c|c|c|}
    \hline
    & n. facce & n. vertici & n. spigoli \\ \hline
    poliedro 1 & 32 & 30 & 60 \\ \hline
poliedro 2 &
14 &
12 &
24 \\ \hline
poliedro 3 &
14 &
24 &
36 \\ \hline
poliedro 4 &
8 &
6 &
12 \\ \hline
poliedro 5  &
12 &
10 &
20  \\ \hline
poliedro 6 &
\dots{}ecc &
\dots{}
     &
\dots{}
   \\ \hline
poliedro 7  &&& \\ \hline
poliedro 8 &&&  \\ \hline
poliedro 9 &&&  \\ \hline
poliedro 10 &&&  \\ \hline
  \end{tabular}
\end{center}
(Per fare questo lavoro è stato necessario numerare i poliedri della
scheda plastificata: potrebbe essere utile fornire la scheda già con i
numeri per differenziare i poliedri.)

Dall'esame di questa tabella, che è stata trascritta in grande alla
lavagna, un alunno ha chiesto se c'era qualche legame tra i tre
numeri, io ho detto ``forse\dots{}'', ``provate a vedere\dots{}''

Così tre ragazzi (i più vispi) hanno trovato la relazione F+V = S + 2 da soli.

È stato fatto presente che una relazione deve valere sempre! Hanno
risposto che sicuramente valeva per otto dei poliedri presi in esame e
che per mezzo di essa avrebbero potuto completare le ultime caselle
dove non riuscivano a fare i conteggi e di questo erano molto
soddisfatti. Nell'altra classe la ricerca della relazione è stata
stimolata dall'insegnante su una base di competizione con l'altra
classe che aveva già trovato la soluzione.

\begin{description}
\item[Classe ``A''] solo un numero ristretto di alunni (5) su 24 ha
  portato avanti il discorso, gli altri ascoltavano;
\item[Classe ``B''] la situazione è stata leggermente diversa con 5, 6
  alunni fortemente interessati e altri che facevano fatica a capire,
  ma non volevano essere tagliati fuori e provavano il calcolo in
  altri solidi della tabella.
\end{description}
In queste situazioni non è facile tenere tutto sotto controllo e si
rischia sempre di non sottolineare a sufficienza l'aspetto geometrico
e/o matematico dell'esperienza.

\subsubsection{Consigli per i colleghi che vogliono proporre le stesse
  attività}
Ci è sembrato utile, dopo le prime esperienze, ritirare il materiale
distribuito nei gruppi, appoggiarlo alla cattedra e far riflettere gli
alunni sull'esperienza fatta senza il materiale in mano per non
distrarsi manipolandolo. Altra osservazione: per facilitare la
produzione di tante osservazioni da parte degli alunni abbiamo
sorvolato molto sulla precisione di linguaggio, che è stato corretto
solo in un secondo tempo.

\subsection{Quarto incontro}

\begin{description}
\item[Alunni presenti:]2 assenti nella classe ``A'', tutti presenti
  nella classe ``B''
\item[Tempo effettivo di lavoro:]2 ore.
\end{description}

\begin{consegna}
  Consegnata la seconda parte della scheda A, ovvero l'attività della
  \attivita{Mosca cieca}.

\materiali{}%
Come previsto dal kit inizialmente non sono stati distribuiti
materiali, i vari capogruppi venivano alla cattedra a richiedere le
forme desiderate per la costruzione.
\end{consegna}

\subsubsection{Osservazioni}

Molto buono l'entusiasmo iniziale di tutti i gruppi, dopo una ventina
di minuti alcuni avevano finito con successo, altri erano ancora al
lavoro con entusiasmo ma altri ancora stavano perdendo interesse
perché con le istruzioni date non riuscivano a costruire il
poliedro%
.
\begin{itemize}
\item Nella classe ``A'' al termine 6 gruppi su 10 avevano lavorato
  con successo utilizzando solo le istruzioni ricevute, 2 erano
  riusciti a terminare solo dopo che l'insegnante ha consegnato loro
  la scheda plastificata e ha ripetuto che il poliedro era uno dei 10
  raffigurati che avevano già costruito la precedente lezione, 2 non
  sono riusciti a costruire il poliedro.
\item Nella classe ``B'' non si è mai mostrata la scheda plastificata,
  ma si è insistito completando o correggendo le informazioni che
  passavano da un gruppo all'altro: è stato necessario ribadire che lo
  scopo era di fornire in modo corretto le indicazioni e non, come in
  molti hanno fatto, quello di confondere e rallentare il lavoro dei
  compagni, che erano visti come ``avversari'' nel gioco.
\end{itemize}
Osservando dall'esterno possiamo dire che alcuni hanno provato
semplicemente a montare il materiale indicato dai compagni nell'unico
modo che sembrava possibile, altri hanno chiesto ulteriori
informazioni, altri non sapevano da dove cominciare.

Il fatto che nella lezione precedente si sia lavorato a contare numero
di facce, vertici e spigoli ha indotto i ragazzi a pensare che queste
fossero le migliori informazioni da dare ma all'atto pratico si sono
accorti che non era così. È stato sicuramente utile far in modo che
fossero gli alunni a chiedere il numero più o meno esatto delle forme
necessarie per la costruzione, anche se nei solidi complessi, le
richieste erano incerte e ripetute.

Frasi degli alunni:
\begin{studente}[ ]
  \begin{itemize}
  \item \bambini{ma adesso tutte queste facce come le metto\dots{}}
  \item \bambini{ma come stanno vicine?}
  \item \bambini{cerca di ricordarti quello dell'altro giorno con gli
      esagoni deve essere quello\dots{}}
  \item \bambini{Ma loro mi hanno scritto solo una riga}
  \item \bambini{Non voglio dirgli di più perché indovinano subito}
  \end{itemize}
\end{studente}
Le difficoltà ci sono state soprattutto per i poliedri
più complessi dove le facce si potevano accostare in modi diversi. Al
termine del lavoro c'è stata da parte degli alunni una breve
discussione su come era meglio dare le istruzioni.

Un alunno ha detto:
\begin{studente}[ ]
  \begin{itemize}
  \item \bambini{Il numero di facce e la forma mi serve solo per
      scegliere nello scatolone i pezzi giusti ma poi non so da dove
      partire\dots{}}
  \item \bambini{io vorrei capire quando prendo in mano la prima
      faccia cosa devo mettere vicino\dots{}}
  \item \bambini{quali forme devo mettere insieme, quelle uguali o
      no?}
  \item \bambini{no, l'altra volta univamo quelle diverse\dots{}}
  \item \bambini{i quadrati sono di più allora incominciamo da quelli
      perché ne abbiamo tanti\dots{}}
  \item \bambini{ti ripeto che io vorrei sapere quali forme mettere
      vicine all'inizio\dots{}}
  \item \bambini{è difficile iniziare\dots{}}
  \end{itemize}
\end{studente}
A questo punto è chiaro che sapere quali figure accostare e in quale
ordine a partire da un vertice è importantissimo e, in particolare,
hanno capito quasi subito che è fondamentale l'informazione di quanti
spigoli escono da ogni vertice. Il gioco \attivita{Mosca cieca} è
stato ripetuto due volte e nella classe ``B'' un gruppo (due alunne
molto deboli in matematica) non sono mai riuscite a costruire il
solido, anche se poi è stata buona e attenta la partecipazione alla
discussione collettiva. Dopo le due partite a \attivita{Mosca cieca} è
stata proposta la scrittura compatta usata nel libretto (4,6,6)
spiegandola. Nella classe ``A'' quattro alunni intuitivamente la
comprendono, altri 10 vengono convinti dai compagni, gli altri 10
restano perplessi. Un alunno chiede ai compagni di scrivere la
``scrittura compatta'' del solido che hanno in mano e dice che sarà
capace subito di ricostruirlo, appare molto sicuro di quello che dice
e i compagni lo mettono alla prova. Alla fine qualcuno mette dei
numeri a caso nella scrittura (3,5,6,4) e chiede allo stesso compagno
di costruire il solido. Il ragazzo prende i pezzi e appena cerca di
accostarli si accorge che non è possibile. Ci fermiamo alla domanda
``ma quando si può?''%
. Il tempo è terminato e l'attenzione di molti è calata, c'è
un gruppo trainante di 5 o 6 alunni che si pone domande ma il resto
della classe è stanco. Le due ore sono state impegnative, un
intervallo di attenzione continuativa, anche se si tratta di un
laboratorio, è un tempo lungo per alunni di prima media. Per la classe
``B'' il giorno successivo, senza ridistribuire tutto il materiale ma
avendo ancora a disposizioni solidi costruiti la mattina precedente,
si sono ripetute le osservazioni fatte, in particolare si è posta
attenzione alla comodità della scrittura sintetica.

\subsection{Quinto incontro}
\begin{description}
\item[Alunni presenti:]tutti gli alunni presenti
\item[Tempo effettivo di lavoro:]Tempo previsto 1 ora.
\end{description}

\begin{consegna}
  Per entrambe le classi la consegna è stata di osservare le diverse
  forme dei poligoni, di arrivare alla definizione di poligono
  regolare. Dopo che è stata data la definizione di [poligono]
  regolare (molti la ricordavano dalle scuole elementari, e altri la
  ritrovavano nelle forme), abbiamo chiesto se fosse possibile
  determinare l'ampiezza degli angoli interni di ciascun poligono.

  \materiali{}%
  Utilizzo delle forme polydron del kit
\end{consegna}

\subsubsection{Osservazioni}
Abbiamo trovato molto utili le forme a disposizione nel kit: abbiamo
distribuito le forme dei poligoni regolari e abbiamo recuperato
dall'osservazione le caratteristiche (angoli e lati uguali); ci siamo
soffermati al calcolo dell'ampiezza degli angoli interni. Alcuni
alunni hanno proposto il goniometro, altri hanno accostato i poligoni
per vedere con quanti si poteva formare l'angolo giro, nel caso fosse
possibile potevano eseguire la divisione e quindi calcolare l'ampiezza
degli angoli (questa osservazione ci verrà utile poi per la
pavimentazione). Dove non era possibile hanno cercato un'altra strada:
alcuni hanno usato le diagonali, suddiviso il poligono in triangoli
per trovare la somma totale e quindi poterla dividere per il numero
degli angoli. Si sono accorti che potevano costruire altri poligoni
regolari oltre a quelli del kit, potevano disegnare qualsiasi poligono
regolare con infiniti lati e quindi hanno pensato alla circonferenza
da loro chiamata cerchio.

\subsubsection{Consigli per i colleghi che vogliono proporre le stesse
  attività}
Questa attività è stata molto utile sia per i poliedri regolari, di
cui conoscevano ormai bene le facce, sia per la pavimentazione che è
venuta spontanea (infatti abbiamo poi anticipato la \attivita{Scheda
  D} alla \attivita{Scheda C}).

\subsection{Sesto incontro}
\begin{description}
\item[Alunni presenti:]24 nella classe ``B'', 24 nella classe ``A''
\item[Tempo effettivo di lavoro:]1 ora e 45 minuti.
\end{description}

\begin{consegna}
  Viene consegnata ai gruppi la \attivita{scheda B} del kit

  \materiali{}%
  Vengono distribuiti ai gruppi i materiali previsti dal kit, ma non
  le fotocopie degli sviluppi, perché la costruzione avrebbe portato
  via troppo tempo e la sola osservazione dello sviluppo piano è
  difficile per gli alunni di prima.
\end{consegna}

\subsubsection{Osservazioni}
A questo punto siamo passati ai poliedri regolari e abbiamo proceduto
secondo la scheda B cartacea, facendo costruire tutti i poliedri della
scheda B, escluso quello azzurro formato da rombi. Tutti avrebbero
voluto anche il cubo sulla scheda e se lo sono costruito lo stesso (la
ragazzina portatrice di handicap ha costruito cubetti per tutti).

In entrambe le classi i gruppi hanno suddiviso i poliedri in modo
diverso e abbiamo fatto un prospetto alla lavagna: sulla base di
queste osservazioni abbiamo verificato che definizioni diverse
portavano a suddivisioni diverse. Le definizioni date erano molto
simili ma sempre incomplete, tutti hanno pensato a facce di poligoni
regolari, in pochi hanno richiesto che a ogni vertice concorresse lo
stesso numero di spigoli, alle facce nessuno ha pensato. Abbiamo
discusso i risultati della scheda e non hanno avuto difficoltà a
compilare la tabella e poi hanno compilato con facilità le risposte
successive sempre osservando la tabella. Sono stati più veloci e
sicuri di quanto avessimo previsto. Anche il numero di possibili
poliedri regolari è stato compreso senza problemi praticamente da
tutti.

La partecipazione è stata buona, i ragazzi erano motivati anche alla
compilazione delle schede; sempre un po' difficoltoso per noi
sistematizzare gli aspetti che emergono nelle discussioni.

\subsubsection{Consigli per i colleghi che vogliono proporre le stesse
  attività}
Ci pare utile l'idea di un prospetto alla lavagna per raccogliere i
dati di ciascun gruppo e poterli poi confrontare. In questo modo tutti
sono coinvolti nella discussione e si ragiona meglio su quei poliedri
che sono stati classificati differentemente dai gruppi.

\subsection{Settimo incontro}
\begin{description}
\item[Alunni presenti:]nella classe ``A'' 2 alunni erano assenti,
  nella classe ``B'' era assente l'alunno dislessico
\item[Tempo effettivo di lavoro:]1 ora e 45 minuti.
\end{description}
\begin{consegna}
  Viene distribuita la \attivita{scheda D} del kit (una per ogni
  alunno)

  \materiali{}%
  Nelle due classi si è proceduto in maniera differente, perché noi
  docenti non ci eravamo accordate (non sono però emerse differenze
  significative nell'attività delle due classi). Nella classe ``A''
  ogni gruppo aveva a disposizione un decina di poligoni per ogni
  forma. Nella classe ``B'' non sono stati distribuiti materiali, ma
  ogni gruppo chiedeva alla docente le forme che riteneva necessarie.
\end{consegna}

\subsubsection{Osservazioni}
\begin{figure}[htp]
  \centering
  \setlength{\fboxrule}{2pt}%
  \fcolorbox{blulogo}{pdcolor3}{
    \begin{minipage}{0.95\linewidth}
      \small{}%
      \vspace*{5pt}%
      \textbf{I POLIEDRI}\label{scheda:studente:poliedri}

      Ogni poliedro ha un'altezza, una lunghezza e una profondità.

      I poliedri sono figure solide, hanno facce uguali e regolari, da
      ogni vertice spuntano lo stesso numero di facce e ne esistono
      solo cinque:
      \begin{multicols}{2}
        \begin{itemize}
        \item 3 formati da i triangoli equilateri(;
        \item 1 formato con i quadrati(cubo);
        \item 1 formato da i pentagoni.
      \end{itemize}
      \end{multicols}
      I poliedri hanno tre punti fondamentali per la costruzione del
      solido:
      \begin{multicols}{3}
        \begin{itemize}
        \item faccia;
        \item spigoli;
        \item vertici.
        \end{itemize}
      \end{multicols}

      \medskip{}%
      \textbf{I POLIGONI}

      I poligoni sono figure piane formate da linee spezzate chiuse,
      hanno gli angoli e i lati uguali,.

      Possono essere infiniti ma i più comuni sono:
      \begin{multicols}{4}
        \begin{itemize}
        \item triangolo(3)
        \item quadrato(4)
        \item pentagono(5)
        \item esagono(6)
        \item ettagono(7)
        \item ottagono(8)
        \item ennagono(9)
        \item decagono(10)
        \end{itemize}
      \end{multicols}

      Ci sono poligoni sia con gli angoli convessi, cioè se si
      prolungano i lati e il prolungamento è fuori alla figura, e gli
      angoli concavi, cioè se si prolungano i lati il prolungamento è
      dentro la figura.

      \medskip{}%
      \textbf{REGOLARI E NON REGOLARI}

      Mi servono minimo 3 lati per disegnare un poligono; mentre non
      esiste un numero massimo di lati perché sono infiniti. Più il
      numero di lati è infinito più si avvicina al cerchio.

      Se parliamo di poligoni convessi li possiamo distinguere in:
      poligoni regolari e non regolari.

      Si dicono regolari i poligoni che hanno lati e angoli uguali.

      \medskip{}%
      \textbf{COME SI MISURANO GLI ANGOLI?}

      In genere si usa il goniometro oppure si possono unire dei
      poligoni. Così:
      \begin{itemize}
      \item 360°:4=90°misura di ogni angolo
      \end{itemize}
      Oppure:
      \begin{itemize}
      \item 180°x3=540° (gradi totali del pentagono)
      \item 540°:5=108° (gradi dell'angolo del pentagono)
      \end{itemize}
      Ci sono delle formule per rendere dei calcoli molto più semplici
      come:
      \begin{itemize}
      \item (n° vertici x n° di spigoli di ogni vertice): 2=n° di
        spigoli in totale
      \item (vertici + facce)-2=spigoli (creatore Eulero)
      \end{itemize}

      \medskip{}%
      \textbf{TASSELLAZIONI REGOLARI}

      Le tassellazioni regolari sono quelle ottenute unendo fra loro
      poligoni regolari dello stesso tipo in modo da restare sul
      piano. Ci sono delle tassellazioni possibili come:
      \begin{multicols}{2}
              \begin{itemize}
      \item 1 con i triangoli equilateri (60°x6= 360°);
      \item 1 con i quadrati (90°x 4=360°);
      \item 1 con gli esagoni (120°x 3=360°).
      \end{itemize}
      \end{multicols}
      E altre non possibili come:
      \begin{multicols}{3}
        \begin{itemize}
        \item pentagoni;
        \item ottagoni;
        \item decagoni.
        \end{itemize}
      \end{multicols}
      \mbox{}
    \end{minipage}
}
\end{figure}

Siamo passate alla \attivita{scheda D} perché il discorso si collegava
bene ai poligoni regolari appena fatti.

Le due alunne che hanno costruito da sole tassellazioni di triangoli e
quadrati secondo le indicazioni del kit e sono state molto orgogliose
mostrando a ciascuna classe il risultato. I ragazzi, forse anche
grazie al lavoro precedente sui poligoni regolari, non hanno trovato
difficoltà neanche nel calcolare l'ampiezza degli angoli interni e
capire quando la loro somma poteva dare 360°. Nella classe ``B'' un
alunno ha posto la domanda se fosse possibile utilizzare due angoli di
180° per ottenere l'angolo giro; i compagni sapevano di no ma solo nel
tentativo di disegnare la tassellazione alla lavagna hanno capito che
era necessaria un semicerchio%
, mancavano i lati per
costruire il poligono. Nella classe ``B'' si sono stupiti nel vedere
quanti vincoli ci siano in una tassellazione: pensavano si potesse
ottenere accostando qualunque poligono. Gli alunni di entrambi le
classi si sono molto divertiti con le tassellazioni uniformi
realizzandone il più possibile. Possiamo sottolineare che entrambe le
classi hanno fatto molta fatica nell'osservare la tassellazione a
destra di pag. 37
 e non abbiamo insistito
più di tanto.

Per prepararsi alla verifica programmata un'alunna della classe ``B''
ha deciso spontaneamente di preparare sul quaderno uno schema
riassuntivo di quanto fatto in classe, schema che, volutamente non
abbiamo corretto, ma che riteniamo utile allegare (si veda
pag.~\pageref{scheda:studente:poliedri}).

\subsubsection{Consigli per i colleghi che vogliono proporre le stesse
  attività}
La suddivisione temporale fra tassellazioni regolari e quelle con
poligoni di diverse forme non è venuta spontanea ma data come
indicazione dalla \attivita{scheda D}; nella successiva discussione
gli alunni hanno trovato più belle le tassellazioni più fantasiose e
con forme diverse.

\subsection{Ottavo incontro}

\begin{description}
\item[Alunni presenti:]classe ``B'' 3 assenti; classe ``A'' un assente
\item[Tempo effettivo di lavoro:] 45 minuti
\end{description}

\begin{figure}[bht]
  \centering
  \begin{dadocente}
    \begin{minipage}{0.98\linewidth}
      \small{}%
      \begin{enumerate}
      \item Prendi uno dei solidi che hai a disposizione. Scrivi la
        lettera che lo contraddistingue e rispondi alle seguenti domande:
        \begin{enumerate}
        \item Scrivi il numeri dei vertici, delle facce e degli spigoli.
        \item Che forma hanno le facce?
        \item Quanti spigoli concorrono in ogni vertice?
        \item Sai scrivere con la notazione abbreviata i numeri che ti
          servono per la costruzione del solido?
        \item Scrivi la relazione di Eulero e controlla che sia
          verificata in questo solido.
        \end{enumerate}
      \item Osserva il solido disegnato:
        \begin{enumerate}
        \item Scrivi il numero dei vertici, delle facce e degli
          spigoli.
        \item Usando le lettere scritte in figura scrivi i nomi
          \raisebox{-3cm}[0pt][0pt]{\makebox[0pt][l]{%
              \hspace*{2cm}\includegraphics{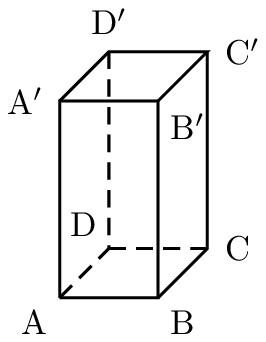}
}}%
          \\ degli
          spigoli, dei vertici e delle facce
        \item Che forma hanno le facce?
        \item Quanti spigoli concorrono in ogni vertice?
        \item Scrivi la relazione di Eulero e controlla che sia
          verificata \\in questo solido.
        \end{enumerate}
      \item Ricorda il lavoro fatto sui poligoni regolari:
        \begin{enumerate}
        \item Scrivi la definizione di poligono regolare.
        \item Scrivi il nome dei primi sei.
        \item Quanti sono i poligoni regolari?
        \item Calcola l'ampiezza dell'angolo interno dell'esagono regolare.
        \item Calcola l'ampiezza dell'angolo interno del pentagono regolare.
        \item Con quali poligoni regolari puoi fare una tassellazione
          regolare?
       \end{enumerate}
      \item Ricorda il lavoro fatto sui poliedri regolari:
        \begin{enumerate}
        \item Scrivi la definizione che abbiamo dato come definitiva.
        \item Quanti sono i poliedri regolari?
        \item Quali poligoni regolari puoi usare come facce?
        \item Ci sono poligono che ti permettono di formare più
          poliedri diversi? Quali?
         \end{enumerate}
      \item Cosa ne pensi dell'esperienza laboratoriale sui poliedri?
        Cosa ti è piaciuto? L'hai trovata interessante? Difficile?
        Quale pensi sia stato il tuo contributo al tuo gruppo? Ti è
        piaciuto lavorare con del materiale per verificare
        concretamente le tue ipotesi?
      \end{enumerate}
    \end{minipage}
  \end{dadocente}
    \caption{Testo della verifica}\label{testo:verifica}
\end{figure}

\begin{consegna}
  Viene consegnata la scheda di verifica (si veda
  Figura~\ref{testo:verifica} pag.~\pageref{testo:verifica}) e viene
  chiesto di eseguirla utilizzando il materiale a disposizione sulla
  cattedra

  \materiali{}%
  Scheda individuale consegnata a ciascun alunno

  Sulla cattedra sono stati messi a disposizione degli alunni 10
  poliedri, alcuni regolari e altri no, etichettati con una
  lettera.
\end{consegna}

\subsubsection{Osservazioni}
Il lavoro è stato individuale, in silenzio, come un normale test. Per
la classe ``B'' la verifica è stata preceduta, nella settimana
precedente, da una mezz'ora di ripasso collettivo, riprendendo le
schede e focalizzando le varie fasi del lavoro. Un'alunna ha preparato
spontaneamente un breve riassunto (cfr. incontro precedente).

Nella classe classe ``A'' non è stato fatto un lavoro analogo, ma è
stato chiesto un ripasso a casa, prima della verifica, delle schede
che gli alunni hanno sul quaderno.

\paragraph{Risultati verifica classe ``A''}
In generale risultati come nelle aspettative: medi.

Difficoltà nelle risposte:
\begin{description}
\item[2b] perché gli alunni di prima non sanno utilizzare nessuna
  notazione per scrivere il nome di una faccia del poliedro, ho
  sbagliato io a dare per scontato che utilizzassero le lettere dei
  vertici!
\item[2 e] solo pochi hanno risposto
\item[3] hanno trovato semplici tutte le domande di questa serie ma
  hanno fatto una gran confusione tra le domande 3 f con 4 d,
  moltissimi hanno messo gli esagoni tra le figure che permettono di
  formare poliedri regolari e i pentagoni tra le figure che permettono
  la tassellazione. Stupita di ciò ho chiesto loro nella lezione
  successiva il perché e molti mi hanno riferito di essersi
  semplicemente confusi con i nomi e le forme (non avevano possibilità
  durante la verifica di maneggiare le forme se non i poliedri già
  costruiti e presenti sulla cattedra). Con in mano le forme capivano
  benissimo che gli esagoni accostati formavano un ``pavimento'' e i
  pentagoni davano ``concavità'' e formavano un solido\dots{} Uno di
  loro mi ha detto \bambini{io ho fatto giusto perché ho pensato alle
    cellette della api nell'alveare!}
\end{description}

Alcune frasi degli alunni
\begin{studente}[ ]
  \begin{itemize}
  \item \bambini{Mi son piaciute molto le parti dove abbiamo costruito
      poliedri; ho dato un buon contributo al gruppo; mi è piaciuto
      molto costruire con il materiale le figure per poi costruire e
      controllare le ipotesi.}
  \item \bambini{Istruttiva ma laboriosa; molto bello il gioco
    }\attivita{Mosca cieca}\bambini{. È molto divertente imparare
      provando. Sbagliando s'impara ma anche provando s'impara.}
  \item \bambini{Secondo me nel mio gruppo ho dato molte idee e ho
      aiutato i miei compagni in difficoltà. Mi è piaciuto verificare
      le mie ipotesi anche se non ne avevo molte. Avere del materiale
      concreto non mi è servito molto perché per me avere materia o no
      è uguale.}
  \item \bambini{Esperienza bellissima, che è servita a integrarmi
      ancora di più nella classe e che mi ha insegnato molto. Mi è
      piaciuto scoprire questi solidi che prima non conoscevo. Ho dato
      un buon contributo al gruppo. Lavorare con il materiale mi ha
      aiutato a capire di più che non con i disegni.}
  \item \bambini{Mi è piaciuto }\attivita{Mosca cieca}\bambini{.}
  \item \bambini{Interessante completare le schede, difficile nel
      contare i vertici, spigoli e facce.}
  \item \bambini{Io ho contribuito molto perché il lavoro di gruppo mi
      piace molto. Ho imparato la geometria quasi giocando.}
  \item \bambini{Nel gruppo abbiamo collaborato tutti insieme. È stato
      molto divertente ma anche un po' preoccupante perché avevo paura
      che i pezzi della Bicocca si rompino.}
\item \bambini{All'inizio un po' difficile perché non avevo mai
      fatto questa esperienza. Penso di aver aiutato e ho fatto il
      possibile. Mi è piaciuto più di tutto verificare con i miei
      occhi. La cosa più bella è stata }\attivita{Mosca
      cieca}\bambini{ perché era basato sul ragionamento.}
  \item \bambini{Avere degli oggetti solidi davanti me l'ha resa più
      facile.}
  \item \bambini{Un po' difficile a volte, ma non sempre.}
  \item \bambini{Costruttiva ma un po' faticosa.}
  \item \bambini{Mi è piaciuto molto perché mi piace costruire e mi
      piacerebbe rifare questa esperienza.}
  \item \bambini{Mi è piaciuto costruire i poliedri. Interessante ma
      un po' difficile. Mi ha insegnato tante cose come quando abbiamo
      fatto }\attivita{Mosca cieca}\bambini{.}
  \item \bambini{Questa esperienza è stata un po' impegnativa ma
      divertente. Mi è piaciuto stare nel gruppo e scoprire le diverse
      forme che potevano saltare fuori. Interessante ma difficile.}
  \item \bambini{Mi è piaciuto perché costruire con il materiale e
      capire le cose logiche mi piace.}
  \item \bambini{Molto utile per capire; da soli potevamo arrivare a
      scoprire le regole senza leggerle e sentirle dalla prof. La
      prima lezione è stata un'improvvisata. Lavorando manualmente
      sono riuscita a capire meglio.}
  \end{itemize}
\end{studente}

\paragraph{Risultati verifica classe ``B''}
Il lavoro degli alunni della classe classe ``B'' non è andato molto
bene. La valutazione ha escluso le due domande sulla formula di Eulero
perché nessuno aveva dato risposte. Errori comuni: poligono al posto
di poliedro, vertice al posto di spigolo, spigolo al posto di lato,
molte difficoltà nell'uso delle lettere della figura per indicare le
facce e/o i vertici e gli spigoli. Le risposte date rispetto
all'andamento e al gradimento dell'attività sono state molto scarne e
spesso semplici monosillabi (SÌ - NO). Eccone una sintesi:
\begin{studente}[ ]
  \begin{itemize}
  \item \bambini{Interessante doversi organizzare in gruppo, e anche
      il fatto di sperimentare con le mani} (l'alunno con disturbi di
    comportamento);
  \item \bambini{ho lavorato tanto cercando di arrivare alle
      conclusioni};
  \item \bambini{mi è piaciuto moltissimo costruire, non è difficile};
  \item \bambini{non mi piace questa materia. Imparo di più senza
      gioco};
  \item \bambini{mi è piaciuta la tassellazione e ascoltare gli
      altri};
  \item \bambini{molto interessante ma a volte complessa};
  \item \bambini{ho imparato giocando e ho conosciuto meglio i
      compagni};
  \item \bambini{mi è piaciuto inventare forme attraverso la logica};
  \item \bambini{piaciuto stare in gruppo, ma difficile};
  \item \bambini{divertente e coinvolgente};
  \item \bambini{esperienza bellissima ma difficile non capisco
      geometria};
  \item \bambini{il mio compito calcolare i numeri e dare
      idee. Bello};
  \item \bambini{imparare giocando, difficile un poco noiosa};
  \item \bambini{in questo modo la geometria è più semplice da
      capire};
  \item \bambini{bella};
  \item \bambini{non servivo più di tanto}.
  \item \bambini{Spero di rifarlo};
  \item \bambini{Non è stato difficile, il gruppo ha lavorato con
      cura} (l'alunna che ha fatto il riassunto finale)
  \end{itemize}
\end{studente}

\subsection{Nono incontro}

Avevamo programmato l'ultima attività per lunedì 22 febbraio nelle
classi classe ``B'' e classe ``A'', ma la professoressa della classe
``A''
è stata ricoverata in ospedale e ne avrà per parecchio tempo. Abbiamo
quindi deciso che per la classe ``A'' il lavoro si può considerare
concluso con l'effettuazione della verifica, mentre la descrizione che
segue è da riferirsi alla sola classe ``B''.

\begin{description}
\item[Alunni presenti:]Solo la classe classe ``B'', 22 alunni, come al
  solito gli alunni si sono divisi in 5 gruppi, mantenendo, per loro
  libera scelta, i gruppi consueti. Due alunne non hanno partecipato
  all'attività perché mi hanno chiesto di poter svolgere la verifica
  che non avevano fatto in quanto assenti
\item[Tempo effettivo di lavoro:] dalle 10.45 alle 11.35
\end{description}

\begin{consegna}

  A ogni alunno è stata consegnata la fotocopia della \attivita{Scheda C}
  \attivita{Oltre i regolari}.

  \materiali{}%
  I sacchetti contenenti le forme geometriche erano liberamente
  disponibili alla cattedra.
\end{consegna}

\subsubsection{Osservazioni}

Nel corso dell'attività sono state effettuate le riprese video. La
\attivita{Scheda C} riportava i dettagli dell'attività da svolgere, ma
ho informato i ragazzi che io avrei fatto il giro di tutti i gruppi
per documentare, attraverso una breve ripresa video, la loro
attività. Ho detto loro che erano liberi di dire qualsiasi cosa
inerente all'attività nel suo complesso: in particolare potevano
esprimere un giudizio su quanto fatto, una precisazione di tipo
geometrico su quanto appreso, una spiegazione dell'attività in
corso\dots{} o altro. Già nella prima giornata di lavoro avevamo
realizzato delle fotografie poi andate perse e tutti loro erano
entusiasti all'idea del filmato. Sottolineo che l'alunno dislessico i
cui genitori non avevano autorizzato riprese o foto e che io pensavo
di oscurare, mi ha chiesto se poteva ``vedersi'' e ha quindi chiesto
in casa la firma per la liberatoria. Avendo a disposizione solo una
fotocamera con microfono incorporato ho dovuto imporre alla classe
molto più silenzio rispetto al solito e l'attività è riuscita
perfettamente. L'analisi della scheda è stato fatto nell'incontro
successivo.

\subsection{Decimo incontro}

\begin{description}
\item[Alunni presenti:] 23 alunni su 24
\item[Tempo effettivo di lavoro:]
  dalle 8.00 alle 8.30
\end{description}

\subsubsection{Osservazioni}

Come previsto l'analisi della scheda è stata fatta due giorni dopo: i
ragazzi non sono più stati divisi in gruppi ma, leggendo e discutendo
insieme si sono riprese le nozioni inerenti ai poliedri non
regolari. I ragazzi hanno detto che l'attività sostanzialmente
riprendeva argomenti già fatti, ma hanno trovato molto interessante
utilizzare la simbologia proposta. Ho ripreso in mano i solidi che la
classe aveva costruito e ho interrogato a caso gli alunni: la quasi
totalità non ha più difficoltà a riconoscere forma, numero delle
facce, numero dei vertici. Una alunna che spesso mostra difficoltà
logiche ancora non riesce a vedere quanti spigoli concorrono in un
vertice. Un alunno ha notato che la simbologia ``non vale'' quando il
poliedro ha un numero di spigoli diverso da vertice a vertice. Mi ha
fatto notare questo con vivo disappunto dicendo \bambini{ma allora non
  serve a niente, cosa ci sforziamo a fare?}. Tutta la classe ha preso
coscienza delle difficoltà incontrate nell'usare termini specifici:
confondono poligono e poliedro, hanno chiaro quando si tratta di un
poliedro regolare ma faticano a definirne le caratteristiche,
confondono vertice e spigolo.

L'alunna con grave handicap non era presente durante il primo incontro
e è entrata a scuola quando la discussione della scheda era già
terminata; avendo molti poliedri montati alla cattedra le ho chiesto
se voleva smontarli e mettere le diverse forme nei vari sacchetti. Ha
svolto il lavoro con molta attenzione, praticamente da sola e senza
errori.

\chapter[Diamo forma alla geometria: Grande o piccolo]{Diamo forma alla geometria: Grande~o~piccolo?}


\section[Sperimentazione \#1: terza primaria]{Sperimentazione \#1:
  classe terza primaria, marzo~2010}

\subsection{Osservazioni generali}

\subsection{Presentazione della classe}

La classe è composta da 24 alunni: 14 ragazze e 10 ragazzi. Sono
presenti due ragazzi stranieri e due bambine che necessitano
dell'insegnante di sostegno.

\subsubsection{Composizione dei gruppi}

5 gruppi eterogenei scelti dall'insegnante, composti da 5/6
bambini. Sono abituati a lavorare in gruppo e a svolgere attività
laboratoriali.

\subsubsection{Insegnanti presenti}

A tutti gli incontri sono presenti le due insegnanti di classe, una
educatrice e la collega della scuola secondaria di primo grado che
partecipa al corso.

\subsubsection{Calendarizzazione degli incontri}

9 incontri di 40 minuti circa.
\begin{calendario}
  \begin{itemize}
  \item 24 febbraio
  \item 25 febbraio
  \item 3 marzo
  \item 4 marzo
  \item 10 marzo
  \item 11 marzo
  \item 24 marzo
  \item 25 marzo
  \item 16 aprile
  \end{itemize}
\end{calendario}

\subsection{Primo incontro}

\begin{description}
\item[Alunni presenti:] 24
\item[Tempo effettivo di lavoro:] 40 minuti
\end{description}

Sono state proposte le attività del quaderno di laboratorio di ``Torri,
serpenti e \dots{} geometria'' della classe terza, modificando alcune
consegne in base al materiale a disposizione. Ad esempio partendo
sempre dal disegno su foglio che veniva ricoperto dalle tessere
triangolo.

\begin{consegna}
  L'insegnante presenta ai bambini il materiale di lavoro, consegnando
  a ogni gruppo un sacchetto. I bambini osservano, toccano e
  iniziano a raggruppare le figure secondo il colore o la forma,
  successivamente qualcuno inizia a creare delle figure astratte.

  L'insegnante richiama l'attenzione dei bambini consegnando dei
  disegni (5 disegni per ogni gruppo) chiedendo loro di ricoprirle con
  le forme di cui dispongono.

  \materiali{} %
  ??
\end{consegna}

\subsubsection{Osservazioni}

I bambini hanno notato subito come uno stesso disegno poteva essere
ricostruito utilizzando forme diverse. L'intento dell'attività è stato
quello di creare un primo approccio con il materiale di lavoro e un
avvio a quello che saranno le richieste della scheda di lavoro.

\subsection{Secondo incontro}
\begin{description}
\item[Alunni presenti:]24
\item[Tempo effettivo di lavoro:] 40 minuti
\end{description}

\begin{consegna}
  A ogni gruppo sono stati mostrati la figure STELLA (pag.1 del
  quaderno di laboratorio) e è stato chiesto di ricostruirla con i
  triangoli a disposizione, concludendo che necessitano 12 tessere
  triangolo.

  Si consegnano le figure LUMACA e DIAMANTE per capire chi delle due
  ha il contorno più lungo e sovrapponendole un bambino ha risposto
  che sono uguali.

  \materiali{}%
  Quelli previsti dal kit.
\end{consegna}

\subsubsection{Osservazioni}

Continuando a sovrapporre le forme ai disegni presentati, alcuni
bambini hanno provato a sovrapporre a una forma rettangolo altre che
riuscissero a ricostruirla, come per esempio due quadrati, oppure un
quadrato e due triangoli rettangoli. La lezione a questo punto ha
preso una piega inaspettata, poiché l'insegnante ha continuato il
gioco e i bambini hanno affrontato a loro insaputa il concetto di
EQUIVALENZA.

\subsection{Terzo incontro}

Continuazione attività con tessere

\begin{description}
\item[Alunni presenti:]24
\item[Tempo effettivo di lavoro:] 40 minuti
\end{description}

\begin{consegna}
  L'insegnante mostra le figure SERPENTE e STELLA (p.3) chiedendo di
  realizzarle con le tessere a loro disposizione, concludendo che ne
  necessitano 12 per la stella e 9 per il serpente. Successivamente
  consegna a ogni gruppo due fili di spago, di due lunghezze diverse,
  e chiede ai bambini di farlo aderire alle figure e capire a chi
  serve il pezzo più lungo. Tutti concludono che per la stella
  necessita lo spago più lungo.
\end{consegna}

\subsubsection{Osservazioni}

Si è riscontrata una evidente difficoltà nel far aderire lo spago alle
figure poiché questo si spostava in continuazione, quindi è stato
necessario l'intervento dell'insegnante. Alcuni bambini avevano
risposto ancora prima di provare con lo spago, che la stella avrebbe
avuto bisogno del pezzo più lungo (senza argomentare l'affermazione,
per alcuni vedendo la stella con più rette e cambi di direzione e
invece il serpente una retta continua, dava l'idea che per la prima
servisse dello spago in più rispetto al serpente, ma l'affermazione
era data più ``a occhio'').

\subsection{Quarto incontro}

\begin{description}
\item[Alunni presenti:] 24
\item[Tempo effettivo di lavoro:] 40 minuti
\end{description}

\begin{consegna}
  L'insegnante presenta le figure CLESSIDRA e LUMACA (pag.~4) e seguendo
  la scheda compila con i bambini la tabella, concludendo che il
  numero di tessere utilizzate per ricostruire le figure non è lo
  stesso, mentre il numero di lati di tessera contenuti nel contorno
  sì. Lavorando sempre sul contorno e la superficie l'insegnante ha
  presentato altri disegni come il pesce, il granchio, il serpente, il
  diamante. A seguito di tutto si è giunti ai concetti di PERIMETRO e
  AREA.
\end{consegna}

\subsubsection{Osservazioni}

Per molti bambini i concetti esposti non sono stati di immediata
comprensione, quindi l'insegnante ha disegnato alla lavagna un
rettangolo andando a calcolare con tutti quello che era il suo
perimetro e la sua area attraverso il semplice conteggio dei
quadretti.

\subsection{Quinto incontro}
\begin{description}
\item[Alunni presenti:] 24
\item[Tempo effettivo di lavoro:] 40 minuti
\end{description}

\begin{consegna}
  L'insegnante consegna a ogni gruppo la scheda 5 di pag.~6.
\end{consegna}

\subsubsection{Osservazioni}

I bambini hanno mostrato serie difficoltà nell'eseguire l'esercizio
proposto (nel trovare una figura che avesse un contorno più lungo di
un'altra inizialmente, poi secondo le indicazioni della scheda,
disegnare una figura di lunghezza 10 e una più corta. La difficoltà
maggiore è stato trovare una figura nuova senza immaginare che
disegnandone una già vista durante gli incontri poteva essere una
buona partenza).

\subsection{Sesto incontro}

\begin{description}
\item[Alunni presenti:]24
\item[Tempo effettivo di lavoro:] 40 minuti
\end{description}
\begin{consegna}
  L'insegnante presenta ai bambini il nuovo materiale di lavoro,
  consegnando a ogni gruppo un quadrato e 4 triangoli equilateri. I
  bambini osservano, toccano e successivamente viene chiesto di
  incastrarle tra di loro e saper dire quale solido hanno creato. Poi
  l'insegnante consegna solo 3 triangoli e i bambini ricostruiscono
  una piramide triangolare. Mettendo a confronto i solidi ottenuti,
  notano subito che nella prima è possibile individuare una ``base di
  appoggio'' (il quadrato), mentre nella seconda
  \bambini{\dots{}comunque la giro è sempre quella\dots{}} (riportando
  l'espressione di un bambino).
\end{consegna}

\subsubsection{Osservazioni}
I bambini non mostrano alcuna difficoltà nella gestione del materiale
e ottenendo in poco tempo due solidi.

\subsection{Settimo incontro}

\begin{description}
\item[Alunni presenti:] 24
\item[Tempo effettivo di lavoro:] 40 minuti
\end{description}

\begin{consegna}
  L'insegnante consegna 6 quadrati a ogni gruppo e chiede loro di
  costruire un solido. Nel frattempo l'insegnante costruisce un cubo
  doppio per metterli a confronto. Ripassa con la classe i concetti di
  FACCE, SPIGOLI, VERTICI e nel conteggio di entrambe le figure i
  bambini notano una sostanziale differenza di grandezza (visibile non
  solo osservandole).

  \materiali{}%
  Tessere Polydron
\end{consegna}

\subsection{Ottavo incontro}

\begin{description}
\item[Alunni presenti:]24
\item[Tempo effettivo di lavoro:] 40 minuti
\end{description}

\begin{consegna}
  L'insegnante consegna 12 pentagoni e i bambini ottengono un
  dodecaedro. Successivamente si chiede loro di contare facce, spigoli
  e vertici e si nota subito come due lati di due pentagoni creano un
  solo spigolo, mentre tre angoli di tre pentagoni unendosi creano un
  solo vertice.
\end{consegna}

\subsubsection{Osservazioni}
Si ottiene così non solo la spiegazione di POLIEDRO ma si arriva a
un'operazione aritmetica (di non facile comprensione da parte di tutta
la classe) per cui è possibile conoscere la quantità di spigoli e
vertici senza contarli uno a uno. 5 (numero lati di un pentagono)
moltiplicato per 12 (numero di pentagoni usati) si ottiene 60 (numero
totale degli spigoli), sapendo però che a due facce corrisponde un
solo spigolo (creato dall'unione di due lati di due pentagoni) 60
viene diviso per 2 e si ottiene 30 (numero effettivo degli spigoli del
dodecaedro)%
. Lo stesso calcolo viene proposto per il conteggio dei vertici (in
questo caso 60 viene diviso per 3).

Inoltre è stato possibile raggruppare in due insiemi separati, figure
come la piramide triangolare e il dodecaedro, il cubo e il
parallelepipedo da una parte, figure come la piramide a base quadrata
dall'altra, poiché, seguendo il ragionamento dei bambini,
\bambini{alcune di queste non cambiano mai a seconda di come le si
  osserva}.

\subsection{Nono incontro}

\begin{description}
\item[Alunni presenti:]24
\item[Tempo effettivo di lavoro:] 40 minuti
\end{description}

\begin{consegna}
  L'insegnante divide in due sottogruppi ogni gruppo di lavoro con i
  separatori e consegna a caso un cubo, una piramide a base quadrata,
  una piramide triangolare, un dodecaedro, un parallelepipedo. Una
  metà del gruppo possiede i solidi e i poliedri, l'altra metà deve
  porre delle domande in modo da capire quale figura possiedono i
  compagni.

  Un'altra proposta: metà gruppo possiede un solido ad esempio il
  cubo, e ordina all'altra metà gruppo di procurarsi 6 quadrati e
  ricreare la figura.
\end{consegna}

\subsection{Conclusioni}

I bambini si sono mostrati sicuramente interessati a ogni attività
proposta. Le forme geometriche di plastica sono risultate il materiale
più efficace, poiché la crepla ha distolto più volte l'attenzione dal
lavoro che i bambini stavano svolgendo, dovendo richiamare più spesso
la loro attenzione.


\section[Sperimentazione \#2: terza secondaria di primo
grado]{Sperimentazione \#2: classe terza secondaria di primo grado,
  marzo~2010}

\subsection{Osservazioni generali }

\subsubsection{Presentazione della classe }
16 alunni 8 ragazze e 8 ragazzi. Sono presenti due ragazzi stranieri
adottati. Un ragazzo con disturbi specifici dell'apprendimento.

\subsubsection{Composizione dei gruppi}
4 gruppi eterogenei scelti dall'insegnante formati da 4 ragazzi (due
maschi e due femmine). Sono abituati a lavorare in gruppo, hanno
partecipato a alcune tappe dei giochi matematici del sito "quaderno a
quadretti".

\subsubsection{Insegnanti presenti}
A tutti gli incontri è presente solo l'insegnante di classe.

\subsubsection{Calendarizzazione degli incontri}
4 incontri di due ore ogni martedì.
\begin{calendario}
  \begin{itemize}
  \item 2 marzo
  \item 9 marzo
  \item 16 marzo
  \item 23 marzo
  \end{itemize}
\end{calendario}

\subsection{Primo incontro }

\begin{datiincontro}
  \begin{description}
  \item[Alunni presenti:] 16
  \item[Tempo effettivo di lavoro:] Scheda A dalle 10,20 alle 11,30
    con all'interno intervallo di 15 minuti. Mosca cieca 11,30-11,55
  \end{description}
\end{datiincontro}
\begin{consegna}
  Sulla cattedra sono posti i sacchetti con le tessere
  \textit{Polydron}.

  Viene consegnata una copia della Scheda A per ogni gruppo. I ragazzi
  sono invitati a leggere autonomamente, interpretare le richieste e
  solo successivamente andare alla cattedra e chiedere all'insegnante
  il materiale che gli necessita per l'attività.

  Ogni gruppo deve rispondere ai quesiti della scheda che verrà poi
  ritirata dall'insegnante per essere valutata. Ogni alunno deve
  appuntare su un foglio del quaderno le fasi dell'attività.

  Durante la compilazione della scheda l'insegnante gira tra i
  gruppi. Quando la maggior parte ha concluso l'attività vengono
  condivisi, discussi e approfonditi i risultati.

  Viene fornita la scheda per il gioco \attivita{Mosca cieca}.
\end{consegna}

\subsubsection{Osservazioni}
I ragazzi sono molto coinvolti, molti giocano con le tessere
costruendo altri solidi oltre a quelli proposti. Il materiale
manipolabile ha sicuramente il pregio di coinvolgerli e mantenerli più
a lunghi attenti all'attività. Necessita però di tempi lunghi per
lasciare ai ragazzi il tempo di prendere dimestichezza col materiale
per poi potersi cimentare nella risoluzione dei quesiti della scheda.

Per un solo insegnante è difficoltoso poter dare la giusta attenzione
a tutti i gruppi. Anche in una classe poco numerosa i quattro gruppi
hanno comunque bisogno di tempo per poter essere guidati alla
riflessione sui risultati ottenuti. La difficoltà da me incontrata è
stata quindi nel passaggio successivo alla compilazione della
scheda. Ritengo infatti necessario che le conclusioni a cui giungono i
ragazzi debbano essere condivise nel gruppo classe e inserite nel
contesto del percorso di studio della geometria già intrapreso. In un
primo tempo ho cercato di guidare i piccoli gruppi ma si è infine
rivelato più utile un momento finale di condivisione insieme. Un'altra
difficoltà è stata decidere quando interrompere l'attività, pur
essendo i gruppi piuttosto equilibrati, uno è stato molto rapido e uno
ha invece necessitato di tempi più lunghi. Ogni gruppo ha comunque
sviluppato interessanti riflessioni soprattutto derivanti dagli errori
fatti. Era quindi importante venissero comunicate ai gruppi che
avevano seguito metodologie diverse.

L'attività di \attivita{Mosca cieca} ha sottolineato come previsto
l'importanza di un linguaggio comune ma è stato interessante osservare
come in un gruppo in particolare, durante la prima fase dell'attività
si erano costruiti un loro gergo particolare (non matematicamente
rigoroso%
)
che li ha resi più veloci degli altri.

\subsection{Secondo incontro }

\begin{datiincontro}
  \begin{description}\item[Alunni presenti:]
    15
  \item[Tempo effettivo di lavoro:] Scheda "Pitagora e la
    similitudine": 10:10 - 11:50 con all'interno l'intervallo di 15
    minuti.
  \end{description}
\end{datiincontro}
\begin{consegna}
  Ogni gruppo riceve una scheda \attivita{Pitagora e la similitudine}
  per le terze medie e un sacchetto con il materiale in gomma crepla.
\end{consegna}

\subsubsection{Osservazioni}
Il materiale si prestava meno al gioco e ha quindi facilitato
l'attività. Gli argomenti richiamati dalla scheda si sono rivelati più
complessi del previsto. È stato importante per i ragazzi confrontarsi
con una attività che richiedesse l'uso di conoscenze precedenti per
poi fare ulteriori approfondimenti.

\subsection{Terzo incontro }

\begin{datiincontro}
  \begin{description}
  \item[Alunni presenti:] 13 (3 assenti giustificati per gare
    sportive). Vengono quindi formati tre gruppi diversi dai
    precedenti.
  \item[Tempo effettivo di lavoro:] Scheda B dalle 10,10 alle 10,40.
    Scheda C dalle 11,05 alle 11,40.
  \end{description}
\end{datiincontro}
\begin{consegna}
  Viene consegnata la \attivita{scheda B} e poi discussa insieme. Dopo
  l'intervallo viene consegnata la \attivita{scheda C}. Il Polydron è
  sulla cattedra a disposizione.
\end{consegna}

\subsubsection{Osservazioni}

La compilazione delle schede B e C ha richiesto meno tempo delle
precedenti. Ciò ha permesso di dedicare più tempo alla condivisione
utilizzando anche il libro di testo per ricercare in esso i concetti
affrontati nelle schede. In questo modo l'attività laboratoriale non è
rimasta una parentesi ma si è inserita pienamente all'interno del
percorso di apprendimento, completandolo e rendendo più efficace
l'acquisizione dei concetti.

\subsection{Quarto incontro }

\begin{datiincontro}
  \begin{description}
  \item[Alunni presenti:] 15 Si ricompongono i gruppi dei primi due
    incontri.
  \item[Tempo effettivo di lavoro:] Scheda D: Dalle 10:10 alle 11:50
    con all'interno un intervallo di 15 minuti.
  \end{description}
\end{datiincontro}
\begin{consegna}
  Viene consegnata a ogni gruppo la \attivita{scheda D} un foglio per
  volta. Il materiale è a disposizione sulla cattedra.
\end{consegna}

\subsubsection{Osservazioni}

La scheda è molto lunga quindi viene consegnata una facciata per volta
per aiutare i ragazzi a concentrarsi e a non disperdere la loro
attenzione. L'insegnante gira nei gruppi per ascoltare le riflessioni
dei ragazzi, chiarire dubbi, valutare le diverse ipotesi di
soluzione. Arrivati alla quarta pagina diviene necessario, viste le
diffuse difficoltà, interrompere il lavoro a gruppi per affrontare
insieme questa parte della scheda:
\begin{dadocente}
  Pensiamo ora al volume dei puzzle. Vi diciamo noi che
  \begin{itemize}
  \item vol(B)=vol(A)
  \item vol(P)=2vol(A)
  \item vol(O)=4vol(A)
  \item vol(A2)=8vol(A)
  \end{itemize}
  Aiutatevi con i puzzle per rendervi conto del fatto che i nostri
  conti qui sopra sono giusti.

  Scrivete qui sotto come convincereste della correttezza di queste
  affermazioni un vostro amico che non ha fatto questo laboratorio.
\end{dadocente}

La mia classe si è trovata in difficoltà nel trovare autonomamente i
passaggi logici che portavano da una affermazione all'altra. È stato
comunque importante per loro scontrarsi con questa difficoltà per poi
comprendere meglio la spiegazione.

In conclusione della attività posso dire che i ragazzi hanno
apprezzato le attività soprattutto per la presenza di materiale
fisicamente utilizzabile e manipolabile. Hanno avuto grandissima
difficoltà nel riportare sul quaderno autonomamente quanto
scoperto. Sarebbe stato più semplice avere una copia ciascuno delle
schede ma oltre all'alto numero di fotocopie necessarie questo li
avrebbe ulteriormente rallentati nell'acquisizione di una competenza
importante come quella del saper prendere appunti. Certo non dover
investire tempo nella formalizzazione insieme sul quaderno avrebbe
lasciato più tempo all'attività vera e propria.

\subsubsection{Considerazioni finali}

Durante il periodo in cui ho avuto a disposizione il kit l'ho inoltre
utilizzato per costruire modelli, farli costruire ai ragazzi,
affrontare altre attività in altre classi ecc. Sarebbe importante
avere sempre a disposizione il materiale a scuola per poter svolgere
le varie attività nel momento più opportuno, per poter utilizzare al
meglio tutte le potenzialità degli strumenti e per rendere i ragazzi
abituati nel tempo al loro utilizzo.

\chapter{Il viaggio segreto: Giochi di aritmetica}


\section[Sperimentazione \#1: quinta primaria]{Sperimentazione \#1:
  classe quinta primaria,
  gennaio/febbraio~2010}

\subsection{Osservazioni generali}

\subsubsection{Presentazione della classe}

Si tratta di una quinta formata da 22 alunni, di cui 11 maschi e 11
femmine. Sono presenti alunni stranieri che sono nati in Italia e che
non presentano difficoltà linguistiche. Ci sono poi tre alunni con
difficoltà di apprendimento: un bambino con ritardo, che ha grosse
difficoltà a livello di ragionamento logico, mentre ha meno difficoltà
a livello tecnico (anche se ha bisogno di un continuo esercizio); una
bambina dislessica e con ritardo, che però dimostra buone capacità nel
ragionamento logico e in generale nelle attività pratiche; una bambina
con problemi psicologici che compromettono l'apprendimento. Tutti e
tre seguono la programmazione della classe ridotta nei contenuti e
lavorano spesso fuori dalla classe con l'insegnante di sostegno,
perché hanno bisogno di spiegazioni individualizzate e di
esercizi/attività mirati.

L'anno scorso è stato affrontato il concetto di frazione e quello di
frazioni equivalenti, ma non quello di somme di frazioni. Abbiamo
scelto di proporre l'attività senza alcun ripasso, come spunto per
riprendere poi l'argomento.

In classe viene spesso utilizzato il lavoro di gruppo per svolgere
ricerche o rielaborazioni delle lezioni di scienze e di
geografia. Meno spesso per proporre problemi logico-matematici (alcuni
presi dal libro ``La formica e il miele'', altri dalle gare Kangourou).

\subsubsection{Composizione dei gruppi}

I bambini si sono suddivisi liberamente in cinque gruppi: 3 da 4
componenti e 2 da 5 componenti (due alunni sono arrivati in ritardo e
si sono inseriti dove c'era posto, cioè non hanno potuto scegliere i
compagni di gruppo).

Quattro gruppi sono abbastanza eterogenei, sia dal punto di vista
delle conoscenze/abilità matematiche che della capacità di attenzione
(anche se due gruppi sono formati da tutti maschi, che in generale in
questa classe sono più agitati e con meno capacità di
concentrazione). Non c'è mescolanza di maschi e femmine perché si sono
scelti liberamente in base all'amicizia e, soprattutto quest'anno, è
più netta la separazione maschi e femmine nei rapporti
interpersonali. In generale però, quando lavorano insieme, riescono a
collaborare in modo proficuo.

Nel quinto gruppo si sono riuniti tutti e tre i bambini con difficoltà
di apprendimento. Anche in altre occasioni, quando i gruppi si sono
formati in base alla libera scelta degli alunni, questi bambini hanno
voluto lavorare insieme. Questo elemento è di vantaggio nel senso che,
quando sono insieme, si sentono più sicuri e lavorano in modo più
spontaneo e attivo, mentre, quando lavorano in gruppo con gli altri
bambini, tendono a delegare, o in alcuni casi sono i compagni che
tendono a sovrastarli.

Di solito quando si lavora in gruppo noi insegnanti tendiamo a mettere
in gruppi diversi i bambini con difficoltà, cercando di inserirli con
i compagni più collaborativi e che fungono da stimolo. In questa
occasione però abbiamo lasciato la scelta dei gruppi ai bambini.

In alcuni incontri i gruppi sono stati modificati dall'insegnante.

\subsubsection{Insegnanti presenti}
L'insegnante di classe è affiancata dall'insegnante di sostegno.

\subsubsection{Calendarizzazione degli incontri}
\begin{calendario}
  \begin{itemize}
  \item 22 gennaio dalle 8.30 alle 10.20
  \item 26 gennaio dalle 8.30 alle 10.20
  \item 8 febbraio dalle 8.30 alle 10.20
  \item 9 febbraio dalle 14.40 alle 16.20
  \end{itemize}
\end{calendario}

\subsection{Primo incontro}

\begin{description}\item[Alunni presenti:] 22 alunni
\item[Tempo effettivo di lavoro] dalle ore 8.30 alle ore 10.20, ovvero
  circa un'ora e mezza. La prima fase manipolativa è durata circa
  mezz'ora
\end{description}

\begin{consegna}
  Entrambe le insegnanti hanno ruotato tra i gruppi, solo all'inizio
  del gioco della tombola l'insegnante di sostegno ha seguito il
  gruppo in difficoltà. Abbiamo estratto le frazioni noi insegnanti e
  abbiamo lasciato che i bambini discutessero e trovassero da soli le
  soluzioni, intervenendo solo per sollecitarli nei ragionamenti.
\begin{enumerate}
\item Prima di iniziare l'attività proposta nel kit abbiamo
  distribuito a ogni gruppo una vetrata con i pezzi corrispondenti,
  lasciando giocare liberamente e senza dare alcuna indicazione.
\item In seguito abbiamo invitato gli alunni a continuare il gioco
  provando a ``spostare i pezzi'' o a ``coprire in modo diverso''.
\item Abbiamo poi distribuito a ogni gruppo le regole del gioco
  (anche in questo caso senza dare indicazioni preliminari). I bambini
  le hanno lette e discusse tra loro e poi ogni gruppo ha spiegato
  agli altri ciò che aveva capito. Noi insegnanti abbiamo ricapitolato
  le regole alla fine del confronto fra i gruppi.
\item Prima di iniziare il gioco della tombola, si stabilisce di
  utilizzare solo la prima modalità di gioco (cioè scegliere un solo
  pezzo per coprire uno spazio vuoto). Questa modalità è stata scelta
  da noi insegnanti per poter controllare meglio la correttezza delle
  mosse e per vedere se tutti i bambini individuavano subito i pezzi
  corrispondenti alla frazione estratta. Si procede poi con il gioco.
\end{enumerate}

\materiali{}%
Il kit comprende tre giochi; noi abbiamo cominciato con
\attivita{Realizziamo una vetrata}, utilizzando le prime cinque
vetrate (frazioni 1/2, 1/4 e 1/8).
\end{consegna}

\subsubsection{Osservazioni}
\begin{enumerate}
\item Tutti i gruppi hanno ricoperto l'intera vetrata, instaurando una
  sorta di gara, utilizzando i pezzi che coincidono perfettamente con
  i vetri delle finestre. Alcuni bambini hanno utilizzato i pezzi
  rimanenti per fare decorazioni.
\item Tutti hanno quindi ricoperto di nuovo la vetrata utilizzando più
  pezzi per coprire uno spazio vuoto e sperimentando varie
  combinazioni.

  Al termine di questa prima fase manipolativa, che è durata circa
  mezz'ora, ogni gruppo ha spiegato agli altri ciò che aveva fatto e
  un bambino ha chiesto se il gioco del kit riguardava la
  geometria. La domanda è stata riposta ai bambini i quali hanno
  risposto che poteva anche riguardare la logica e le operazioni,
  perché avevano diviso o sommato le figure.
\item Dopo la lettura delle regole del gioco, ogni gruppo ha
  illustrato la sua interpretazione e alcuni bambini hanno precisato
  che il gioco riguardava le frazioni. Il punto che ha suscitato
  domande e confronti fra i gruppi è quello che riguarda le due
  diverse possibilità di gioco: coprire uno spazio vuoto usando un
  unico vetro o più vetri. Tutti i bambini hanno capito la prima
  modalità di gioco (un bambino ha detto \bambini{è come la
    tombola!}), ma non sono riusciti subito a spiegare la seconda
  modalità (cioè usare più vetri per coprire un unico spazio). Li
  abbiamo allora invitati a ripensare a quanto fatto nella fase
  manipolativa, quando cioè avevano usato più pezzi per coprire uno
  spazio vuoto, sperimentando varie combinazioni. Attraverso vari
  interventi dei bambini di diversi gruppi si è arrivati a spiegare
  tutte e due le modalità di gioco.  Un bambino:
  \begin{studente}[]
    Si può usare un solo pezzo oppure più pezzi e si
  formano figure uguali
  \end{studente}
  un altro precisa:
  \begin{studente}[]
    Non sono proprio uguali ma hanno lo stesso valore
  \end{studente}
  E una bambina ricorda:
  \begin{studente}[]
    Ah sì, sono le frazioni equivalenti
  \end{studente}
\item All'inizio due gruppi (quello con gli alunni con difficoltà e un
  gruppo formato da quattro femmine, di cui tre con buone abilità
  matematiche e buona preparazione e una con molte lacune e poca
  capacità di concentrazione) hanno difficoltà a riconoscere a quale
  frazione corrisponde ogni pezzo, sia perché le finestre non sono
  divise in parti uguali, sia perché a volte la finestra è già coperta
  da altri pezzi che li mandano in confusione. Per capire devono
  coprire tutta la finestra con pezzi uguali, individuando così la
  frazione corrispondente. Usando tutti pezzi uguali per coprire la
  finestra, riescono infatti a individuare l'unità frazionaria. È una
  sorta di verifica per il gruppo delle quattro bambine, mentre è una
  necessità per il gruppo con difficoltà: i tre bambini segnalati
  hanno bisogno di disporre tutti i pezzi uguali e di contarli per
  capire in quante parti risulta frazionata la finestra.

  Dopo le prime incertezze il gioco prosegue speditamente e si arriva
  al gruppo vincitore.

  Tutti i bambini hanno partecipato attivamente, anche se alcuni solo
  nella parte del gioco vero e proprio.
\end{enumerate}

Nel prossimo incontro vogliamo riproporre il gioco, lasciandoli liberi
di combinare più tessere per coprire gli spazi vuoti e, se c'è tempo,
utilizzando tutte le vetrate (anche quelle con le frazioni 1/3 1/6 e
1/9).

Al termine del gioco vorremmo far scrivere a ogni gruppo le
riflessioni sull'attività svolta, ma non sappiamo se sia meglio dar
loro una traccia con delle domande precise per guidarli o se lasciarli
esprimere liberamente, dando indicazioni generiche.

\subsubsection{Consigli per i colleghi che vogliono proporre le stesse
  attività}

A noi è sembrato molto utile far precedere al gioco vero e proprio la
fase manipolativa con gioco libero, sia per far entrare gli alunni in
un'atmosfera più ludica della normale lezione, sia perché in questo
modo gli alunni hanno cominciato a conoscere i pezzi del kit e nella
fase successiva sono riusciti a trovare più facilmente e più
velocemente i vetri corrispondenti alle frazioni estratte.

Per quanto riguarda il gioco della tombola, prima di iniziare è
indispensabile chiarire le due possibilità di gioco e scegliere se
ammetterne solo una o tutte e due. A noi è sembrato più utile iniziare
con la prima modalità (cioè scegliere un solo pezzo per coprire uno
spazio vuoto), perché, avendo affrontato l'anno scorso il concetto di
frazione e non avendo fatto un ripasso, volevamo essere sicure che
tutti i bambini riuscissero almeno a individuare la frazione estratta,
prima di passare alla somma di frazioni. Nel corso del gioco però
alcuni bambini si sono lamentati proprio perché non potevano combinare
i vetri.

Secondo noi questo gioco può essere proposto anche in una classe
quarta, una volta che sia stato introdotto in linea generale il
concetto di frazione.

\subsection{Secondo incontro}
\begin{description}
\item[Alunni presenti:] Erano assenti cinque alunni (non abbiamo
  rimandato l'attività perché è difficile avere due ore consecutive in
  compresenza e non volevamo che i due incontri fossero troppo
  distanti nel tempo).

  Perché i gruppi non fossero troppo squilibrati nel numero, abbiamo
  apportato delle variazioni, provando a spostare in altri gruppi i
  bambini con difficoltà di apprendimento. Si sono quindi creati tre
  gruppi da tre componenti (uno formato da tre maschi con buone
  conoscenze matematiche e con discreta capacità di attenzione; uno
  con un bambino e una bambina con difficoltà di apprendimento e un
  bambino con buone conoscenze matematiche; un gruppo con due maschi e
  una femmina con discrete conoscenze matematiche e scarsa capacità di
  attenzione) e due da quattro (formati da tutte femmine, ma
  abbastanza eterogenei sia dal punto di vista delle
  conoscenze/abilità matematiche che della capacità di attenzione. In
  uno si è inserita una bambina con difficoltà).
\item[Tempo effettivo di lavoro] Dalle 8.30 alle 10.30, anche se
  l'attività vera e propria è durata circa un'ora e mezza
\end{description}
\begin{consegna}
  Entrambe le insegnanti hanno ruotato tra i gruppi e i bambini sono
  stati lasciati liberi di discutere, affinché trovassero da soli le
  soluzioni.
  \begin{enumerate}
  \item Gioco con le prime cinque vetrate, lasciando la libertà di
    utilizzare più pezzi per formare la frazione estratta.
  \item Trovare le nuove frazioni che potevano essere estratte con le
    vetrate da 6 a 10.
  \end{enumerate}

  \materiali{}%
  Gioco con le prime cinque vetrate e successivamente con le vetrate
  da 6 a 10.
\end{consegna}

\subsubsection{Osservazioni}

1. I bambini hanno chiesto se dovevano rispettare gli spazi vuoti o se
potevano mettere i pezzi liberamente sulle finestre senza rispettare
la suddivisione.

Dopo un breve confronto fra insegnanti e alunni abbiamo deciso insieme
di provare questa seconda modalità di gioco, ma abbiamo scoperto che,
così facendo, tutti i gruppi vincevano nello stesso momento. In
generale non ci sono stati problemi nell'individuare le frazioni e nel
combinare i pezzi: i bambini hanno dimostrato di aver dimestichezza
con i vetri da 1/2, 1/4 e 1/8. Inoltre c'è stata più organizzazione:
tutti i gruppi hanno stabilito al proprio interno turni per ricercare
i pezzi, in modo da far partecipare tutti e da evitare litigi.

2. Abbiamo quindi ritirato le vetrate e riconsegnato ai gruppi le
nuove vetrate, dalla 6 alla 10, chiedendo di trovare le nuove frazioni
che potevano essere estratte. Dopo qualche minuto, in cui i bambini
hanno manipolato e provato i nuovi vetri, tutti hanno individuato le
frazioni 1/3, 1/6 e 1/9.

Prima di iniziare la seconda partita, abbiamo stabilito (le insegnanti
insieme ai bambini) che si potevano usare più vetri per formare una
frazione, ma anche che si dovevano rispettare gli spazi vuoti.

Questa modalità di gioco è quella che, secondo il parere dei bambini,
si è dimostrata più facile (\bambini{perché se non hai un pezzo lo
  puoi formare}) e è anche quella che richiede più abilità
(\bambini{perché non conta solo la fortuna, come quando non puoi
  combinare i pezzi}).

Aggiungo che questa modalità crea anche più confronto nel gruppo.

È capitato a esempio che un gruppo dicesse di non avere la frazione
estratta; sollecitati dall'insegnante a ``cercare bene'' i bambini
cominciavano allora a provare varie combinazioni e a discutere fra di
loro, trovando poi il modo per formare la frazione.

C'è da dire anche però che richiede più tempo, sia per la ricerca da
parte dei bambini dei pezzi giusti, sia per il controllo delle mosse
da parte dell'insegnante. Ma ne vale la pena perché è più stimolante,
crea maggior discussione e confronto all'interno del gruppo.

In questo incontro non c'è stato tempo per far scrivere le riflessioni
sull'attività; pensavamo di farlo in un altro momento e di dedicare
poi gli altri incontri al gioco del domino sulle frazioni equivalenti.

\subsubsection{Consigli per i colleghi che vogliono proporre le stesse attività}

A noi è sembrato utile andare per gradi, cioè fare più partite e
giocare prima con le cartelle più facili (dalla prima alla quinta) e
poi con le altre, perché in questo modo tutti i bambini sono riusciti
a prendere dimestichezza con il gioco e a divertirsi. È anche utile,
se si fanno più partite, far scambiare le cartelle ai gruppi, perché,
essendo diverse, presentano difficoltà diverse.

Su sollecitazione dei bambini abbiamo provato in una partita a far
mettere i pezzi sulla finestra senza rispettare gli spazi vuoti. Ma
questa modalità di gioco è sconsigliabile perché, non ponendo limiti,
crea meno difficoltà e quindi meno discussione; inoltre tutti i gruppi
hanno vinto nello stesso momento.

\subsection{{Terzo incontro}}
\begin{description}\item[Alunni presenti:]
  20 alunni presenti e due assenti.  Gli alunni si sono divisi nei
  cinque gruppi del primo incontro. I due bambini assenti fanno parte
  dello stesso gruppo (quello in cui sono presenti i bambini con
  difficoltà, uno è l'assente) che risulta quindi formato da tre
  componenti. Gli altri rimangono invariati.
\item[Tempo effettivo di lavoro:] dalle ore 8.30 alle ore 10.20. Un'ora
  e mezza circa di effettiva attività.
\end{description}

\begin{consegna}
  \begin{enumerate}
  \item Viene riproposto il gioco con le vetrate dalla 6 alla 10. Si
  consegna a ogni gruppo una vetrata, diversa da quelle usate nelle
  partite precedenti. Prima di iniziare i bambini ricapitolano le
  regole del gioco e insieme si stabilisce di utilizzare la modalità
  di gioco che si è dimostrata più proficua e più divertente: si
  lascia la possibilità di usare più vetri per formare una frazione,
  ma devono essere rispettati gli spazi vuoti.
\item Al termine del gioco, che è durato complessivamente tre quarti
  d'ora (compresa la condivisione delle regole), a ogni gruppo sono
  stati consegnati dei fogli sui cui scrivere le riflessioni sulle
  attività svolte.
\end{enumerate}

Entrambe le insegnanti hanno ruotato tra i gruppi e i bambini sono
stati lasciati liberi di discutere, affinché trovassero da soli delle
soluzioni.  %
\materiali{}%
Gioco con le vetrate da 6 a 10. Fogli su cui scrivere le riflessioni.
\end{consegna}

\subsubsection{Osservazioni}
1. Il gioco si è svolto senza particolari problemi, perché i bambini,
anche se era passata più di una settimana dall'incontro precedente,
ricordavano le combinazioni sperimentate nelle partite precedenti. Due
gruppi, gli stessi che hanno avuto bisogno delle sollecitazioni delle
insegnanti negli incontri precedenti, sono riusciti a combinare le
frazioni più ``facili'' (ad esempio 1/4 e 1/4 per fare 1/2), ma non
hanno trovato la frazione quando si trattava di unire 1/6 e 1/6 per
formare 1/3.

2. Per quanto riguarda le riflessioni, si è chiesto ai bambini di
esprimersi sul lavoro di gruppo (come si sono trovati con i compagni,
se volevano lavorare in un altro gruppo, se hanno trovato
un'organizzazione interna) e sul gioco (se hanno imparato cose nuove,
se hanno ripassato o capito meglio concetti già appresi in precedenza,
se hanno trovato difficoltà e come le hanno superate, se si sono
divertiti). Si è lasciato poi uno spazio per le ``riflessioni libere''.

Questo lavoro è stato fatto in gruppo, ma alcuni bambini hanno
espresso opinioni personali che potevano essere in disaccordo con il
resto del gruppo.

Ogni gruppo ha poi scelto un portavoce per leggere le riflessioni, ma,
per mancanza di tempo, si è rimandato il confronto fra i gruppi
all'incontro successivo.

\paragraph{Sintesi delle riflessioni dei vari gruppi:}
\begin{itemize}
\item tutti hanno trovato piacevole e divertente il lavorare in
  gruppo; solo alcune bambine dello stesso gruppo hanno scritto che in
  alcuni momenti si sono annoiate. Uno dei motivi di divertimento è
  stato il fatto di aver potuto scegliere i propri amici come compagni
  di squadra, anche se questo li ha portati a distrarsi e a perdere la
  concentrazione. Due gruppi hanno anche segnalato difficoltà:
  \bambini{lavorare in gruppo è stato anche complicato perché dovevamo
    trovare un'intesa comune e con un po' di pazienza ci siamo
    riusciti}; \bambini{a volte è stato difficile riunire le idee}.
\item Tutti hanno detto di aver trovato un'organizzazione all'interno
  del gruppo stabilendo dei turni, anche se a volte ci sono stati
  ``litigi''.
\item Tutti hanno trovato il gioco divertente e un gruppo ha
  specificato \bambini{perché dovevamo combinare i pezzi in modo
    strategico}. Hanno aggiunto che il gioco è servito per ripassare
  in generale le frazioni e le forme geometriche e \bambini{in
    particolare per capire meglio alcune cose, come ad esempio che 2/6
    equivale a 1/3}. Hanno anche trovato difficoltà \bambini{a formare
    alcune figure o a trovare il posto adatto per appoggiare la
    figura} e le hanno superate discutendo fra di loro o chiedendo
  aiuto alle maestre.
\end{itemize}

\subsubsection{Consigli per i colleghi che vogliono proporre le stesse
  attività}
Nelle varie partite noi abbiamo lasciato più o meno invariati i
gruppi, per rispettare le scelte dei bambini, ma così facendo abbiamo
notato che permanevano le stesse difficoltà negli stessi gruppi. Può
quindi essere utile cambiare i componenti del gruppo.

\subsection{Quarto incontro}
\begin{description}\item[Alunni presenti:] 17 alunni presenti e 5
  assenti. I gruppi sono stati scelti dalle insegnanti, perché nel
  pomeriggio i bambini che hanno più difficoltà di concentrazione sono
  più agitati e per creare gruppi più equilibrati per quanto riguarda
  la suddivisione maschi e femmine. Si sono creati 4 gruppi abbastanza
  eterogenei sia dal punto di vista del comportamento/attenzione, che
  da quello delle competenze/conoscenze matematiche:

  GRUPPO 1: due maschi e due femmine. Una bambina con molte lacune in
  matematica, gli altri con buone capacità.

  GRUPPO 2: due maschi e due femmine. Una bambina con il sostegno e
  gli altri con buone/ottime capacità e competenze matematiche.

  GRUPPO 3: tre femmine e un maschio molto agitato. Tutti con discrete
  o buone capacità e competenze matematiche.

  GRUPPO 4: un bambino con difficoltà di apprendimento, un bambino con
  buone capacità matematiche e molto disponibile a aiutare i compagni
  in difficoltà e due bambine con buone capacità e competenze
  matematiche.
\item[Tempo effettivo di lavoro:] Dalle 14.40 alle 16.20. Un'ora e
  mezza circa di effettiva attività.
\end{description}

\begin{consegna}
  \begin{enumerate}
  \item Ogni portavoce legge agli altri le riflessioni del proprio
    gruppo. Al termine si sollecitano gli alunni a esprimere opinioni
    su quanto espresso dai compagni. L'attività è durata circa tre
    quarti d'ora.
  \item Gioco del domino:
    \begin{itemize}
    \item Viene consegnato a ogni gruppo un numero uguale di carte, senza spiegare qual è il gioco.
    \item Viene lasciato del tempo (circa dieci minuti) per osservare le carte e capire il gioco.
    \item Ogni gruppo spiega agli altri le sue ipotesi sulle regole del gioco.
    \item Vengono distribuite a ogni gruppo le regole e vengono lette dall'insegnante.
    \item Si procede con il gioco.
    \end{itemize}
  \end{enumerate}
  \materiali{}%
  Gioco del domino.
\end{consegna}

\subsubsection{Osservazioni}

1) La lettura delle riflessioni dei vari gruppi non ha provocato
discussioni e confronti. Noi insegnanti abbiamo quindi cercato di
stimolare degli interventi chiedendo se questo gioco poteva essere
proposto a altre classi. Queste le risposte:
\begin{studente}[ ]
  \begin{itemize}
  \item \bambini{per fare questo gioco bisogna conoscere le frazioni};
  \item \bambini{si può proporre a tutte le classi, anche in prima, ma
      dando spiegazioni più semplificate delle regole, dicendo ad
      esempio metà figura invece di 1/2};
  \item \bambini{si può proporre a tutte le classi perché con questo
      gioco si possono imparare le frazioni e anche i nomi delle
      figure geometriche};
  \item \bambini{in una classe prima si può proporre facendo delle
      vetrate più semplici o comunque usando solo le prime cinque
      vetrate};
  \item \bambini{per i bambini più piccoli è più facile capire le
      frazioni attraverso un gioco, piuttosto che attraverso una
      lezione normale};
  \item \bambini{per gli alunni delle medie questo gioco può essere un
      ripasso o un approfondimento e può servire per capire meglio
      alcune cose. E poi è divertente}.
  \end{itemize}
\end{studente}
2) Gioco del domino: i gruppi 1 e 2 hanno capito da soli, discutendo
fra i componenti, che il gioco consisteva nell'unire le tessere con
frazioni equivalenti e il gruppo 1 ha specificato che si poteva
trattare del gioco del domino, spiegandone le regole generali. Il
gruppo 4 ha interpellato le insegnanti per chiedere conferma circa le
ipotesi fatte, ma ha comunque compreso la regola per unire le
tessere. Il gruppo 3 invece non ha pensato al gioco del domino e ha
cercato di trovare una coerenza fra le frazioni presenti su una stessa
tessera (ad esempio sommando i numeratori e i denominatori).

Dopo la condivisione e il confronto delle ipotesi di ogni gruppo
abbiamo letto insieme le regole del gioco.

Per cominciare il gioco noi insegnanti abbiamo posto su un tavolo al
centro dell'aula due tessere già collegate e abbiamo dato il via alla
partita. All'inizio tutti i gruppi hanno mostrato difficoltà, perché,
pur avendo compreso il concetto di frazioni equivalenti, faticavano a
trovare le combinazioni fra le frazioni con numeri più alti. Io e la
mia collega siamo dovute intervenire suggerendo alcune mosse e
spiegando le strategie usate (ad esempio moltiplicare o dividere per
uno stesso numero i termini della frazione). Sono bastati pochi
suggerimenti per sbloccare la situazione e, dopo le prime difficoltà,
si è instaurata la competizione fra le squadre e la partita è
proseguita speditamente (abbiamo dovuto stabilire un tempo massimo per
dare le risposte).

La partita è finita perché nessun gruppo aveva tessere da attaccare e
sono stati dichiarati vincitori due gruppi ai quali erano rimaste solo
tre tessere.

I bambini, nonostante la partenza un po' titubante, hanno apprezzato
molto il gioco e hanno chiesto di poter costruire un loro domino per
poter giocare anche nei momenti di intervallo.

Non essendo riuscite a organizzare in tempi brevi un altro incontro in
compresenza, questa attività di costruzione è stata svolta da alcuni
alunni nei ritagli di tempo all'interno delle lezioni, utilizzando le
stesse frazioni del kit in dotazione.

\subsubsection{Consigli per i colleghi che vogliono proporre le stesse
  attività}
Noi abbiamo ritenuto opportuno proporre il gioco del domino dopo
quello delle vetrate, per vedere se, dopo aver manipolato e
``costruito'' alcune frazioni equivalenti, queste venivano ricordate e
riconosciute. In effetti con le frazioni più semplici questo è
avvenuto.

Probabilmente c'è bisogno di fare più partite perché i bambini
riconoscano velocemente anche le equivalenze con frazioni più
difficili.


\section[Sperimentazione \#2: terza primaria]{Sperimentazione \#2:
  classe terza primaria, gennaio/febbraio~2010}

\subsection{Osservazioni generali}

Le due docenti sperimentatrici
hanno progettato insieme il percorso da proporre agli alunni delle
classi terze (classe ``A'' e classe ``B'')
e hanno raccolto e sintetizzato il lavoro svolto dalle due classi
sulle prime cinque vetrate come di seguito riportato.

\subsubsection{Presentazione della classe}
La classe ``A'' è formata da 26 alunni, 19 maschi e 7 femmine.
Sono presenti alunni stranieri che sono nati in Italia e non
presentano difficoltà linguistiche, un'alunna con deficit di
apprendimento, che ha seguito la proposta con i colleghi di sostegno,
e un alunno con scarse capacità logiche e deficit di apprendimento non
ancora certificati.  Dall'anno scorso in classe viene utilizzato il
lavoro di gruppo per svolgere attività di laboratorio teatrale.  Gli
alunni non hanno ancora affrontato il concetto di frazione

La classe ``B'' è formata da 26 alunni,15 maschi e 11 femmine.  Sono
presenti alunni stranieri che sono nati in Italia e non presentano
difficoltà linguistiche, tre alunni con difficoltà specifiche di
apprendimento, i quali hanno seguito la proposta con i colleghi di
sostegno (cfr.~sperimentazione\#3).  Anche in questa classe gli alunni
non hanno ancora affrontato il concetto di frazione

\subsubsection{Composizione dei gruppi}

Suddivisione dei bambini in cinque gruppi omogenei per sesso e
competenze
\begin{itemize}
\item Nella classe ``A'' il bambino con alcune difficoltà non
  certificate è stato inserito in un gruppo con compagni più disposti
  alla collaborazione
\item Nella classe ``B'' i due bambini con difficoltà specifiche di
  apprendimento sono stati inseriti in un gruppo con compagni più
  disposti alla collaborazione
\end{itemize}

\subsubsection{Insegnanti presenti}
Per gli incontri della classe ``A'' sono presenti le due insegnanti di
classe e l'insegnante di sostegno.
Per gli incontri della classe ``B'' sono presenti la docente
sperimentatrice, l'insegnante di sostegno e una educatrice.

\subsubsection{Calendarizzazione degli incontri}
\begin{calendario}
  \begin{itemize}
  \item 25 gennaio (classe ``A'') e 26 gennaio (classe ``B'')
  \item 1 e 2 febbraio
  \item 8 e 9 febbraio
  \item 15 e 16 febbraio
  \end{itemize}
\end{calendario}

\subsection{Primo incontro}

\begin{description}\item[Alunni presenti:] Nella classe ``A'' erano
  presenti tutti gli alunni tranne la bambina con deficit di
  apprendimento che ha seguito l'attività con il gruppo dei colleghi
  di sostegno. Nella classe ``B'' era assente solo un alunno con
  deficit di apprendimento
\item[Tempo effettivo di lavoro:] 1 ora: classe ``A'' dalle 15.15 alle
  16.15; classe ``B'' dalle 15.15 alle 16.15
\end{description}

\begin{consegna}
  \begin{itemize}
  \item elezione dei capigruppo
  \item Spiegazioni riguardanti il materiale, rispetto dello stesso e
    ``regole del gioco'' (attività da svolgere).
\end{itemize}
\materiali{}%
Nel primo incontro all'inizio non è stato distribuito alcun materiale,
solo dopo la formazione dei gruppi è stata presentata la scatola del
kit Viaggio segreto, e in particolare sono stati mostrati i sacchetti
con le figure geometriche.
\end{consegna}

\subsubsection{Osservazioni}

I bambini hanno mostrato interesse per l'incontro, anche se in questa
prima fase non hanno manipolato il materiale.

Non sono stati raccolti materiali.

L'incontro è stato necessario in quanto le due classi sono piuttosto
esuberanti e è stato necessario puntare sul rispetto del materiale e
sulle regole del lavoro di gruppo.

\subsubsection{Consigli per i colleghi che vogliono proporre le stesse attività}

In questo primo incontro non abbiamo indicazioni particolari da
proporre. Favorevole è stata la possibilità di lavorare durante le ore
di compresenza.
\subsection{Secondo incontro}
\begin{description}\item[Alunni presenti:]
  Nella classe ``A'' erano presenti tutti gli alunni tranne la
    bambina con deficit di apprendimento che ha seguito l'attività con
    il gruppo dei colleghi di sostegno.
   Nella classe ``B'' erano presenti tutti gli alunni a
    eccezione di un alunno con deficit di apprendimento.
\item[Tempo effettivo di lavoro:] 1h 30 (classe ``A'' dalle 14.45 alle
  16.15; classe ``B'' dalle 14.45 alle 16.15)
\end{description}

\begin{consegna}
  Sono state utilizzate le prime cinque vetrate con frazioni ½, ¼,
  1/8.

  A ogni gruppo è stata consegnata una vetrata con i pezzi
  corrispondenti e i bambini sono stati lasciati liberi di manipolare
  il materiale per ricoprire le finestre senza dare alcuna spiegazione

  A ogni gruppo è stata consegnata la seguente griglia di
  osservazione:

  \begin{dadocente}
      OSSERVATE LA SEGUENTE VETRATA E COMPILATE LA TABELLA

  Quante finestre ci sono?

  Che forma ha ogni finestra?

 Utilizzando le forme colorate, ricoprite le vetrate.

 Quali forme sono state utilizzate per ricoprire i vetri di ogni finestra?

 C'è un solo modo per ricoprire le vetrate?

 Perché?

 Quante tessere possono servire per ricoprire una finestra?

 Per ogni finestra avete usato la stessa forma delle tessere?
  \end{dadocente}

 \materiali{}%
 Distribuzione del materiale ``kit viaggio segreto''
 \attivita{Realizziamo una vetrata}.
\end{consegna}

\subsubsection{Osservazioni}

Dopo la compilazione discussione e socializzazione dei risultati di
ogni gruppo, con esternazione di eventuali difficoltà:

Vetrata 1:

\begin{itemize}
\item I rettangoli hanno la stessa forma ma diverse dimensioni
\item Forme uguali con suddivisioni diverse
\item Una vetrata è formata da pezzi uguali
\item Non hanno avuto particolari difficoltà nel ricoprire le vetrate.
\end{itemize}

Vetrata 2:

\begin{itemize}
\item I triangoli hanno la stessa forma ma diverse dimensioni
\item Le finestre sono tutte uguali
\item Un triangolo è la metà dell'altro
\item Difficoltà nell'orientare in modo adeguato i pezzi
\end{itemize}

Vetrata 3:

\begin{itemize}
  \item Sono tutti divisi a metà
  \item 2 quadrati piccoli equivalgono a un rettangolo
  \item Non hanno avuto difficoltà
\end{itemize}

Vetrata4:

\begin{itemize}
\item C'è un rettangolo con forme tutte uguali;
\item 4 triangoli sono metà finestra, cioè un quadrato
\item Tutti i quadrati sono formati da figure diverse
\item Hanno avuto difficoltà nel posizionare il triangolo rettangolo.
\end{itemize}

Vetrata 5:

\begin{itemize}
\item I rettangoli hanno la stessa forma ma diverse dimensioni
\item forme uguali per coprire le vetrate
\item i triangoli, uno è la metà dell'altro
\item una finestra è divisa in rettangoli, in una metà sono messi in
  orizzontale e nell'alta in verticale, è diviso in otto parti
\item una metà è un quadrato, tutta la figura cioè un rettangolo è
  formata da due quadrati
\item Nessuna difficoltà nel posizionare le forme
\end{itemize}

Alla fine della discussione abbiamo concluso che tutte le figure
divise in parti uguali sono state frazionate e abbiamo sintetizzato le
risposte in un'unica tabella:
\begin{center}
  \begin{tabular}{|p{0.46\linewidth}|p{0.46\linewidth}|}
    \hline
    Quante finestre ci sono?&
    Ci sono sei finestre \\ \hline
    Che forma ha ogni finestra? &
    Ha la forma di un rettangolo \\ \hline
    Utilizzando le forme colorate, ricoprite le vetrate. &
    Sì, abbiamo ricoperto tutte le vetrate \\ \hline
    Quali forme sono state utilizzate per ricoprire i vetri di ogni
    finestra? &
    Triangoli, quadrati, rettangoli, pentagono \\ \hline
    C'è un solo modo per ricoprire le vetrate?
    Perché? &
    No, posso sostituire le tessere con altre che occupano lo stesso
    spazio \\ \hline
    Quante tessere possono servire per
    ricoprire una finestra? &
    Da una a otto tessere \\ \hline
    Per ogni finestra avete usato la stessa forma delle tessere? &
    In alcune sono state usate le stesse forme, in altre forme
    diverse. \\ \hline
  \end{tabular}
\end{center}

\subsubsection{Consigli per i colleghi che vogliono proporre le stesse attività}

Il lavoro svolto non ha presentato particolari difficoltà operative,
l'attenzione maggiore è stata posta nell'osservare che ciascun alunno
partecipasse attivamente al lavoro di gruppo. Anche in questo caso è
stato fondamentale prevedere l'attività durante le ore di
``compresenza''. L'attività proposta può a nostro parere essere
proposta anche in classi precedenti

\subsection{Terzo incontro}

\begin{description}
\item[Alunni presenti:]
  Nella classe ``A'' erano presenti tutti gli alunni tranne la
    bambina con deficit di apprendimento che ha seguito l'attività con
    il gruppo dei colleghi di sostegno.
  Nella classe ``B'' erano presenti tutti gli alunni a
    eccezione di un alunno con deficit di apprendimento.
  \item[Tempo effettivo di lavoro:] 1 ora e 30 minuti (classe ``A''
    14.45-16.15; classe ``B'' 8.45-10.15).
\end{description}

\begin{consegna}
  Ricoprire in modo diverso le vetrate, cioè senza tenere conto delle
  suddivisioni date.

  Distribuzione delle vetrate ai singoli gruppi e manipolazione del
  materiale.

  \materiali{}%
  Kit: Vetrate n 1-2-3-4-5 e relativi sacchetti con le forme

  Fotocopie delle vetrate proposte da utilizzare singolarmente sul
  quaderno
\end{consegna}

\subsubsection{Osservazioni}
Dopo la manipolazione e la discussione nel piccolo gruppo,
condivisione nel grande gruppo:
\begin{itemize}
\item Alcune finestre erano divise in parti uguali (un mezzo, un
  quarto e un ottavo) e altre no.
\item I pezzi potevano essere sostituiti con altrettanti più piccoli o
  più grandi
\item Lo stesso spazio poteva essere occupato con pezzi diversi
  (equiestensione/equivalenza)
\end{itemize}

Al termine della discussione a ogni singolo alunno è stata data la
fotocopia della vetrata che aveva precedentemente utilizzato; ogni
alunno ha ritagliato i pezzi, li ha colorati e ha riprodotto la
vetrata sul quaderno. Purtroppo non abbiamo foto da allegare

\subsubsection{Consigli per i colleghi che vogliono proporre le stesse attività}

La presenza di più insegnanti ha permesso di svolgere le attività in
modo efficace, soprattutto durante il lavoro individuale sul quaderno
con le fotocopie delle vetrate. Il lavoro con la carta è stato più
difficoltoso in quanto il materiale è risultato più fragile e quindi
le forme ritagliate non erano precise come quelle proposte dal
kit. Sicuramente l'avere a disposizione più kit all'interno della
classe durante l'attività avrebbe giovato.

\subsection{Quarto incontro}
\begin{description}
\item[Alunni presenti:] Nella classe ``A'' erano presenti tutti gli
  alunni tranne la bambina con deficit di apprendimento che ha seguito
  l'attività con il gruppo dei colleghi di sostegno.  Nella classe
  ``B'' erano presenti tutti gli alunni a eccezione di un alunno con
  deficit di apprendimento.
\item[Tempo effettivo di lavoro:] 1 ora e 30 minuti
\end{description}

\begin{consegna}
  Ogni gruppo ha potuto utilizzare in un primo tempo solo i pezzi
  indicati dalla frazione estratta. In un secondo momento i bambini
  hanno coperto utilizzando uno o più pezzi corrispondenti alla
  frazione estratta.
  \materiali{}%
  Gioco della tombola

  Per il gioco abbiamo utilizzato solo frazioni 1/2, 1/4, 1/8
\end{consegna}

\subsubsection{Osservazioni}

Nella seconda fase di lavoro, quando i bambini hanno coperto
utilizzando uno o più pezzi corrispondenti alla frazione estratta, gli
alunni più ``attenti'' hanno notato che la frazione 1/8 non poteva
essere sostituita, in quanto era la parte più piccola.

\subsection{Osservazioni finali}
Sicuramente il materiale è stato molto utile per avviare i bambini ai
concetti di frazione e di equivalenza. Sarebbe stato ancora più
vantaggioso se avessimo potuto avere il kit a disposizione per più
tempo. Il manipolare ha suscitato negli alunni grandissimo interesse e
entusiasmo, chiedendo anche di poterlo acquistare. Una maggiore
verifica di ciò che è stato appreso sarà possibile nel prossimo anno,
quando affronteremo ampliamente la frazione.


\section[Sperimentazione \#3: sostegno primaria]{Sperimentazione \#3:
  sostegno primaria, febbraio~2010}

\subsection{Osservazioni generali}
La sperimentazione non si svolge in classe, ma con gruppi di alunni
diversamente abili.

\subsubsection{Presentazione della classe}

Gli insegnanti di sostegno decidono di somministrare le attività sulle
vetrate del Kit ``Viaggio Segreto'' agli alunni da loro seguiti (quattro
alunni di quarta e due di terza). Con tale attività si vuole
introdurre il concetto di frazione attraverso la manipolazione di
materiale sperimentale, separatamente dal gruppo classe che ha già
consolidato tale nozione. Infatti, quando le insegnanti curricolari
avevano affrontato l'argomento, i tempi di apprendimento rispetto alle
abilità e alle potenzialità dei bambini (secondo il loro Piano
Educativo Individualizzato) non erano ancora tali da consentire loro
un primo approccio alle frazioni. Diversamente, sarebbe stato
possibile sperimentare il kit in classe, consentendo anche ai bambini
in situazione di handicap una maggiore integrazione e collaborazione
con il gruppo classe attraverso il cooperative-learning.

Tenendo conto delle disabilità e, di conseguenza, dei tempi di
attenzione e concentrazione ridotti, si prevedono due soli incontri
con durata non superiore a un'ora e mezza, per volta.

\subsubsection{Composizione dei gruppi}

Tre femmine e tre maschi: cinque alunni diversamente abili (due con
disturbi dell'apprendimento, uno con ADHD, uno con ritardo mentale,
uno con ritardo evolutivo) e una bambina normodotata con svantaggio
socio-culturale.

I bambini saranno divisi in due gruppi da tre, la cui composizione è
determinata dagli insegnanti secondo i diversi livelli cognitivi.

\subsubsection{Insegnanti presenti}
Il \emph{team} di insegnanti di sostegno (3 docenti)

\subsubsection{Calendarizzazione degli incontri}
\begin{calendario}
  \begin{itemize}
  \item 3 febbraio
  \item 10 febbraio
  \end{itemize}
\end{calendario}

\subsection{Primo incontro}

\begin{description}\item[Alunni presenti:] tutti
\item[Luogo:] aula di sostegno
\item[Tempo effettivo di lavoro:] un'ora e mezza, dalle 11,00 alle
  12,30
\end{description}

\begin{consegna}
  Stabiliti i gruppi, i maestri consegnano il materiale agli alunni
  senza dare alcuna indicazione, lasciandoli liberi di sperimentare e
  giocare insieme.  %
  \materiali{}%
  vetrata 1 primo gruppo, vetrata 2 secondo gruppo
\end{consegna}

\subsubsection{Osservazioni}

Stabiliti i gruppi, i maestri consegnano il materiale agli alunni
(vetrata 1 primo gruppo, vetrata 2 secondo gruppo) senza dare alcuna
indicazione, lasciandoli liberi di sperimentare e giocare insieme.

Entrambi i gruppi osservano con curiosità e manipolano con entusiasmo
le tessere, cercando in modo istintivo di ricoprire con le stesse le
finestre delle due vetrate.

\begin{figure}[hbtp]
  \centering
  \label{pic:vetrate:3}
  \begin{tabular}{cc}
    \includegraphics[width=0.47\textwidth]{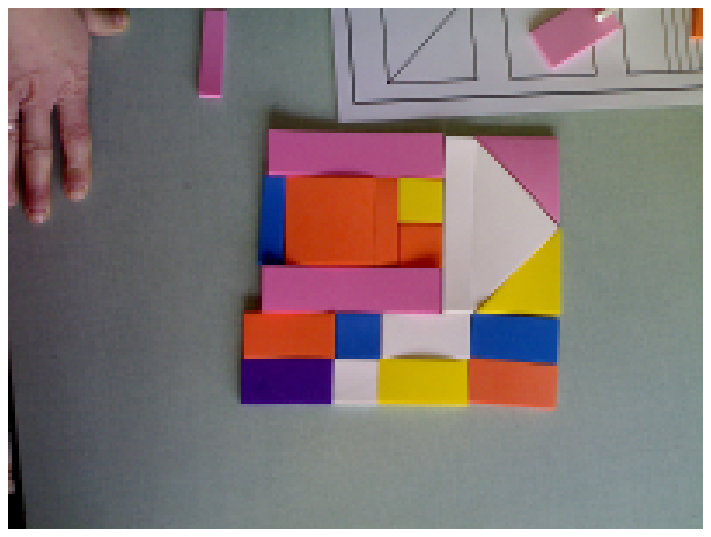} &
    \includegraphics[width=0.47\textwidth]{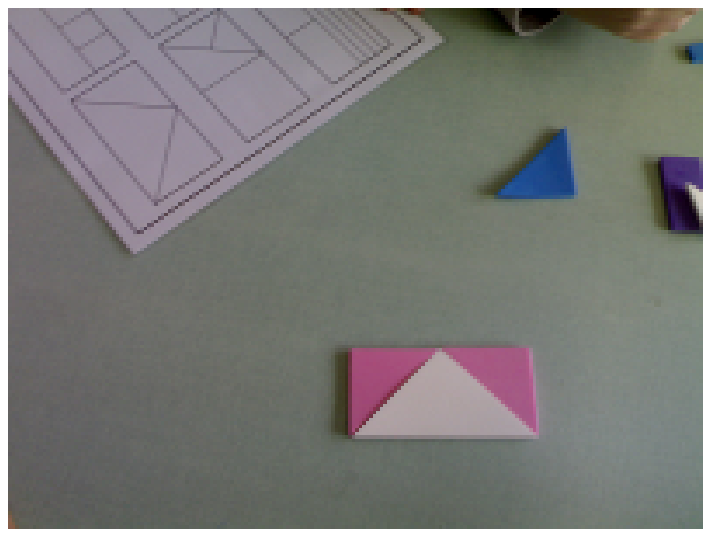} \\
    \includegraphics[width=0.47\textwidth]{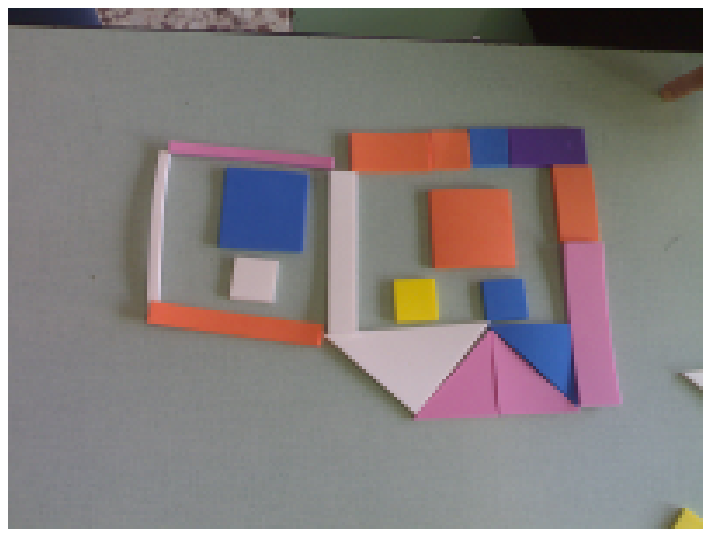} &
    \includegraphics[width=0.47\textwidth]{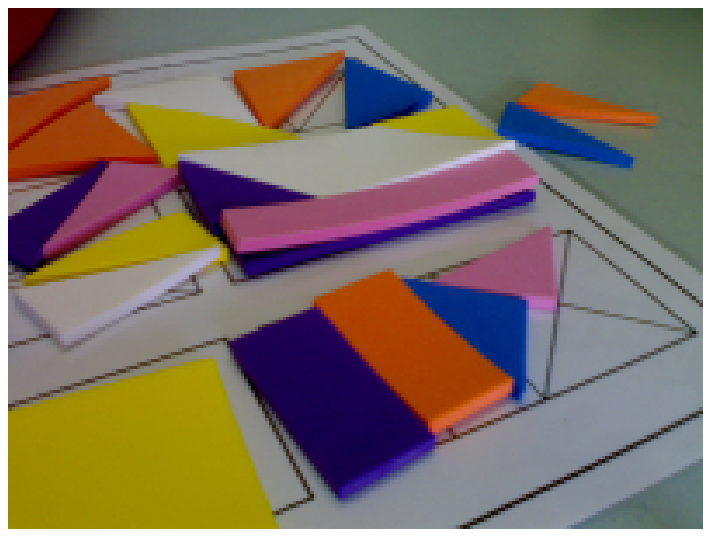} \\
    \includegraphics[width=0.47\textwidth]{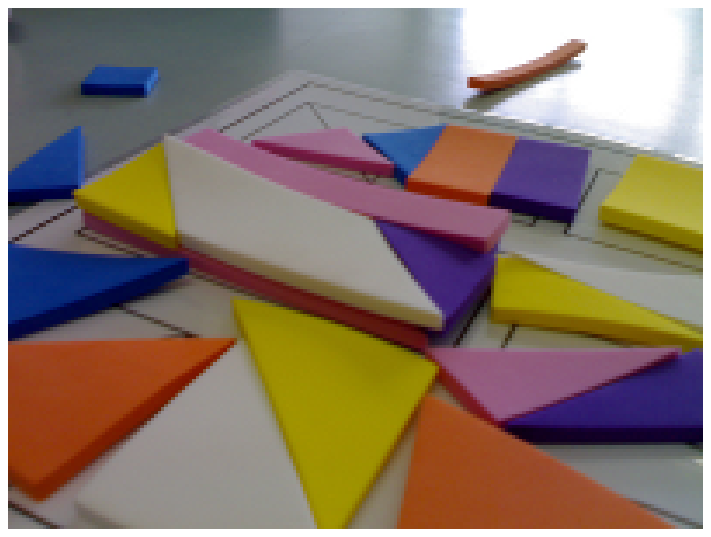} &
    \includegraphics[width=0.47\textwidth]{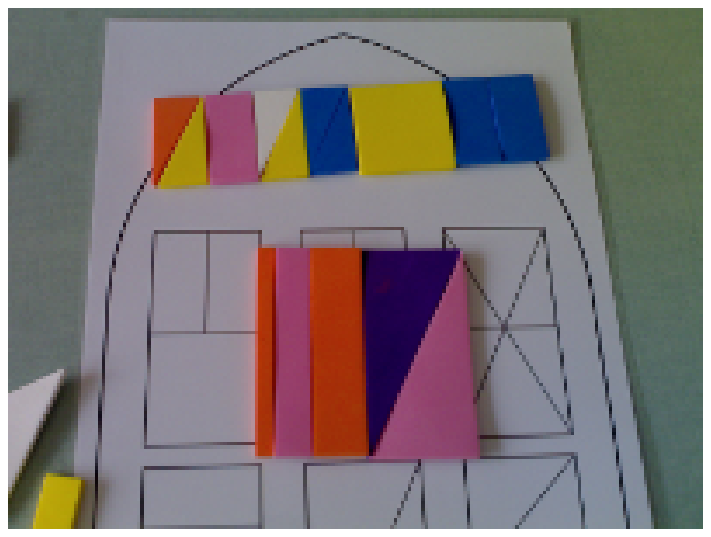} \\
  \end{tabular}
\end{figure}
Si osserva che il \textit{primo gruppo}, composto da soli alunni di
genere maschile, esegue con poche difficoltà la copertura delle
finestre (vetrata 2) con le rispettive tessere. Durante la
manipolazione, prima di trovare la giusta collocazione della tessera
presa, i bambini sperimentano più rotazioni, arrivando a comprendere
solo successivamente che la forma presa in esame ha una sola posizione
nella finestra.
\begin{studente}[Bb]
  Maestra ma sei sicura che non manca qualche pezzo?
\end{studente}
\begin{studente}[Bb]
  Alcuni pezzi mi sembrano troppo grandi, altri troppo piccoli.
\end{studente}
\begin{tutor}[Ins]
  Fai più tentativi nel mettere le tessere!
\end{tutor}
Una volta raggiunto l'obiettivo di ricoprire le finestre,
spontaneamente decidono di utilizzare i pezzi per costruire altre
forme a piacere (casa, pinocchio, palazzo).

I bambini scoprono che con i pezzi dati di varie forme e grandezze
geometriche è possibile realizzare oggetti da loro immaginati. In
sintesi utilizzano le tessere come se appartenessero a un tangram.

Il \textit{secondo gruppo, }composto da sole bambine, ricopre le
finestre della vetrata 1 richiedendo più tempo per l'esecuzione
rispetto all'altro gruppo, a causa di maggiori difficoltà. Le alunne,
infatti, non riescono a individuare facilmente la giusta disposizione
delle tessere, l'incastro e il combaciare dei lati fra i triangoli più
piccoli.

\subsection{Secondo incontro}

\begin{description}
\item[Alunni presenti:] tutti
\item[Luogo:] aula di sostegno
\item[Tempo effettivo di lavoro:] un'ora e mezza, dalle 11,00 alle
  12,30
\end{description}

\begin{consegna}
  Durante il secondo incontro, si invertono le vetrate con la consegna
  di ricoprire le finestre con le tessere corrispondenti. %
  \materiali{}%
  vetrata 2 primo gruppo, vetrata 1 secondo gruppo
\end{consegna}

\subsubsection{Osservazioni}
Il \textit{gruppo dei maschi}, ancora una volta, esegue senza
difficoltà il compito. Successivamente, avendo memorizzato la
disposizione delle forme delle finestre, riproduce volontariamente le
stesse sul tavolo senza guardare la vetrata.

Il \textit{gruppo delle bambine}, dopo una discussione sulla divisione
delle finestre, esegue la stessa consegna data ai compagni, ma questa
volta con minore difficoltà rispetto al primo incontro.

Gli insegnanti pongono a tutti i bambini delle domande-stimolo sulla
forma delle tessere, chiedendo loro quali sono le forme che conoscono.
\begin{studente}[Bb]
  Quadrato, triangolo, rettangolo\dots{}
\end{studente}
\begin{studente}[Bb]
  Ma anche questo assomiglia a un rettangolo!
\end{studente}
(riferendosi al parallelogrammo)
\begin{tutor}[Ins]
  Hai ragione!
\end{tutor}

L'insegnante disegna e ritaglia su un foglio un parallelogrammo; dopo
aver tracciato le altezze, piega i due triangoli e chiede agli alunni
di ritagliarli.
\begin{tutor}[Ins]
  Che forma hanno i pezzettini di carta che avete ritagliato?
\end{tutor}
\begin{studente}[Bb]
  Triangoli
\end{studente}
\begin{tutor}[Ins]
  Cosa rimane?
\end{tutor}
\begin{studente}[Bb]
  Un rettangolo
\end{studente}
A questo punto, una bambina (normodotata) prende tra le tessere due
triangoli uguali a quelli ritagliati e li sovrappone al
parallelogrammo
\begin{studente}[Bb]
  Guarda maestra!
\end{studente}
Gli altri compagni cominciano a sovrapporre pezzi più piccoli su
quelli più grandi. Gli insegnanti ritengo opportuno orientare il
pensiero dei bambini verso l'obiettivo prefissato, quindi, si mostra
una tessera a forma di rettangolo chiedendo loro se è possibile
ricoprirla con altre tessere.

Un bambino, dopo alcuni tentativi, ricopre il rettangolo con due
quadratini. Gli altri, imitano il compagno. Gli insegnanti fanno
notare che con un quadratino è stata ricoperta solo una parte del
rettangolo e chiedono quale.
\begin{studente}[Bb]
  Metà
\end{studente}
A questo punto, per poter ricoprire il nostro rettangolo (intero), si
chiede quanti quadratini necessitano per formare la figura
considerata.
\begin{studente}[Bb]
  Due quadratini per creare il rettangolo intero
\end{studente}
Si introduce così il concetto di intero e della sua metà.

Le insegnanti chiedono ai bambini se altre forme possono essere
costituite da più pezzi della stessa dimensione. Alcuni alunni si
orientano a considerare un quadrato grande sovrapponendo due
rettangoli uguali, altri partendo sempre dal quadrato grande, invece,
lo ricoprono con quattro quadratini di ugual misura.

L'insegnante riferendosi alla prima sovrapposizione chiede ai bambini
di togliere un rettangolo e domanda quante parti ha preso su quelle
considerate.
\begin{studente}[Bb]
  Una pezzo \underline{su} due pezzi
\end{studente}
Il maestro riporta sulla lavagna la frazione corrispondente a quanto
appena detto dal bambino, spiegando che la linea (di frazione) indica
la parola ``\underline{su}'', che divide il quadrato in due parti
uguali, di cui ne è stata considerata una (un rettangolo).

La riflessione è fatta anche sul quadrato ricoperto da quattro
quadratini più piccoli.

\subsection{Osservazioni finali}

L'equipe di lavoro ha ritenuto valido e interessante il materiale
utilizzato in quanto ha permesso ai bambini diversamente abili di
avere un primo approccio al concetto di frazione sperimentando
concretamente e personalmente attraverso la manipolazione di
forme. Per consentire ai bambini un consolidamento dei contenuti
emersi, tenendo conto dei livelli cognitivi e della memoria a breve
termine, sarebbe opportuno riproporre il materiale attraverso
ulteriori attività, ripetute nel tempo.

Dalle attività è stato osservato che il gruppo dei maschi è risultato
più abile nella costruzione delle finestre e più creativo
nell'immaginare e riprodurre altre forme.

Le bambine, invece, hanno mostrato più attenzione alla grandezza e ai
colori delle varie figure, notando sin dall'inizio che fra tante
tessere vi sono forme uguali ma di diversa dimensione; di contro, non
hanno mostrato le stesse abilità dei maschi nell'incastrare e
nell'individuare facilmente la giusta disposizione delle tessere.

È possibile che questa diversità di competenze sia determinata dalla
differenza di genere dei nostri alunni, in quanto nella formazione dei
due gruppi si è cercato di equilibrare i diversi livelli cognitivi.

Tra i risultati non attesi, non era prevista l'osservazione di una
bambina sulla somiglianza del parallelogramma a un rettangolo, ciò è
stato per gli insegnanti spunto per introdurre a grandi linee
l'appartenenza di questa figura alla famiglia dei quadrilateri. Questo
argomento verrà affrontato e approfondito il prossimo anno,
inserendolo tra gli obiettivi didattici del PEI.

Al termine di queste osservazioni il gruppo di lavoro si ritiene
soddisfatto dei risultati ottenuti e per la collaborazione e l'intesa
raggiunta fra i colleghi di sostegno.

Si propone di sperimentare il kit in classe, progettando un percorso
con le colleghe curricolari che preveda delle attività tali da
consentire la partecipazione dei bambini diversamente abili, che
attraverso il tutoraggio dei compagni consentire un apprendimento tra
pari (peer education ovvero peer to peer).


\section[Sperimentazione \#4: seconda secondaria di primo
grado]{Sperimentazione \#4: classe seconda secondaria di primo grado,
  febbraio~2010}

\subsection{Osservazioni generali}

La sperimentazione si svolge in due classi seconde (classe ``A'' e
classe ``B'').

\subsubsection{Presentazione della classe}

La classe ``A'' è composta da 17 allievi, di cui 6 femmine e 11
maschi. La classe ``B'' è composta da 17, di cui 8 femmine e 9 maschi.

Il livello delle due classi è medio-basso in entrambe le
sezioni. Parecchi studenti provengono da un contesto socio-culturale
molto basso e presentano difficoltà di comprensione sia orale che
scritta e hanno tempi di attenzione e concentrazione molto
ridotti. Nella classe ``A'' è inserito un ragazzo russo di recente
immigrazione, un nomade in difficoltà e una ragazzina con disabilità
cognitiva. Nella ``B'' sono inseriti 3 studenti con disabilità
medio-grave di tipo cognitivo.

\subsubsection{Composizione dei gruppi}

Il primo laboratorio \attivita{Realizziamo una vetrata} è stato svolto
separatamente dalle due classi, mentre gli altri due a classi unite
per ovviare al problema della presenza di un unico insegnante a
seguire le attività. Ciò è stato possibile anche perché le due classi
sono poco numerose e caratterizzate quotidianamente da un elevato
numero di assenti. Si sono potuti così costituire 5 gruppi ciascuno
formato da 5-6 alunni.

\subsubsection{Insegnanti presenti}

Al primo incontro in ognuna delle due classi è presente solo
l'insegnante di classe. Negli incontri a classi unite sono presenti
entrambe le docenti.

\subsubsection{Calendarizzazione degli incontri}
\begin{calendario}
  \begin{itemize}
  \item 5 febbraio (classe ``A'') e 8 febbraio (classe ``B'')
  \item 10 febbraio (classi unite)
  \item 12 febbraio (classi unite)
  \end{itemize}
\end{calendario}

\subsection{Primo incontro}

\begin{description}
\item[Alunni presenti:] Classe ``A'': 14 alunni, classe ``B'': 13
  alunni
\item[Tempo effettivo di lavoro:]
 1h 40min
\end{description}

\begin{consegna}
  A ogni gruppo è stata distribuita una vetrata con le frazioni da
  1/2, 1/4, 1/8 e le relative istruzioni di gioco:
  \begin{itemize}
  \item gli alunni hanno liberamente manipolato il materiale a
    disposizione e sono state ascoltate le loro osservazioni
  \item successivamente a ogni gruppo è stata consegnata una vetrata
    con le frazioni da 1/2, 1/3, 1/4, 1/6, 1/9
  \item di nuovo gli alunni hanno manipolato il materiale e discusso
    sia all'interno del piccolo gruppo, sia a classe intera
  \item è stato svolto il gioco della Tombola seguito da discussione
  \end{itemize}

  \materiali{}%
  quelli previsti dal Kit in dotazione
\end{consegna}

\subsubsection{Osservazioni}

Il primo gioco \attivita{Realizziamo una vetrata} si è svolto a classi
separate ma secondo le stesse modalità: in ogni classe gli studenti si
sono divisi in 3 gruppi da 3-4 alunni.

In entrambe le classi gli alunni, lasciati liberi di scegliere i
compagni, hanno costituito due gruppi omogenei, uno di livello buono e
l'altro di alunni più in difficoltà, e un gruppo eterogeneo. Ogni
gruppo ha letto le istruzioni reagendo in maniera differente. Gli
studenti di buon livello hanno capito velocemente e si sono messi a
lavorare collaborando, nei gruppi misti i più attenti hanno trascinato
il gruppo riuscendo a coinvolgerli tutti, il gruppo degli alunni con
minori capacità ha avuto un'iniziale difficoltà perché impaurito di
fronte alla novità e l'insegnante ha dovuto in qualche caso
intervenire per dare suggerimenti e incoraggiare. L'atteggiamento di
ogni alunno non è stato molto dissimile da quello abitualmente
evidenziato durante le ore di lezione, cioè i più capaci seguono le
indicazioni e lavorano da subito autonomamente, gli altri hanno
bisogno di supporto e sollecitazioni.

Comunque in ciascun gruppo gli alunni hanno provato e riprovato con
entrambe le tipologie di frazioni e sono riusciti a arrivare al
concetto di equivalenza di aree ma non di frazioni, non sono riusciti
a identificare le varie tessere come frazioni perché nella loro
immaginazione rimanevano sempre e solo comunque superfici di
poligoni. A questo sono giunti, non tutti, solo dopo avere ricordato
insieme il percorso svolto in precedenza su tale concetto. I diversi
gruppi hanno interagito tra loro scambiando opinioni e i più brillanti
si sono spesi con i meno brillanti in spiegazioni e dimostrazioni,
cioè hanno mostrato, provando e riprovando, a ricoprire le singole
finestre o parti di finestra con pezzi diversi. L'attività è risultata
utile per i ragazzi con maggiori difficoltà che hanno incominciato a
acquisire, senza tuttavia prenderne consapevolezza, il concetto, molto
ostico, delle frazioni equivalenti. Successivamente si è giocato alla
Tombola utilizzando vetrate con le frazioni da 1/2, 1/3, 1/4, 1/6, 1/9
e scegliendo di comporre due o più tessere per volta per ottenere la
stessa frazione. L'interesse al gioco è stato apprezzabile e è
scattato un meccanismo di collaborazione all'interno del gruppo e di
competizione tra i diversi gruppi. Tutti hanno commesso l'errore di
iniziare con più pezzi a comporre una frazione restando alla fine
senza le giuste tessere per completare la vetrata. Dalla discussione
seguita, breve per il poco tempo a disposizione, e dalle riflessioni
scritte da parte di ogni studente è emerso che il laboratorio è
servito a ricordare un concetto già affrontato ma dimenticato, senza
però consolidarne la sua assimilazione.

Riportiamo alcune riflessioni dei ragazzi:
\begin{studente}[ ]
  \begin{itemize}
  \item \bambini{Mi è piaciuto questo gioco che ci ha fatto usare in
      pratica le frazioni}
  \item \bambini{Mi è servito per ripassare le frazioni equivalenti}
  \item \bambini{Non ho ben capito cosa c'entra l'area con le
      frazioni}
  \item \bambini{Mi è piaciuto fare un lavoro di matematica con gli
      altri}
  \end{itemize}
\end{studente}

\subsubsection{Consigli per i colleghi che vogliono proporre le stesse attività}

È necessario che i vari kit siano più ricchi di pezzi e che ci siano
più kit in modo che una stessa classe o scuola possa avere a
disposizione il materiale per più tempo.

Non sono richieste attenzioni particolari da parte dei docenti,
escluse quelle abituali che i docenti devono avere quando si lavora in
gruppo.

È positivo lavorare in gruppo e fare della matematica un'attività
laboratoriale però questa metodologia è difficilmente gestibile da un
solo insegnante che non può, nello stesso momento, organizzare
l'attività, rendere e mantenere il clima di lavoro adeguato, fare
osservazioni sull'atteggiamento e sull'apprendimento degli studenti,
intervenire nelle problematiche dei gruppi, rispondere a domande,
ecc\dots{}

In un futuro preferiremmo proporre il gioco \attivita{Realizziamo una
  vetrata} in classi nelle quali non abbiamo ancora affrontato
l'argomento frazioni nelle quali la scoperta può diventare più
facilmente un'acquisizione personale e pertanto rimane
significativamente nel proprio bagaglio di conoscenze.

\subsection{Secondo incontro}

Il secondo laboratorio è stato \attivita{Effetto domino}, svolto a
classi unite con 5 gruppi di 5 alunni ciascuno.

\begin{description}
\item[Alunni presenti:]25 alunni
\item[Tempo effettivo di lavoro:] 1h 40min
\end{description}
\begin{consegna}
  Giocare al gioco del domino con le frazioni %
  \materiali{}%
  quelli previsti dal Kit in dotazione
\end{consegna}

\subsubsection{Osservazioni}

I gruppi si sono costituiti spontaneamente tra allievi di classi
diverse senza problemi perché le nostre fasce di classi provenendo
dalla scuola primaria sottostante sono abituate a lavorare in più
occasioni insieme: ore opzionali pomeridiane, laboratori aperti,
uscite insieme, progetti didattici comuni.

L'approccio al gioco, già conosciuto dalla maggior parte degli
studenti, è avvenuto con più entusiasmo e con la sicurezza di essere
in grado di poterlo svolgere senza difficoltà, tranne poi accorgersi
in seguito di non essere sempre in grado di riconoscere le frazioni
equivalenti. Rispetto al gioco \attivita{Realizziamo una vetrata},
all'interno di ogni gruppo i più brillanti, dopo un primo momento di
incertezza nel riconoscimento delle frazioni equivalenti, hanno
monopolizzato la scelta delle tessere e la competizione tra gruppi è
stata più evidente. Gli studenti di livello intermedio hanno avuto
bisogno di un tempo più lungo per individuare le frazioni equivalenti
ma alla fine quasi tutti ci sono riusciti, senza però essere in grado
di esprimere a parole le conoscenze possedute e il ragionamento
seguito per individuare le tessere giuste. Gli alunni con maggiori
difficoltà, dopo un primo momento di entusiasmo, hanno continuato il
gioco con il timore di non essere all'altezza e, bloccati dalla paura
dell'insuccesso, non sono riusciti a individuare alcuna frazione
equivalente per cui, da un certo momento in poi, hanno assistito al
gioco senza più parteciparvi. Tutti coloro che invece hanno avuto la
percezione di avere compreso il concetto in questione, hanno
continuato con entusiasmo e con la voglia di ripetere il gioco una
volta conclusa la partita. Dalla discussione e dalle riflessioni
scritte
\begin{studente}[]
  \begin{itemize}
  \item \bambini{Con questo gioco ho imparato a riconoscere le
      frazioni equivalenti}
  \item \bambini{Ho capito bene le frazioni equivalenti perché ho
      potuto vederle}
  \item \bambini{Ho imparato a essere più veloce nei calcoli con le
      frazioni}
  \item \bambini{Ho ripassato le frazioni}
  \end{itemize}
\end{studente}
È emersa chiaramente la distinzione tra chi ha assimilato il concetto
e se ne è appropriato e chi invece ancora non è riuscito concretamente
a comprenderlo. L'appropriazione è risultata evidente nel momento in
cui, poco tempo dopo, è stato affrontato il nuovo argomento relativo
alle proporzioni. Gli alunni stessi hanno segnalato che in una
proporzione l'uguaglianza è tra due frazioni equivalenti. Partendo da
questi presupposti l'argomento è risultato agli alunni di più facile
comprensione e applicazione%
.

\subsection{Terzo incontro}

\begin{description}\item[Alunni presenti:] 26 alunni
\item[Tempo effettivo di lavoro:] 1h 40min
\end{description}
\begin{consegna}
  Giocare al \attivita{Viaggio segreto}
\materiali{}%
quelli previsti dal Kit in dotazione

\end{consegna}

\subsubsection{Osservazioni}

Quest'ultimo laboratorio è stato poco efficace e produttivo. A ogni
gruppo costituito come nel secondo incontro, è stato consegnato il
foglio delle istruzioni e il materiale necessario allo svolgimento del
gioco. L'atteggiamento dei gruppi è stato sostanzialmente diverso: i
gruppi con gli alunni più brillanti hanno messo in luce la volontà di
comprendere le istruzioni e trovare la soluzione. Hanno dovuto
comunque essere aiutati nella corretta comprensione di qualche
passaggio (ad esempio hanno capito solo dopo suggerimento
dell'insegnante che si dovesse procedere provando le varie
combinazioni). Gli studenti degli altri gruppi, dopo la lettura delle
prime righe di istruzioni, hanno rinunciato a ogni tentativo di
capire quali operazioni dovessero fare per eseguire il gioco. Non
hanno accettato le sollecitazioni dei docenti e neppure l'aiuto
offerto dai compagni che sono riusciti a portare a termine il
gioco. Di fatto due soli gruppi hanno svolto l'attività, un gruppo a
fatica l'ha iniziata ma il tempo limitato a disposizione non gli ha
permesso di concluderla, tra gli studenti degli altri due gruppi è
prevalso un atteggiamento rinunciatario di fronte alle prime
difficoltà. Poiché le indicazioni non venivano offerte dal docente ma
implicavano la fatica del comprendere un testo scritto, anziché
chiedere e continuare hanno reputato l'attività troppo difficile e
quindi, come al solito, hanno preferito rinunciare. Neppure i
compagni che hanno compreso le regole del gioco sono riusciti a
interagire con i compagni in modo tale da stimolare la volontà e il
desiderio di provare a capire e a mettersi in gioco. Quest'ultimo
laboratorio è stato vissuto dalla maggior parte degli studenti come un
momento noioso, per nulla proficuo, e ha attivato una sorta di
demotivazione all'attività. D'altra parte, quegli studenti che sono
riusciti a trovare la soluzione non ne hanno compreso il senso e cosa
avrebbero dovuto imparare. Il gioco, molto meno accattivante del
domino anche perché più complicato nelle istruzioni, non è riuscito a
catturare la curiosità iniziale che era stata il motore che ha
supportato le altre attività e soprattutto il gioco \attivita{Effetto
  domino}.


\section[Sperimentazione \#5: prima secondaria di primo grado]{Sperimentazione \#5:
  classe prima secondaria di primo grado,
  marzo~2010}

\subsection{Osservazioni generali}
\subsubsection{Presentazione della classe}

Classe composta da 23 alunni (12 maschi e 11 femmine) di cui 3 (2
maschi e 1 femmina) che necessitano del supporto dell'insegnante di
sostegno per problemi di apprendimento cognitivo. Il gruppo classe si
presenta vivace sia sotto l'aspetto comportamentale che sotto quello
propriamente didattico, dando il meglio di sé in quelle attività in
cui il loro coinvolgimento è diretto (es. esperimenti scientifici,
giochi matematici a squadre, esercitazioni a computer).

\subsubsection{Composizione dei gruppi}

I gruppi sono variati nel corso delle attività e sono quindi descritti
via via.

\subsubsection{Insegnanti presenti}
Agli incontri è presente solo l'insegnante sperimentatore. Si
auspicava la presenza dell'insegnante di sostegno, ma non è stato
possibile organizzare la compresenza.

\subsubsection{Calendarizzazione degli incontri}
\begin{calendario}
  \begin{itemize}
  \item 22 marzo (\attivita{Realizziamo una vetrata})
  \item 24 marzo (\attivita{Il gioco del domino})
  \item 26 marzo (\attivita{Il viaggio segreto})
  \end{itemize}
\end{calendario}

\subsection{Primo incontro}
\begin{description}
\item[Alunni presenti:] Tutti gli alunni erano presenti, compresi
  quelli con disabilità
\item[Tempo effettivo di lavoro:] due ore
\end{description}
\begin{consegna}
  Nel primo incontro si è svolta l'attività \attivita{Realizziamo una
    vetrata}.

  Ho iniziato a distribuire le vetrate (N. 1 - 5) con le relative
  superfici in gommapiuma + alcune tessere contenenti frazioni.

  Si è fatta notare l'importanza della conoscenza dei concetti base
  rispetto alla loro semplice esecuzione tecnica.

\materiali{}%
\begin{itemize}
\item le vetrate (N. 1 - 5) con le relative superfici in gommapiuma +
  alcune tessere contenenti frazioni
\item il foglio-istruzioni
\item lavagnetta per risposte di gruppo
\end{itemize}
\end{consegna}

\subsubsection{Osservazioni}

La classe è stata suddivisa in 5 gruppi chiedendo agli alunni più
riservati, indicati dal sottoscritto, di fare loro da capisquadra. La
proposta è stata accolta positivamente dagli studenti. Si sono formati
gruppi omogenei tra loro con un'equa distribuzione dei diversi livelli
di competenze e abilità matematiche al proprio interno, compresi i
portatori di disabilità anch'essi uniformemente suddivisi nei gruppi.

Ogni gruppo ha scelto come logo/distintivo un fumetto contenente un
simbolo matematico (potenze, triangoli, piramidi, sfere, radici
quadrate) da predisporre individualmente per l'incontro successivo, e
un proprio nome d'appartenenza scelto tra quello degli scienziati da
loro più conosciuti (nomi scelti: galileiani, newtoniani,
einsteiniani, leonardeschi, archimedei). Al termine di questo primo
momento aggregativo (15') si è passati alla spiegazione
dell'attività. Da subito è stato chiarito che la modalità di approccio
è sì di tipo ludico ma l'obiettivo finale da raggiungere è il
potenziamento sia delle conoscenze che delle competenze matematiche
individuali. È stato inoltre rimarcato il fatto che il kit didattico è
stato predisposto da degli esperti universitari e ciò ha creato nei
ragazzi un senso di coinvolgimento ancora maggiore verso l'attività
proposta. Come consigliato dalle maestre delle sperimentazioni
precedenti,
iniziando la distribuzione delle vetrate (N. 1 - 5) e le relative
superfici in gommapiuma assieme a alcune tessere frazionarie, ho
lasciato liberi i ragazzi di ricoprire tutti gli spazi predisposti,
permettendogli inoltre di osservare, nella rotazione presso gli altri
gruppi, come ognuno di essi avesse adottato modalità diverse di
esecuzione sia dal punto di vista cromatico che di copertura delle
superfici. Alla mia domanda
\begin{tutor}[Ins]
  Come mai ciò è stato possibile?
\end{tutor}
ogni gruppo ha cercato al proprio interno una risposta plausibile,
scrivendola poi su una lavagnetta, in dotazione come proprio materiale
scolastico. Due gruppi hanno citato in particolare i concetti di unità
frazionaria e di somma di frazioni corrispondenti a superfici da
ricoprire e, dopo rapida votazione per alzata di mano, è stata scelta
tra esse la risposta/definizione più completa e chiara la quale è
stata trascritta dal capogruppo sulla lavagna della classe con 3
esempi, scelti al momento, di frazioni (diverse da quelle considerate
nel gioco) scomposte in somme di unità frazionarie. Nel frattempo è
stato distribuito il foglio-istruzioni e si sono lasciati 5' per la
lettura di esso all'interno di ogni squadra. Ho notato che qualche
componente tendeva a privilegiare la fase manipolativa rispetto a
quella di conoscenza delle regole del gioco. Trascorso questo tempo,
un componente per gruppo, propostosi in modo volontario, ha riesposto
a memoria, sia ai propri compagni che a quelli degli altri gruppi, le
istruzioni apprese precedentemente (10').

Si è passati quindi alla vera e propria fase di realizzazione delle
vetrate, dopo una prima estrazione di prova e relativa spiegazione
pratica di posizionamento della relativa superficie. Già alla terza
estrazione di una tessera frazionaria, qualcuno ha fatto notare il
parallelismo con la classica tombola dei numeri interi aggiungendo
altresì la commistione di elementi aritmetici con quelli riferibili a
grandezze geometriche. Inizialmente tutti i gruppi cercavano di usare
pezzi esattamente corrispondenti ai contorni delle superfici
predisposte, quando però un gruppo ha fatto notare che non aveva la
figura esattamente corrispondente a quella frazione, si è aperto il
confronto con tutta la classe se non era possibile trovare una
soluzione alternativa, giungendo così alla conclusione che si possono
sommare più superfici per ottenere quella richiesta (es. 2 superfici
da ¼ sono equivalenti a 1 da ½). Dopo questo chiarimento tutti gli
studenti hanno capito meglio il concetto di frazione e relativa
corrispondenza con parti di figure geometriche, procedendo in modo più
spedito nell'avanzamento della composizione del proprio mosaico. Il
coinvolgimento emotivo dei ragazzi è stato sempre elevato, in ciò
testimoniato dal coro di esultanza dei singoli gruppi all'atto del
completamento di ogni singola finestra. La vittoria contemporanea di 2
formazioni ha comportato infine la verifica dei risultati la quale ha
richiesto una durata abbastanza lunga e, per alcuni studenti meno
riflessivi, un'attesa alquanto noiosa. Al termine del gioco si è fatta
una breve sintesi delle modalità di esecuzione dell'attività e, in
particolare, si è rimarcata la possibilità di poter usare più sagome
per ottenere una stessa frazione estratta. Tutti gli alunni hanno
mostrato di aver apprezzato questa nuova modalità di ``fare
matematica''. Per il giorno successivo è stato richiesto a ognuno di
comporre una breve relazione dell'attività svolta, corredandola degli
aspetti positivi e negativi che sono emersi durante il suo
svolgimento. Durata complessiva dell'attività (prima parte): 75'.

Dopo un momento di intervallo-merenda, si è passati alla realizzazione
delle vetrate 6 - 10 (vista come seconda manche/rivincita), mantenendo
gli stessi gruppi precedenti, ma ponendo quale elemento di novità
l'introduzione di frazioni a denominatore multiplo di 3 (e quindi
corrispondentemente un numero maggiore di pezzi di vetrata, subito
notato dai ragazzi, e quindi la richiesta di competenze matematiche
leggermente superiori rispetto a quelle della prima parte del
laboratorio). A gran richiesta ho dovuto lasciare un tempo iniziale
per la fase manipolativa libera. Ciò ha permesso di osservare e
prendere dimestichezza con le nuove frazioni (1/3, 1/6 e 1/9). Per
evitare il prolungato controllo finale, a ogni estrazione giravo tra
i gruppi per verificare l'esattezza della corrispondenza
frazione-superficie. Anche in questo modo però i tempi si sono
dilatati eccessivamente riconfermando quindi la necessità di un secondo
insegnante collaboratore (non mi è stato possibile averlo in questa
settimana per esigenze molto particolari all'interno dell'Istituto
Comprensivo). In questa seconda manche, ogni tessera estratta veniva
segnata sulla lavagna e si lasciava un breve momento di riflessione
per lasciar calcolare, a partire da tale valore, quale somma di unità
frazionarie poteva generarla e quindi facendo interagire tutti i
componenti del gruppo alla ricerca di variabili sia aritmetiche
(calcoli) che geometriche (corrispondenza di superfici). Grazie a una
organizzazione puntuale all'interno di ogni gruppo, senza l'intervento
del docente a imporlo dall'esterno, i tempi di esecuzione di questa
seconda manche si sono leggermente ridotti rispetto a quella
precedente. Al termine, l'osservazione verbale più interessante che
hanno fatto alcuni alunni è l'aver constatato che non è stato
fondamentale avere raggiunto l'obiettivo della vittoria finale del
proprio gruppo, bensì aver capito bene le modalità di svolgimento del
gioco, mediante un'azione di rinforzo dei concetti precedenti, e al
contempo essersi maggiormente divertiti in gruppo per completare la
realizzazione della vetrata a mosaico più complessa. Durata
complessiva dell'attività (seconda parte): 50'.

Osservazioni scritte degli alunni, riportate sinteticamente ma
fedelmente nella loro spontaneità:

\underline{aspetti positivi}:
\begin{studente}[ ]
  \begin{itemize}
  \item \bambini{il gioco è molto divertente e coinvolgente};
  \item \bambini{è stato un altro modo per imparare le frazioni,
      certamente molto divertente};
  \item \bambini{questo gioco è stato molto divertente, mi ha fatto
      comprendere molte cose sulla matematica e ci voleva molta
      astuzia e furbizia per cercare di vincere};
  \item \bambini{mi sono divertito perché parlavamo e ci confrontavamo
      come quando bisognava imparare le regole, e credo che ci siamo
      divertiti tutti quanti};
  \item \bambini{è un gioco impegnativo che ti fa capire le frazioni e
      saperle mettere a posto, ma era anche molto divertente};
  \item \bambini{la cosa che mi è piaciuta di più è quando dovevamo
      sistemare le figure nelle vetrate};
  \item \bambini{mi è piaciuto molto questo gioco perché è un modo
      diverso per imparare la matematica};
  \item \bambini{mi è piaciuto stare in compagnia e riuscire a
      completare pezzo per pezzo tutte le vetrate anche se a volte era
      un po' difficile perché non sempre usciva la frazione che ci
      serviva per finire una finestra};
  \item \bambini{questo gioco mi è piaciuto molto e spero di rifarlo};
  \item \bambini{è bello e è un gioco di logica che mette insieme le
      frazioni, le somme di frazioni e le aree di figure geometriche
      uguali alla frazione iniziale};
  \item \bambini{abbiamo fatto un'attività diversa dal solito e tutti
      insieme};
  \item \bambini{è stata un'esperienza entusiasmante perché ci siamo
      divertiti};
  \item \bambini{mi sono divertita molto perché è stato un nuovo modo
      di lavorare};
  \item \bambini{è un bel gioco d'intelligenza};
  \item \bambini{è stato bello perché ci siamo divertiti};
  \item \bambini{mi sono divertita molto e spero che si ripeta};
  \item \bambini{per me è stato molto divertente e ho imparato,
      finalmente, bene le frazioni}.
  \end{itemize}
\end{studente}
Si può ben notare che per aspetti positivi gli alunni privilegiano, in
genere, l'esternazione della loro emotività tenendo in sottordine
l'aspetto propriamente didattico-concettuale.

\underline{aspetti negativi}:
\begin{studente}[ ]
  \begin{itemize}
  \item \bambini{nella prima partita sono stato quasi l'unico del
      gruppo a ragionare};
  \item \bambini{dopo un po' di tempo mi sono annoiato, soprattutto
      nella prima partita, perché non uscivano le frazioni giuste};
  \item \bambini{non ci sono stati aspetti negativi (confermato da
      parecchi alunni)};
  \item \bambini{all'inizio era un po' confusionario perché non avevo
      capito bene le regole};
  \item \bambini{alla fine, quando una persona perde, gli altri ridono
      senza un perché (scritto da una ragazza un po' permalosa per
      natura)};
  \item \bambini{è stato un po' lungo perché non uscivano le frazioni
      che servivano al nostro gruppo per completare la finestra e
      quindi mi sono annoiato un po'};
  \item \bambini{all'inizio mi sembrava un gioco da piccoli perché si
      dovevano riempire le finestre, ma poi ho capito che è da grandi
      perché bisognava riempirle seguendo le frazioni e le operazioni
      con le frazioni};
  \item \bambini{a volte capitava che sbagliavamo a abbinare le figure
      o non capivamo bene lo svolgimento};
  \item \bambini{i pezzi sono abbastanza facili da perdere e da
      rompere e si spostano facilmente dalla vetrata}.
  \end{itemize}
\end{studente}

Gli alunni con disabilità hanno partecipato attivamente al gioco, in
particolare diventando attori fisici nella deposizione delle sagome
geometriche nei riquadri proposti dal gruppo e al tempo stesso
collaborando nella ricerca delle zone di inserimento. Assieme a altri
hanno però mal sopportato il conteggio finale per la proclamazione
ufficiale dei gruppi vincitori, momento importante per la regolarità
del gioco ma da snellire in futuro, se possibile.

\subsubsection{Consigli per i colleghi che vogliono proporre le stesse attività}

Per il futuro è il caso di pensare di lasciare a un altro docente
l'incarico del controllo dell'esito finale di questo gioco o, in
alternativa, il controllo sistematico della correttezza manipolativa
a ogni estrazione effettuata.

Conclusione: forse dovrebbe essere predisposto da parte del team
docenti somministratori un questionario di soddisfazione e
rielaborazione concettuale ben più articolato rispetto a quello da me
predisposto (spiegazione del gioco; aspetti positivi; aspetti
negativi) in modo che dalle risposte degli alunni si possano rendere
più oggettivi gli obiettivi metacognitivi raggiunti nelle
sperimentazioni laboratoriali praticate.

\subsection{Secondo incontro}

\begin{description}
\item[Alunni presenti:] Tutti gli alunni con disabilità erano
  presenti, mancava una alunna normodotata
\item[Tempo effettivo di lavoro:] due ore
\end{description}
\begin{consegna}
  Nel secondo incontro si è svolta l'attività \attivita{Effetto
    domino}.

  Presentazione dell'attività: si è presentata brevemente la nuova
  attività laboratoriale (frazioni equivalenti) come un ampliamento
  dei concetti appresi in quella precedente sulle vetrate, in modo da
  rendere esplicito il quadro generale all'interno del quale si andava
  a collocare ogni singola tappa del percorso formativo del progetto.

  Calcolo delle tessere: ho fatto calcolare dai ragazzi quante
  sarebbero state le tessere a cui aveva diritto ognuno di essi
  conoscendo il totale delle stesse e il numero complessivo dei
  partecipanti

  Scelta delle tessere: È stato fatto girare il mazzo completo con le
  tessere coperte e ognuno ha scelto le proprie

  Riportare definizioni alla lavagna: scrittura alla lavagna della
  definizione prescelta contenente i concetti di frazione ridotta ai
  minimi termini, frazioni equivalenti e classi di equivalenza con 3
  esempi relativi scelti al momento, indipendentemente dai contenuti
  delle tessere

  Lettura delle regole: distribuzione del foglio-regole e lettura
  collettiva delle stesse

  Regole del gioco: A turno ogni squadra doveva cercare di attaccare
  una sua tessera equivalente alle estremità di quelle presenti sul
  tavolo e riportarne i valori in lavagna, nel tempo massimo di 30-40''
  \materiali{}%
  \begin{itemize}
  \item tessere del domino
  \item lavagnetta per le risposte di gruppo
  \item foglio-regole
  \end{itemize}
\end{consegna}

\subsubsection{Osservazioni}

A livello organizzativo, gli alunni hanno presentato il loro
logo/distintivo (si veda l'attività precedente) e si è deciso quindi,
per giustificare il lavoro preparatorio così ben eseguito, di operare
per la formazione delle squadre con la seguente modalità: per la prima
manche si sarebbero mantenuti i gruppi già formati la volta scorsa
mentre per quella successiva si sarebbe dato mandato a altri alunni,
scelti dal sottoscritto tra quelli del livello intermedio, di fare da
capisquadra scambiando quindi tali distintivi tra gli alunni, gesto
significativo di unitarietà del gruppo classe \textit{in toto}.

Il concetto di frazioni equivalenti era stato sviluppato in classe la
settimana precedente, partendo dal libro di testo ``Sistema matematica
vol. 1 - Anna Montemurro - De Agostini''. Tale libro è strutturato
ponendo in successione la parte teorica di ogni argomento che forma
un'unità didattica con la parte prettamente esercitativa, in ciò
innovativo in quanto chiarisce rapidamente agli studenti, a partire
dalle definizioni e dai concetti base, quali contenuti applicativi
devono essere immediatamente acquisiti.

Ciò detto, si è presentata brevemente la nuova attività laboratoriale
(frazioni equivalenti) come un ampliamento dei concetti appresi in
quella precedente sulle vetrate, in modo da rendere esplicito il
quadro generale all'interno del quale si andava a collocare ogni
singola tappa del percorso formativo del progetto. Ho fatto calcolare
dai ragazzi quante sarebbero state le tessere a cui aveva diritto
ognuno di essi conoscendo il totale delle stesse e il numero
complessivo dei partecipanti (i più capaci hanno dedotto che il resto
della divisione indicava quanti erano coloro che avrebbero ricevuto
una tessera in più, che però nel computo finale per determinare il
gruppo vincitore doveva comunque essere tenuto in debita
considerazione e quindi risottratto). A questo punto è stato fatto
girare il mazzo completo con le tessere coperte e ognuno ha scelto le
proprie, condividendone poi i valori frazionari con gli altri
componenti del gruppo. Curiosamente ho notato che la scelta cadeva su
una composizione cromatica delle tessere che fosse la più varia
possibile, con nessun gruppo che si era orientato a acquisire carte
di uno stesso o al più di 2 colori, quasi a sottolineare che nella
variabilità fosse implicita l'unitarietà dei componenti (concetto
troppo filosofico o semplice casualità nelle scelte dei
ragazzi?). Sono stati lasciati 5' per cercare di capire il senso di
frazioni così diverse tra loro. Anche qui, analogamente all'attività
precedente, lavagnetta con risposta e condivisione delle
stesse. Successiva distribuzione del foglio-regole e lettura
collettiva delle stesse. Due studenti (un ragazzo e una ragazza) hanno
poi volontariamente rielaborato tali regole di fronte ai compagni che
hanno dimostrato di aver compreso le regole generali del gioco (il più
chiacchierone della classe ha proposto di chiamarlo \bambini{il colmo
  del prof di mate} perché sapeva della battuta \dots{}  \bambini{abitare
  in una frazione}, aggiungendoci questa volta, il
corollario\dots{} \bambini{equivalente}!  Velo pietoso, sigh). Si è dato
inizio al gioco dopo che ho posto su un tavolo al centro della classe
una tessera estratta in precedenza e contenete da entrambe le parti
frazioni irriducibili, segnandone sulla lavagna i valori, così come
quelli delle tessere successive, per aiutare i ragazzi nella
visualizzazione del percorso del domino. A turno ogni squadra doveva
cercare di attaccare una sua tessera equivalente alle estremità di
quelle presenti sul tavolo e riportarne i valori in lavagna, nel tempo
massimo di 30-40''. Dopo pochi giri ci si è resi conto però che
qualche gruppo non riusciva a procedere nel gioco in quanto in
possesso di frazioni che non si attaccavano facilmente. L'interesse
iniziale è quindi un po' scemato, e comunque si è riusciti a
concludere la gara con la proclamazione dei vincitori (2 gruppi, in
base al numero di tessere rimaste) nel momento in cui più nessuna
tessera poteva più essere aggiunta. Tempo necessario: 70'. Nel
riconsiderare i concetti affrontati, qualche gruppo si è accorto di
non aver posto sufficiente attenzione nel calcolo delle possibili
frazioni equivalenti tra quelle in proprio possesso. Per rendere
quindi decisamente più semplice e partecipata la manche successiva,
dopo il momento di ricomposizione dei gruppi coi nuovi capisquadra e
relativo scambio dei loghi/distintivi, si è optato per l'utilizzo di
un certo numero di giri totali (per la precisione 3) in cui era
consentito anche il collegamento con punti interni del domino in modo
da rendere doppia la possibilità di usare gli stessi e, al contempo,
l'opportunità di esplicitare ogni volta sulla lavagna qual era la
nuova frazione ridotta ai minimi termini che si aggiungeva nel
percorso. Ogni gruppo doveva fare tale calcolo segnandolo sulla
propria lavagnetta personale, con conseguente bonus di ½ punto per
ogni risposta esatta individuata, da usarsi nel conteggio a fine
gara. Queste 2 varianti hanno consentito una maggiore partecipazione
di tutti allo svolgimento del gioco e, al contempo, a un gruppo di
concludere il domino con l'inserimento di tutte le sue tessere nel
percorso. I ragazzi hanno esplicitato che questo secondo gioco
richiedeva, rispetto al primo, una maggiore padronanza delle abilità
di calcolo e dei concetti di equivalenza aritmetica e quindi una
maggior concentrazione per farsi trovare pronti al proprio turno, e lo
hanno definito più da scuola media rispetto al precedente. È piaciuta
molto anche l'idea di trasferire nuovi concetti aritmetici sulle
tessere scaricabili dal sito \matemm{} e infatti un gruppo ha proposto
di creare un nuovo domino di classe con brevi espressioni contenenti
elevamenti a potenza e relative proprietà, + risultati finali delle
stesse ``mini-espressioni''. Svilupperemo la nuova proposta
applicativa non appena, terminato il programma di prima media,
affronteremo il ripasso degli argomenti fondamentali appresi
(indicativamente nelle ultime 2 settimane di scuola). Per esigenze di
verifiche e interrogazioni programmate nei giorni successivi ho
evitato di richiedere ulteriori commenti scritti sul kit didattico
testato anche se la percezione diretta è sovrapponibile a quella
relativa alla prima attività laboratoriale effettuata. Tempo impiegato
per la seconda manche: 55' circa.

\subsubsection{Consigli per i colleghi che vogliono proporre le stesse attività}

Consiglio di anticipare l'attività con un'efficace spiegazione del
concetto di frazione equivalente, supportata da una verifica in
itinere sull'argomento. Valida anche la proposta, al termine del
laboratorio, di un sintetico questionario di soddisfazione + relative
proposte di estensione del gioco a altri argomenti precedenti o
successivi (in base alle loro competenze acquisite nel periodo delle
Elementari). Attenzione nel mantenere sempre alto l'interesse della
classe mediante un tono di voce e gestualità coinvolgenti e
contemporaneamente miranti all'acquisizione dei concetti basilari
proposti come obiettivi didattici.

Eventuali utilizzi di videocamere per filmare, rivedere all'interno
del gruppo e/o riproporlo a classi in orizzontale/verticale piuttosto
che al corpo docente durante un collegio unitario possono essere
ulteriori modalità di confronto e rinforzo dei rispettivi ruoli
studente/docente.

\subsection{Terzo incontro}
\begin{description}\item[Alunni presenti:]Tutti gli alunni erano presenti, compresi quelli con disabilità
\item[Tempo effettivo di lavoro:] Tempo complessivo: poco meno di 2 ore.
\end{description}
\begin{consegna}
  Nel terzo incontro si è svolta l'attività \attivita{Viaggio segreto}

  Data la complessità dell'aspetto crittografico del gioco, ho letto
  assieme a loro le regole generali e alla lavagna ho riportato, passo
  passo, l'esempio contenuto nel libretto delle istruzioni, facendo
  prendere appunti ragionati ai ragazzi sul loro quaderno di
  matematica. %
  \materiali{}%
  Kit didattico
\end{consegna}

\subsubsection{Osservazioni}

Per questa terza attività laboratoriale ho deciso di porre a
capisquadra i 5 migliori elementi della classe: 3 ragazzi e 2
ragazze. I gruppi risultanti sono apparsi omogenei tra loro e motivati
alla competizione.

Apparentemente il più complesso, in realtà il più coinvolgente per
quanto riguarda i ragazzi con buone capacità di ragionamento. Alla
presentazione del kit ha subito incuriosito i ragazzi il collegamento
della matematica con concetti inerenti la geografia
(multidisciplinarietà). L'aspetto motivazionale è stato quindi
raggiunto col semplice consegnare la cartina in ogni gruppo e
chiarendo che, per loro, questa era come la mappa del tesoro
all'interno della quale tracciare il percorso circolare di un
misterioso viaggio segreto. L'algoritmo risolutivo è stato quindi
ribadito con calma una seconda volta e, al termine della spiegazione,
ho fatto eseguire un ulteriore esempio chiarificatore da una ragazza
con capacità medie. Appurato che la maggior parte degli alunni avesse
colto il punto chiave del discorso (scomposizione in fattori primi e
traduzione dei prodotti, escludendo di volta in volta uno dei fattori
primi dal computo totale, per tradurli poi in termini combinatori
alfabetici), sono passato alla distribuzione delle coppie di numeri e
relativi codici cifrati, comprensivi di tabelle dei numeri primi e
calcolatrici. In realtà le calcolatrici inserite nel kit non sono
state utilizzate in quanto il giorno precedente avevo comunicato ai
ragazzi che avrebbero potuto utilizzare quella dei propri cellulari
(risparmio energetico per il provveditorato agli studi vicino al
100\%!). Al via tutti i gruppi hanno cominciato a macinare calcoli ma,
purtroppo, sono dovuto intervenire da lì a poco presso tutti i gruppi
per consentirgli, soluzioni alla mano, di conoscere divisori primi non
comuni come ad esempio il 41 o altri valori similari a doppia
cifra. Infatti su valori alti facevano fatica a procedere nel calcolo
(i criteri di divisibilità affrontati alle medie riguardano i numeri
2, 3, 4, 5, 6, 9, 10, 11, 25, 100, 1000). Ciò in quanto ritenevo
prioritaria la fase di traduzione della scomposizione finale in
combinazioni di lettere rispetto a quella eccessivamente meccanica
della ricerca di numeri primi non di ``uso comune'' e conseguente
abbandono della fase di ricerca per ``logoramento da calcolo''. Questo
anche per mantenere alto l'interesse del maggior numero di componenti
in ogni gruppo. Probabilmente la stessa attività proposta in una
seconda o in una terza media non avrebbe comportato lo stesso problema
nella ricerca di tutti i fattori primi! Con mia grande sorpresa,
invece, una volta superato lo scoglio della scomposizione, la fase
successiva della traduzione in nomi comuni di città è stata risolta
brillantemente (usando al di là dell'aspetto combinatorio anche il
buon senso per la ricerca di nomi di città abbastanza note), anche se
in tempi leggermente diversi da gruppo a gruppo, con completamento da
parte di tutti i partecipanti del percorso nella sua totalità. Anche
qui entusiasmo alle stelle quando un gruppo individuava una tappa del
viaggio misterioso. L'idea di un premio alla classe, consigliato nel
libretto istruzioni, non è stato portato avanti
da me
per una
dimenticanza ma i ragazzi mi sono sembrati comunque pienamente
soddisfatti dal risultato matematico-geografico raggiunto con le loro
forze mentali. Aggiungo il particolare che più di un gruppo, al
termine della soluzione del proprio percorso parziale mi ha chiesto di
poter risolvere anche quello di altri gruppi, segno inequivocabile di
predisposizione per questa attività che unisce capacità di calcolo,
mente logica, creatività, multidisciplinarietà e spirito di gruppo.

\subsubsection{Consigli per i colleghi che vogliono proporre le stesse attività}

L'utilizzo dei cellulari/calcolatrici personali rendono gli studenti
più protagonisti e ciò favorisce un approccio corretto a una attività
di una certa complessità e che richiede un livello superiore di
concentrazione rispetto a quelle precedenti.

A maggior ragione in questo laboratorio, rispetto al precedente,
bisognerebbe sforzarsi di escogitare o lasciar escogitare da un gruppo
di studenti particolarmente motivato, altri usi di argomenti
matematici e non matematici (arte; storia; letteratura; sport;\dots{})
che uniti assieme siano di rinforzo didattico reciproco. La curiosità
e il sano protagonismo stimolano sempre l'acquisizione di nuove
modalità di apprendimento nei giovani studenti.

Durante una supplenza nella medesima settimana, ho proposto questa
attività in una seconda media particolarmente problematica dal punto
di vista disciplinare, e i risultati sono stati esattamente
sovrapponibili a quelli sperimentati nella mia prima classe.

\subsection{Conclusioni finali}

La modalità full-immersion settimanale mi pare abbia dato risultati
positivi, la cosa interessante per il futuro sarà poter avere questi
e/o altri kit analoghi (perché non anche in ambito prettamente
scientifico sperimentale?) come propria dotazione standard o scambiati
tra scuole aderenti all'iniziativa + scambio continuo di esiti
sperimentali personalizzati (sullo stile della creazione di un nuovo
kit-domino di scuola) partendo da materiali scaricabili in rete da
internet/piattaforma
dell'università.

\backmatter{}
\tableofcontents{}
\end{document}